\documentclass[12pt]{article}
\usepackage{amsmath,amssymb,amsfonts}     % AMS LateX is standard nowadays.
\usepackage{graphics}                     % For picture environment.
\usepackage{graphicx}                     % For picture environment.
\date{}
\newtheorem{lemma}{Lemma}[section]
\newtheorem{theorem}{Theorem}[section]

\newcommand{\adfQED}{\hfill $\square$}           % QED.

\newcommand{\adfnull}[1]{}                       % Comment.

\newcommand{\adfmod}[1]{~(\mathrm{mod}~#1)}      % Attempt to do (mod #1) properly.

\newcommand{\adfTgap}{\vskip 2mm}                % Start of a part of a theorem or proof.

\newcommand{\adfPgap}{\vskip 2mm}                % Start of a paragraph after a proof.

\newcommand{\adfInd}{}                     % Indentation

\begin{document}
\title{\bf Group divisible designs with block size 4 and type $g^u m^1$ -- II\\}
\author{A.~D.~Forbes\\
        School of Mathematics and Statistics,\\
        The Open University,\\
        Walton Hall, Milton Keynes MK7 6AA, UK\\
        {\sf anthony.d.forbes@gmail.com}
%       \\ Version 9.15, 19 Jun 2018  % HIDE FOR FINAL SUBMISSION.
        }
\maketitle
\begin{abstract}
\noindent We show that the necessary conditions for the existence of 4-GDDs of type
$g^u m^1$ are sufficient for
$g \equiv 0$ (mod $h$), $h$ = 39, 51, 57, 69, 87, 93, and for
$g$ = 13, 17, 19, 23, 25, 29, 31 and 35.
More generally, we show that for all $g \equiv 3 \adfmod{6}$, the possible exceptions occur only when
$u = 8$ and $g$ is not divisible by any of 9, 15, 21, 33, 39, 51, 57, 69, 87 or 93.
Consequently we are able to extend the known spectrum for $g \equiv 1$ and $g \equiv 5$ (mod 6).
Also we complete the spectrum for 4-GDDs of type $(3a)^4 (6a)^1 (3b)^1$.
\end{abstract}



%%%%%%%%%%%%%%%%%%%%%%%%%%%%%%%%%%%%%%%%%%%%%%%%%%%%%%%%%%%%%%%%%%%%%%%%%%%%%%%%%%%%%%%%%%
%%%%%%%%%%%%%%%%%%%%%%%%%%%%%%%%%%%%%%%%%%%%%%%%%%%%%%%%%%%%%%%%%%%%%%%%%%%%%%%%%%%%%%%%%%
%%%%%%%%%%%%%%%%%%%%%%%%%%%%%%%%%%%%%%%%%%%%%%%%%%%%%%%%%%%%%%%%%%%%%%%%%%%%%%%%%%%%%%%%%%

\section{Introduction}\label{sec:Introduction}

For the purpose of this paper, a {\em group divisible design}, $K$-GDD, of type $g_1^{u_1} \dots g_r^{u_r}$ is an ordered triple ($V,\mathcal{G},\mathcal{B}$) such that:
\begin{enumerate}
\item[(i)]{$V$
is a base set of cardinality $u_1 g_1 + \dots + u_r g_r$;}
\item[(ii)]{$\mathcal{G}$
is a partition of $V$ into $u_i$ subsets of cardinality $g_i$, $i = 1, \dots, r$, called \textit{groups};}
\item[(iii)]{$\mathcal{B}$
is a non-empty collection of subsets of $V$ with cardinalities $k \in K$, called \textit{blocks}; and}
\item[(iv)]{each pair of elements from distinct groups occurs in precisely one block but no pair of
elements from the same group occurs in any block.}
\end{enumerate}
A $\{k\}$-GDD is also called a $k$-GDD, and
a $k$-GDD of type $q^k$ is a {\em transversal design}, TD$(k,q)$.

A {\em parallel class} in a group divisible design is a subset of the block set in which each element of the base set appears exactly once.
A $k$-GDD is called {\em resolvable}, and is denoted by $k$-RGDD, if the entire set of blocks can be partitioned into parallel classes.
A $k$-RGDD of type $q^k$ is also called a {\em resolvable transversal design}, RTD$(k,q)$.
If there exist $k$ mutually orthogonal Latin squares (MOLS) of side $q$, then there exists
a $(k+2)$-GDD of type $q^{k+2}$ and
a $(k+1)$-RGDD of type $q^{k+1}$, \cite[Theorem III.3.18]{AbelColbournDinitz2007}.
Furthermore, as is well known, there exist $q - 1$ MOLS of side $q$ if $q$ is a prime power.

Group divisible designs are useful and important structures.
They provide the main ingredients for establishing the existence of infinite classes of various combinatorial objects by means of
a standard technique known as Wilson's Fundamental Construction, \cite{WilsonRM1972}, \cite[Theorem IV.2.5]{GreigMullen2007}.
To make use of Wilson's Fundamental Construction one needs group divisible designs, and therefore the problem of existence naturally arises.
However the general existence problem appears to be a long way from being solved.
Nevertheless, considerable progress has been achieved for GDDs where the block size is a small constant and the
group type is restricted to either $g^u$ or $g^u m^1$.
In particular, the existence spectrum for $k$-GDDs has been completely determined in the following cases.
\begin{enumerate}
\item[(i)]{
A 2-GDD of type $g_1^{u_1} \dots g_r^{u_r}$ is essentially a complete multipartite graph $K_{g_1^{u_1} \dots\ g_r^{u_r}}$.}
\item[(ii)]{
A 3-GDD of type $g^u$ exists if and only if
$u \ge 3$, $g(u - 1) \equiv 0 \adfmod{2}$ and $g^2 u (u - 1) \equiv 0 \adfmod{6}$,
\cite{Hanini1975}, \cite[Theorem IV.4.1]{Ge2007}.}
\item[(iii)]{
A 3-GDD of type $g^u m^1$ with $m > 0$ and $m \neq g$ exists if and only if
$u \ge 3$, $gu \equiv g(u - 1) + m \equiv 0 \adfmod{2}$, $g^2 u (u - 1) + 2 g u m \equiv 0 \adfmod{6}$ and $m \le g(u - 1)$,
\cite{ColbournHoffmanRees1992}, \cite[Theorem IV.4.4]{Ge2007}.}
\item[(iv)]{
A 4-GDD of type $g^u$ exists if and only if
$u \ge 4$, $g(u - 1) \equiv 0 \adfmod{3}$ and $g^2 u (u - 1) \equiv 0 \adfmod{12}$,
except that there is no 4-GDD of type $2^4$ or $6^4$,
\cite{BrouwerSchrijverHanani1977}, \cite[Theorem IV.4.6]{Ge2007}.}
\end{enumerate}
The next case one would naturally consider is that of 4-GDDs of type $g^u m^1$, and here a partial solution has been achieved
in the sense that for each $g$ there are at most a small number of $u$ where existence remains undecided,
\cite{GeRees2002}, \cite{GeReesZhu2002}, \cite{GeLing2004}, \cite{GeRees2004}, \cite{GeLing2005}, \cite{Schuster2010},
\cite{WeiGe2013}, \cite{WeiGe2014}, \cite{Schuster2014}, \cite{WeiGe2015}, \cite{ForbesForbes2018}.
For theorems concerning the existence of GDDs of other types, including 5-GDDs of type $g^u$,
see \cite{Ge2007} and the references therein.

Assuming $m > 0$ and $m \neq g$,
we review the necessary conditions for the existence of a 4-GDD of type $g^u m^1$:
\begin{align*}
&u \ge 4,~~ m \le g(u - 1)/2,\\
&gu \equiv g (u - 1) + m \equiv 0 \adfmod{3}, \textrm{~and~}\\
&\dfrac{g^2 u (u-1)}{12} + \dfrac{g u m}{6}, \textrm{~the~number~of~blocks,~is~a~positive~integer}.
\end{align*}
They are illustrated in Table~\ref{tab:4-GDD g^u m^1 necessary}.
Given $g$, we say that $u$ and $m$ are {\em admissible} if they satisfy these conditions.
\begin{table}[h]
\begin{center}
\label{tab:4-GDD g^u m^1 necessary}
\caption{Necessary conditions for 4-GDDs of type $g^u m^1$, $m > 0$, $m \neq g$}
\begin{tabular}{l@{~~~~}l@{~~~~}l}
 & & \\
$g \equiv 0 \adfmod{6}$ & -                        & $m \equiv 0 \adfmod{3}$\\
\hline
$g \equiv 1 \adfmod{6}$ & $u \equiv 0 \adfmod{12}$ & $m \equiv 1 \adfmod{3}$\\
                        & $u \equiv 3 \adfmod{12}$ & $m \equiv 1 \adfmod{6}$\\
                        & $u \equiv 9 \adfmod{12}$ & $m \equiv 4 \adfmod{6}$\\
\hline
$g \equiv 2 \adfmod{6}$ & $u \equiv 0 \adfmod{3}$  & $m \equiv 2 \adfmod{3}$\\
\hline
$g \equiv 3 \adfmod{6}$ & $u \equiv 0 \adfmod{4}$  & $m \equiv 0 \adfmod{3}$\\
                        & $u \equiv 1 \adfmod{4}$  & $m \equiv 0 \adfmod{6}$\\
                        & $u \equiv 3 \adfmod{4}$  & $m \equiv 3 \adfmod{6}$\\
\hline
$g \equiv 4 \adfmod{6}$ & $u \equiv 0 \adfmod{3}$  & $m \equiv 1 \adfmod{3}$\\
\hline
$g \equiv 5 \adfmod{6}$ & $u \equiv 0 \adfmod{12}$ & $m \equiv 2 \adfmod{3}$\\
                        & $u \equiv 3 \adfmod{12}$ & $m \equiv 5 \adfmod{6}$\\
                        & $u \equiv 9 \adfmod{12}$ & $m \equiv 2 \adfmod{6}$
\end{tabular}
\end{center}
\end{table}
As for sufficiency, we summarize here the main conclusions of \cite{Schuster2010}, \cite{WeiGe2015} and \cite{ForbesForbes2018} as a theorem.
\begin{theorem}
\label{thm:4-GDD g^u m^1 existence}
The necessary conditions for the existence of a 4-GDD of type $g^u m^1$ with $m > 0$ and $m \neq g$ are sufficient
except that there is no 4-GDD of type $2^6 5^1$ and except possibly for the following:
\begin{enumerate}
\item[]{for $g \equiv 3 \adfmod{6}$, $g \ge 39$, $g$ not divisible by $15$, $21$, $27$ or $33$: \\
      $u \in \{7, 11\}$ and $0 < m < g$;}
\item[]{for $g \equiv 3 \adfmod{6}$, $g \ge 39$, $g$ not divisible by $9$, $15$, $21$ or $33$: \\
      $u = 8$ and $3g < m < (7g - 3)/2$;}
\item[]{for $g \equiv 1, 5 \adfmod{6}$, $g \ge 13$: \\
      $u = 9$,  \\
      $u \in \{12, 24, 15, 27, 39, 51\}$ and $0 < m < g$,  \\
      $u \in \{21, 33\}$ and $0 < m < 4g$, \\
      $u = 24$ and $10g < m < (23g - 3)/2$;}
\item[]{for $g \equiv 2, 4 \adfmod{6}$, $g \ge 14$, $g \neq 16$:\\
      $u \in \{6, 9\}$, \\
      $u \in \{12, 15, 18, 21, 27\}$ and $0 < m < g$.}
\end{enumerate}
\end{theorem}
\noindent {\bf Proof~} See \cite[Theorem 3.5]{Schuster2010} for $g = 16$, \cite[Theorem 7.1]{WeiGe2015} and \cite{ForbesForbes2018}.
~\adfQED

\adfPgap
The purpose of this paper is to extend the known spectrum of 4-GDDs of type $g^u m^1$ by
eliminating some of the possible exceptions of Theorem~\ref{thm:4-GDD g^u m^1 existence} for $g \equiv 1, 3, 5 \adfmod{6}$.
We clear all possible exceptions for $g \equiv 3 \adfmod{6}$, $u \in$ \{7, 11\}.
We achieve a significant reduction for $g \equiv 3 \adfmod{30}$, $u = 8$, and consequently
for $g \equiv 1, 11 \adfmod{30}$, $u = 24$.
Also, the spectrum of 4-GDDs of type $(3rq)^u m^1$ is now determined for $r \ge 1$ and all primes $q \le 31$.
Moreover, we reduce the number of possible exceptions for $g \equiv 1,5 \adfmod{6}$, $u \in$ \{21, 33\},
and we clear them all when $g \in$ \{13, 17, 19, 23, 25, 29, 31, 35\}.

Section~\ref{sec:new 4-GDDs} is a catalogue of new 4-GDDs which we use in
Sections~\ref{sec:Theorems g = 3 (mod 6)} and \ref{sec:Theorems g = 1, 5 (mod 6)}, where we obtain our main results,
Theorems~\ref{thm:g = 3 (mod 6)}, \ref{thm:g = 39r 51r 57r 69r}, \ref{thm:g = 1, 5 (mod 6), u = 21, 33} and \ref{thm:g = 13 17 19 23 25}.
These are combined, together with Theorem~\ref{thm:4-GDD g^u m^1 existence}, in the following, a concise summary of,
as far as we are aware, the latest situation regarding the solution of the existence problem for 4-GDDs of type $g^u m^1$.
\begin{theorem}
\label{thm:4-GDD new g^u m^1 existence}
The necessary conditions for the existence of a 4-GDD of type $g^u m^1$ with $m > 0$ and $m \neq g$ are sufficient
except that there is no 4-GDD of type $2^6 5^1$ and except possibly for the following:
\begin{enumerate}
\item[]{for $g \equiv 3 \adfmod{6}$, $g \ge 111$, $g$ not divisible by $9$, $15$, $21$, $33$, $39$, $51$, $57$, $69$, $87$ or $93$: \\
      $u = 8$ and $m' < m < (7g - 3)/2$, where $m' = (17g - 6)/5$ if $g \equiv 3 \adfmod{30}$, $m' = 3g$ otherwise;}
\item[]{for $g \equiv 1, 5 \adfmod{6}$, $g \ge 37$: \\
      $u = 9$, \\
      $u \in \{12, 24, 15, 27, 39, 51, 21, 33\}$ and $0 < m < g$, \\
      $u = 24$ and $m'' < m < (23g - 3)/2$, where $m'' = (56g - 6)/5$ if $g \equiv 1$ or $11 \adfmod{30}$, $m'' = 10g$ otherwise;}
\item[]{for $g \equiv 2, 4 \adfmod{6}$, $g \ge 14$, $g \neq 16$:\\
      $u \in \{6, 9\}$, \\
      $u \in \{12, 15, 18, 21, 27\}$ and $0 < m < g$.}
\end{enumerate}
\end{theorem}
\noindent{\bf Proof}~ This follows from
Theorems~\ref{thm:4-GDD g^u m^1 existence},
\ref{thm:g = 3 (mod 6)}, \ref{thm:g = 39r 51r 57r 69r}, \ref{thm:g = 1, 5 (mod 6), u = 21, 33} and \ref{thm:g = 13 17 19 23 25}.
~\adfQED

\adfPgap
Of course, Theorem~\ref{thm:4-GDD new g^u m^1 existence} is not a complete statement of what is known.
It omits some instances of \cite[Theorem IV.4.11]{Ge2007} concerning the minimum and maximum values of $m$,
and it does not include various other minor improvements.
For example, as indicated by the remarks following the proof of Theorem~\ref{thm:g = 3 (mod 6)},
it is possible to increase $m'$ and $m''$ for some values of $g$.
We have also omitted new results concerning $g \equiv 2, 4 \adfmod{6}$.
Observe that the spectrum is completely determined for each $g \le 36$ except 14, 20, 22, 26, 28, 32 and 34.
However, these seven cases have actually been settled; the details will appear in a future paper.

%%%%%%%%%%%%%%%%%%%%%%%%%%%%%%%%%%%%%%%%%%%%%%%%%%%%%%%%%%%%%%%%%%%%%%%%%%%%%%%%%%%%%%%%%%
%%%%%%%%%%%%%%%%%%%%%%%%%%%%%%%%%%%%%%%%%%%%%%%%%%%%%%%%%%%%%%%%%%%%%%%%%%%%%%%%%%%%%%%%%%
%%%%%%%%%%%%%%%%%%%%%%%%%%%%%%%%%%%%%%%%%%%%%%%%%%%%%%%%%%%%%%%%%%%%%%%%%%%%%%%%%%%%%%%%%%
%%%%%%%%%%%%%%%%%%%%%%%%%%%%%%%%%%%%%%%%%%%%%%%%%%%%%%%%%%%%%%%%%%%%%%%%%%%%%%%%%%%%%%%%%%

\section{New 4-GDDs}\label{sec:new 4-GDDs}

The 4-GDDs that we use to prove our theorems are given here, grouped into lemmas for easy reference.
The blocks of the designs are generated from sets of base blocks by appropriate mappings.
The base blocks themselves and the instructions for expanding them are collected together in the Appendix.
If the point set of the 4-GDD has $v$ elements, it is represented by $Z_v = \{0, 1, \dots, v-1\}$ partitioned into groups as indicated.
The expression $a\adfmod{b}$ denotes the integer $n$ such that $0 \le n < b$ and $b \,|\, n - a$.
Where applicable, addition in $\mathbb{Z}_{gu/3} \times \mathbb{Z}_3$ is denoted by $\oplus$ and for brevity we represent element $(a,b)$ of this group by the number $3a+b$.
As an aid to checking the correctness of the designs we also provide coded forms of the block generation instructions.
% In case of inconsistency please ignore the English form.
It is a possibility that some of the designs presented in this section might be constructible by other means.

% % FOR THE SUBMITTED VERSION:
% To save space, sections~C--O of the rather lengthy Appendix have been omitted.
% However, in the unabridged version of this paper, which is available on {\sf ArXiv.org},
% [[ArXiv location]],
% all sections of the Appendix are present and hence the proofs of all the lemmas in this section are complete.
% % NOTE THAT THE REFERENCES TO APPENDIXES C-O WILL HAVE TO BE HARD-CODED.

% 4-GDDs a^u b^1 c^1:       5   a^u b^1 c^1    5
% Designs for g = 39:       5
% Designs for g = 51:       6
% Designs for g = 57:       8
% Designs for g = 69:       6
% Designs for g = 87:       2
% Designs for g = 93:       2   3 (mod 6)     29
% Designs for g = 13:       9
% Designs for g = 17:      16
% Designs for g = 19:      19
% Designs for g = 23:      17
% Designs for g = 25:      22
% Designs for g = 29:      26
% Designs for g = 31:      19
% Designs for g = 35:      16   1,5 (mod 6)  144   TOTAL 178
%
% ALL DESIGNS ARE IN THE APPENDIX.

%%%%%%%%%%%%%%%%%%%%%%%%%%%%%%%%%%%%%%%%%%%%%%%
%%%%%% a^u b^1 c^1
%%%%%%%%%%%%%%%%%%%%%%%%%%%%%%%%%%%%%%%%%%%%%%%
\begin{lemma}
\label{lem:4-GDD a^u b^1 c^1}
There exist 4-GDDs of types
$ 9^4 18^1 15^1 $,
$ 12^4 24^1 21^1 $,
$ 18^4 36^1 21^1 $,
$ 18^4 36^1 33^1 $ and
$ 9^6 27^1 15^1 $.
\end{lemma}

\noindent{\bf Proof}~ The designs are presented in Appendix~\ref{app:4-GDD a^u b^1 c^1}.
~\adfQED

%%%%%%%%%%%%%%%%%%%%%%%%%%%%%%%%%%%%%%%%%%%%%%%
%%%%%% 39^8 m^1
%%%%%%%%%%%%%%%%%%%%%%%%%%%%%%%%%%%%%%%%%%%%%%%
\begin{lemma}
\label{lem:4-GDD 39^8 m^1}
There exist 4-GDDs of types
$ 39^8 120^1 $,
$ 39^8 123^1 $,
$ 39^8 126^1 $,
$ 39^8 129^1 $ and
$ 39^8 132^1 $.
\end{lemma}

\noindent{\bf Proof}~ The designs are presented in Appendix~\ref{app:4-GDD 39^8 m^1}.
~\adfQED

%%%%%%%%%%%%%%%%%%%%%%%%%%%%%%%%%%%%%%%%%%%%%%%
%%%%%% 51^8 m^1
%%%%%%%%%%%%%%%%%%%%%%%%%%%%%%%%%%%%%%%%%%%%%%%
\begin{lemma}
\label{lem:4-GDD 51^8 m^1}
There exist 4-GDDs of types
$ 51^8 159^1 $,
$ 51^8 162^1 $,
$ 51^8 165^1 $,
$ 51^8 168^1 $,
$ 51^8 171^1 $ and
$ 51^8 174^1 $.
\end{lemma}

\noindent{\bf Proof}~ The designs are presented in Appendix~\ref{app:4-GDD 51^8 m^1}.
~\adfQED

%%%%%%%%%%%%%%%%%%%%%%%%%%%%%%%%%%%%%%%%%%%%%%%
%%%%%% 57^8 m^1
%%%%%%%%%%%%%%%%%%%%%%%%%%%%%%%%%%%%%%%%%%%%%%%
\begin{lemma}
\label{lem:4-GDD 57^8 m^1}
There exist 4-GDDs of types
$ 57^8 174^1 $,
$ 57^8 177^1 $,
$ 57^8 180^1 $,
$ 57^8 183^1 $,
$ 57^8 186^1 $,
$ 57^8 189^1 $,
$ 57^8 192^1 $ and
$ 57^8 195^1 $.
\end{lemma}

\noindent{\bf Proof}~ The designs are presented in Appendix~\ref{app:4-GDD 57^8 m^1}.
~\adfQED

%%%%%%%%%%%%%%%%%%%%%%%%%%%%%%%%%%%%%%%%%%%%%%%
%%%%%% 69^8 m^1
%%%%%%%%%%%%%%%%%%%%%%%%%%%%%%%%%%%%%%%%%%%%%%%
\begin{lemma}
\label{lem:4-GDD 69^8 m^1}
There exist 4-GDDs of types
$ 69^8 222^1 $,
$ 69^8 225^1 $,
$ 69^8 228^1 $,
$ 69^8 231^1 $,
$ 69^8 234^1 $ and
$ 69^8 237^1 $.
\end{lemma}

\noindent{\bf Proof}~ The designs are presented in Appendix~\ref{app:4-GDD 69^8 m^1}.
~\adfQED

%%%%%%%%%%%%%%%%%%%%%%%%%%%%%%%%%%%%%%%%%%%%%%%
%%%%%% 87^8 m^1
%%%%%%%%%%%%%%%%%%%%%%%%%%%%%%%%%%%%%%%%%%%%%%%
\begin{lemma}
\label{lem:4-GDD 87^8 m^1}
There exist 4-GDDs of types
$ 87^8 297^1 $ and
$ 87^8 300^1 $.
\end{lemma}

\noindent{\bf Proof}~ The designs are presented in Appendix~\ref{app:4-GDD 87^8 m^1}.
~\adfQED

%%%%%%%%%%%%%%%%%%%%%%%%%%%%%%%%%%%%%%%%%%%%%%%
%%%%%% 93^8 m^1
%%%%%%%%%%%%%%%%%%%%%%%%%%%%%%%%%%%%%%%%%%%%%%%
\begin{lemma}
\label{lem:4-GDD 93^8 m^1}
There exist 4-GDDs of types
$ 93^8 318^1 $ and
$ 93^8 321^1 $.
\end{lemma}

\noindent{\bf Proof}~ The designs are presented in Appendix~\ref{app:4-GDD 93^8 m^1}.
~\adfQED

%%%%%%%%%%%%%%%%%%%%%%%%%%%%%%%%%%%%%%%%%%%%%%%
%%%%%% 13^u m^1
%%%%%%%%%%%%%%%%%%%%%%%%%%%%%%%%%%%%%%%%%%%%%%%
\begin{lemma}
\label{lem:4-GDD 13^u m^1}
There exist 4-GDDs of types
$13^{12}  7^1$,
$13^{12} 10^1$,
$13^{9} 10^1$,
$13^{9} 16^1$,
$13^{9} 22^1$,
$13^{9} 28^1$,
$13^{9} 34^1$,
$13^{9} 40^1$ and
$13^{9} 46^1$.
\end{lemma}

\noindent{\bf Proof}~ The designs are presented in Appendix~\ref{app:4-GDD 13^u m^1}.
~\adfQED

%%%%%%%%%%%%%%%%%%%%%%%%%%%%%%%%%%%%%%%%%%%%%%%
%%%%%% 17^u m^1
%%%%%%%%%%%%%%%%%%%%%%%%%%%%%%%%%%%%%%%%%%%%%%%
\begin{lemma}
\label{lem:4-GDD 17^u m^1}
There exist 4-GDDs of types
$17^{12}  2^1$,
$17^{12}  5^1$,
$17^{12}  8^1$,
$17^{12} 11^1$,
$17^{12} 14^1$,
$17^{15} 11^1$,
$17^{9}  8^1$,
$17^{9} 14^1$,
$17^{9} 20^1$,
$17^{9} 26^1$,
$17^{9} 32^1$,
$17^{9} 38^1$,
$17^{9} 44^1$,
$17^{9} 50^1$,
$17^{9} 56^1$ and
$17^{9} 62^1$.
\end{lemma}

\noindent{\bf Proof}~ The designs are presented in Appendix~\ref{app:4-GDD 17^u m^1}.
\adfQED

%%%%%%%%%%%%%%%%%%%%%%%%%%%%%%%%%%%%%%%%%%%%%%%
%%%%%% 19^u m^1
%%%%%%%%%%%%%%%%%%%%%%%%%%%%%%%%%%%%%%%%%%%%%%%
\begin{lemma}
\label{lem:4-GDD 19^u m^1}
There exist 4-GDDs of types
$19^{12}  7^1$,
$19^{12} 10^1$,
$19^{12} 13^1$,
$19^{12} 16^1$,
$19^{24} 13^1$,
$19^{24} 16^1$,
$19^{15} 13^1$,
$19^{9} 10^1$,
$19^{9} 16^1$,
$19^{9} 22^1$,
$19^{9} 28^1$,
$19^{9} 34^1$,
$19^{9} 40^1$,
$19^{9} 46^1$,
$19^{9} 52^1$,
$19^{9} 58^1$,
$19^{9} 64^1$,
$19^{9} 70^1$ and
$19^{21} 16^1$.
\end{lemma}

\noindent{\bf Proof}~ The designs are presented in Appendix~\ref{app:4-GDD 19^u m^1}.
\adfQED

%%%%%%%%%%%%%%%%%%%%%%%%%%%%%%%%%%%%%%%%%%%%%%%
%%%%%% 23^u m^1
%%%%%%%%%%%%%%%%%%%%%%%%%%%%%%%%%%%%%%%%%%%%%%%
\begin{lemma}
\label{lem:4-GDD 23^u m^1}
There exist 4-GDDs of types
$23^{12} 14^1$,
$23^{12} 17^1$,
$23^{12} 20^1$,
$23^{15} 17^1$,
$23^{9} 14^1$,
$23^{9} 20^1$,
$23^{9} 26^1$,
$23^{9} 32^1$,
$23^{9} 38^1$,
$23^{9} 44^1$,
$23^{9} 50^1$,
$23^{9} 56^1$,
$23^{9} 62^1$,
$23^{9} 68^1$,
$23^{9} 74^1$,
$23^{9} 80^1$ and
$23^{9} 86^1$.
\end{lemma}

\noindent{\bf Proof}~ The designs are presented in Appendix~\ref{app:4-GDD 23^u m^1}.
\adfQED

%%%%%%%%%%%%%%%%%%%%%%%%%%%%%%%%%%%%%%%%%%%%%%%
%%%%%% 25^u m^1
%%%%%%%%%%%%%%%%%%%%%%%%%%%%%%%%%%%%%%%%%%%%%%%
\begin{lemma}
\label{lem:4-GDD 25^u m^1}
There exist 4-GDDs of types
$25^{12}  7^1$,
$25^{12} 13^1$,
$25^{12} 16^1$,
$25^{12} 19^1$,
$25^{12} 22^1$,
$25^{15} 13^1$,
$25^{15} 19^1$,
$25^{27} 19^1$,
$25^{9} 16^1$,
$25^{9} 22^1$,
$25^{9} 28^1$,
$25^{9} 34^1$,
$25^{9} 46^1$,
$25^{9} 52^1$,
$25^{9} 58^1$,
$25^{9} 64^1$,
$25^{9} 76^1$,
$25^{9} 82^1$,
$25^{9} 88^1$,
$25^{9} 94^1$,
$25^{21} 16^1$ and
$25^{21} 22^1$.
\end{lemma}

\noindent{\bf Proof}~ The designs are presented in Appendix~\ref{app:4-GDD 25^u m^1}.
\adfQED

%%%%%%%%%%%%%%%%%%%%%%%%%%%%%%%%%%%%%%%%%%%%%%%
%%%%%% 29^u m^1
%%%%%%%%%%%%%%%%%%%%%%%%%%%%%%%%%%%%%%%%%%%%%%%
\begin{lemma}
\label{lem:4-GDD 29^u m^1}
There exist 4-GDDs of types
$29^{12} 14^1$,
$29^{12} 17^1$,
$29^{12} 20^1$,
$29^{12} 23^1$,
$29^{12} 26^1$,
$29^{24} 26^1$,
$29^{15} 17^1$,
$29^{15} 23^1$,
$29^{9} 14^1$,
$29^{9} 20^1$,
$29^{9} 26^1$,
$29^{9} 32^1$,
$29^{9} 38^1$,
$29^{9} 44^1$,
$29^{9} 50^1$,
$29^{9} 56^1$,
$29^{9} 62^1$,
$29^{9} 68^1$,
$29^{9} 74^1$,
$29^{9} 80^1$,
$29^{9} 86^1$,
$29^{9} 92^1$,
$29^{9} 98^1$,
$29^{9} 104^1$,
$29^{9} 110^1$ and
$29^{21} 26^1$.
\end{lemma}

\noindent{\bf Proof}~ The designs are presented in Appendix~\ref{app:4-GDD 29^u m^1}.
\adfQED

%%%%%%%%%%%%%%%%%%%%%%%%%%%%%%%%%%%%%%%%%%%%%%%
%%%%%% 31^u m^1
%%%%%%%%%%%%%%%%%%%%%%%%%%%%%%%%%%%%%%%%%%%%%%%
\begin{lemma}
\label{lem:4-GDD 31^u m^1}
There exist 4-GDDs of types
$31^{12} 25^1$,
$31^{12} 28^1$,
$31^{9} 22^1$,
$31^{9} 28^1$,
$31^{9} 34^1$,
$31^{9} 40^1$,
$31^{9} 46^1$,
$31^{9} 52^1$,
$31^{9} 58^1$,
$31^{9} 64^1$,
$31^{9} 70^1$,
$31^{9} 76^1$,
$31^{9} 82^1$,
$31^{9} 88^1$,
$31^{9} 94^1$,
$31^{9} 100^1$,
$31^{9} 106^1$,
$31^{9} 112^1$ and
$31^{9} 118^1$.
\end{lemma}

\noindent{\bf Proof}~ The designs are presented in Appendix~\ref{app:4-GDD 31^u m^1}.
\adfQED

%%%%%%%%%%%%%%%%%%%%%%%%%%%%%%%%%%%%%%%%%%%%%%%
%%%%%% 35^u m^1
%%%%%%%%%%%%%%%%%%%%%%%%%%%%%%%%%%%%%%%%%%%%%%%
\begin{lemma}
\label{lem:4-GDD 35^u m^1}
There exist 4-GDDs of types
$35^{12} 29^1$,
$35^{12} 32^1$,
$35^{9} 26^1$,
$35^{9} 32^1$,
$35^{9} 38^1$,
$35^{9} 44^1$,
$35^{9} 62^1$,
$35^{9} 68^1$,
$35^{9} 74^1$,
$35^{9} 86^1$,
$35^{9} 92^1$,
$35^{9} 104^1$,
$35^{9} 116^1$,
$35^{9} 122^1$,
$35^{9} 128^1$ and
$35^{9} 134^1$.
\end{lemma}

\noindent{\bf Proof}~ The designs are presented in Appendix~\ref{app:4-GDD 35^u m^1}.
\adfQED

%%%%%%%%%%%%%%%%%%%%%%%%%%%%%%%%%%%%%%%%%%%%%%%%%%%%%%%%%%%%%%%%%%%%%%%%%%%%%%%%%%%%%%%%%%
%%%%%%%%%%%%%%%%%%%%%%%%%%%%%%%%%%%%%%%%%%%%%%%%%%%%%%%%%%%%%%%%%%%%%%%%%%%%%%%%%%%%%%%%%%
%%%%%%%%%%%%%%%%%%%%%%%%%%%%%%%%%%%%%%%%%%%%%%%%%%%%%%%%%%%%%%%%%%%%%%%%%%%%%%%%%%%%%%%%%%
%%%%%%%%%%%%%%%%%%%%%%%%%%%%%%%%%%%%%%%%%%%%%%%%%%%%%%%%%%%%%%%%%%%%%%%%%%%%%%%%%%%%%%%%%%

\section{4-GDDs of type $g^u m^1$ with $g \equiv 3 \adfmod{6}$}\label{sec:Theorems g = 3 (mod 6)}
A {\em double group divisible design}, $k$-DGDD, of type
$(g_1, h_1^w)^{u_1} (g_2, h_2^w)^{u_2} \dots (g_r, h_r^w)^{u_r}$, $g_i = w h_i$, $i = 1, 2, \dots, r$, is an ordered quadruple
($V,\mathcal{G},\mathcal{H},\mathcal{B}$) such that:
\begin{enumerate}
\item[(i)]{$V$
is a base set of cardinality $u_1 g_1 + \dots + u_r g_r$;}
\item[(ii)]{$\mathcal{G}$
is a partition of $V$ into $u_i$ subsets of cardinality $g_i$, $i = 1, 2, \dots, r$, the \textit{groups};}
\item[(iii)]{$\mathcal{H}$
is a partition of $V$ into $w$ subsets of cardinality $h_1 u_1 + h_2 u_2 + \dots + h_r u_r$, the \textit{holes},
where each group of size $g_i$ intersects each hole in $h_i$ points, $i = 1, 2, \dots, r$;}
\item[(iv)]{$\mathcal{B}$
is a non-empty collection of subsets of $V$ of cardinality $k$, the \textit{blocks};}
\item[(v)]{%
each pair of elements of $V$ not in the same group and not in the same hole occurs in precisely one block; and}
\item[(vi)]{%
a pair of elements from the same group or from the same hole does not occur in any block.}
\end{enumerate}
We quote the following result concerning the existence of 4-DGDDs.
\begin{theorem}
\label{thm:4-DGDD existence}
There exists a 4-DGDD of type $(h n, h^n)^t$ if and only if $t, n \ge 4$ and
$(t - 1)(n - 1) h \equiv 0 \adfmod{3}$ except for $(h,n,t) \in \{(1,4,6), (1,6,4)\}$.
\end{theorem}
\noindent {\bf Proof~} See \cite{GeWei2004} and \cite{CaoWangWei2009}.
~\adfQED

\adfPgap
\noindent Some of our constructions require 4-GDDs of the form $(3a)^4 (6a)^1 (3b)^1$,
and so we extend a result of Wang \& Shen, \cite{WangShen2008}, with the following theorem.
%%%%%%%%%%%%%%%%%%%%%%%%%%%%%%%%%%%%%%%%%%%%%%%%%%%%%%%%%%%%%%%%%%%%%%%%%%%%%%%%%%%%%%%%%%
%%%%%%%%%%%%%%%%%%%%%%%%%%%%%%%%%%%%%%%%%%%%%%%%%%%%%%%%%%%%%%%%%%%%%%%%%%%%%%%%%%%%%%%%%%
\begin{theorem}
\label{thm:4-GDD a^u b^1 c^1}
Let $a$ and $b$ be positive integers. Then there exists a 4-GDD of type $(3a)^4 (6a)^1 (3b)^1$ if and only if $a \le b \le 2a$.
\end{theorem}
\noindent {\bf Proof~}
The four cases left unresolved by Theorem~2.13 of \cite{WangShen2008}, namely $(a,b) \in$ \{(3,5), (4,7), (6,7), (6,11)\},
are given by Lemma~\ref{lem:4-GDD a^u b^1 c^1}.
~\adfQED

\adfPgap
\noindent We also make use of the following construction.
%
%%%%%%%%%%%%%%%%%%%%%%%%%%%%%%%%%%%%%%%%%%%%%%%%%%%%%%%%%%%%%%%%%%%%%%%%%%%%%%%%%%%%%%%%%%
%%%%%%%%%%%%%%%%%%%%%%%%%%%%%%%%%%%%%%%%%%%%%%%%%%%%%%%%%%%%%%%%%%%%%%%%%%%%%%%%%%%%%%%%%%
\begin{theorem}
\label{thm:4-GDD (va + b)^u, (ct + m)^1}
Let $a$, $b$, $c$, $d$, $t$, $u$, $v$ be non-negative integers such that
$2 \le v \le u - 1$ and $0 \le t \le u - 1$.
Suppose there exists a $(v + 1)$-RGDD of type $u^{v + 1}$.
Suppose also that there exist
4-GDDs of types $a^u d^1$ and $b^u d^1$ as well as
type $a^v b^1 c^1$ if $t c > 0$ and type $a^v b^1$ if $c = 0$ or $t < u - 1$.
Then there exists a 4-GDD of type $(v a + b)^u (c t + d)^1$.
\end{theorem}
\noindent {\bf Proof~}
Take the $(v+1)$-RGDD of type $u^{v+1}$ and note that it has $u$ parallel classes.
Remove all blocks of a single parallel class.
Inflate all points in one group by a factor of $b$ and inflate all other points by $a$.
Adjoin $t c + d$ points.
If $tc > 0$, assign $tc$ points to $t$ of the remaining parallel classes, $c$ points to each, and
overlay each inflated block of these parallel classes plus the $c$ points assigned to it with a 4-GDD of type $a^v b^1 c^1$.
Overlay all other inflated blocks, if any, with 4-GDDs of type $a^v b^1$.
Overlay the inflated groups plus the remaining $d$ points with 4-GDDs of type $a^u d^1$ or $b^u d^1$, as appropriate.
The result is a 4-GDD of type $(v a + b)^u (t c + d)^1$, where
the $u$ groups of size $v a + b$ correspond to the $u$ (inflated) blocks of the removed parallel class.~\adfQED

\adfPgap
\noindent We are now ready to prove the main theorems of this section.

%%%%%%%%%%%%%%%%%%%%%%%%%%%%%%%%%%%%%%%%%%%%%%%%%%%%%%%%%%%%%%%%%%%%%%%%%%%%%%%%%%%%%%%%%%
%%%%%%%%%%%%%%%%%%%%%%%%%%%%%%%%%%%%%%%%%%%%%%%%%%%%%%%%%%%%%%%%%%%%%%%%%%%%%%%%%%%%%%%%%%
\begin{theorem}
\label{thm:g = 3 (mod 6)}
Let $g$ be a positive integer such that $g \equiv 3 \adfmod{6}$.
Then there exists a 4-GDD of type $g^u m^1$ if
\begin{enumerate}
\item[]{$u \equiv 0 \adfmod{4}$, $ u \ge 4$, $m \equiv 0 \adfmod{3}$ and $0 \le m \le g (u - 1)/2$, or}
\item[]{$u \equiv 1 \adfmod{4}$, $ u \ge 5$, $m \equiv 0 \adfmod{6}$ and $0 \le m \le g (u - 1)/2$, or}
\item[]{$u \equiv 3 \adfmod{4}$, $ u \ge 7$, $m \equiv 3 \adfmod{6}$ and $3 \le m \le g (u - 1)/2$,}
\end{enumerate}
except possibly when $g \ge 39$, $g$ is not divisible by any of $9$, $15$, $21$ or $33$, $u = 8$ and
$m' < m < (7g - 3)/2$, where $m' = (17g - 6)/5$ if $g \equiv 3 \adfmod{30}$, $m' = 3g$ otherwise.
\end{theorem}
\noindent{\bf Proof~}
By Theorem~\ref{thm:4-GDD g^u m^1 existence}, the theorem is true for $g \le 33$ and for $g$ divisible by any of 15, 21, 27, 33.
Suppose $g \ge 39$, $g \equiv 3 \adfmod{6}$ is given and assume the theorem holds for all smaller $g \equiv 3 \adfmod{6}$.

Again by Theorem~\ref{thm:4-GDD g^u m^1 existence}, we may assume either
(i) $u \in$ \{7, 11\}, $g$ is not divisible by any of 15, 21, 27 or 33, and $m \le g - 6$, or
(ii) $u = 8$, $g$ is not divisible by any of 9, 15, 21 or 33, and $3g < m < (7g - 3)/2$.

First, we assume $u \in$ \{7, 11\}. Choose $a$, $v$ and $b$, depending on $g$, according to the following scheme:
\begin{center}
\begin{tabular}{l@{~}l@{~~~~~~}l@{~}l}
 $g = 39:$ & $a =  6$,  $v = 5$,  $b =  9$;  &
 $g = 51:$ & $a = 12$,  $v = 4$,  $b =  3$; \\
 $g = 57:$ & $a = 12$,  $v = 4$,  $b =  9$;  &
 $g = 69:$ & $a = 12$,  $v = 5$,  $b =  9$; \\
 $g = 87:$ & $a = 18$,  $v = 4$,  $b = 15$;  &
 $g = 93:$ & $a = 18$,  $v = 4$,  $b = 21$; \\
 $g \ge 111:$ & \multicolumn{3}{l}{$a = 6 \left\lfloor \dfrac{g + 18}{36} \right\rfloor$,  $v = 5$, $b = g - 5a$.}
\end{tabular}
\end{center}
It is easily verified that
\begin{align*}
& g = v a + b,~~~ v \le u - 1, \\
& a \equiv 0 \adfmod{6},~~~ b \equiv 3 \adfmod{6},~~~ 0 < b \le \dfrac{a(v - 1)}{2} < g, \\
& \dfrac{b (u - 1)}{2} \ge a - 3,~~~ \dfrac{a(u - 1)}{2} \ge a - 3,~~~ 7a - 3 \ge g - 6.
\end{align*}
% Checked.
So, by Theorem~\ref{thm:4-GDD g^u m^1 existence} and this theorem, there exist 4-GDDs of types
$a^v b^1 a^1$ and $a^v b^1$ as well as
$a^u d^1$ and $b^u d^1$ for $d = 3, 9, \dots, a - 3$.
Then we use Theorem~\ref{thm:4-GDD (va + b)^u, (ct + m)^1} with $c = a$, $d = 3, 9, \dots, a - 3$ and $t = 0$, 1, \dots, 6, noting that $6 \le u - 1$.
Hence there exists a 4-GDD of type $g^u (t a + d)^1$.
Moreover, $t a + d$ takes values 3, 9, \dots, $7a - 3 \ge g - 6$, as required.

When $u = 8$, by Theorem~\ref{thm:4-GDD g^u m^1 existence}, there is nothing to prove unless $g \ge 93$ and $g \equiv 3 \adfmod{30}$.

So now we assume $g = 30h + 3$, $h \ge 3$ and $u = 8$.
Let $a = 6h$, $b = 6h + 3$, $c = 12h$, $v = 4$ and let $t = 0$, 1, \dots, 7.
By Theorem~\ref{thm:4-GDD a^u b^1 c^1} there exists a 4-GDD of type $a^v b^1 c^1$.
Also, by Theorem~\ref{thm:4-GDD g^u m^1 existence}, there exist 4-GDDs of types
$a^v b^1$ as well as $a^u d^1$ and $b^u d^1$ for $d = 0, 3, \dots, d_\mathrm{max}$, where $d_\mathrm{max} = 3b \ge c$.
Therefore, by Theorem~\ref{thm:4-GDD (va + b)^u, (ct + m)^1}, there exists a 4-GDD of type
$(v a + b)^u (t c + d)^1$, that is, $g^8 (t c + d)^1$.
Moreover, $t c + d$ takes values 0, 3, \dots, $7c + d_\mathrm{max}$, and
$7c + d_\mathrm{max} = 102h + 9 = (17g - 6)/5$, as required.
~\adfQED

It is perhaps worth pointing out that if $g = 30h + 3$, $h \ge 3$, and $b = 6h + 3 = (g + 12)/5$ is such that
there exists a 4-GDD of type $b^8 d^1$ for all admissible $d \le (7b - 21)/2 = 21h$,
then $7c + d_\mathrm{max}$ in the proof of Theorem~\ref{thm:g = 3 (mod 6)} can be increased to $105h = (7g - 21)/2$ and consequently
the number of possible exceptions for $g$ is reduced to just two, namely
$g^8 ((7g - 15)/2)^1$ and $g^8 ((7g - 9)/2)^1$.

%%%%%%%%%%%%%%%%%%%%%%%%%%%%%%%%%%%%%%%%%%%%%%%%%%%%%%%%%%%%%%%%%%%%%%%%%%%%%%%%%%%%%%%%%%
%%%%%%%%%%%%%%%%%%%%%%%%%%%%%%%%%%%%%%%%%%%%%%%%%%%%%%%%%%%%%%%%%%%%%%%%%%%%%%%%%%%%%%%%%%
\begin{theorem}
\label{thm:g = 39r 51r 57r 69r}
Let $g \in \{39, 51, 57, 69, 87, 93\}$ and let $r$ be a positive integer.
Then there exists a 4-GDD of type $(rg)^u m^1$ if
\begin{enumerate}
\item[]{$u \equiv 0 \adfmod{4}$, $ u \ge 4$, $m \equiv 0 \adfmod{3}$ and $0 \le m \le rg (u - 1)/2$, or}
\item[]{$u \equiv 1 \adfmod{4}$, $ u \ge 5$, $m \equiv 0 \adfmod{6}$ and $0 \le m \le rg (u - 1)/2$, or}
\item[]{$u \equiv 3 \adfmod{4}$, $ u \ge 7$, $m \equiv 3 \adfmod{6}$ and $3 \le m \le rg (u - 1)/2$.}
\end{enumerate}
\end{theorem}
\noindent{\bf Proof~}
Let $g \in \{39, 51, 57, 69, 87, 93\}$. We first prove that there exists a 4-GDD of type $g^u m^1$ for all admissible $u$ and $m$.

By Theorems~\ref{thm:4-GDD g^u m^1 existence} and \ref{thm:g = 3 (mod 6)}, we need to consider only 
$u = 8$, $m \equiv 0 \adfmod{3}$, $3g < m < (7g - 3)/2$.
We deal with each value of $g$ in turn.

{\boldmath $g = 39$~}
The 4-GDDs, corresponding to $m \in$ \{120, 123, \dots, 132\}, are given by Lemma~\ref{lem:4-GDD 39^8 m^1}.

{\boldmath $g = 51$~}
If $m = 156$, we use Theorem~\ref{thm:4-GDD (va + b)^u, (ct + m)^1} with
$a = 9$, $v = 4$, $b = 15$, $c = 18$, $d = 30$, $t = 7$ and $u = 8$;
the relevant 4-GDD types are $9^4 15^1 18^1$ from Theorem~\ref{thm:4-GDD a^u b^1 c^1} as well as $9^8 30^1$ and
$15^8 30^1$ from Theorem~\ref{thm:4-GDD g^u m^1 existence}.
The remaining 4-GDDs, corresponding to $m \in$ \{159, 162, \dots, 174\}, are given by Lemma~\ref{lem:4-GDD 51^8 m^1}.

{\boldmath $g = 57$~} The 4-GDDs, corresponding to $m \in$ \{174, 177, \dots, 195\}, are given by Lemma~\ref{lem:4-GDD 57^8 m^1}.
% We can't use 6^8 24^1 9^1 because v = 8 > u - 1.

{\boldmath $g = 69$~}
If $m \in \{210, 213, 216, 219\}$, we use Theorem~\ref{thm:4-GDD (va + b)^u, (ct + m)^1} with
$a = 9$, $v = 6$, $b = 15$, $c = 27$, $d \in \{21, 24, 27, 30\}$, $t = 7$ and $u = 8$;
the relevant 4-GDD types are $9^6 15^1 27^1$ from Lemma~\ref{lem:4-GDD a^u b^1 c^1} as well as
$9^8 d^1$ and $15^8 d^1$ from Theorem~\ref{thm:4-GDD g^u m^1 existence}.
Hence there exists a 4-GDD of type $(v a + b)^u (t c + d)^1$, that is, $69^8 (189 + d)^1$.
The remaining 4-GDDs, corresponding to $m \in$ \{222, 225, \dots, 237\}, are given by Lemma~\ref{lem:4-GDD 69^8 m^1}.
% We can't use 6^10 30^1 9^1 because v = 10 > u - 1.

{\boldmath $g = 87$~}
If $m \in$ \{264, 267, \dots, 294\}, we use Theorem~\ref{thm:4-GDD (va + b)^u, (ct + m)^1} with
$a = 12$, $v = 6$, $b = 15$, $c = 36$, $d \in \{12, 15, \dots, 42\}$, $t = 7$ and $u = 8$;
the relevant 4-GDD types are $12^6 15^1 36^1$  from \cite[Proposition IV.4.15]{Ge2007} as well as
$12^8 d^1$ and $15^8 d^1$ from Theorem~\ref{thm:4-GDD g^u m^1 existence}.
The remaining 4-GDDs, corresponding to $m \in$ \{297, 300\}, are given by Lemma~\ref{lem:4-GDD 87^8 m^1}.

{\boldmath $g = 93$~}
If $m \in$ \{282, 285, \dots, 315\}, we use Theorem~\ref{thm:4-GDD (va + b)^u, (ct + m)^1} with
$a = 18$, $v = 4$, $b = 21$, $c = 36$, $d \in \{30, 33, \dots, 63\}$, $t = 7$ and $u = 8$;
the relevant 4-GDD types are $18^4 21^1 36^1$  from Theorem~\ref{thm:4-GDD a^u b^1 c^1} as well as
$18^8 d^1$ and $21^8 d^1$ from Theorem~\ref{thm:4-GDD g^u m^1 existence}.
The remaining 4-GDDs, corresponding to $m \in$ \{{318, 321\}, are given by Lemma~\ref{lem:4-GDD 87^8 m^1}.

Hence the theorem is proved when $r = 1$.

If $r$ is even, then $rg \equiv 0 \adfmod{6}$ and the result follows from Theorem~\ref{thm:4-GDD g^u m^1 existence}.
If $r$ is odd and $r \ge 5$, the result follows from \cite[Lemma 5.24]{WeiGe2015} (with $v = r$).

Finally, when $r = 3$,
by Theorem~\ref{thm:g = 3 (mod 6)}, there are no possible exceptions since $3g$ is divisible by 9.
~\adfQED

%%%%%%%%%%%%%%%%%%%%%%%%%%%%%%%%%%%%%%%%%%%%%%%%%%%%%%%%%%%%%%%%%%%%%%%%%%%%%%%%%%%%%%%%%%
%%%%%%%%%%%%%%%%%%%%%%%%%%%%%%%%%%%%%%%%%%%%%%%%%%%%%%%%%%%%%%%%%%%%%%%%%%%%%%%%%%%%%%%%%%
%%%%%%%%%%%%%%%%%%%%%%%%%%%%%%%%%%%%%%%%%%%%%%%%%%%%%%%%%%%%%%%%%%%%%%%%%%%%%%%%%%%%%%%%%%
%%%%%%%%%%%%%%%%%%%%%%%%%%%%%%%%%%%%%%%%%%%%%%%%%%%%%%%%%%%%%%%%%%%%%%%%%%%%%%%%%%%%%%%%%%

\section{4-GDDs of type $g^u m^1$ with $g \equiv 1, 5 \adfmod{6}$}\label{sec:Theorems g = 1, 5 (mod 6)}

We prove two theorems concerning 4-GDDs of type $g^u m^1$ where $g \equiv 1, 5 \adfmod{6}$.
Theorem~\ref{thm:g = 1, 5 (mod 6), u = 21, 33} eliminates $m \ge g$ from the possible exceptions for $u \in \{21, 33\}$ and
reduces the number of possible exceptions when $g \equiv 1,11 \adfmod{30}$ and $u = 24$.
Theorem~\ref{thm:g = 13 17 19 23 25} completes the spectrum for 4-GDDs of type $g^u m^1$ when $g \equiv 1, 5 \adfmod{6}$ and $g \le 35$.
Before proving Theorems~\ref{thm:g = 1, 5 (mod 6), u = 21, 33} and \ref{thm:g = 13 17 19 23 25}
we describe four simple constructions the first of which enables us to use the results of Section~\ref{sec:Theorems g = 3 (mod 6)}.

%%%%%%%%%%%%%%%%%%%%%%%%%%%%%%%%%%%%%%%%%%%%%%%%%%%%%%%%%%%%%%%%%%%%%%%%%%%%%%%%%%%%%%%%%%
%%%%%%%%%%%%%%%%%%%%%%%%%%%%%%%%%%%%%%%%%%%%%%%%%%%%%%%%%%%%%%%%%%%%%%%%%%%%%%%%%%%%%%%%%%
\begin{theorem}
\label{thm:4-GDD (3g)^(u/3) (m - g)^1 + g^4 -> g^u m^1}
Suppose there exists a 4-GDD of type $(3g)^{u/3} (m - g)^1$ and $g \neq 2, 6$.
Then there exists a 4-GDD of type $g^u m^1$.
\end{theorem}
\noindent {\bf Proof~}
Take the 4-GDD of type $(3g)^{u/3} (m - g)^1$,
adjoin $g$ points and overlay the groups plus the new points with 4-GDDs of type $g^4$.
~\adfQED

%%%%%%%%%%%%%%%%%%%%%%%%%%%%%%%%%%%%%%%%%%%%%%%%%%%%%%%%%%%%%%%%%%%%%%%%%%%%%%%%%%%%%%%%%%
%%%%%%%%%%%%%%%%%%%%%%%%%%%%%%%%%%%%%%%%%%%%%%%%%%%%%%%%%%%%%%%%%%%%%%%%%%%%%%%%%%%%%%%%%%
\begin{theorem}
\label{thm:4-GDD gi^ti + gi^u m^1}
Suppose there exists a 4-GDD of type $g_1^{t_1} g_2^{t_2} \dots g_r^{t_r}$.
Let $u$ and $m$ be non-negative integers such that $u \ge 4$ and $u \neq 6$.
Suppose there exists a 4-GDD of type $g_i^u m^1$ for $i = 1$, $2$, \dots, $r$.
Then there exists a 4-GDD of type $(g_1 t_1 + g_2 t_2 + \dots + g_r t_r)^u m^1$.
\end{theorem}
\noindent {\bf Proof~}
This is a special case of \cite[Constructions 1.8--1.9]{GeReesZhu2002}.
By Theorem~\ref{thm:4-DGDD existence}, there exists a 4-DGDD of type $(u, 1^u)^4$.
Take the 4-GDD of type $g_1^{t_1} g_2^{t_2} \dots g_r^{t_r}$,
inflate each point by a factor of $u$ and overlay each inflated block with the 4-DGDD
to create a 4-DGDD of type $(u g_1, g_1^u)^{t_1} (u g_2, g_2^u)^{t_2} \dots (u g_r, g_r^u)^{t_r}$.
Adjoin $m$ points and
overlay the groups plus the $m$ points with 4-GDDs of types $g_1^u m^1$, $g_2^u m^1$, \dots, $g_r^u m^1$, as appropriate.
~\adfQED

%%%%%%%%%%%%%%%%%%%%%%%%%%%%%%%%%%%%%%%%%%%%%%%%%%%%%%%%%%%%%%%%%%%%%%%%%%%%%%%%%%%%%%%%%%
%%%%%%%%%%%%%%%%%%%%%%%%%%%%%%%%%%%%%%%%%%%%%%%%%%%%%%%%%%%%%%%%%%%%%%%%%%%%%%%%%%%%%%%%%%
\begin{theorem}
\label{thm:4-DGDD (u, 1^u)^g + 4-GDD 1^u m^1}
Suppose $g, u \ge 4$, $g \equiv 1 \adfmod{3}$ and
$m \le (u - 1)/2$.
Suppose also that
$u \equiv 0 \adfmod{12}$ and $m \equiv 1 \adfmod{3}$, or
$u \equiv 3 \adfmod{12}$ and $m \equiv 1 \adfmod{6}$, or
$u \equiv 9 \adfmod{12}$ and $m \equiv 4 \adfmod{6}$.
Then there exists a 4-GDD of type $g^u m^1$.
\end{theorem}
% Note that this is not a special case of the previous theorem.
\noindent {\bf Proof~}
By Theorem~\ref{thm:4-DGDD existence}, since $\{g, u\} \neq \{4, 6\}$, there exists a 4-DGDD of type $(g, 1^g)^u$, and
by Theorem~\ref{thm:4-GDD g^u m^1 existence}, there exists a 4-GDD of type $1^u m^1$.
To construct the 4-GDD of type $g^u m^1$, take the 4-DGDD,
adjoin $m$ points and overlay the holes plus the $m$ points with 4-GDDs of type $1^u m^1$.
~\adfQED

%%%%%%%%%%%%%%%%%%%%%%%%%%%%%%%%%%%%%%%%%%%%%%%%%%%%%%%%%%%%%%%%%%%%%%%%%%%%%%%%%%%%%%%%%%
%%%%%%%%%%%%%%%%%%%%%%%%%%%%%%%%%%%%%%%%%%%%%%%%%%%%%%%%%%%%%%%%%%%%%%%%%%%%%%%%%%%%%%%%%%
\begin{theorem}
\label{thm:4-GDD (rg)^u (rm)^1}
Suppose there exists a 4-GDD of type $g^u m^1$ and let $r \ge 3$, $r \neq 6$ be an integer.
Then there exists a 4-GDD of type $(rg)^u (rm)^1$.
\end{theorem}
\noindent {\bf Proof~}
Take the 4-GDD of type $g^u m^1$, inflate its points by $r$ and overlay the inflated blocks with 4-GDDs of type $r^4$.
~\adfQED

\adfPgap
We are now ready to prove the main theorems of this section.

%%%%%%%%%%%%%%%%%%%%%%%%%%%%%%%%%%%%%%%%%%%%%%%%%%%%%%%%%%%%%%%%%%%%%%%%%%%%%%%%%%%%%%%%%%
%%%%%%%%%%%%%%%%%%%%%%%%%%%%%%%%%%%%%%%%%%%%%%%%%%%%%%%%%%%%%%%%%%%%%%%%%%%%%%%%%%%%%%%%%%
\begin{theorem}
\label{thm:g = 1, 5 (mod 6), u = 21, 33}
$\mathrm{(i)}$ Suppose $g \equiv 1 \adfmod{6}$ and $m \equiv 4 \adfmod{6}$, or
            $g \equiv 5 \adfmod{6}$ and $m \equiv 2 \adfmod{6}$.
Suppose also that $u \in$ \{$21$, $33$\} and $g \le m \le g (u - 1)/2$.
Then there exists a 4-GDD of type $g^u m^1$.

\adfTgap
$\mathrm{(ii)}$  Suppose $g \equiv 1$ or $11 \adfmod{30}$ and $m \equiv g \adfmod{3}$.
Then there exists a 4-GDD of type $g^{24} m^1$ for $g \le m \le (56g - 6)/5$.
\end{theorem}
\noindent{\bf Proof~}
To prove (i), observe that
$3g \equiv 3 \adfmod{6}$, $u/3 \in$ \{7, 11\}, $m - g \equiv 3 \adfmod{6}$ and $m - g \le 3g(u/3 - 1)/2$.
So, by Theorem~\ref{thm:g = 3 (mod 6)}, there exists a 4-GDD of type $(3g)^{u/3} (m - g)^1$.
Hence, by Theorem~\ref{thm:4-GDD (3g)^(u/3) (m - g)^1 + g^4 -> g^u m^1}, there exists a 4-GDD of type $g^u m^1$.

\adfTgap
For (ii), we have
$3g \equiv 3 \adfmod{30}$, $m - g \equiv 0 \adfmod{3}$ and $m - g \le (51g - 6)/5$.
Hence, by Theorem~\ref{thm:g = 3 (mod 6)}, there exists a 4-GDD of type $(3g)^{8} (m - g)^1$, and
it follows from Theorem~\ref{thm:4-GDD (3g)^(u/3) (m - g)^1 + g^4 -> g^u m^1} that there exists a 4-GDD of type $g^{24} m^1$.
~\adfQED

%%%%%%%%%%%%%%%%%%%%%%%%%%%%%%%%%%%%%%%%%%%%%%%%%%%%%%%%%%%%%%%%%%%%%%%%%%%%%%%%%%%%%%%%%%
%%%%%%%%%%%%%%%%%%%%%%%%%%%%%%%%%%%%%%%%%%%%%%%%%%%%%%%%%%%%%%%%%%%%%%%%%%%%%%%%%%%%%%%%%%

\begin{theorem}
\label{thm:g = 13 17 19 23 25}
Suppose $g \in \{13, 19, 25, 31\}$ and either
\begin{enumerate}
\item[]{$u \equiv 0 \adfmod{12}$, $u \ge 12$, $m \equiv 1 \adfmod{3}$ and $1 \le m \le g (u - 1)/2$, or}
\item[]{$u \equiv 3 \adfmod{12}$, $u \ge 15$, $m \equiv 1 \adfmod{6}$ and $1 \le m \le g (u - 1)/2$, or}
\item[]{$u \equiv 9 \adfmod{12}$, $u \ge 9$,  $m \equiv 4 \adfmod{6}$ and $4 \le m \le g (u - 1)/2$.}
\end{enumerate}
Then there exists a 4-GDD of type $g^u m^1$.

\adfTgap
Suppose $g \in \{17, 23, 29, 35\}$ and either
\begin{enumerate}
\item[]{$u \equiv 0 \adfmod{12}$, $u \ge 12$, $m \equiv 2 \adfmod{3}$ and $2 \le m \le g (u - 1)/2$, or}
\item[]{$u \equiv 3 \adfmod{12}$, $u \ge 15$, $m \equiv 5 \adfmod{6}$ and $5 \le m \le g (u - 1)/2$, or}
\item[]{$u \equiv 9 \adfmod{12}$, $u \ge 9$,  $m \equiv 2 \adfmod{6}$ and $2 \le m \le g (u - 1)/2$.}
\end{enumerate}
Then there exists a 4-GDD of type $g^u m^1$.
\end{theorem}
\noindent{\bf Proof~}
Let $g \in \{13, 17, 19, 23, 25, 29, 31, 35\}$.
By Theorems~\ref{thm:4-GDD g^u m^1 existence} and \ref{thm:g = 1, 5 (mod 6), u = 21, 33}, and \cite{GeReesZhu2002} (or see \cite[Theorem IV.4.11]{Ge2007}),
we need to consider only:
\begin{enumerate}
\item[]{$u \in$ \{12, 15, 21, 24, 27, 33, 39, 51\} and $0 < m < g$;}
\item[]{$u = 9$ and $0 < m < 4g$;}
\item[]{$u = 24$ and $10g < m < (23g - 3)/2$.}
\end{enumerate}
Moreover, the necessary conditions imply the following:
\begin{enumerate}
\item[]{for $g \in \{13, 19, 25, 31\}$,
$m \equiv 1 \adfmod{3}$ when $ u \in \{12, 24\}$,
$m \equiv 1 \adfmod{6}$ when $ u \in \{15, 27, 39, 51\}$, and
$m \equiv 4 \adfmod{6}$ when $ u \in \{9, 21, 33\}$;}
\item[]{for $g \in \{17, 23, 29, 35\}$,
$m \equiv 2 \adfmod{3}$ when $ u \in \{12, 24\}$,
$m \equiv 5 \adfmod{6}$ when $ u \in \{15, 27, 39, 51\}$, and
$m \equiv 2 \adfmod{6}$ when $ u \in \{9, 21, 33\}$.}
\end{enumerate}
We address each $g$ separately for $u \in \{9, 12, 15, 21, 24, 27\}$.
Then we use the results to deal with $u \in \{33, 39, 51\}$.

% There exists a k=RGDD u^k for (k,u) = (9,9), (6,12), (5,15), (6,21), (8,24), (27,27).

\adfTgap
\noindent {\bf\boldmath $g = 13$~}
For type $13^u m^1$, where

\adfInd
$u = 12$, $m \in \{1, 4\}$, or
$u = 24$, $m \in \{1, 4, 7, 10\}$, or

\adfInd
$u = 15$, $m \in \{1, 7\}$, or
$u = 27$, $m \in \{1, 7\}$, or

\adfInd
$u = 9$, $m = 4$, or
$u = 21$, $m \in \{4, 10\}$,

\noindent use Theorem~\ref{thm:4-DGDD (u, 1^u)^g + 4-GDD 1^u m^1}.

For type $13^{12} m^1$, $m \in \{7, 10\}$,
see Lemma~\ref{lem:4-GDD 13^u m^1}.

For type $13^{24} m^1$, $m \in \{133, 136, 139, 142, 145\}$,
use Theorem~\ref{thm:4-GDD (3g)^(u/3) (m - g)^1 + g^4 -> g^u m^1}
with a 4-GDD of type $39^8 (m - 13)^1$ from Theorem~\ref{thm:g = 39r 51r 57r 69r}.

For type $13^{9} m^1$, $m \in \{10, 16, 22, 28, 34, 40, 46\}$,
see Lemma~\ref{lem:4-GDD 13^u m^1}.

%%%%%%%%%%%%%%%%%%%%%%%%%%%%%%%%%%%%%%%%%%%%%%%%%%%%%%%%%%%%%%%%%%%%%%%%%%%%%%%%%%%%%%%%%%
\adfTgap
\noindent {\bf\boldmath $g = 17$~}
For type $17^{a} 5^1$, $a \in \{15, 27\}$ and $17^{b} 2^1$, $b \in \{9, 21, 24\}$,
see \cite{GeLing2005} or \cite[Theorem IV.4.11]{Ge2007}.

For type $17^{12} m^1$, $m \in \{2, 5, 8, 11, 14\}$,
see Lemma~\ref{lem:4-GDD 17^u m^1}.
Incidentally, $17^{12} 2^1$ fills the last remaining gap in \cite[Theorem IV.4.11]{Ge2007}.

For type $17^{u} m^1$, $u \in \{21, 24, 27\}$, $m < 17$,
we argue similarly to \cite[Lemma 6.2]{WeiGe2015}.
For $n \in \{4, 5, 6\}$, there exists a 4-DGDD, $D_n$  say, of type $(2n, 2^n)^6 (5n, 5^n)^1$, \cite[Lemma 2.14]{GeLing2005}.
Take a \{4,5,6\}-GDD of type $1^u$, \cite[Table IV.3.23]{AbelBennettGreig2007},
inflate each point by a factor of 17 and overlay each inflated block of size $n$ with a copy of $D_n$,
aligning the holes of $D_n$ with the inflated points in a consistent manner.
The result is a 4-DGDD of type $(2u, 2^u)^6 (5u, 5^u)^1$.
Now adjoin $m$ new points and overlay the groups plus the new points with 4-GDDs of types $2^u m^1$ and $5^u m^1$, as appropriate.
Thus we create 4-GDDs of types
$17^{24}d^1$, $d \in \{5, 8, 11, 14\}$, $17^{27} 11^1$, $17^{21} 8^1$ and $17^{21} 14^1$.

For type $17^{24} m^1$, $m \in \{173, 176, \dots, 191\}$,
use Theorem~\ref{thm:4-GDD (3g)^(u/3) (m - g)^1 + g^4 -> g^u m^1}
with a 4-GDD of type $51^8 (m - 17)^1$ from Theorem~\ref{thm:g = 39r 51r 57r 69r}.

For type $17^{15} 11^1$,
see Lemma~\ref{lem:4-GDD 17^u m^1}.

For type $17^{9} m^1$, $m \in \{8, 14, \dots, 62\}$,
see Lemma~\ref{lem:4-GDD 17^u m^1}.

%%%%%%%%%%%%%%%%%%%%%%%%%%%%%%%%%%%%%%%%%%%%%%%%%%%%%%%%%%%%%%%%%%%%%%%%%%%%%%%%%%%%%%%%%%
\adfTgap
\noindent {\bf\boldmath $g = 19$~}
For type $19^u m^1$, where

\adfInd
$u = 12$, $m \in \{1, 4\}$, or
$u = 24$, $m \in \{1, 4, 7, 10\}$, or

\adfInd
$u = 15$, $m \in \{1, 7\}$, or
$u = 27$, $m \in \{1, 7, 13\}$, or

\adfInd
$u = 9$, $m = 4$, or
$u = 21$, $m \in \{4, 10\}$,

\noindent use Theorem~\ref{thm:4-DGDD (u, 1^u)^g + 4-GDD 1^u m^1}.

For type $19^{12} m^1$, $m \in \{7, 10, 13, 16\}$,
see Lemma~\ref{lem:4-GDD 19^u m^1}.

For type $19^{24} m^1$, $m \in \{13, 16\}$,
see Lemma~\ref{lem:4-GDD 19^u m^1}.

For type $19^{24} m^1$, $m \in \{193, 196, \dots, 214\}$,
use Theorem~\ref{thm:4-GDD (3g)^(u/3) (m - g)^1 + g^4 -> g^u m^1}
with a 4-GDD of type $57^8 (m - 19)^1$ from Theorem~\ref{thm:g = 39r 51r 57r 69r}.

For type $19^{15} 13^1$,
see Lemma~\ref{lem:4-GDD 19^u m^1}.

For type $19^{9} m^1$, $m \in \{10, 16, \dots, 70\}$,
see Lemma~\ref{lem:4-GDD 19^u m^1}.

For type $19^{21} 16^1$,
see Lemma~\ref{lem:4-GDD 19^u m^1}.

%%%%%%%%%%%%%%%%%%%%%%%%%%%%%%%%%%%%%%%%%%%%%%%%%%%%%%%%%%%%%%%%%%%%%%%%%%%%%%%%%%%%%%%%%%
\adfTgap
\noindent {\bf\boldmath $g = 23$~}
For type $23^u m^1$, where

\adfInd
$u = 12$, $m \in \{2, 5, 8, 11\}$, or
$u = 24$, $m \in \{2, 5, 8, 11, 14, 17, 20\}$, or

\adfInd
$u = 15$, $m \in \{5, 11\}$, or
$u = 27$, $m \in \{5, 11, 17\}$, or

\adfInd
$u = 9$, $m \in \{2, 8\}$, or
$u = 21$, $m \in \{2, 8, 14, 20\}$,

\noindent use Theorem~\ref{thm:4-GDD gi^ti + gi^u m^1} with 4-GDDs of types $2^9 5^1$, $2^u m^1$ and $5^u m^1$ from Theorem~\ref{thm:4-GDD g^u m^1 existence}.

For type $23^{12} m^1$, $m \in \{14, 17, 20\}$,
see Lemma~\ref{lem:4-GDD 23^u m^1}.

For type $23^{24} m^1$, $m \in \{233, 236, \dots, 260\}$,
use Theorem~\ref{thm:4-GDD (3g)^(u/3) (m - g)^1 + g^4 -> g^u m^1}
with a 4-GDD of type $69^8 (m - 23)^1$ from Theorem~\ref{thm:g = 39r 51r 57r 69r}.

For type $23^{15} 17^1$,
see Lemma~\ref{lem:4-GDD 23^u m^1}.

For type $23^{9} m^1$, $m \in \{14, 20, \dots, 86\}$,
see Lemma \ref{lem:4-GDD 23^u m^1}.

%%%%%%%%%%%%%%%%%%%%%%%%%%%%%%%%%%%%%%%%%%%%%%%%%%%%%%%%%%%%%%%%%%%%%%%%%%%%%%%%%%%%%%%%%%
\adfTgap
\noindent {\bf\boldmath $g = 25$~}
For type $25^u m^1$, where

\adfInd
$u = 12$, $m \in \{1, 4\}$, or

\adfInd
$u = 15$, $m \in \{1, 7\}$, or
$u = 27$, $m \in \{1, 7, 13\}$, or

\adfInd
$u = 9$, $m = 4$, or
$u = 21$, $m \in \{4, 10\}$,

\noindent use Theorem~\ref{thm:4-DGDD (u, 1^u)^g + 4-GDD 1^u m^1}.

For type $25^{12} 10^1$,
use Theorem~\ref{thm:4-GDD (rg)^u (rm)^1} with a 4-GDD of type $5^{12} 2^1$ (Theorem~\ref{thm:4-GDD g^u m^1 existence}).

For type $25^{12} m^1$, $m \in \{7, 13, 16, 19, 22\}$,
see Lemma~\ref{lem:4-GDD 25^u m^1}.

For type $25^{24} m^1$, $m \in \{1, 4, 7, 10, 13, 16, 19, 22\}$,
take a 4-DGDD of type $(60, 10^{6})^{10}$ (Theorem~\ref{thm:4-DGDD existence}),
adjoin $m$ points and overlay the groups plus the new points with 4-GDDs of type $10^{6} m^1$ (Theorem~\ref{thm:4-GDD g^u m^1 existence})
to construct a 4-GDD of type $100^6 m^1$.
Then overlay the groups with 4-GDDs of type $25^4$.

For type $25^{24} m^1$, $m \in \{253, 256, \dots, 283\}$,
use Theorem~\ref{thm:4-GDD (3g)^(u/3) (m - g)^1 + g^4 -> g^u m^1}
with a 4-GDD of type $75^8 (m - 25)^1$ from Theorem~\ref{thm:4-GDD g^u m^1 existence}.

For types $25^{15} 13^1$, $25^{15} 19^1$ and $25^{27} 19^1$,
see Lemma~\ref{lem:4-GDD 25^u m^1}.

For type $25^{9} m^1$, $m \in \{10, 40, 70\}$,
use Theorem~\ref{thm:4-GDD (rg)^u (rm)^1} with a 4-GDD of type $5^9 (m/5)^1$ (Theorem~\ref{thm:4-GDD g^u m^1 existence}).

For type $25^{9} m^1$, $m \in \{16, 22, \dots, 94\} \setminus \{10, 40, 70\}$,
see Lemma~\ref{lem:4-GDD 25^u m^1}.

For types $25^{21} 16^1$ and $25^{21} 22^1$,
see Lemma~\ref{lem:4-GDD 25^u m^1}.

%%%%%%%%%%%%%%%%%%%%%%%%%%%%%%%%%%%%%%%%%%%%%%%%%%%%%%%%%%%%%%%%%%%%%%%%%%%%%%%%%%%%%%%%%%
\adfTgap
\noindent {\bf\boldmath $g = 29$~}
For type $29^u m^1$, where

\adfInd
$u = 12$, $m \in \{2, 5, 8, 11\}$, or
$u = 24$, $m \in \{2, 5, 8, 11, 14, 17, 20, 23\}$, or

\adfInd
$u = 15$, $m \in \{5, 11\}$, or
$u = 27$, $m \in \{5, 11, 17, 23\}$, or

\adfInd
$u = 9$, $m \in \{2, 8\}$, or
$u = 21$, $m \in \{2, 8, 14, 20\}$,

\noindent use Theorem~\ref{thm:4-GDD gi^ti + gi^u m^1} with 4-GDDs of types $2^{12} 5^1$, $2^u m^1$ and $5^u m^1$ from Theorem~\ref{thm:4-GDD g^u m^1 existence}.

For type $29^{12} m^1$, $m \in \{14, 17, 20, 23, 26\}$,
see Lemma~\ref{lem:4-GDD 29^u m^1}.

For type $29^{24} 26^1$,
see Lemma~\ref{lem:4-GDD 29^u m^1}.

For type $29^{24} m^1$, $m \in \{293, 296, \dots, 329\}$,
use Theorem~\ref{thm:4-GDD (3g)^(u/3) (m - g)^1 + g^4 -> g^u m^1}
with a 4-GDD of type $87^8 (m - 29)^1$ from Theorem~\ref{thm:g = 39r 51r 57r 69r}.

For type $29^{15} m^1$, $m \in \{17, 23\}$,
see Lemma~\ref{lem:4-GDD 29^u m^1}.

For type $29^{9} m^1$, $m \in \{14, 20, \dots, 110\}$,
see Lemma \ref{lem:4-GDD 29^u m^1}.

For type $29^{21} 26^1$,
see Lemma~\ref{lem:4-GDD 29^u m^1}.

%%%%%%%%%%%%%%%%%%%%%%%%%%%%%%%%%%%%%%%%%%%%%%%%%%%%%%%%%%%%%%%%%%%%%%%%%%%%%%%%%%%%%%%%%%
\adfTgap
\noindent {\bf\boldmath $g = 31$~}
For type $31^u m^1$, where

\adfInd
$u = 12$, $m \in \{1, 4, 7, 10, 13, 16, 19, 22\}$, or

\adfInd
$u = 24$, $m \in \{1, 4, 7, 10, 13, 16, 19, 22, 25, 28\}$, or

\adfInd
$u = 15$, $m \in \{1, 7, 13, 19, 25\}$, or
$u = 27$, $m \in \{1, 7, 13, 19, 25\}$, or

\adfInd
$u = 9$, $m \in \{4, 10, 16\}$, or
$u = 21$, $m \in \{4, 10, 16, 16, 22, 28\}$,

\noindent use Theorem~\ref{thm:4-GDD gi^ti + gi^u m^1} with 4-GDDs of types $4^6 7^1$, $4^u m^1$ and $7^u m^1$ from Theorem~\ref{thm:4-GDD g^u m^1 existence}.

For type $31^{12} m^1$, $m \in \{25, 28\}$,
see Lemma~\ref{lem:4-GDD 31^u m^1}.

For type $31^{24} m^1$, $m \in \{313, 316, \dots, 352\}$,
use Theorem~\ref{thm:4-GDD (3g)^(u/3) (m - g)^1 + g^4 -> g^u m^1}
with a 4-GDD of type $93^8 (m - 31)^1$ from Theorem~\ref{thm:g = 39r 51r 57r 69r}.

For type $31^{9} m^1$, $m \in \{22, 28, \dots, 118\}$,
see Lemma~\ref{lem:4-GDD 31^u m^1}.

%%%%%%%%%%%%%%%%%%%%%%%%%%%%%%%%%%%%%%%%%%%%%%%%%%%%%%%%%%%%%%%%%%%%%%%%%%%%%%%%%%%%%%%%%%
\adfTgap
\noindent {\bf\boldmath $g = 35$~}
For type $35^u m^1$, where

\adfInd
$u = 12$, $m \in \{2, 5, \dots, 26\}$, or
$u = 24$, $m \in \{2, 5, \dots, 32\}$, or

\adfInd
$u \in \{15, 27\}$, $m \in \{5, 11, 17, 23, 29\}$, or

\adfInd
$u = 9$, $m \in \{2, 8, 14, 20\}$, or
$u = 21$, $m \in \{2, 8, 14, 20, 26, 32\}$,

\noindent take a 4-DGDD of type $(35, 5^7)^{u}$, adjoin $m$ points, and overlay each hole together with the new points with a 4-GDD of type $5^{u} m^1$.

For type $35^{12} m^1$, $m \in \{29, 32\}$,
see Lemma~\ref{lem:4-GDD 35^u m^1}.

For type $35^{24} m^1$, $m \in \{353, 356, \dots, 398\}$,
use Theorem~\ref{thm:4-GDD (3g)^(u/3) (m - g)^1 + g^4 -> g^u m^1}
with a 4-GDD of type $105^8 (m - 35)^1$ from Theorem~\ref{thm:4-GDD g^u m^1 existence}.

For type $35^{9} m^1$, $m \in \{50, 80, 110\}$,
use Theorem~\ref{thm:4-GDD (rg)^u (rm)^1} with a 4-GDD of type $7^9 (m/5)^1$ from Theorem~\ref{thm:4-GDD g^u m^1 existence}.

For type $35^{9} m^1$, $m \in \{56, 98\}$,
use Theorem~\ref{thm:4-GDD (rg)^u (rm)^1} with a 4-GDD of type $5^9 (m/7)^1$ from Theorem~\ref{thm:4-GDD g^u m^1 existence}.

For type $35^{9} m^1$, $m \in \{26, 32, \dots, 134\} \setminus \{50, 56, 80, 98, 110\}$,
see Lemma~\ref{lem:4-GDD 35^u m^1}.

%%%%%%%%%%%%%%%%%%%%%%%%%%%%%%%%%%%%%%%%%%%%%%%%%%%%%%%%%%%%%%%%%%%%%%%%%%%%%%%%%%%%%%%%%%
\adfTgap
\noindent {\bf\boldmath $u = 33, 39, 51$~}
For type $g^{u} m^1$, $u \in \{39, 51\}$, $m < g$,
take a 4-GDD of type $g^{u - 12} (12g + m)^1$ and overlay the big group with a 4-GDD of type $g^{12} m^1$.

For type $g^{33} m^1$, $m < g$,
take a 4-GDD of type $g^{24} (9g + m)^1$ and overlay the big group with a 4-GDD of type $g^{9} m^1$.
~\adfQED

%%%%%%%%%%%%%%%%%%%%%%%%%%%%%%%%%%%%%%%%%%%%%%%%%%%%%%%%%%%%%%%%%%%%%%%%%%%%%%%%%%%%%%%%%%
%%%%%%%%%%%%%%%%%%%%%%%%%%%%%%%%%%%%%%%%%%%%%%%%%%%%%%%%%%%%%%%%%%%%%%%%%%%%%%%%%%%%%%%%%%
%%%%%%%%%%%%%%%%%%%%%%%%%%%%%%%%%%%%%%%%%%%%%%%%%%%%%%%%%%%%%%%%%%%%%%%%%%%%%%%%%%%%%%%%%%
%%%%%%%%%%%%%%%%%%%%%%%%%%%%%%%%%%%%%%%%%%%%%%%%%%%%%%%%%%%%%%%%%%%%%%%%%%%%%%%%%%%%%%%%%%

%%%%%%%%%%%%%%%%%%%%%%%%%%%%%%%%%%%%%%%%%%%%%%%%%%%%%%%%%%%%%%%%%%%%%%%%%%%%%%%%%%%%%%%%%%
%%%%%%%%%%%%%%%%%%%%%%%%%%%%%%%%%%%%%%%%%%%%%%%%%%%%%%%%%%%%%%%%%%%%%%%%%%%%%%%%%%%%%%%%%%
%%%%%%%%%%%%%%%%%%%%%%%%%%%%%%%%%%%%%%%%%%%%%%%%%%%%%%%%%%%%%%%%%%%%%%%%%%%%%%%%%%%%%%%%%%
%%%%%%%%%%%%%%%%%%%%%%%%%%%%%%%%%%%%%%%%%%%%%%%%%%%%%%%%%%%%%%%%%%%%%%%%%%%%%%%%%%%%%%%%%%

\newpage
\appendix
%
% THE DESIGN VERIFICATION MECHANISM.
%
% %XXXvfySectionStart  Start of a section (XXX represents ADF). The first item {x}
%                      is the title of the section.
%
% %XXXvfyBlocksStart   Start of the blocks of the designs. The first item {x} defines the
%                      graph that is being decomposed: K_{x}.
%                      A design begins with \adfLgap or %\adfLgap.
%
% %XXXvfyDesignStart:  Start of the blocks of a single design.
%
% %XXXvfyBlocksEnd     End of the blocks of the designs.
%
% \XXXvfyParStart      Start of the design parameters:
%                      (v, (..., (b_i, o_i, (..., (v_ij, s_ij), ...)), ...), (..., (g_k, u_k), ...)).
%
% %XXXvfyParEnd        End of the design parameters.
%
% %XXXvfySectionEnd    End of a section
%
% When activated, this command prints the parameters of the design in angle brackets.
\newcommand{\ADFvfyParStart}[1]{{\par\noindent#1}}
% \newcommand{\ADFvfyParStart}[1]{}
%
% Design spacings.
\newcommand{\adfDgap}{\vskip 1.75mm}      % Start of a design
\newcommand{\adfLgap}{\vskip 0.75mm}      % Start of the actual blocks of a design
\newcommand{\adfsplit}{\par}              % Split the list of blocks at these points

% FULL VERSION:  include the whole of GDD4-1-3-5-mod-6-TeX-gen-A.TEX
% SHORT VERSION: include only Sections A and B of GDD4-1-3-5-mod-6-TeX-gen-A.TEX

% \input{GDD4-1-3-5-mod-6-TeX-gen-A.TEX}

%ADFvfySectionStart {4-GDDs}

%%%%%%%%%%%%%%%%%%%%%%%%%%%%%%%%%%%%%%%%%%%%%%%%%%%%%%%%%%%%%%%%%%%%%%%%%%%%%%%%%%%%%%%%%%
%%%%%%%%%%%%%%%%%%%%%%%%%%%%%%%%%%%%%%%%%%%%%%%%%%%%%%%%%%%%%%%%%%%%%%%%%%%%%%%%%%%%%%%%%%
\section{4-GDDs for the proof of Lemma \ref{lem:4-GDD a^u b^1 c^1}}
\label{app:4-GDD a^u b^1 c^1}
\adfnull{
$ 9^4 18^1 15^1 $,
$ 12^4 24^1 21^1 $,
$ 18^4 36^1 21^1 $,
$ 18^4 36^1 33^1 $ and
$ 9^6 27^1 15^1 $.
}

% Charlotte:GDD4-1-3-5-mod-6-TeX-gen-A:HITS-fun:4.10
\adfDgap
%ADFvfyBlocksStart {9,9,9,9,18,15}
\noindent{\boldmath $ 9^{4} 18^{1} 15^{1} $}~
With the point set $Z_{69}$ partitioned into
 residue classes modulo $4$ for $\{0, 1, \dots, 35\}$,
 $\{36, 37, \dots, 53\}$, and
 $\{54, 55, \dots, 68\}$,
 the design is generated from

\adfLgap %ADFvfyDesignStart
$(0, 1, 36, 61)$,
$(0, 2, 38, 56)$,
$(0, 3, 40, 66)$,
$(1, 2, 37, 66)$,\adfsplit
$(1, 3, 41, 56)$,
$(2, 3, 39, 61)$,
$(0, 5, 42, 65)$,
$(0, 6, 44, 60)$,\adfsplit
$(0, 7, 46, 55)$,
$(1, 6, 43, 55)$,
$(1, 7, 49, 65)$,
$(2, 19, 43, 60)$,\adfsplit
$(0, 9, 19, 51)$,
$(0, 10, 17, 53)$,
$(0, 13, 39, 54)$,
$(0, 11, 41, 64)$,\adfsplit
$(0, 18, 37, 59)$,
$(0, 35, 43, 57)$,
$(0, 22, 45, 67)$,
$(1, 19, 39, 57)$,\adfsplit
$(1, 18, 53, 58)$,
$(0, 31, 47, 63)$,
$(0, 34, 49, 58)$,
$(0, 15, 21, 30)$,\adfsplit
$(0, 33, 48, 62)$,
$(0, 23, 26, 52)$,
$(0, 25, 50, 68)$,
$(0, 14, 27, 29)$,\adfsplit
$(2, 11, 52, 67)$,
$(1, 14, 52, 62)$,
$(2, 7, 48, 64)$,
$(1, 26, 48, 59)$,\adfsplit
$(1, 23, 34, 42)$,
$(2, 31, 50, 58)$,
$(1, 27, 46, 64)$,
$(1, 15, 44, 63)$

%ADFvfyBlocksEnd
\adfLgap \noindent by the mapping:
$x \mapsto x + 4 j \adfmod{36}$ for $x < 36$,
$x \mapsto (x + 2 j \adfmod{18}) + 36$ for $36 \le x < 54$,
$x \mapsto (x - 54 + 5 j \adfmod{15}) + 54$ for $x \ge 54$,
$0 \le j < 9$.
\ADFvfyParStart{(69, ((36, 9, ((36, 4), (18, 2), (15, 5)))), ((9, 4), (18, 1), (15, 1)))} %ADFvfyParEnd
% End of 9^4 18^1 15^1
%%%%%%%%%%%%%%%%%%%%%%%%%%%%%%%%%%%%%%%%%%%%%%%%%%%%%%%%%%%%%%%%%%%%%%%%%%%%%%%%%%%%%%%%%%
%%%%%%%%%%%%%%%%%%%%%%%%%%%%%%%%%%%%%%%%%%%%%%%%%%%%%%%%%%%%%%%%%%%%%%%%%%%%%%%%%%%%%%%%%%

% Charlotte:GDD4-1-3-5-mod-6-TeX-gen-A:HITS-fun:4.10
\adfDgap
%ADFvfyBlocksStart {12,12,12,12,24,21}
\noindent{\boldmath $ 12^{4} 24^{1} 21^{1} $}~
With the point set $Z_{93}$ partitioned into
 residue classes modulo $4$ for $\{0, 1, \dots, 47\}$,
 $\{48, 49, \dots, 71\}$, and
 $\{72, 73, \dots, 92\}$,
 the design is generated from

\adfLgap %ADFvfyDesignStart
$(70, 0, 41, 78)$,
$(58, 32, 7, 34)$,
$(60, 14, 76, 24)$,
$(16, 13, 75, 70)$,\adfsplit
$(64, 9, 74, 28)$,
$(14, 83, 29, 67)$,
$(48, 79, 47, 18)$,
$(25, 16, 80, 49)$,\adfsplit
$(53, 19, 8, 9)$,
$(2, 67, 81, 40)$,
$(12, 27, 83, 69)$,
$(66, 30, 4, 17)$,\adfsplit
$(8, 91, 27, 60)$,
$(17, 53, 80, 31)$,
$(17, 60, 74, 26)$,
$(0, 3, 54, 90)$,\adfsplit
$(1, 2, 55, 90)$,
$(1, 7, 48, 83)$,
$(0, 27, 52, 81)$,
$(1, 3, 51, 74)$,\adfsplit
$(2, 3, 53, 88)$,
$(0, 6, 58, 87)$,
$(1, 19, 66, 87)$,
$(0, 18, 57, 80)$,\adfsplit
$(2, 15, 50, 87)$,
$(1, 6, 62, 72)$,
$(0, 7, 68, 79)$,
$(0, 34, 60, 75)$,\adfsplit
$(1, 35, 54, 89)$,
$(2, 39, 70, 92)$,
$(1, 22, 68, 84)$,
$(1, 18, 70, 85)$,\adfsplit
$(1, 43, 56, 91)$,
$(0, 17, 42, 47)$,
$(2, 27, 59, 89)$,
$(1, 23, 38, 65)$,\adfsplit
$(0, 22, 61, 72)$,
$(0, 37, 51, 86)$,
$(1, 27, 71, 79)$,
$(0, 25, 55, 82)$,\adfsplit
$(0, 21, 69, 91)$,
$(0, 39, 67, 85)$,
$(0, 14, 31, 33)$,
$(0, 5, 43, 46)$,\adfsplit
$(2, 11, 69, 75)$,
$(1, 46, 63, 78)$,
$(0, 30, 53, 84)$,
$(2, 43, 57, 84)$,\adfsplit
$(0, 35, 59, 92)$

%ADFvfyBlocksEnd
\adfLgap \noindent by the mapping:
$x \mapsto x + 4 j \adfmod{48}$ for $x < 48$,
$x \mapsto (x + 2 j \adfmod{24}) + 48$ for $48 \le x < 72$,
$x \mapsto (x - 72 + 7 j \adfmod{21}) + 72$ for $x \ge 72$,
$0 \le j < 12$.
\ADFvfyParStart{(93, ((49, 12, ((48, 4), (24, 2), (21, 7)))), ((12, 4), (24, 1), (21, 1)))} %ADFvfyParEnd
% End of 12^4 24^1 21^1
%%%%%%%%%%%%%%%%%%%%%%%%%%%%%%%%%%%%%%%%%%%%%%%%%%%%%%%%%%%%%%%%%%%%%%%%%%%%%%%%%%%%%%%%%%
%%%%%%%%%%%%%%%%%%%%%%%%%%%%%%%%%%%%%%%%%%%%%%%%%%%%%%%%%%%%%%%%%%%%%%%%%%%%%%%%%%%%%%%%%%

% Charlotte:GDD4-1-3-5-mod-6-TeX-gen-A:HITS-fun:4.10
\adfDgap
%ADFvfyBlocksStart {18,18,18,18,36,21}
\noindent{\boldmath $ 18^{4} 36^{1} 21^{1} $}~
With the point set $Z_{129}$ partitioned into
 residue classes modulo $4$ for $\{0, 1, \dots, 71\}$,
 $\{72, 73, \dots, 107\}$, and
 $\{108, 109, \dots, 128\}$,
 the design is generated from

\adfLgap %ADFvfyDesignStart
$(108, 106, 50, 5)$,
$(108, 78, 15, 30)$,
$(108, 101, 56, 13)$,
$(108, 80, 21, 40)$,\adfsplit
$(108, 73, 71, 48)$,
$(108, 81, 58, 55)$,
$(109, 83, 44, 21)$,
$(109, 81, 60, 22)$,\adfsplit
$(109, 84, 41, 66)$,
$(109, 85, 7, 1)$,
$(109, 86, 3, 14)$,
$(109, 106, 64, 11)$,\adfsplit
$(110, 82, 4, 31)$,
$(110, 89, 26, 59)$,
$(110, 72, 22, 25)$,
$(110, 105, 66, 53)$,\adfsplit
$(110, 92, 12, 21)$,
$(110, 79, 3, 20)$,
$(111, 91, 22, 24)$,
$(111, 75, 11, 42)$,\adfsplit
$(111, 88, 61, 40)$,
$(111, 95, 8, 2)$,
$(111, 78, 5, 19)$,
$(111, 86, 15, 69)$,\adfsplit
$(112, 92, 0, 14)$,
$(112, 99, 49, 6)$,
$(112, 72, 31, 45)$,
$(112, 88, 71, 32)$,\adfsplit
$(112, 107, 39, 28)$,
$(1, 6, 77, 126)$,
$(0, 1, 3, 99)$,
$(0, 2, 37, 83)$,\adfsplit
$(0, 6, 61, 101)$,
$(0, 5, 31, 95)$,
$(0, 7, 33, 103)$,
$(0, 10, 47, 107)$,\adfsplit
$(1, 10, 35, 83)$,
$(0, 13, 54, 75)$,
$(0, 42, 63, 89)$,
$(1, 39, 91, 113)$,\adfsplit
$(1, 55, 103, 127)$,
$(1, 43, 89, 114)$,
$(0, 26, 77, 113)$,
$(1, 22, 23, 97)$,\adfsplit
$(1, 50, 81, 121)$,
$(0, 51, 81, 114)$,
$(2, 55, 104, 128)$,
$(0, 15, 17, 82)$,\adfsplit
$(0, 22, 41, 71)$,
$(0, 30, 43, 84)$,
$(1, 51, 58, 84)$,
$(0, 59, 65, 76)$,\adfsplit
$(0, 50, 57, 90)$,
$(0, 18, 35, 98)$,
$(1, 34, 63, 92)$,
$(2, 47, 74, 120)$,\adfsplit
$(1, 2, 11, 78)$,
$(0, 62, 67, 88)$,
$(0, 58, 69, 100)$,
$(0, 46, 106, 128)$,\adfsplit
$(0, 38, 78, 127)$,
$(0, 25, 102, 120)$,
$(0, 45, 94, 121)$

%ADFvfyBlocksEnd
\adfLgap \noindent by the mapping:
$x \mapsto x + 4 j \adfmod{72}$ for $x < 72$,
$x \mapsto (x + 2 j \adfmod{36}) + 72$ for $72 \le x < 108$,
$x \mapsto (x - 108 + 7 j \adfmod{21}) + 108$ for $x \ge 108$,
$0 \le j < 18$.
\ADFvfyParStart{(129, ((63, 18, ((72, 4), (36, 2), (21, 7)))), ((18, 4), (36, 1), (21, 1)))} %ADFvfyParEnd
% End of 18^4 36^1 21^1
%%%%%%%%%%%%%%%%%%%%%%%%%%%%%%%%%%%%%%%%%%%%%%%%%%%%%%%%%%%%%%%%%%%%%%%%%%%%%%%%%%%%%%%%%%
%%%%%%%%%%%%%%%%%%%%%%%%%%%%%%%%%%%%%%%%%%%%%%%%%%%%%%%%%%%%%%%%%%%%%%%%%%%%%%%%%%%%%%%%%%

% Charlotte:GDD4-1-3-5-mod-6-TeX-gen-A:HITS-fun:4.10
\adfDgap
%ADFvfyBlocksStart {18,18,18,18,36,33}
\noindent{\boldmath $ 18^{4} 36^{1} 33^{1} $}~
With the point set $Z_{141}$ partitioned into
 residue classes modulo $4$ for $\{0, 1, \dots, 71\}$,
 $\{72, 73, \dots, 107\}$, and
 $\{108, 109, \dots, 140\}$,
 the design is generated from

\adfLgap %ADFvfyDesignStart
$(77, 48, 30, 15)$,
$(90, 21, 58, 3)$,
$(77, 14, 61, 52)$,
$(90, 63, 36, 49)$,\adfsplit
$(25, 32, 2, 63)$,
$(11, 32, 5, 62)$,
$(68, 53, 46, 7)$,
$(32, 19, 61, 22)$,\adfsplit
$(14, 0, 1, 43)$,
$(108, 89, 7, 16)$,
$(108, 74, 9, 54)$,
$(108, 105, 35, 17)$,\adfsplit
$(108, 102, 13, 24)$,
$(108, 76, 50, 63)$,
$(108, 73, 46, 8)$,
$(109, 72, 8, 43)$,\adfsplit
$(109, 100, 24, 45)$,
$(109, 105, 25, 15)$,
$(109, 103, 40, 38)$,
$(109, 77, 23, 46)$,\adfsplit
$(109, 74, 66, 5)$,
$(110, 100, 48, 2)$,
$(110, 81, 42, 35)$,
$(110, 95, 41, 8)$,\adfsplit
$(110, 74, 69, 19)$,
$(110, 78, 49, 58)$,
$(110, 79, 28, 3)$,
$(111, 96, 24, 43)$,\adfsplit
$(111, 89, 41, 62)$,
$(111, 79, 8, 63)$,
$(111, 105, 57, 58)$,
$(111, 74, 16, 13)$,\adfsplit
$(111, 94, 42, 23)$,
$(112, 81, 36, 51)$,
$(112, 104, 30, 44)$,
$(112, 91, 26, 35)$,\adfsplit
$(112, 89, 34, 5)$,
$(112, 78, 9, 28)$,
$(1, 67, 82, 112)$,
$(0, 3, 95, 113)$,\adfsplit
$(0, 17, 107, 114)$,
$(0, 2, 93, 125)$,
$(0, 23, 74, 136)$,
$(1, 6, 84, 114)$,\adfsplit
$(1, 3, 94, 125)$,
$(2, 7, 75, 114)$,
$(0, 46, 92, 115)$,
$(0, 22, 86, 116)$,\adfsplit
$(0, 71, 78, 117)$,
$(0, 6, 96, 118)$,
$(0, 25, 98, 124)$,
$(2, 43, 92, 113)$,\adfsplit
$(0, 7, 94, 128)$,
$(2, 39, 80, 118)$,
$(2, 47, 96, 115)$,
$(0, 10, 106, 138)$,\adfsplit
$(1, 47, 100, 116)$,
$(1, 59, 98, 137)$,
$(1, 42, 88, 113)$,
$(0, 18, 80, 139)$,\adfsplit
$(1, 51, 92, 140)$,
$(1, 70, 103, 139)$,
$(2, 27, 103, 138)$,
$(1, 71, 95, 138)$,\adfsplit
$(1, 35, 85, 118)$,
$(1, 11, 81, 128)$,
$(1, 54, 107, 117)$,
$(0, 67, 77, 135)$,\adfsplit
$(0, 66, 85, 140)$,
$(2, 3, 83, 137)$,
$(1, 18, 73, 135)$,
$(0, 5, 103, 126)$,\adfsplit
$(0, 41, 73, 137)$,
$(0, 37, 97, 127)$,
$(0, 49, 79, 129)$

%ADFvfyBlocksEnd
\adfLgap \noindent by the mapping:
$x \mapsto x + 4 j \adfmod{72}$ for $x < 72$,
$x \mapsto (x + 2 j \adfmod{36}) + 72$ for $72 \le x < 108$,
$x \mapsto (x - 108 + 11 j \adfmod{33}) + 108$ for $x \ge 108$,
$0 \le j < 18$.
\ADFvfyParStart{(141, ((75, 18, ((72, 4), (36, 2), (33, 11)))), ((18, 4), (36, 1), (33, 1)))} %ADFvfyParEnd
% End of 18^4 36^1 33^1
%%%%%%%%%%%%%%%%%%%%%%%%%%%%%%%%%%%%%%%%%%%%%%%%%%%%%%%%%%%%%%%%%%%%%%%%%%%%%%%%%%%%%%%%%%
%%%%%%%%%%%%%%%%%%%%%%%%%%%%%%%%%%%%%%%%%%%%%%%%%%%%%%%%%%%%%%%%%%%%%%%%%%%%%%%%%%%%%%%%%%

% Charlotte:GDD4-1-3-5-mod-6-TeX-gen-A:HITS-fun:4.10
\adfDgap
%ADFvfyBlocksStart {9,9,9,9,9,9,27,15}
\noindent{\boldmath $ 9^{6} 27^{1} 15^{1} $}~
With the point set $Z_{96}$ partitioned into
 residue classes modulo $6$ for $\{0, 1, \dots, 53\}$,
 $\{54, 55, \dots, 80\}$, and
 $\{81, 82, \dots, 95\}$,
 the design is generated from

\adfLgap %ADFvfyDesignStart
$(81, 54, 15, 49)$,
$(82, 54, 0, 43)$,
$(83, 54, 25, 47)$,
$(84, 54, 36, 19)$,\adfsplit
$(85, 54, 51, 34)$,
$(86, 54, 52, 44)$,
$(87, 54, 24, 17)$,
$(88, 54, 22, 29)$,\adfsplit
$(0, 1, 62, 86)$,
$(0, 3, 60, 92)$,
$(0, 2, 68, 88)$,
$(0, 4, 14, 73)$,\adfsplit
$(0, 5, 9, 79)$,
$(0, 13, 22, 61)$,
$(0, 41, 74, 89)$,
$(0, 51, 76, 94)$,\adfsplit
$(0, 33, 57, 95)$,
$(0, 25, 65, 85)$,
$(1, 11, 39, 71)$,
$(0, 23, 28, 72)$,\adfsplit
$(0, 15, 16, 56)$,
$(0, 27, 29, 78)$,
$(0, 11, 19, 34)$,
$(0, 21, 35, 80)$

%ADFvfyBlocksEnd
\adfLgap \noindent by the mapping:
$x \mapsto x + 2 j \adfmod{54}$ for $x < 54$,
$x \mapsto (x +  j \adfmod{27}) + 54$ for $54 \le x < 81$,
$x \mapsto (x - 81 + 5 j \adfmod{15}) + 81$ for $x \ge 81$,
$0 \le j < 27$.
\ADFvfyParStart{(96, ((24, 27, ((54, 2), (27, 1), (15, 5)))), ((9, 6), (27, 1), (15, 1)))} %ADFvfyParEnd
% End of 9^6 27^1 15^1
%%%%%%%%%%%%%%%%%%%%%%%%%%%%%%%%%%%%%%%%%%%%%%%%%%%%%%%%%%%%%%%%%%%%%%%%%%%%%%%%%%%%%%%%%%
%%%%%%%%%%%%%%%%%%%%%%%%%%%%%%%%%%%%%%%%%%%%%%%%%%%%%%%%%%%%%%%%%%%%%%%%%%%%%%%%%%%%%%%%%%

%%%%%%%%%%%%%%%%%%%%%%%%%%%%%%%%%%%%%%%%%%%%%%%%%%%%%%%%%%%%%%%%%%%%%%%%%%%%%%%%%%%%%%%%%%
%%%%%%%%%%%%%%%%%%%%%%%%%%%%%%%%%%%%%%%%%%%%%%%%%%%%%%%%%%%%%%%%%%%%%%%%%%%%%%%%%%%%%%%%%%
\section{4-GDDs for the proof of Lemma \ref{lem:4-GDD 39^8 m^1}}
\label{app:4-GDD 39^8 m^1}
\adfnull{
$ 39^8 120^1 $,
$ 39^8 123^1 $,
$ 39^8 126^1 $,
$ 39^8 129^1 $ and
$ 39^8 132^1 $.
}

% Charlotte:GDD4-1-3-5-mod-6-TeX-gen-A:HITS-fun:4.10
\adfDgap
%ADFvfyBlocksStart {39,39,39,39,39,39,39,39,120}
\noindent{\boldmath $ 39^{8} 120^{1} $}~
With the point set $Z_{432}$ partitioned into
 residue classes modulo $8$ for $\{0, 1, \dots, 311\}$, and
 $\{312, 313, \dots, 431\}$,
 the design is generated from

\adfLgap %ADFvfyDesignStart
$(312, 130, 256, 227)$,
$(312, 236, 243, 285)$,
$(312, 306, 168, 225)$,\adfsplit
$(312, 26, 241, 28)$,
$(312, 37, 262, 239)$,
$(312, 271, 276, 113)$,\adfsplit
$(312, 128, 159, 6)$,
$(312, 283, 293, 14)$,
$(313, 6, 108, 125)$,\adfsplit
$(313, 89, 229, 43)$,
$(313, 278, 18, 256)$,
$(313, 165, 263, 220)$,\adfsplit
$(313, 146, 240, 142)$,
$(313, 193, 175, 308)$,
$(313, 123, 154, 153)$,\adfsplit
$(313, 255, 35, 56)$,
$(314, 201, 208, 262)$,
$(314, 21, 267, 60)$,\adfsplit
$(314, 223, 50, 161)$,
$(314, 259, 47, 61)$,
$(314, 104, 265, 138)$,\adfsplit
$(314, 231, 48, 53)$,
$(314, 302, 11, 172)$,
$(314, 92, 250, 174)$,\adfsplit
$(315, 29, 257, 82)$,
$(315, 244, 286, 45)$,
$(315, 240, 177, 194)$,\adfsplit
$(315, 151, 139, 14)$,
$(315, 191, 36, 59)$,
$(315, 49, 15, 253)$,\adfsplit
$(315, 291, 208, 258)$,
$(315, 126, 236, 8)$,
$(316, 195, 74, 286)$,\adfsplit
$(316, 182, 114, 37)$,
$(316, 43, 15, 164)$,
$(316, 78, 185, 79)$,\adfsplit
$(316, 120, 130, 225)$,
$(316, 309, 23, 80)$,
$(316, 16, 241, 292)$,\adfsplit
$(316, 12, 155, 101)$,
$(317, 142, 97, 271)$,
$(317, 30, 48, 59)$,\adfsplit
$(317, 306, 45, 20)$,
$(317, 283, 8, 87)$,
$(317, 122, 208, 13)$,\adfsplit
$(317, 178, 191, 292)$,
$(317, 14, 171, 209)$,
$(317, 33, 101, 36)$,\adfsplit
$(318, 167, 189, 152)$,
$(318, 277, 204, 34)$,
$(318, 158, 221, 139)$,\adfsplit
$(318, 289, 99, 112)$,
$(318, 17, 120, 159)$,
$(318, 199, 122, 270)$,\adfsplit
$(318, 273, 203, 212)$,
$(318, 138, 244, 22)$,
$(319, 164, 53, 94)$,\adfsplit
$(319, 199, 106, 272)$,
$(319, 165, 121, 35)$,
$(319, 64, 17, 198)$,\adfsplit
$(319, 132, 159, 0)$,
$(319, 37, 273, 28)$,
$(319, 187, 167, 2)$,\adfsplit
$(319, 186, 206, 171)$,
$(0, 19, 108, 320)$,
$(0, 3, 62, 321)$,\adfsplit
$(0, 6, 41, 252)$,
$(0, 12, 59, 205)$,
$(0, 14, 301, 400)$,\adfsplit
$(0, 44, 99, 135)$,
$(0, 45, 140, 331)$,
$(0, 67, 69, 341)$,\adfsplit
$(0, 30, 169, 219)$,
$(0, 28, 103, 390)$,
$(0, 58, 181, 350)$,\adfsplit
$(0, 38, 162, 361)$,
$(0, 81, 133, 381)$,
$(0, 227, 285, 431)$,\adfsplit
$(0, 115, 263, 360)$,
$(0, 109, 233, 340)$,
$(0, 85, 187, 247)$,\adfsplit
$(0, 92, 237, 410)$,
$(0, 141, 147, 401)$,
$(0, 53, 171, 421)$,\adfsplit
$(1, 5, 95, 320)$,
$(0, 78, 156, 234)$,
$(1, 79, 157, 235)$

%ADFvfyBlocksEnd
\adfLgap \noindent by the mapping:
$x \mapsto x + 2 j \adfmod{312}$ for $x < 312$,
$x \mapsto (x - 312 + 10 j \adfmod{120}) + 312$ for $x \ge 312$,
$0 \le j < 156$
 for the first 85 blocks,
$0 \le j < 39$
 for the last two blocks.
\ADFvfyParStart{(432, ((85, 156, ((312, 2), (120, 10))), (2, 39, ((312, 2), (120, 10)))), ((39, 8), (120, 1)))} %ADFvfyParEnd
% End of 39^8 120^1
%%%%%%%%%%%%%%%%%%%%%%%%%%%%%%%%%%%%%%%%%%%%%%%%%%%%%%%%%%%%%%%%%%%%%%%%%%%%%%%%%%%%%%%%%%
%%%%%%%%%%%%%%%%%%%%%%%%%%%%%%%%%%%%%%%%%%%%%%%%%%%%%%%%%%%%%%%%%%%%%%%%%%%%%%%%%%%%%%%%%%

% Charlotte:GDD4-1-3-5-mod-6-TeX-gen-A:HITS-fun:4.10
\adfDgap
%ADFvfyBlocksStart {39,39,39,39,39,39,39,39,123}
\noindent{\boldmath $ 39^{8} 123^{1} $}~
With the point set $Z_{435}$ partitioned into
 residue classes modulo $8$ for $\{0, 1, \dots, 311\}$, and
 $\{312, 313, \dots, 434\}$,
 the design is generated from

\adfLgap %ADFvfyDesignStart
$(432, 240, 241, 206)$,
$(312, 274, 80, 260)$,
$(312, 194, 43, 23)$,\adfsplit
$(312, 151, 83, 100)$,
$(312, 180, 142, 105)$,
$(312, 126, 21, 288)$,\adfsplit
$(312, 37, 162, 121)$,
$(312, 279, 14, 75)$,
$(312, 221, 257, 16)$,\adfsplit
$(313, 273, 164, 160)$,
$(313, 31, 86, 217)$,
$(313, 34, 150, 251)$,\adfsplit
$(313, 143, 309, 100)$,
$(313, 293, 240, 211)$,
$(313, 228, 113, 80)$,\adfsplit
$(313, 94, 159, 218)$,
$(313, 133, 138, 195)$,
$(314, 171, 263, 273)$,\adfsplit
$(314, 278, 115, 16)$,
$(314, 26, 240, 164)$,
$(314, 32, 253, 295)$,\adfsplit
$(314, 82, 261, 15)$,
$(314, 28, 203, 257)$,
$(314, 73, 54, 42)$,\adfsplit
$(314, 293, 166, 300)$,
$(315, 49, 77, 3)$,
$(315, 238, 261, 19)$,\adfsplit
$(315, 157, 2, 88)$,
$(315, 258, 198, 177)$,
$(315, 76, 155, 161)$,\adfsplit
$(315, 252, 82, 79)$,
$(315, 104, 86, 215)$,
$(0, 63, 140, 315)$,\adfsplit
$(0, 9, 39, 198)$,
$(0, 90, 190, 316)$,
$(0, 25, 202, 371)$,\adfsplit
$(0, 27, 121, 356)$,
$(0, 22, 119, 416)$,
$(0, 44, 239, 391)$,\adfsplit
$(0, 2, 89, 147)$,
$(0, 52, 158, 381)$,
$(0, 11, 26, 386)$,\adfsplit
$(0, 13, 143, 401)$,
$(0, 78, 156, 234)$

%ADFvfyBlocksEnd
\adfLgap \noindent by the mapping:
$x \mapsto x +  j \adfmod{312}$ for $x < 312$,
$x \mapsto (x - 312 + 5 j \adfmod{120}) + 312$ for $312 \le x < 432$,
$x \mapsto (x +  j \adfmod{3}) + 432$ for $x \ge 432$,
$0 \le j < 312$
 for the first 43 blocks,
$0 \le j < 78$
 for the last block.
\ADFvfyParStart{(435, ((43, 312, ((312, 1), (120, 5), (3, 1))), (1, 78, ((312, 1), (120, 5), (3, 1)))), ((39, 8), (123, 1)))} %ADFvfyParEnd
% End of 39^8 123^1
%%%%%%%%%%%%%%%%%%%%%%%%%%%%%%%%%%%%%%%%%%%%%%%%%%%%%%%%%%%%%%%%%%%%%%%%%%%%%%%%%%%%%%%%%%
%%%%%%%%%%%%%%%%%%%%%%%%%%%%%%%%%%%%%%%%%%%%%%%%%%%%%%%%%%%%%%%%%%%%%%%%%%%%%%%%%%%%%%%%%%

% Charlotte:GDD4-1-3-5-mod-6-TeX-gen-A:HITS-fun:4.10
\adfDgap
%ADFvfyBlocksStart {39,39,39,39,39,39,39,39,126}
\noindent{\boldmath $ 39^{8} 126^{1} $}~
With the point set $Z_{438}$ partitioned into
 residue classes modulo $7$ for $\{0, 1, \dots, 272\}$,
 $\{273, 274, \dots, 311\}$, and
 $\{312, 313, \dots, 437\}$,
 the design is generated from

\adfLgap %ADFvfyDesignStart
$(70, 124, 153, 366)$,
$(24, 25, 258, 387)$,
$(47, 50, 157, 433)$,\adfsplit
$(106, 252, 311, 365)$,
$(29, 149, 308, 319)$,
$(54, 70, 92, 318)$,\adfsplit
$(82, 244, 304, 389)$,
$(19, 91, 256, 324)$,
$(12, 62, 157, 401)$,\adfsplit
$(15, 47, 182, 356)$,
$(80, 111, 169, 393)$,
$(47, 150, 267, 411)$,\adfsplit
$(26, 120, 167, 429)$,
$(18, 103, 302, 356)$,
$(117, 162, 179, 418)$,\adfsplit
$(60, 259, 299, 326)$,
$(35, 220, 278, 422)$,
$(12, 42, 146, 406)$,\adfsplit
$(123, 174, 305, 407)$,
$(52, 224, 304, 376)$,
$(20, 197, 207, 337)$,\adfsplit
$(30, 189, 303, 342)$,
$(181, 208, 254, 382)$,
$(101, 226, 251, 352)$,\adfsplit
$(227, 233, 303, 421)$,
$(62, 191, 260, 434)$,
$(135, 237, 308, 376)$,\adfsplit
$(112, 124, 303, 409)$,
$(40, 223, 232, 372)$,
$(3, 81, 163, 406)$,\adfsplit
$(0, 15, 295, 364)$,
$(0, 2, 43, 318)$,
$(0, 19, 281, 348)$,\adfsplit
$(0, 48, 298, 390)$,
$(0, 4, 209, 311)$,
$(0, 18, 303, 327)$,\adfsplit
$(0, 92, 302, 333)$,
$(0, 20, 306, 423)$,
$(0, 5, 57, 405)$,\adfsplit
$(0, 13, 80, 313)$,
$(0, 23, 116, 391)$,
$(0, 24, 176, 349)$,\adfsplit
$(0, 8, 130, 356)$,
$(0, 76, 155, 335)$,
$(0, 26, 87, 320)$,\adfsplit
$(0, 33, 142, 350)$,
$(0, 59, 124, 419)$,
$(0, 37, 136, 323)$,\adfsplit
$(0, 34, 100, 347)$,
$(0, 11, 55, 213)$

%ADFvfyBlocksEnd
\adfLgap \noindent by the mapping:
$x \mapsto x +  j \adfmod{273}$ for $x < 273$,
$x \mapsto (x +  j \adfmod{39}) + 273$ for $273 \le x < 312$,
$x \mapsto (x - 312 + 6 j \adfmod{126}) + 312$ for $x \ge 312$,
$0 \le j < 273$.
\ADFvfyParStart{(438, ((50, 273, ((273, 1), (39, 1), (126, 6)))), ((39, 7), (39, 1), (126, 1)))} %ADFvfyParEnd
% End of 39^8 126^1
%%%%%%%%%%%%%%%%%%%%%%%%%%%%%%%%%%%%%%%%%%%%%%%%%%%%%%%%%%%%%%%%%%%%%%%%%%%%%%%%%%%%%%%%%%
%%%%%%%%%%%%%%%%%%%%%%%%%%%%%%%%%%%%%%%%%%%%%%%%%%%%%%%%%%%%%%%%%%%%%%%%%%%%%%%%%%%%%%%%%%

% Charlotte:GDD4-1-3-5-mod-6-TeX-gen-A:HITS-fun:4.10
\adfDgap
%ADFvfyBlocksStart {39,39,39,39,39,39,39,39,129}
\noindent{\boldmath $ 39^{8} 129^{1} $}~
With the point set $Z_{441}$ partitioned into
 residue classes modulo $8$ for $\{0, 1, \dots, 311\}$, and
 $\{312, 313, \dots, 440\}$,
 the design is generated from

\adfLgap %ADFvfyDesignStart
$(432, 236, 175, 237)$,
$(433, 300, 266, 13)$,
$(434, 191, 133, 177)$,\adfsplit
$(312, 270, 163, 263)$,
$(312, 193, 79, 158)$,
$(312, 0, 3, 12)$,\adfsplit
$(312, 305, 253, 207)$,
$(312, 11, 116, 237)$,
$(312, 248, 269, 42)$,\adfsplit
$(312, 28, 122, 9)$,
$(312, 160, 214, 250)$,
$(313, 156, 23, 154)$,\adfsplit
$(313, 136, 245, 9)$,
$(313, 93, 76, 262)$,
$(313, 144, 260, 193)$,\adfsplit
$(313, 229, 199, 75)$,
$(313, 8, 110, 185)$,
$(313, 282, 251, 174)$,\adfsplit
$(313, 290, 303, 139)$,
$(314, 131, 199, 194)$,
$(314, 274, 117, 246)$,\adfsplit
$(314, 305, 263, 40)$,
$(314, 153, 133, 171)$,
$(314, 128, 260, 5)$,\adfsplit
$(314, 36, 206, 96)$,
$(314, 73, 4, 166)$,
$(314, 39, 306, 211)$,\adfsplit
$(315, 281, 39, 50)$,
$(315, 175, 16, 134)$,
$(315, 3, 36, 13)$,\adfsplit
$(315, 28, 193, 120)$,
$(315, 56, 251, 294)$,
$(315, 260, 245, 70)$,\adfsplit
$(0, 66, 167, 326)$,
$(0, 22, 283, 386)$,
$(0, 65, 228, 426)$,\adfsplit
$(0, 71, 201, 385)$,
$(0, 26, 141, 330)$,
$(0, 50, 259, 321)$,\adfsplit
$(0, 6, 125, 401)$,
$(0, 37, 215, 416)$,
$(0, 27, 82, 173)$,\adfsplit
$(0, 39, 138, 361)$,
$(0, 4, 87, 336)$,
$(0, 78, 156, 234)$

%ADFvfyBlocksEnd
\adfLgap \noindent by the mapping:
$x \mapsto x +  j \adfmod{312}$ for $x < 312$,
$x \mapsto (x - 312 + 5 j \adfmod{120}) + 312$ for $312 \le x < 432$,
$x \mapsto (x + 3 j \adfmod{9}) + 432$ for $x \ge 432$,
$0 \le j < 312$
 for the first 44 blocks,
$0 \le j < 78$
 for the last block.
\ADFvfyParStart{(441, ((44, 312, ((312, 1), (120, 5), (9, 3))), (1, 78, ((312, 1), (120, 5), (9, 3)))), ((39, 8), (129, 1)))} %ADFvfyParEnd
% End of 39^8 129^1
%%%%%%%%%%%%%%%%%%%%%%%%%%%%%%%%%%%%%%%%%%%%%%%%%%%%%%%%%%%%%%%%%%%%%%%%%%%%%%%%%%%%%%%%%%
%%%%%%%%%%%%%%%%%%%%%%%%%%%%%%%%%%%%%%%%%%%%%%%%%%%%%%%%%%%%%%%%%%%%%%%%%%%%%%%%%%%%%%%%%%

% Charlotte:GDD4-1-3-5-mod-6-TeX-gen-A:HITS-fun:4.10
\adfDgap
%ADFvfyBlocksStart {39,39,39,39,39,39,39,39,132}
\noindent{\boldmath $ 39^{8} 132^{1} $}~
With the point set $Z_{444}$ partitioned into
 residue classes modulo $8$ for $\{0, 1, \dots, 311\}$, and
 $\{312, 313, \dots, 443\}$,
 the design is generated from

\adfLgap %ADFvfyDesignStart
$(312, 16, 5, 139)$,
$(312, 158, 277, 4)$,
$(312, 248, 118, 33)$,\adfsplit
$(312, 106, 305, 75)$,
$(312, 30, 144, 295)$,
$(312, 44, 135, 114)$,\adfsplit
$(312, 50, 300, 119)$,
$(312, 169, 213, 179)$,
$(313, 292, 41, 128)$,\adfsplit
$(313, 99, 103, 154)$,
$(313, 258, 78, 68)$,
$(313, 245, 278, 0)$,\adfsplit
$(313, 203, 119, 276)$,
$(313, 13, 63, 217)$,
$(313, 67, 170, 262)$,\adfsplit
$(313, 225, 93, 112)$,
$(314, 213, 273, 251)$,
$(314, 286, 151, 92)$,\adfsplit
$(314, 294, 29, 56)$,
$(314, 95, 170, 113)$,
$(314, 229, 139, 192)$,\adfsplit
$(314, 204, 1, 202)$,
$(314, 90, 196, 112)$,
$(314, 207, 243, 302)$,\adfsplit
$(315, 0, 58, 269)$,
$(315, 166, 189, 68)$,
$(315, 73, 62, 263)$,\adfsplit
$(315, 201, 150, 59)$,
$(315, 172, 243, 272)$,
$(315, 186, 136, 252)$,\adfsplit
$(315, 39, 187, 157)$,
$(315, 2, 305, 55)$,
$(316, 102, 203, 16)$,\adfsplit
$(316, 47, 14, 140)$,
$(316, 271, 194, 297)$,
$(316, 211, 28, 32)$,\adfsplit
$(316, 253, 82, 190)$,
$(316, 303, 165, 264)$,
$(316, 252, 89, 210)$,\adfsplit
$(316, 27, 193, 245)$,
$(317, 235, 170, 47)$,
$(317, 197, 278, 210)$,\adfsplit
$(317, 154, 102, 244)$,
$(317, 41, 141, 286)$,
$(317, 243, 156, 241)$,\adfsplit
$(317, 231, 308, 56)$,
$(317, 96, 301, 81)$,
$(317, 79, 59, 64)$,\adfsplit
$(318, 287, 310, 275)$,
$(318, 24, 68, 113)$,
$(318, 133, 222, 252)$,\adfsplit
$(318, 91, 63, 177)$,
$(318, 242, 241, 171)$,
$(318, 149, 292, 130)$,\adfsplit
$(318, 136, 42, 271)$,
$(318, 32, 165, 206)$,
$(319, 22, 168, 103)$,\adfsplit
$(319, 277, 167, 184)$,
$(319, 123, 278, 169)$,
$(319, 141, 4, 138)$,\adfsplit
$(319, 163, 156, 128)$,
$(319, 236, 250, 65)$,
$(319, 207, 102, 26)$,\adfsplit
$(319, 129, 197, 203)$,
$(320, 29, 239, 156)$,
$(320, 274, 121, 303)$,\adfsplit
$(320, 19, 186, 232)$,
$(320, 142, 271, 196)$,
$(0, 1, 300, 342)$,\adfsplit
$(0, 9, 195, 320)$,
$(0, 20, 173, 419)$,
$(0, 27, 143, 375)$,\adfsplit
$(0, 38, 163, 321)$,
$(0, 57, 207, 322)$,
$(0, 17, 31, 332)$,\adfsplit
$(0, 5, 26, 365)$,
$(0, 110, 275, 443)$,
$(0, 6, 249, 431)$,\adfsplit
$(0, 55, 309, 343)$,
$(0, 41, 147, 420)$,
$(0, 139, 215, 387)$,\adfsplit
$(0, 102, 263, 442)$,
$(0, 49, 115, 432)$,
$(0, 25, 230, 421)$,\adfsplit
$(0, 36, 79, 399)$,
$(0, 73, 188, 355)$,
$(0, 95, 140, 388)$,\adfsplit
$(0, 18, 251, 305)$,
$(1, 43, 215, 443)$,
$(0, 78, 156, 234)$,\adfsplit
$(1, 79, 157, 235)$

%ADFvfyBlocksEnd
\adfLgap \noindent by the mapping:
$x \mapsto x + 2 j \adfmod{312}$ for $x < 312$,
$x \mapsto (x - 312 + 11 j \adfmod{132}) + 312$ for $x \ge 312$,
$0 \le j < 156$
 for the first 89 blocks,
$0 \le j < 39$
 for the last two blocks.
\ADFvfyParStart{(444, ((89, 156, ((312, 2), (132, 11))), (2, 39, ((312, 2), (132, 11)))), ((39, 8), (132, 1)))} %ADFvfyParEnd
% End of 39^8 132^1
%%%%%%%%%%%%%%%%%%%%%%%%%%%%%%%%%%%%%%%%%%%%%%%%%%%%%%%%%%%%%%%%%%%%%%%%%%%%%%%%%%%%%%%%%%
%%%%%%%%%%%%%%%%%%%%%%%%%%%%%%%%%%%%%%%%%%%%%%%%%%%%%%%%%%%%%%%%%%%%%%%%%%%%%%%%%%%%%%%%%%

%%%%%%%%%%%%%%%%%%%%%%%%%%%%%%%%%%%%%%%%%%%%%%%%%%%%%%%%%%%%%%%%%%%%%%%%%%%%%%%%%%%%%%%%%%
%%%%%%%%%%%%%%%%%%%%%%%%%%%%%%%%%%%%%%%%%%%%%%%%%%%%%%%%%%%%%%%%%%%%%%%%%%%%%%%%%%%%%%%%%%
\section{4-GDDs for the proof of Lemma \ref{lem:4-GDD 51^8 m^1}}
\label{app:4-GDD 51^8 m^1}
\adfnull{
$ 51^8 159^1 $,
$ 51^8 162^1 $,
$ 51^8 165^1 $,
$ 51^8 168^1 $,
$ 51^8 171^1 $ and
$ 51^8 174^1 $.
}

% Charlotte:GDD4-1-3-5-mod-6-TeX-gen-A:HITS-fun:4.10
\adfDgap
%ADFvfyBlocksStart {51,51,51,51,51,51,51,51,159}
\noindent{\boldmath $ 51^{8} 159^{1} $}~
With the point set $Z_{567}$ partitioned into
 residue classes modulo $8$ for $\{0, 1, \dots, 407\}$, and
 $\{408, 409, \dots, 566\}$,
 the design is generated from

\adfLgap %ADFvfyDesignStart
$(552, 363, 88, 17)$,
$(553, 218, 118, 297)$,
$(554, 264, 316, 113)$,\adfsplit
$(555, 253, 228, 71)$,
$(556, 318, 250, 311)$,
$(408, 76, 43, 81)$,\adfsplit
$(408, 213, 276, 26)$,
$(408, 271, 318, 289)$,
$(408, 44, 32, 306)$,\adfsplit
$(408, 112, 209, 101)$,
$(408, 230, 397, 226)$,
$(408, 203, 120, 383)$,\adfsplit
$(408, 147, 190, 279)$,
$(409, 15, 250, 60)$,
$(409, 290, 152, 51)$,\adfsplit
$(409, 239, 288, 182)$,
$(409, 140, 283, 162)$,
$(409, 245, 7, 275)$,\adfsplit
$(409, 261, 366, 177)$,
$(409, 277, 113, 148)$,
$(409, 406, 400, 313)$,\adfsplit
$(410, 40, 135, 37)$,
$(410, 259, 81, 148)$,
$(410, 386, 251, 144)$,\adfsplit
$(410, 123, 31, 286)$,
$(410, 54, 239, 389)$,
$(410, 116, 17, 93)$,\adfsplit
$(410, 178, 300, 374)$,
$(410, 186, 73, 56)$,
$(411, 256, 201, 142)$,\adfsplit
$(411, 396, 41, 91)$,
$(411, 172, 247, 53)$,
$(411, 327, 325, 306)$,\adfsplit
$(411, 215, 78, 0)$,
$(411, 164, 49, 165)$,
$(411, 106, 326, 339)$,\adfsplit
$(411, 347, 296, 362)$,
$(412, 110, 144, 309)$,
$(412, 321, 90, 54)$,\adfsplit
$(412, 196, 407, 368)$,
$(412, 26, 227, 173)$,
$(412, 205, 400, 219)$,\adfsplit
$(0, 58, 318, 533)$,
$(0, 91, 222, 485)$,
$(0, 31, 362, 443)$,\adfsplit
$(0, 20, 266, 479)$,
$(0, 26, 86, 436)$,
$(0, 10, 149, 514)$,\adfsplit
$(0, 42, 326, 550)$,
$(0, 70, 155, 413)$,
$(0, 28, 109, 234)$,\adfsplit
$(0, 41, 339, 539)$,
$(0, 27, 252, 527)$,
$(0, 65, 159, 282)$,\adfsplit
$(0, 9, 127, 461)$,
$(0, 37, 154, 198)$,
$(0, 102, 204, 306)$

%ADFvfyBlocksEnd
\adfLgap \noindent by the mapping:
$x \mapsto x +  j \adfmod{408}$ for $x < 408$,
$x \mapsto (x - 408 + 6 j \adfmod{144}) + 408$ for $408 \le x < 552$,
$x \mapsto (x - 552 + 5 j \adfmod{15}) + 552$ for $x \ge 552$,
$0 \le j < 408$
 for the first 56 blocks,
$0 \le j < 102$
 for the last block.
\ADFvfyParStart{(567, ((56, 408, ((408, 1), (144, 6), (15, 5))), (1, 102, ((408, 1), (144, 6), (15, 5)))), ((51, 8), (159, 1)))} %ADFvfyParEnd
% End of 51^8 159^1
%%%%%%%%%%%%%%%%%%%%%%%%%%%%%%%%%%%%%%%%%%%%%%%%%%%%%%%%%%%%%%%%%%%%%%%%%%%%%%%%%%%%%%%%%%
%%%%%%%%%%%%%%%%%%%%%%%%%%%%%%%%%%%%%%%%%%%%%%%%%%%%%%%%%%%%%%%%%%%%%%%%%%%%%%%%%%%%%%%%%%

% Charlotte:GDD4-1-3-5-mod-6-TeX-gen-A:HITS-fun:4.10
\adfDgap
%ADFvfyBlocksStart {51,51,51,51,51,51,51,51,162}
\noindent{\boldmath $ 51^{8} 162^{1} $}~
With the point set $Z_{570}$ partitioned into
 residue classes modulo $8$ for $\{0, 1, \dots, 407\}$, and
 $\{408, 409, \dots, 569\}$,
 the design is generated from

\adfLgap %ADFvfyDesignStart
$(552, 18, 80, 253)$,
$(552, 309, 203, 190)$,
$(553, 325, 180, 358)$,\adfsplit
$(553, 21, 314, 71)$,
$(554, 392, 263, 94)$,
$(554, 355, 216, 255)$,\adfsplit
$(555, 228, 41, 355)$,
$(555, 69, 94, 92)$,
$(556, 261, 210, 175)$,\adfsplit
$(556, 316, 347, 350)$,
$(557, 197, 324, 170)$,
$(557, 39, 298, 91)$,\adfsplit
$(408, 87, 249, 254)$,
$(408, 23, 339, 22)$,
$(408, 222, 360, 109)$,\adfsplit
$(408, 247, 116, 355)$,
$(408, 400, 244, 314)$,
$(408, 317, 233, 250)$,\adfsplit
$(408, 272, 156, 59)$,
$(408, 282, 73, 141)$,
$(409, 269, 295, 1)$,\adfsplit
$(409, 270, 228, 219)$,
$(409, 200, 91, 110)$,
$(409, 264, 263, 209)$,\adfsplit
$(409, 138, 325, 356)$,
$(409, 347, 76, 105)$,
$(409, 286, 250, 375)$,\adfsplit
$(409, 218, 189, 232)$,
$(410, 364, 351, 120)$,
$(410, 325, 332, 98)$,\adfsplit
$(410, 152, 126, 187)$,
$(410, 387, 393, 127)$,
$(410, 25, 378, 304)$,\adfsplit
$(410, 189, 12, 347)$,
$(410, 310, 245, 209)$,
$(410, 202, 383, 302)$,\adfsplit
$(411, 295, 97, 261)$,
$(411, 253, 94, 332)$,
$(411, 179, 24, 226)$,\adfsplit
$(411, 230, 176, 283)$,
$(411, 77, 348, 87)$,
$(411, 225, 246, 330)$,\adfsplit
$(411, 112, 51, 292)$,
$(411, 257, 338, 383)$,
$(412, 162, 13, 206)$,\adfsplit
$(412, 44, 187, 370)$,
$(412, 286, 72, 151)$,
$(412, 273, 252, 95)$,\adfsplit
$(412, 73, 347, 327)$,
$(412, 136, 333, 2)$,
$(412, 246, 125, 27)$,\adfsplit
$(412, 41, 364, 344)$,
$(413, 151, 336, 5)$,
$(413, 171, 237, 30)$,\adfsplit
$(413, 196, 98, 41)$,
$(413, 119, 36, 114)$,
$(413, 61, 105, 59)$,\adfsplit
$(413, 200, 82, 92)$,
$(413, 283, 25, 70)$,
$(413, 230, 351, 88)$,\adfsplit
$(414, 317, 74, 28)$,
$(414, 299, 65, 224)$,
$(414, 202, 208, 374)$,\adfsplit
$(414, 94, 27, 9)$,
$(414, 380, 399, 205)$,
$(414, 175, 12, 357)$,\adfsplit
$(414, 241, 211, 23)$,
$(414, 264, 294, 18)$,
$(415, 179, 110, 170)$,\adfsplit
$(415, 273, 391, 163)$,
$(415, 76, 171, 293)$,
$(415, 239, 329, 322)$,\adfsplit
$(415, 111, 94, 396)$,
$(415, 164, 258, 136)$,
$(415, 325, 385, 120)$,\adfsplit
$(415, 237, 150, 128)$,
$(416, 29, 375, 200)$,
$(416, 348, 225, 295)$,\adfsplit
$(416, 282, 381, 385)$,
$(416, 202, 254, 4)$,
$(416, 190, 43, 215)$,\adfsplit
$(416, 266, 277, 195)$,
$(416, 126, 305, 192)$,
$(416, 92, 227, 16)$,\adfsplit
$(417, 285, 94, 76)$,
$(417, 349, 386, 19)$,
$(417, 227, 352, 302)$,\adfsplit
$(417, 401, 154, 336)$,
$(417, 20, 63, 121)$,
$(417, 383, 99, 198)$,\adfsplit
$(417, 271, 248, 345)$,
$(417, 108, 402, 269)$,
$(0, 33, 171, 212)$,\adfsplit
$(0, 37, 153, 418)$,
$(0, 59, 245, 284)$,
$(0, 3, 15, 419)$,\adfsplit
$(0, 12, 269, 502)$,
$(0, 93, 188, 305)$,
$(0, 321, 397, 491)$,\adfsplit
$(0, 47, 203, 526)$,
$(0, 63, 193, 515)$,
$(0, 291, 319, 442)$,\adfsplit
$(0, 57, 71, 550)$,
$(0, 293, 315, 466)$,
$(0, 49, 255, 539)$,\adfsplit
$(0, 73, 111, 349)$,
$(0, 339, 381, 503)$,
$(0, 133, 148, 431)$,\adfsplit
$(0, 140, 359, 455)$,
$(0, 91, 222, 490)$,
$(0, 4, 150, 301)$,\adfsplit
$(0, 38, 130, 527)$,
$(0, 58, 340, 538)$,
$(0, 102, 204, 306)$,\adfsplit
$(1, 103, 205, 307)$

%ADFvfyBlocksEnd
\adfLgap \noindent by the mapping:
$x \mapsto x + 2 j \adfmod{408}$ for $x < 408$,
$x \mapsto (x - 408 + 12 j \adfmod{144}) + 408$ for $408 \le x < 552$,
$x \mapsto (x - 552 + 6 j \adfmod{18}) + 552$ for $x \ge 552$,
$0 \le j < 204$
 for the first 113 blocks,
$0 \le j < 51$
 for the last two blocks.
\ADFvfyParStart{(570, ((113, 204, ((408, 2), (144, 12), (18, 6))), (2, 51, ((408, 2), (144, 12), (18, 6)))), ((51, 8), (162, 1)))} %ADFvfyParEnd
% End of 51^8 162^1
%%%%%%%%%%%%%%%%%%%%%%%%%%%%%%%%%%%%%%%%%%%%%%%%%%%%%%%%%%%%%%%%%%%%%%%%%%%%%%%%%%%%%%%%%%
%%%%%%%%%%%%%%%%%%%%%%%%%%%%%%%%%%%%%%%%%%%%%%%%%%%%%%%%%%%%%%%%%%%%%%%%%%%%%%%%%%%%%%%%%%

% Charlotte:GDD4-1-3-5-mod-6-TeX-gen-A:HITS-fun:4.10
\adfDgap
%ADFvfyBlocksStart {51,51,51,51,51,51,51,51,165}
\noindent{\boldmath $ 51^{8} 165^{1} $}~
With the point set $Z_{573}$ partitioned into
 residue classes modulo $8$ for $\{0, 1, \dots, 407\}$, and
 $\{408, 409, \dots, 572\}$,
 the design is generated from

\adfLgap %ADFvfyDesignStart
$(552, 38, 148, 51)$,
$(553, 56, 174, 229)$,
$(554, 126, 44, 64)$,\adfsplit
$(555, 216, 364, 329)$,
$(556, 163, 6, 23)$,
$(557, 91, 24, 2)$,\adfsplit
$(558, 348, 73, 176)$,
$(408, 114, 405, 108)$,
$(408, 212, 240, 35)$,\adfsplit
$(408, 102, 55, 153)$,
$(408, 99, 110, 325)$,
$(408, 146, 233, 53)$,\adfsplit
$(408, 224, 97, 382)$,
$(408, 239, 154, 364)$,
$(408, 91, 136, 15)$,\adfsplit
$(409, 297, 76, 43)$,
$(409, 351, 206, 84)$,
$(409, 391, 97, 406)$,\adfsplit
$(409, 78, 261, 266)$,
$(409, 250, 216, 359)$,
$(409, 245, 75, 353)$,\adfsplit
$(409, 232, 155, 157)$,
$(409, 128, 404, 138)$,
$(410, 287, 371, 97)$,\adfsplit
$(410, 223, 198, 261)$,
$(410, 190, 43, 402)$,
$(410, 168, 226, 182)$,\adfsplit
$(410, 386, 221, 224)$,
$(410, 148, 17, 184)$,
$(410, 231, 300, 243)$,\adfsplit
$(410, 133, 297, 380)$,
$(411, 329, 229, 15)$,
$(411, 339, 285, 380)$,\adfsplit
$(411, 162, 341, 70)$,
$(411, 6, 288, 265)$,
$(411, 225, 232, 251)$,\adfsplit
$(411, 52, 218, 176)$,
$(411, 119, 187, 346)$,
$(411, 223, 324, 398)$,\adfsplit
$(412, 213, 299, 303)$,
$(412, 138, 109, 294)$,
$(412, 193, 236, 184)$,\adfsplit
$(0, 27, 329, 503)$,
$(0, 1, 207, 497)$,
$(0, 21, 279, 484)$,\adfsplit
$(0, 37, 107, 509)$,
$(0, 46, 245, 412)$,
$(0, 81, 234, 496)$,\adfsplit
$(0, 66, 303, 418)$,
$(0, 18, 71, 508)$,
$(0, 78, 273, 473)$,\adfsplit
$(0, 31, 222, 413)$,
$(0, 39, 292, 479)$,
$(0, 50, 115, 539)$,\adfsplit
$(0, 30, 91, 269)$,
$(0, 59, 119, 270)$,
$(0, 73, 219, 467)$,\adfsplit
$(0, 102, 204, 306)$

%ADFvfyBlocksEnd
\adfLgap \noindent by the mapping:
$x \mapsto x +  j \adfmod{408}$ for $x < 408$,
$x \mapsto (x - 408 + 6 j \adfmod{144}) + 408$ for $408 \le x < 552$,
$x \mapsto (x - 552 + 7 j \adfmod{21}) + 552$ for $x \ge 552$,
$0 \le j < 408$
 for the first 57 blocks,
$0 \le j < 102$
 for the last block.
\ADFvfyParStart{(573, ((57, 408, ((408, 1), (144, 6), (21, 7))), (1, 102, ((408, 1), (144, 6), (21, 7)))), ((51, 8), (165, 1)))} %ADFvfyParEnd
% End of 51^8 165^1
%%%%%%%%%%%%%%%%%%%%%%%%%%%%%%%%%%%%%%%%%%%%%%%%%%%%%%%%%%%%%%%%%%%%%%%%%%%%%%%%%%%%%%%%%%
%%%%%%%%%%%%%%%%%%%%%%%%%%%%%%%%%%%%%%%%%%%%%%%%%%%%%%%%%%%%%%%%%%%%%%%%%%%%%%%%%%%%%%%%%%

% Charlotte:GDD4-1-3-5-mod-6-TeX-gen-A:HITS-fun:4.10
\adfDgap
%ADFvfyBlocksStart {51,51,51,51,51,51,51,51,168}
\noindent{\boldmath $ 51^{8} 168^{1} $}~
With the point set $Z_{576}$ partitioned into
 residue classes modulo $7$ for $\{0, 1, \dots, 356\}$,
 $\{357, 358, \dots, 407\}$, and
 $\{408, 409, \dots, 575\}$,
 the design is generated from

\adfLgap %ADFvfyDesignStart
$(94, 258, 407, 453)$,
$(203, 222, 346, 451)$,
$(173, 231, 286, 552)$,\adfsplit
$(206, 208, 226, 505)$,
$(27, 71, 135, 472)$,
$(165, 222, 387, 539)$,\adfsplit
$(13, 110, 116, 557)$,
$(5, 67, 387, 408)$,
$(1, 94, 394, 443)$,\adfsplit
$(5, 177, 403, 562)$,
$(8, 68, 186, 481)$,
$(158, 262, 296, 513)$,\adfsplit
$(178, 193, 205, 469)$,
$(2, 68, 147, 434)$,
$(50, 132, 355, 434)$,\adfsplit
$(87, 96, 375, 524)$,
$(7, 262, 334, 545)$,
$(272, 312, 406, 458)$,\adfsplit
$(131, 279, 308, 528)$,
$(14, 212, 405, 527)$,
$(129, 198, 266, 570)$,\adfsplit
$(196, 232, 405, 462)$,
$(46, 138, 268, 510)$,
$(321, 322, 401, 467)$,\adfsplit
$(179, 201, 406, 545)$,
$(2, 41, 322, 483)$,
$(44, 82, 255, 572)$,\adfsplit
$(204, 209, 357, 562)$,
$(48, 330, 365, 510)$,
$(5, 112, 179, 469)$,\adfsplit
$(53, 66, 301, 564)$,
$(320, 343, 381, 526)$,
$(84, 169, 335, 438)$,\adfsplit
$(146, 156, 404, 455)$,
$(158, 199, 383, 565)$,
$(121, 153, 297, 554)$,\adfsplit
$(26, 284, 362, 533)$,
$(87, 91, 208, 548)$,
$(217, 267, 356, 466)$,\adfsplit
$(19, 64, 385, 480)$,
$(46, 245, 325, 564)$,
$(44, 353, 392, 428)$,\adfsplit
$(141, 256, 327, 560)$,
$(41, 57, 333, 551)$,
$(89, 296, 304, 437)$,\adfsplit
$(288, 305, 406, 524)$,
$(0, 25, 51, 533)$,
$(0, 3, 246, 468)$,\adfsplit
$(0, 11, 163, 520)$,
$(0, 120, 383, 456)$,
$(0, 31, 405, 553)$,\adfsplit
$(0, 59, 261, 481)$,
$(0, 110, 363, 569)$,
$(0, 129, 370, 479)$,\adfsplit
$(0, 87, 257, 361)$,
$(0, 61, 267, 499)$,
$(0, 33, 262, 555)$,\adfsplit
$(0, 46, 241, 483)$,
$(0, 74, 234, 422)$,
$(0, 83, 169, 534)$,\adfsplit
$(0, 43, 192, 518)$,
$(0, 53, 141, 439)$,
$(0, 54, 127, 221)$,\adfsplit
$(0, 47, 204, 543)$,
$(0, 24, 156, 471)$,
$(0, 101, 226, 495)$

%ADFvfyBlocksEnd
\adfLgap \noindent by the mapping:
$x \mapsto x +  j \adfmod{357}$ for $x < 357$,
$x \mapsto (x +  j \adfmod{51}) + 357$ for $357 \le x < 408$,
$x \mapsto (x - 408 + 8 j \adfmod{168}) + 408$ for $x \ge 408$,
$0 \le j < 357$.
\ADFvfyParStart{(576, ((66, 357, ((357, 1), (51, 1), (168, 8)))), ((51, 7), (51, 1), (168, 1)))} %ADFvfyParEnd
% End of 51^8 168^1
%%%%%%%%%%%%%%%%%%%%%%%%%%%%%%%%%%%%%%%%%%%%%%%%%%%%%%%%%%%%%%%%%%%%%%%%%%%%%%%%%%%%%%%%%%
%%%%%%%%%%%%%%%%%%%%%%%%%%%%%%%%%%%%%%%%%%%%%%%%%%%%%%%%%%%%%%%%%%%%%%%%%%%%%%%%%%%%%%%%%%

% Charlotte:GDD4-1-3-5-mod-6-TeX-gen-A:HITS-fun:4.10
\adfDgap
%ADFvfyBlocksStart {51,51,51,51,51,51,51,51,171}
\noindent{\boldmath $ 51^{8} 171^{1} $}~
With the point set $Z_{579}$ partitioned into
 residue classes modulo $8$ for $\{0, 1, \dots, 407\}$, and
 $\{408, 409, \dots, 578\}$,
 the design is generated from

\adfLgap %ADFvfyDesignStart
$(576, 220, 231, 389)$,
$(408, 229, 386, 100)$,
$(408, 118, 36, 312)$,\adfsplit
$(408, 158, 271, 320)$,
$(408, 131, 63, 402)$,
$(408, 93, 184, 20)$,\adfsplit
$(408, 123, 215, 353)$,
$(408, 331, 197, 73)$,
$(408, 10, 126, 201)$,\adfsplit
$(409, 310, 103, 138)$,
$(409, 215, 313, 126)$,
$(409, 322, 100, 375)$,\adfsplit
$(409, 387, 89, 317)$,
$(409, 170, 336, 109)$,
$(409, 156, 59, 368)$,\adfsplit
$(409, 278, 260, 256)$,
$(409, 91, 45, 129)$,
$(410, 190, 343, 220)$,\adfsplit
$(410, 174, 200, 225)$,
$(410, 75, 60, 383)$,
$(410, 38, 183, 90)$,\adfsplit
$(410, 241, 202, 328)$,
$(410, 61, 216, 137)$,
$(410, 35, 218, 77)$,\adfsplit
$(410, 21, 115, 116)$,
$(411, 267, 193, 149)$,
$(411, 116, 318, 128)$,\adfsplit
$(411, 201, 376, 139)$,
$(411, 103, 371, 288)$,
$(411, 119, 122, 117)$,\adfsplit
$(411, 327, 300, 130)$,
$(411, 253, 190, 76)$,
$(411, 378, 89, 350)$,\adfsplit
$(412, 371, 306, 342)$,
$(412, 29, 48, 175)$,
$(412, 328, 108, 385)$,\adfsplit
$(412, 262, 181, 33)$,
$(412, 333, 135, 34)$,
$(412, 287, 172, 152)$,\adfsplit
$(412, 161, 291, 170)$,
$(412, 134, 68, 91)$,
$(413, 4, 17, 371)$,\adfsplit
$(413, 163, 372, 14)$,
$(0, 10, 55, 161)$,
$(0, 31, 265, 420)$,\adfsplit
$(0, 6, 297, 441)$,
$(0, 14, 249, 532)$,
$(0, 59, 300, 427)$,\adfsplit
$(0, 7, 348, 574)$,
$(0, 33, 189, 455)$,
$(0, 34, 139, 414)$,\adfsplit
$(0, 77, 163, 456)$,
$(0, 71, 266, 477)$,
$(0, 21, 58, 498)$,\adfsplit
$(0, 47, 154, 540)$,
$(0, 17, 243, 526)$,
$(0, 78, 203, 442)$,\adfsplit
$(0, 90, 193, 561)$,
$(0, 102, 204, 306)$

%ADFvfyBlocksEnd
\adfLgap \noindent by the mapping:
$x \mapsto x +  j \adfmod{408}$ for $x < 408$,
$x \mapsto (x - 408 + 7 j \adfmod{168}) + 408$ for $408 \le x < 576$,
$x \mapsto (x +  j \adfmod{3}) + 576$ for $x \ge 576$,
$0 \le j < 408$
 for the first 58 blocks,
$0 \le j < 102$
 for the last block.
\ADFvfyParStart{(579, ((58, 408, ((408, 1), (168, 7), (3, 1))), (1, 102, ((408, 1), (168, 7), (3, 1)))), ((51, 8), (171, 1)))} %ADFvfyParEnd
% End of 51^8 171^1
%%%%%%%%%%%%%%%%%%%%%%%%%%%%%%%%%%%%%%%%%%%%%%%%%%%%%%%%%%%%%%%%%%%%%%%%%%%%%%%%%%%%%%%%%%
%%%%%%%%%%%%%%%%%%%%%%%%%%%%%%%%%%%%%%%%%%%%%%%%%%%%%%%%%%%%%%%%%%%%%%%%%%%%%%%%%%%%%%%%%%

% Charlotte:GDD4-1-3-5-mod-6-TeX-gen-A:HITS-fun:4.10
\adfDgap
%ADFvfyBlocksStart {51,51,51,51,51,51,51,51,174}
\noindent{\boldmath $ 51^{8} 174^{1} $}~
With the point set $Z_{582}$ partitioned into
 residue classes modulo $8$ for $\{0, 1, \dots, 407\}$, and
 $\{408, 409, \dots, 581\}$,
 the design is generated from

\adfLgap %ADFvfyDesignStart
$(576, 97, 162, 179)$,
$(576, 243, 184, 386)$,
$(577, 164, 371, 306)$,\adfsplit
$(577, 91, 63, 118)$,
$(408, 183, 339, 402)$,
$(408, 133, 377, 248)$,\adfsplit
$(408, 400, 260, 251)$,
$(408, 403, 28, 143)$,
$(408, 153, 149, 222)$,\adfsplit
$(408, 261, 58, 158)$,
$(408, 265, 382, 266)$,
$(408, 264, 271, 180)$,\adfsplit
$(409, 364, 278, 203)$,
$(409, 375, 277, 122)$,
$(409, 392, 228, 225)$,\adfsplit
$(409, 64, 404, 293)$,
$(409, 298, 99, 311)$,
$(409, 406, 18, 216)$,\adfsplit
$(409, 401, 163, 366)$,
$(409, 265, 151, 357)$,
$(410, 48, 124, 234)$,\adfsplit
$(410, 61, 372, 49)$,
$(410, 211, 16, 87)$,
$(410, 272, 383, 317)$,\adfsplit
$(410, 254, 380, 155)$,
$(410, 30, 2, 247)$,
$(410, 357, 57, 315)$,\adfsplit
$(410, 178, 406, 377)$,
$(411, 365, 200, 130)$,
$(411, 216, 379, 39)$,\adfsplit
$(411, 213, 44, 343)$,
$(411, 257, 86, 364)$,
$(411, 358, 372, 3)$,\adfsplit
$(411, 208, 57, 277)$,
$(411, 198, 354, 287)$,
$(411, 169, 218, 179)$,\adfsplit
$(412, 227, 122, 197)$,
$(412, 40, 95, 28)$,
$(412, 350, 103, 261)$,\adfsplit
$(412, 105, 102, 327)$,
$(412, 235, 401, 214)$,
$(412, 258, 1, 277)$,\adfsplit
$(412, 24, 171, 300)$,
$(412, 356, 394, 296)$,
$(413, 359, 322, 312)$,\adfsplit
$(413, 69, 145, 46)$,
$(413, 63, 57, 13)$,
$(413, 194, 67, 7)$,\adfsplit
$(413, 101, 76, 315)$,
$(413, 398, 112, 162)$,
$(413, 84, 209, 176)$,\adfsplit
$(413, 188, 83, 126)$,
$(414, 353, 156, 407)$,
$(414, 293, 370, 364)$,\adfsplit
$(414, 393, 320, 20)$,
$(414, 87, 336, 222)$,
$(414, 91, 69, 266)$,\adfsplit
$(414, 38, 275, 13)$,
$(414, 18, 166, 112)$,
$(414, 271, 1, 219)$,\adfsplit
$(415, 102, 195, 376)$,
$(415, 186, 335, 107)$,
$(415, 44, 49, 127)$,\adfsplit
$(415, 19, 76, 166)$,
$(415, 13, 33, 264)$,
$(415, 399, 398, 12)$,\adfsplit
$(415, 389, 104, 338)$,
$(415, 322, 185, 45)$,
$(416, 320, 262, 196)$,\adfsplit
$(416, 336, 233, 251)$,
$(416, 253, 86, 195)$,
$(416, 393, 287, 306)$,\adfsplit
$(416, 211, 337, 88)$,
$(416, 255, 380, 29)$,
$(416, 36, 10, 357)$,\adfsplit
$(416, 338, 150, 31)$,
$(417, 195, 74, 85)$,
$(417, 319, 177, 200)$,\adfsplit
$(417, 238, 202, 204)$,
$(417, 116, 131, 222)$,
$(417, 137, 135, 101)$,\adfsplit
$(417, 232, 405, 311)$,
$(417, 235, 374, 73)$,
$(417, 168, 124, 306)$,\adfsplit
$(418, 57, 347, 284)$,
$(418, 256, 39, 1)$,
$(418, 53, 387, 326)$,\adfsplit
$(418, 100, 262, 104)$,
$(418, 2, 405, 264)$,
$(418, 343, 198, 394)$,\adfsplit
$(418, 276, 43, 258)$,
$(418, 325, 17, 191)$,
$(419, 41, 134, 277)$,\adfsplit
$(419, 234, 151, 316)$,
$(419, 366, 10, 53)$,
$(0, 49, 242, 433)$,\adfsplit
$(0, 9, 327, 503)$,
$(0, 29, 150, 447)$,
$(0, 179, 295, 531)$,\adfsplit
$(0, 113, 139, 475)$,
$(0, 77, 275, 420)$,
$(0, 31, 377, 421)$,\adfsplit
$(0, 46, 313, 434)$,
$(0, 27, 74, 449)$,
$(0, 78, 263, 519)$,\adfsplit
$(0, 133, 255, 491)$,
$(0, 30, 401, 575)$,
$(0, 118, 395, 560)$,\adfsplit
$(0, 117, 214, 476)$,
$(0, 41, 363, 463)$,
$(0, 170, 393, 561)$,\adfsplit
$(0, 219, 397, 505)$,
$(0, 154, 355, 532)$,
$(0, 137, 178, 391)$,\adfsplit
$(0, 42, 349, 504)$,
$(0, 81, 127, 518)$,
$(1, 15, 339, 546)$,\adfsplit
$(0, 102, 204, 306)$,
$(1, 103, 205, 307)$

%ADFvfyBlocksEnd
\adfLgap \noindent by the mapping:
$x \mapsto x + 2 j \adfmod{408}$ for $x < 408$,
$x \mapsto (x - 408 + 14 j \adfmod{168}) + 408$ for $408 \le x < 576$,
$x \mapsto (x + 2 j \adfmod{6}) + 576$ for $x \ge 576$,
$0 \le j < 204$
 for the first 117 blocks,
$0 \le j < 51$
 for the last two blocks.
\ADFvfyParStart{(582, ((117, 204, ((408, 2), (168, 14), (6, 2))), (2, 51, ((408, 2), (168, 14), (6, 2)))), ((51, 8), (174, 1)))} %ADFvfyParEnd
% End of 51^8 174^1
%%%%%%%%%%%%%%%%%%%%%%%%%%%%%%%%%%%%%%%%%%%%%%%%%%%%%%%%%%%%%%%%%%%%%%%%%%%%%%%%%%%%%%%%%%
%%%%%%%%%%%%%%%%%%%%%%%%%%%%%%%%%%%%%%%%%%%%%%%%%%%%%%%%%%%%%%%%%%%%%%%%%%%%%%%%%%%%%%%%%%

%%%%%%%%%%%%%%%%%%%%%%%%%%%%%%%%%%%%%%%%%%%%%%%%%%%%%%%%%%%%%%%%%%%%%%%%%%%%%%%%%%%%%%%%%%
%%%%%%%%%%%%%%%%%%%%%%%%%%%%%%%%%%%%%%%%%%%%%%%%%%%%%%%%%%%%%%%%%%%%%%%%%%%%%%%%%%%%%%%%%%
\section{4-GDDs for the proof of Lemma \ref{lem:4-GDD 57^8 m^1}}
\label{app:4-GDD 57^8 m^1}
\adfnull{
$ 57^8 174^1 $,
$ 57^8 177^1 $,
$ 57^8 180^1 $,
$ 57^8 183^1 $,
$ 57^8 186^1 $,
$ 57^8 189^1 $,
$ 57^8 192^1 $ and
$ 57^8 195^1 $.
}

% Charlotte:GDD4-1-3-5-mod-6-TeX-gen-A:HITS-fun:4.10
\adfDgap
%ADFvfyBlocksStart {57,57,57,57,57,57,57,57,174}
\noindent{\boldmath $ 57^{8} 174^{1} $}~
With the point set $Z_{630}$ partitioned into
 residue classes modulo $8$ for $\{0, 1, \dots, 455\}$, and
 $\{456, 457, \dots, 629\}$,
 the design is generated from

\adfLgap %ADFvfyDesignStart
$(624, 63, 94, 221)$,
$(625, 104, 390, 403)$,
$(456, 412, 90, 441)$,\adfsplit
$(456, 191, 98, 45)$,
$(456, 355, 230, 145)$,
$(456, 411, 406, 365)$,\adfsplit
$(456, 328, 78, 305)$,
$(456, 183, 301, 323)$,
$(456, 199, 20, 408)$,\adfsplit
$(456, 348, 152, 442)$,
$(457, 434, 385, 440)$,
$(457, 419, 113, 375)$,\adfsplit
$(457, 333, 256, 147)$,
$(457, 190, 383, 441)$,
$(457, 221, 420, 342)$,\adfsplit
$(457, 68, 307, 250)$,
$(457, 302, 268, 79)$,
$(457, 312, 138, 229)$,\adfsplit
$(458, 411, 94, 450)$,
$(458, 279, 326, 164)$,
$(458, 112, 129, 149)$,\adfsplit
$(458, 268, 449, 223)$,
$(458, 432, 443, 49)$,
$(458, 130, 325, 438)$,\adfsplit
$(458, 119, 170, 331)$,
$(458, 189, 372, 344)$,
$(459, 353, 454, 61)$,\adfsplit
$(459, 111, 1, 146)$,
$(459, 380, 237, 195)$,
$(459, 154, 249, 276)$,\adfsplit
$(459, 48, 211, 127)$,
$(459, 167, 340, 232)$,
$(459, 422, 330, 77)$,\adfsplit
$(459, 342, 152, 299)$,
$(460, 294, 303, 384)$,
$(460, 10, 367, 433)$,\adfsplit
$(460, 290, 86, 12)$,
$(460, 263, 160, 41)$,
$(460, 253, 147, 234)$,\adfsplit
$(460, 115, 165, 356)$,
$(460, 33, 347, 317)$,
$(460, 416, 148, 286)$,\adfsplit
$(461, 123, 399, 41)$,
$(461, 374, 241, 367)$,
$(461, 215, 384, 213)$,\adfsplit
$(461, 109, 176, 378)$,
$(461, 124, 413, 451)$,
$(461, 418, 94, 332)$,\adfsplit
$(0, 1, 221, 552)$,
$(0, 4, 319, 615)$,
$(0, 10, 397, 553)$,\adfsplit
$(0, 25, 61, 567)$,
$(0, 52, 259, 497)$,
$(0, 107, 332, 588)$,\adfsplit
$(0, 26, 97, 609)$,
$(0, 14, 116, 131)$,
$(0, 3, 21, 258)$,\adfsplit
$(0, 54, 297, 581)$,
$(0, 12, 135, 291)$,
$(0, 75, 380, 546)$,\adfsplit
$(0, 60, 149, 302)$,
$(0, 70, 245, 490)$,
$(0, 114, 228, 342)$

%ADFvfyBlocksEnd
\adfLgap \noindent by the mapping:
$x \mapsto x +  j \adfmod{456}$ for $x < 456$,
$x \mapsto (x - 456 + 7 j \adfmod{168}) + 456$ for $456 \le x < 624$,
$x \mapsto (x + 2 j \adfmod{6}) + 624$ for $x \ge 624$,
$0 \le j < 456$
 for the first 62 blocks,
$0 \le j < 114$
 for the last block.
\ADFvfyParStart{(630, ((62, 456, ((456, 1), (168, 7), (6, 2))), (1, 114, ((456, 1), (168, 7), (6, 2)))), ((57, 8), (174, 1)))} %ADFvfyParEnd
% End of 57^8 174^1
%%%%%%%%%%%%%%%%%%%%%%%%%%%%%%%%%%%%%%%%%%%%%%%%%%%%%%%%%%%%%%%%%%%%%%%%%%%%%%%%%%%%%%%%%%
%%%%%%%%%%%%%%%%%%%%%%%%%%%%%%%%%%%%%%%%%%%%%%%%%%%%%%%%%%%%%%%%%%%%%%%%%%%%%%%%%%%%%%%%%%

% Charlotte:GDD4-1-3-5-mod-6-TeX-gen-A:HITS-fun:4.10
\adfDgap
%ADFvfyBlocksStart {57,57,57,57,57,57,57,57,177}
\noindent{\boldmath $ 57^{8} 177^{1} $}~
With the point set $Z_{633}$ partitioned into
 residue classes modulo $8$ for $\{0, 1, \dots, 455\}$, and
 $\{456, 457, \dots, 632\}$,
 the design is generated from

\adfLgap %ADFvfyDesignStart
$(624, 312, 41, 214)$,
$(624, 213, 302, 271)$,
$(625, 154, 53, 435)$,\adfsplit
$(625, 386, 379, 294)$,
$(626, 208, 385, 180)$,
$(626, 3, 137, 446)$,\adfsplit
$(456, 172, 283, 365)$,
$(456, 171, 254, 191)$,
$(456, 224, 205, 265)$,\adfsplit
$(456, 18, 111, 54)$,
$(456, 214, 34, 179)$,
$(456, 242, 424, 209)$,\adfsplit
$(456, 21, 308, 153)$,
$(456, 360, 36, 7)$,
$(457, 101, 346, 107)$,\adfsplit
$(457, 160, 309, 73)$,
$(457, 307, 268, 302)$,
$(457, 157, 249, 27)$,\adfsplit
$(457, 247, 94, 170)$,
$(457, 320, 87, 18)$,
$(457, 71, 324, 294)$,\adfsplit
$(457, 113, 260, 24)$,
$(458, 227, 20, 274)$,
$(458, 432, 381, 354)$,\adfsplit
$(458, 438, 152, 335)$,
$(458, 427, 160, 172)$,
$(458, 137, 252, 230)$,\adfsplit
$(458, 77, 51, 385)$,
$(458, 271, 201, 382)$,
$(458, 386, 109, 111)$,\adfsplit
$(459, 356, 283, 86)$,
$(459, 453, 46, 132)$,
$(459, 96, 322, 389)$,\adfsplit
$(459, 447, 258, 217)$,
$(459, 208, 81, 383)$,
$(459, 296, 387, 2)$,\adfsplit
$(459, 209, 174, 295)$,
$(459, 13, 4, 203)$,
$(460, 165, 360, 92)$,\adfsplit
$(460, 109, 395, 370)$,
$(460, 171, 62, 287)$,
$(460, 32, 132, 55)$,\adfsplit
$(460, 265, 426, 220)$,
$(460, 399, 377, 328)$,
$(460, 283, 342, 245)$,\adfsplit
$(460, 201, 74, 118)$,
$(461, 104, 217, 142)$,
$(461, 24, 19, 438)$,\adfsplit
$(461, 228, 75, 225)$,
$(461, 134, 89, 35)$,
$(461, 197, 424, 359)$,\adfsplit
$(461, 434, 327, 172)$,
$(461, 13, 426, 260)$,
$(461, 295, 178, 285)$,\adfsplit
$(462, 51, 88, 333)$,
$(462, 292, 37, 118)$,
$(462, 144, 409, 231)$,\adfsplit
$(462, 156, 38, 202)$,
$(462, 283, 247, 368)$,
$(462, 347, 153, 452)$,\adfsplit
$(462, 269, 330, 359)$,
$(462, 390, 185, 410)$,
$(463, 348, 2, 325)$,\adfsplit
$(463, 80, 393, 390)$,
$(463, 361, 452, 15)$,
$(463, 52, 311, 58)$,\adfsplit
$(463, 65, 315, 384)$,
$(463, 214, 16, 5)$,
$(463, 187, 186, 45)$,\adfsplit
$(463, 371, 446, 367)$,
$(464, 301, 329, 115)$,
$(464, 314, 158, 212)$,\adfsplit
$(464, 280, 276, 447)$,
$(464, 439, 358, 345)$,
$(464, 296, 244, 154)$,\adfsplit
$(464, 21, 264, 275)$,
$(464, 243, 119, 450)$,
$(464, 173, 361, 414)$,\adfsplit
$(465, 153, 320, 319)$,
$(465, 405, 136, 6)$,
$(465, 192, 207, 61)$,\adfsplit
$(465, 89, 426, 214)$,
$(465, 265, 410, 53)$,
$(465, 284, 291, 335)$,\adfsplit
$(465, 418, 59, 196)$,
$(465, 134, 307, 228)$,
$(466, 13, 81, 346)$,\adfsplit
$(466, 339, 280, 121)$,
$(466, 108, 379, 406)$,
$(466, 52, 122, 395)$,\adfsplit
$(466, 152, 30, 317)$,
$(466, 134, 295, 189)$,
$(466, 159, 336, 18)$,\adfsplit
$(466, 449, 95, 116)$,
$(467, 297, 76, 62)$,
$(467, 432, 92, 414)$,\adfsplit
$(467, 146, 91, 103)$,
$(467, 286, 69, 147)$,
$(467, 421, 359, 82)$,\adfsplit
$(467, 203, 221, 137)$,
$(467, 361, 375, 232)$,
$(467, 248, 258, 108)$,\adfsplit
$(468, 168, 62, 252)$,
$(468, 167, 29, 443)$,
$(468, 414, 91, 106)$,\adfsplit
$(0, 31, 139, 469)$,
$(0, 21, 214, 496)$,
$(0, 26, 374, 497)$,\adfsplit
$(0, 60, 238, 525)$,
$(0, 65, 74, 237)$,
$(0, 2, 68, 126)$,\adfsplit
$(0, 53, 105, 172)$,
$(0, 62, 305, 580)$,
$(0, 115, 431, 623)$,\adfsplit
$(0, 187, 439, 622)$,
$(0, 47, 196, 246)$,
$(0, 61, 219, 417)$,\adfsplit
$(0, 33, 131, 608)$,
$(0, 327, 361, 553)$,
$(0, 159, 285, 581)$,\adfsplit
$(0, 17, 135, 291)$,
$(0, 63, 227, 524)$,
$(0, 95, 252, 511)$,\adfsplit
$(0, 101, 151, 347)$,
$(1, 31, 411, 581)$,
$(0, 114, 228, 342)$,\adfsplit
$(1, 115, 229, 343)$

%ADFvfyBlocksEnd
\adfLgap \noindent by the mapping:
$x \mapsto x + 2 j \adfmod{456}$ for $x < 456$,
$x \mapsto (x - 456 + 14 j \adfmod{168}) + 456$ for $456 \le x < 624$,
$x \mapsto (x - 624 + 3 j \adfmod{9}) + 624$ for $x \ge 624$,
$0 \le j < 228$
 for the first 125 blocks,
$0 \le j < 57$
 for the last two blocks.
\ADFvfyParStart{(633, ((125, 228, ((456, 2), (168, 14), (9, 3))), (2, 57, ((456, 2), (168, 14), (9, 3)))), ((57, 8), (177, 1)))} %ADFvfyParEnd
% End of 57^8 177^1
%%%%%%%%%%%%%%%%%%%%%%%%%%%%%%%%%%%%%%%%%%%%%%%%%%%%%%%%%%%%%%%%%%%%%%%%%%%%%%%%%%%%%%%%%%
%%%%%%%%%%%%%%%%%%%%%%%%%%%%%%%%%%%%%%%%%%%%%%%%%%%%%%%%%%%%%%%%%%%%%%%%%%%%%%%%%%%%%%%%%%

% Charlotte:GDD4-1-3-5-mod-6-TeX-gen-A:HITS-fun:4.10
\adfDgap
%ADFvfyBlocksStart {57,57,57,57,57,57,57,57,180}
\noindent{\boldmath $ 57^{8} 180^{1} $}~
With the point set $Z_{636}$ partitioned into
 residue classes modulo $8$ for $\{0, 1, \dots, 455\}$, and
 $\{456, 457, \dots, 635\}$,
 the design is generated from

\adfLgap %ADFvfyDesignStart
$(624, 93, 143, 244)$,
$(625, 125, 202, 294)$,
$(626, 92, 304, 246)$,\adfsplit
$(627, 391, 342, 386)$,
$(456, 357, 201, 63)$,
$(456, 187, 230, 364)$,\adfsplit
$(456, 112, 332, 293)$,
$(456, 83, 402, 32)$,
$(456, 58, 233, 223)$,\adfsplit
$(456, 180, 6, 13)$,
$(456, 239, 310, 73)$,
$(456, 242, 360, 363)$,\adfsplit
$(457, 77, 433, 274)$,
$(457, 421, 105, 72)$,
$(457, 52, 99, 270)$,\adfsplit
$(457, 450, 176, 204)$,
$(457, 327, 237, 2)$,
$(457, 151, 308, 160)$,\adfsplit
$(457, 430, 59, 65)$,
$(457, 355, 71, 86)$,
$(458, 25, 387, 164)$,\adfsplit
$(458, 5, 438, 47)$,
$(458, 235, 254, 237)$,
$(458, 416, 273, 255)$,\adfsplit
$(458, 280, 397, 170)$,
$(458, 251, 196, 450)$,
$(458, 48, 295, 226)$,\adfsplit
$(458, 17, 166, 180)$,
$(459, 151, 346, 292)$,
$(459, 349, 313, 324)$,\adfsplit
$(459, 427, 224, 39)$,
$(459, 422, 112, 209)$,
$(459, 167, 357, 282)$,\adfsplit
$(459, 92, 275, 362)$,
$(459, 54, 147, 317)$,
$(459, 360, 238, 33)$,\adfsplit
$(460, 69, 275, 10)$,
$(460, 314, 190, 116)$,
$(460, 282, 183, 229)$,\adfsplit
$(460, 48, 151, 121)$,
$(460, 200, 180, 102)$,
$(460, 424, 383, 29)$,\adfsplit
$(460, 449, 187, 412)$,
$(460, 123, 249, 278)$,
$(461, 20, 153, 259)$,\adfsplit
$(461, 365, 192, 62)$,
$(461, 333, 4, 183)$,
$(461, 150, 226, 195)$,\adfsplit
$(0, 1, 394, 609)$,
$(0, 26, 324, 490)$,
$(0, 52, 119, 623)$,\adfsplit
$(0, 12, 125, 321)$,
$(0, 105, 347, 580)$,
$(0, 81, 373, 468)$,\adfsplit
$(0, 4, 211, 475)$,
$(0, 22, 252, 538)$,
$(0, 60, 374, 567)$,\adfsplit
$(0, 27, 84, 581)$,
$(0, 38, 108, 616)$,
$(0, 21, 116, 525)$,\adfsplit
$(0, 13, 79, 268)$,
$(0, 35, 276, 588)$,
$(0, 34, 123, 345)$,\adfsplit
$(0, 114, 228, 342)$

%ADFvfyBlocksEnd
\adfLgap \noindent by the mapping:
$x \mapsto x +  j \adfmod{456}$ for $x < 456$,
$x \mapsto (x - 456 + 7 j \adfmod{168}) + 456$ for $456 \le x < 624$,
$x \mapsto (x + 4 j \adfmod{12}) + 624$ for $x \ge 624$,
$0 \le j < 456$
 for the first 63 blocks,
$0 \le j < 114$
 for the last block.
\ADFvfyParStart{(636, ((63, 456, ((456, 1), (168, 7), (12, 4))), (1, 114, ((456, 1), (168, 7), (12, 4)))), ((57, 8), (180, 1)))} %ADFvfyParEnd
% End of 57^8 180^1
%%%%%%%%%%%%%%%%%%%%%%%%%%%%%%%%%%%%%%%%%%%%%%%%%%%%%%%%%%%%%%%%%%%%%%%%%%%%%%%%%%%%%%%%%%
%%%%%%%%%%%%%%%%%%%%%%%%%%%%%%%%%%%%%%%%%%%%%%%%%%%%%%%%%%%%%%%%%%%%%%%%%%%%%%%%%%%%%%%%%%

% Charlotte:GDD4-1-3-5-mod-6-TeX-gen-A:HITS-fun:4.10
\adfDgap
%ADFvfyBlocksStart {57,57,57,57,57,57,57,57,183}
\noindent{\boldmath $ 57^{8} 183^{1} $}~
With the point set $Z_{639}$ partitioned into
 residue classes modulo $8$ for $\{0, 1, \dots, 455\}$, and
 $\{456, 457, \dots, 638\}$,
 the design is generated from

\adfLgap %ADFvfyDesignStart
$(624, 359, 96, 433)$,
$(624, 194, 87, 22)$,
$(625, 84, 106, 85)$,\adfsplit
$(625, 321, 104, 269)$,
$(626, 143, 24, 283)$,
$(626, 254, 357, 154)$,\adfsplit
$(627, 193, 150, 364)$,
$(627, 437, 392, 129)$,
$(628, 276, 145, 322)$,\adfsplit
$(628, 2, 207, 257)$,
$(456, 226, 102, 171)$,
$(456, 177, 197, 307)$,\adfsplit
$(456, 375, 155, 84)$,
$(456, 286, 61, 450)$,
$(456, 31, 233, 392)$,\adfsplit
$(456, 422, 16, 52)$,
$(456, 313, 92, 335)$,
$(456, 362, 72, 69)$,\adfsplit
$(457, 110, 80, 453)$,
$(457, 447, 266, 33)$,
$(457, 305, 92, 282)$,\adfsplit
$(457, 294, 173, 240)$,
$(457, 439, 99, 157)$,
$(457, 348, 136, 383)$,\adfsplit
$(457, 274, 145, 443)$,
$(457, 283, 46, 148)$,
$(458, 239, 113, 427)$,\adfsplit
$(458, 408, 330, 63)$,
$(458, 412, 157, 26)$,
$(458, 101, 382, 271)$,\adfsplit
$(458, 140, 80, 446)$,
$(458, 261, 222, 84)$,
$(458, 424, 27, 25)$,\adfsplit
$(458, 393, 202, 227)$,
$(459, 297, 100, 403)$,
$(459, 161, 378, 23)$,\adfsplit
$(459, 339, 46, 200)$,
$(459, 317, 422, 228)$,
$(459, 6, 399, 217)$,\adfsplit
$(459, 16, 371, 298)$,
$(459, 48, 205, 164)$,
$(459, 165, 79, 386)$,\adfsplit
$(460, 242, 127, 184)$,
$(460, 306, 203, 272)$,
$(460, 353, 92, 166)$,\adfsplit
$(460, 333, 423, 134)$,
$(460, 95, 270, 325)$,
$(460, 241, 269, 276)$,\adfsplit
$(460, 363, 216, 442)$,
$(460, 57, 388, 19)$,
$(461, 382, 380, 387)$,\adfsplit
$(461, 103, 230, 240)$,
$(461, 165, 263, 226)$,
$(461, 417, 12, 231)$,\adfsplit
$(461, 272, 115, 125)$,
$(461, 253, 193, 378)$,
$(461, 41, 222, 74)$,\adfsplit
$(461, 172, 184, 107)$,
$(462, 43, 415, 421)$,
$(462, 245, 249, 116)$,\adfsplit
$(462, 112, 385, 285)$,
$(462, 32, 131, 230)$,
$(462, 450, 216, 436)$,\adfsplit
$(462, 99, 290, 238)$,
$(462, 431, 204, 449)$,
$(462, 270, 178, 255)$,\adfsplit
$(463, 181, 66, 145)$,
$(463, 327, 75, 281)$,
$(463, 347, 392, 78)$,\adfsplit
$(463, 190, 201, 218)$,
$(463, 101, 4, 280)$,
$(463, 453, 360, 215)$,\adfsplit
$(463, 404, 163, 298)$,
$(463, 175, 180, 398)$,
$(464, 0, 346, 323)$,\adfsplit
$(464, 454, 19, 7)$,
$(464, 363, 455, 436)$,
$(464, 423, 181, 49)$,\adfsplit
$(464, 425, 258, 342)$,
$(464, 352, 201, 356)$,
$(464, 29, 410, 228)$,\adfsplit
$(464, 182, 368, 309)$,
$(465, 439, 179, 70)$,
$(465, 137, 20, 327)$,\adfsplit
$(465, 384, 226, 302)$,
$(465, 124, 242, 153)$,
$(465, 42, 389, 136)$,\adfsplit
$(465, 78, 215, 37)$,
$(465, 51, 145, 309)$,
$(465, 344, 300, 451)$,\adfsplit
$(466, 0, 295, 348)$,
$(466, 328, 310, 359)$,
$(466, 45, 438, 425)$,\adfsplit
$(466, 33, 80, 419)$,
$(466, 265, 373, 163)$,
$(466, 114, 134, 363)$,\adfsplit
$(466, 159, 106, 436)$,
$(466, 194, 188, 221)$,
$(467, 428, 143, 240)$,\adfsplit
$(467, 396, 147, 181)$,
$(467, 111, 405, 406)$,
$(467, 34, 179, 320)$,\adfsplit
$(467, 151, 398, 305)$,
$(467, 197, 194, 340)$,
$(467, 138, 270, 321)$,\adfsplit
$(467, 451, 328, 25)$,
$(468, 294, 172, 219)$,
$(0, 17, 206, 468)$,\adfsplit
$(0, 68, 246, 539)$,
$(0, 38, 371, 552)$,
$(0, 13, 425, 469)$,\adfsplit
$(0, 26, 326, 510)$,
$(0, 66, 431, 483)$,
$(0, 121, 294, 595)$,\adfsplit
$(0, 15, 42, 358)$,
$(0, 185, 196, 447)$,
$(0, 155, 204, 566)$,\adfsplit
$(0, 61, 195, 524)$,
$(0, 85, 417, 623)$,
$(1, 67, 223, 594)$,\adfsplit
$(0, 81, 143, 419)$,
$(1, 55, 123, 580)$,
$(0, 91, 105, 134)$,\adfsplit
$(0, 159, 254, 609)$,
$(0, 62, 287, 437)$,
$(0, 141, 385, 581)$,\adfsplit
$(1, 27, 173, 483)$,
$(0, 114, 228, 342)$,
$(1, 115, 229, 343)$

%ADFvfyBlocksEnd
\adfLgap \noindent by the mapping:
$x \mapsto x + 2 j \adfmod{456}$ for $x < 456$,
$x \mapsto (x - 456 + 14 j \adfmod{168}) + 456$ for $456 \le x < 624$,
$x \mapsto (x - 624 + 5 j \adfmod{15}) + 624$ for $x \ge 624$,
$0 \le j < 228$
 for the first 127 blocks,
$0 \le j < 57$
 for the last two blocks.
\ADFvfyParStart{(639, ((127, 228, ((456, 2), (168, 14), (15, 5))), (2, 57, ((456, 2), (168, 14), (15, 5)))), ((57, 8), (183, 1)))} %ADFvfyParEnd
% End of 57^8 183^1
%%%%%%%%%%%%%%%%%%%%%%%%%%%%%%%%%%%%%%%%%%%%%%%%%%%%%%%%%%%%%%%%%%%%%%%%%%%%%%%%%%%%%%%%%%
%%%%%%%%%%%%%%%%%%%%%%%%%%%%%%%%%%%%%%%%%%%%%%%%%%%%%%%%%%%%%%%%%%%%%%%%%%%%%%%%%%%%%%%%%%

% Charlotte:GDD4-1-3-5-mod-6-TeX-gen-A:HITS-fun:4.10
\adfDgap
%ADFvfyBlocksStart {57,57,57,57,57,57,57,57,186}
\noindent{\boldmath $ 57^{8} 186^{1} $}~
With the point set $Z_{642}$ partitioned into
 residue classes modulo $8$ for $\{0, 1, \dots, 455\}$, and
 $\{456, 457, \dots, 641\}$,
 the design is generated from

\adfLgap %ADFvfyDesignStart
$(624, 75, 236, 1)$,
$(625, 98, 312, 412)$,
$(626, 407, 348, 349)$,\adfsplit
$(627, 260, 205, 57)$,
$(628, 207, 361, 251)$,
$(629, 353, 432, 349)$,\adfsplit
$(456, 235, 120, 157)$,
$(456, 439, 261, 390)$,
$(456, 183, 365, 406)$,\adfsplit
$(456, 58, 302, 409)$,
$(456, 100, 227, 2)$,
$(456, 328, 41, 243)$,\adfsplit
$(456, 129, 320, 60)$,
$(456, 236, 191, 114)$,
$(457, 266, 100, 219)$,\adfsplit
$(457, 226, 33, 139)$,
$(457, 325, 127, 84)$,
$(457, 382, 203, 72)$,\adfsplit
$(457, 20, 32, 197)$,
$(457, 354, 169, 167)$,
$(457, 232, 222, 165)$,\adfsplit
$(457, 230, 305, 351)$,
$(458, 333, 436, 3)$,
$(458, 259, 325, 375)$,\adfsplit
$(458, 160, 299, 174)$,
$(458, 201, 300, 170)$,
$(458, 199, 289, 118)$,\adfsplit
$(458, 8, 401, 164)$,
$(458, 335, 154, 192)$,
$(458, 389, 42, 230)$,\adfsplit
$(459, 304, 59, 285)$,
$(459, 23, 236, 393)$,
$(459, 217, 411, 128)$,\adfsplit
$(459, 148, 288, 139)$,
$(459, 330, 12, 414)$,
$(459, 230, 98, 245)$,\adfsplit
$(459, 394, 61, 39)$,
$(459, 103, 142, 401)$,
$(460, 439, 85, 155)$,\adfsplit
$(460, 260, 234, 285)$,
$(460, 455, 81, 248)$,
$(460, 318, 173, 144)$,\adfsplit
$(460, 268, 257, 310)$,
$(460, 280, 348, 139)$,
$(460, 15, 10, 267)$,\adfsplit
$(460, 26, 97, 62)$,
$(461, 216, 237, 334)$,
$(461, 371, 196, 105)$,\adfsplit
$(461, 259, 232, 49)$,
$(0, 52, 395, 560)$,
$(0, 20, 383, 595)$,\adfsplit
$(0, 65, 303, 510)$,
$(0, 92, 255, 490)$,
$(0, 30, 180, 531)$,\adfsplit
$(0, 34, 220, 608)$,
$(0, 6, 449, 587)$,
$(0, 117, 251, 503)$,\adfsplit
$(0, 3, 348, 546)$,
$(0, 137, 292, 567)$,
$(0, 124, 286, 497)$,\adfsplit
$(0, 60, 321, 504)$,
$(0, 18, 94, 323)$,
$(0, 33, 394, 476)$,\adfsplit
$(0, 17, 206, 234)$,
$(0, 114, 228, 342)$

%ADFvfyBlocksEnd
\adfLgap \noindent by the mapping:
$x \mapsto x +  j \adfmod{456}$ for $x < 456$,
$x \mapsto (x - 456 + 7 j \adfmod{168}) + 456$ for $456 \le x < 624$,
$x \mapsto (x - 624 + 6 j \adfmod{18}) + 624$ for $x \ge 624$,
$0 \le j < 456$
 for the first 64 blocks,
$0 \le j < 114$
 for the last block.
\ADFvfyParStart{(642, ((64, 456, ((456, 1), (168, 7), (18, 6))), (1, 114, ((456, 1), (168, 7), (18, 6)))), ((57, 8), (186, 1)))} %ADFvfyParEnd
% End of 57^8 186^1
%%%%%%%%%%%%%%%%%%%%%%%%%%%%%%%%%%%%%%%%%%%%%%%%%%%%%%%%%%%%%%%%%%%%%%%%%%%%%%%%%%%%%%%%%%
%%%%%%%%%%%%%%%%%%%%%%%%%%%%%%%%%%%%%%%%%%%%%%%%%%%%%%%%%%%%%%%%%%%%%%%%%%%%%%%%%%%%%%%%%%

% Charlotte:GDD4-1-3-5-mod-6-TeX-gen-A:HITS-fun:4.10
\adfDgap
%ADFvfyBlocksStart {57,57,57,57,57,57,57,57,189}
\noindent{\boldmath $ 57^{8} 189^{1} $}~
With the point set $Z_{645}$ partitioned into
 residue classes modulo $7$ for $\{0, 1, \dots, 398\}$,
 $\{399, 400, \dots, 455\}$, and
 $\{456, 457, \dots, 644\}$,
 the design is generated from

\adfLgap %ADFvfyDesignStart
$(37, 306, 309, 573)$,
$(10, 238, 288, 565)$,
$(86, 157, 334, 531)$,\adfsplit
$(0, 249, 388, 483)$,
$(11, 350, 450, 501)$,
$(92, 240, 259, 475)$,\adfsplit
$(139, 184, 417, 548)$,
$(2, 375, 448, 626)$,
$(360, 385, 437, 572)$,\adfsplit
$(41, 61, 331, 586)$,
$(135, 270, 346, 495)$,
$(112, 129, 411, 534)$,\adfsplit
$(7, 130, 243, 609)$,
$(175, 369, 393, 458)$,
$(242, 330, 440, 463)$,\adfsplit
$(23, 325, 448, 469)$,
$(242, 386, 415, 520)$,
$(17, 358, 437, 605)$,\adfsplit
$(42, 363, 452, 590)$,
$(83, 122, 205, 495)$,
$(160, 274, 312, 598)$,\adfsplit
$(119, 121, 321, 572)$,
$(48, 193, 241, 479)$,
$(87, 259, 387, 462)$,\adfsplit
$(84, 222, 379, 456)$,
$(208, 357, 370, 580)$,
$(66, 67, 407, 514)$,\adfsplit
$(210, 265, 375, 602)$,
$(52, 261, 265, 628)$,
$(278, 378, 384, 587)$,\adfsplit
$(48, 120, 454, 472)$,
$(164, 201, 223, 479)$,
$(144, 233, 264, 473)$,\adfsplit
$(227, 263, 268, 531)$,
$(71, 291, 324, 586)$,
$(90, 260, 355, 498)$,\adfsplit
$(1, 263, 327, 581)$,
$(89, 329, 455, 619)$,
$(153, 245, 255, 464)$,\adfsplit
$(146, 239, 271, 526)$,
$(91, 358, 421, 595)$,
$(38, 106, 341, 498)$,\adfsplit
$(66, 153, 354, 582)$,
$(248, 338, 403, 603)$,
$(8, 69, 229, 511)$,\adfsplit
$(176, 384, 414, 549)$,
$(215, 244, 297, 586)$,
$(5, 280, 288, 626)$,\adfsplit
$(40, 83, 405, 458)$,
$(36, 252, 355, 466)$,
$(61, 329, 398, 629)$,\adfsplit
$(56, 130, 388, 506)$,
$(28, 373, 409, 485)$,
$(49, 340, 435, 497)$,\adfsplit
$(0, 85, 427, 565)$,
$(0, 15, 101, 521)$,
$(0, 34, 246, 494)$,\adfsplit
$(0, 27, 79, 577)$,
$(0, 46, 204, 496)$,
$(0, 143, 402, 612)$,\adfsplit
$(0, 40, 284, 468)$,
$(0, 44, 434, 506)$,
$(0, 18, 65, 497)$,\adfsplit
$(0, 57, 233, 634)$,
$(0, 107, 243, 490)$,
$(0, 174, 417, 562)$,\adfsplit
$(0, 12, 226, 498)$,
$(0, 23, 117, 192)$,
$(0, 142, 403, 606)$,\adfsplit
$(0, 16, 416, 474)$,
$(0, 30, 348, 418)$,
$(0, 66, 184, 552)$,\adfsplit
$(0, 9, 450, 492)$,
$(0, 180, 433, 480)$

%ADFvfyBlocksEnd
\adfLgap \noindent by the mapping:
$x \mapsto x +  j \adfmod{399}$ for $x < 399$,
$x \mapsto (x +  j \adfmod{57}) + 399$ for $399 \le x < 456$,
$x \mapsto (x - 456 + 9 j \adfmod{189}) + 456$ for $x \ge 456$,
$0 \le j < 399$.
\ADFvfyParStart{(645, ((74, 399, ((399, 1), (57, 1), (189, 9)))), ((57, 7), (57, 1), (189, 1)))} %ADFvfyParEnd
% End of 57^8 189^1
%%%%%%%%%%%%%%%%%%%%%%%%%%%%%%%%%%%%%%%%%%%%%%%%%%%%%%%%%%%%%%%%%%%%%%%%%%%%%%%%%%%%%%%%%%
%%%%%%%%%%%%%%%%%%%%%%%%%%%%%%%%%%%%%%%%%%%%%%%%%%%%%%%%%%%%%%%%%%%%%%%%%%%%%%%%%%%%%%%%%%

% Charlotte:GDD4-1-3-5-mod-6-TeX-gen-A:HITS-fun:4.10
\adfDgap
%ADFvfyBlocksStart {57,57,57,57,57,57,57,57,192}
\noindent{\boldmath $ 57^{8} 192^{1} $}~
With the point set $Z_{648}$ partitioned into
 residue classes modulo $8$ for $\{0, 1, \dots, 455\}$, and
 $\{456, 457, \dots, 647\}$,
 the design is generated from

\adfLgap %ADFvfyDesignStart
$(456, 357, 194, 368)$,
$(456, 347, 144, 143)$,
$(456, 328, 42, 67)$,\adfsplit
$(456, 385, 10, 156)$,
$(456, 366, 75, 172)$,
$(456, 401, 310, 103)$,\adfsplit
$(456, 380, 269, 201)$,
$(456, 325, 182, 255)$,
$(457, 94, 221, 132)$,\adfsplit
$(457, 78, 284, 128)$,
$(457, 181, 218, 347)$,
$(457, 385, 450, 381)$,\adfsplit
$(457, 110, 123, 271)$,
$(457, 9, 52, 455)$,
$(457, 34, 88, 375)$,\adfsplit
$(457, 449, 67, 144)$,
$(458, 302, 351, 325)$,
$(458, 188, 143, 8)$,\adfsplit
$(458, 448, 54, 417)$,
$(458, 74, 363, 172)$,
$(458, 83, 300, 394)$,\adfsplit
$(458, 49, 21, 120)$,
$(458, 354, 53, 113)$,
$(458, 7, 451, 22)$,\adfsplit
$(459, 40, 70, 212)$,
$(459, 234, 453, 97)$,
$(459, 449, 416, 414)$,\adfsplit
$(459, 195, 62, 297)$,
$(459, 63, 173, 348)$,
$(459, 83, 199, 349)$,\adfsplit
$(459, 274, 143, 235)$,
$(459, 146, 360, 100)$,
$(460, 289, 59, 166)$,\adfsplit
$(460, 295, 26, 373)$,
$(460, 209, 317, 71)$,
$(460, 274, 15, 20)$,\adfsplit
$(460, 451, 304, 12)$,
$(460, 224, 177, 141)$,
$(460, 216, 28, 182)$,\adfsplit
$(460, 30, 219, 402)$,
$(461, 454, 445, 379)$,
$(461, 164, 447, 216)$,\adfsplit
$(461, 366, 130, 148)$,
$(461, 315, 41, 439)$,
$(461, 9, 359, 272)$,\adfsplit
$(461, 158, 213, 36)$,
$(461, 234, 275, 77)$,
$(461, 314, 328, 433)$,\adfsplit
$(462, 360, 219, 38)$,
$(462, 440, 124, 337)$,
$(0, 21, 395, 551)$,\adfsplit
$(0, 59, 126, 582)$,
$(0, 6, 57, 543)$,
$(0, 117, 326, 559)$,\adfsplit
$(0, 7, 125, 278)$,
$(0, 29, 251, 566)$,
$(0, 90, 211, 567)$,\adfsplit
$(0, 85, 186, 510)$,
$(0, 19, 343, 487)$,
$(0, 95, 294, 503)$,\adfsplit
$(0, 22, 255, 519)$,
$(0, 42, 86, 614)$,
$(0, 3, 79, 574)$,\adfsplit
$(0, 63, 212, 599)$,
$(0, 20, 159, 606)$,
$(0, 114, 228, 342)$

%ADFvfyBlocksEnd
\adfLgap \noindent by the mapping:
$x \mapsto x +  j \adfmod{456}$ for $x < 456$,
$x \mapsto (x - 456 + 8 j \adfmod{192}) + 456$ for $x \ge 456$,
$0 \le j < 456$
 for the first 65 blocks,
$0 \le j < 114$
 for the last block.
\ADFvfyParStart{(648, ((65, 456, ((456, 1), (192, 8))), (1, 114, ((456, 1), (192, 8)))), ((57, 8), (192, 1)))} %ADFvfyParEnd
% End of 57^8 192^1
%%%%%%%%%%%%%%%%%%%%%%%%%%%%%%%%%%%%%%%%%%%%%%%%%%%%%%%%%%%%%%%%%%%%%%%%%%%%%%%%%%%%%%%%%%
%%%%%%%%%%%%%%%%%%%%%%%%%%%%%%%%%%%%%%%%%%%%%%%%%%%%%%%%%%%%%%%%%%%%%%%%%%%%%%%%%%%%%%%%%%

% Charlotte:GDD4-1-3-5-mod-6-TeX-gen-A:HITS-fun:4.10
\adfDgap
%ADFvfyBlocksStart {57,57,57,57,57,57,57,57,195}
\noindent{\boldmath $ 57^{8} 195^{1} $}~
With the point set $Z_{651}$ partitioned into
 residue classes modulo $8$ for $\{0, 1, \dots, 455\}$, and
 $\{456, 457, \dots, 650\}$,
 the design is generated from

\adfLgap %ADFvfyDesignStart
$(648, 354, 53, 56)$,
$(648, 346, 295, 453)$,
$(456, 121, 400, 35)$,\adfsplit
$(456, 173, 257, 380)$,
$(456, 253, 200, 242)$,
$(456, 174, 367, 442)$,\adfsplit
$(456, 422, 23, 148)$,
$(456, 231, 81, 358)$,
$(456, 192, 123, 285)$,\adfsplit
$(456, 19, 426, 396)$,
$(457, 298, 27, 65)$,
$(457, 8, 63, 438)$,\adfsplit
$(457, 424, 132, 245)$,
$(457, 25, 155, 350)$,
$(457, 240, 47, 146)$,\adfsplit
$(457, 405, 81, 46)$,
$(457, 139, 188, 349)$,
$(457, 282, 295, 292)$,\adfsplit
$(458, 210, 80, 148)$,
$(458, 453, 188, 27)$,
$(458, 115, 373, 384)$,\adfsplit
$(458, 108, 63, 302)$,
$(458, 179, 136, 77)$,
$(458, 321, 394, 119)$,\adfsplit
$(458, 414, 265, 410)$,
$(458, 271, 89, 310)$,
$(459, 68, 16, 235)$,\adfsplit
$(459, 95, 10, 285)$,
$(459, 127, 290, 368)$,
$(459, 190, 289, 181)$,\adfsplit
$(459, 195, 329, 183)$,
$(459, 278, 297, 364)$,
$(459, 66, 323, 365)$,\adfsplit
$(459, 228, 246, 432)$,
$(460, 160, 391, 166)$,
$(460, 322, 289, 299)$,\adfsplit
$(460, 272, 163, 209)$,
$(460, 172, 129, 66)$,
$(460, 349, 62, 116)$,\adfsplit
$(460, 411, 288, 429)$,
$(460, 255, 222, 101)$,
$(460, 84, 23, 194)$,\adfsplit
$(461, 105, 446, 101)$,
$(461, 67, 39, 96)$,
$(461, 204, 184, 138)$,\adfsplit
$(461, 305, 299, 284)$,
$(461, 241, 146, 151)$,
$(461, 253, 246, 274)$,\adfsplit
$(461, 123, 261, 142)$,
$(461, 32, 407, 340)$,
$(462, 413, 273, 244)$,\adfsplit
$(462, 60, 327, 107)$,
$(462, 96, 26, 97)$,
$(462, 186, 133, 195)$,\adfsplit
$(462, 418, 353, 278)$,
$(462, 368, 391, 286)$,
$(462, 150, 184, 211)$,\adfsplit
$(462, 119, 380, 453)$,
$(463, 328, 412, 327)$,
$(463, 359, 10, 270)$,\adfsplit
$(463, 237, 242, 139)$,
$(463, 419, 193, 110)$,
$(463, 223, 293, 105)$,\adfsplit
$(463, 450, 70, 253)$,
$(463, 368, 92, 281)$,
$(463, 348, 147, 360)$,\adfsplit
$(464, 362, 273, 304)$,
$(464, 440, 211, 174)$,
$(464, 141, 383, 258)$,\adfsplit
$(464, 52, 373, 305)$,
$(464, 343, 35, 70)$,
$(464, 300, 1, 207)$,\adfsplit
$(464, 219, 446, 202)$,
$(464, 260, 245, 384)$,
$(465, 398, 376, 293)$,\adfsplit
$(465, 11, 183, 445)$,
$(465, 311, 120, 430)$,
$(465, 7, 226, 363)$,\adfsplit
$(465, 18, 260, 1)$,
$(465, 67, 429, 177)$,
$(465, 60, 362, 152)$,\adfsplit
$(465, 294, 244, 353)$,
$(466, 443, 53, 417)$,
$(466, 210, 112, 163)$,\adfsplit
$(466, 39, 190, 250)$,
$(466, 237, 411, 32)$,
$(466, 121, 412, 157)$,\adfsplit
$(466, 270, 215, 420)$,
$(466, 192, 20, 233)$,
$(466, 415, 314, 422)$,\adfsplit
$(467, 356, 370, 238)$,
$(467, 208, 325, 113)$,
$(467, 96, 295, 218)$,\adfsplit
$(467, 227, 351, 80)$,
$(467, 261, 302, 211)$,
$(467, 49, 222, 5)$,\adfsplit
$(467, 407, 441, 268)$,
$(467, 195, 156, 66)$,
$(468, 102, 427, 205)$,\adfsplit
$(468, 8, 215, 82)$,
$(468, 81, 258, 384)$,
$(468, 386, 221, 156)$,\adfsplit
$(468, 257, 237, 424)$,
$(468, 238, 127, 412)$,
$(468, 351, 171, 326)$,\adfsplit
$(468, 121, 356, 203)$,
$(469, 171, 312, 169)$,
$(469, 284, 151, 429)$,\adfsplit
$(469, 146, 246, 12)$,
$(0, 44, 247, 630)$,
$(0, 151, 209, 599)$,\adfsplit
$(0, 87, 166, 517)$,
$(1, 55, 197, 534)$,
$(0, 211, 327, 487)$,\adfsplit
$(0, 2, 431, 471)$,
$(0, 343, 419, 486)$,
$(0, 38, 159, 597)$,\adfsplit
$(0, 135, 241, 549)$,
$(0, 31, 102, 629)$,
$(1, 53, 271, 565)$,\adfsplit
$(0, 97, 175, 503)$,
$(1, 127, 291, 631)$,
$(0, 45, 220, 614)$,\adfsplit
$(1, 15, 397, 470)$,
$(0, 127, 286, 567)$,
$(0, 69, 206, 355)$,\adfsplit
$(0, 163, 278, 551)$,
$(0, 36, 238, 519)$,
$(0, 143, 443, 502)$,\adfsplit
$(0, 138, 294, 598)$,
$(0, 116, 314, 534)$,
$(0, 114, 228, 342)$,\adfsplit
$(1, 115, 229, 343)$

%ADFvfyBlocksEnd
\adfLgap \noindent by the mapping:
$x \mapsto x + 2 j \adfmod{456}$ for $x < 456$,
$x \mapsto (x - 456 + 16 j \adfmod{192}) + 456$ for $456 \le x < 648$,
$x \mapsto (x +  j \adfmod{3}) + 648$ for $x \ge 648$,
$0 \le j < 228$
 for the first 131 blocks,
$0 \le j < 57$
 for the last two blocks.
\ADFvfyParStart{(651, ((131, 228, ((456, 2), (192, 16), (3, 1))), (2, 57, ((456, 2), (192, 16), (3, 1)))), ((57, 8), (195, 1)))} %ADFvfyParEnd
% End of 57^8 195^1
%%%%%%%%%%%%%%%%%%%%%%%%%%%%%%%%%%%%%%%%%%%%%%%%%%%%%%%%%%%%%%%%%%%%%%%%%%%%%%%%%%%%%%%%%%
%%%%%%%%%%%%%%%%%%%%%%%%%%%%%%%%%%%%%%%%%%%%%%%%%%%%%%%%%%%%%%%%%%%%%%%%%%%%%%%%%%%%%%%%%%

%%%%%%%%%%%%%%%%%%%%%%%%%%%%%%%%%%%%%%%%%%%%%%%%%%%%%%%%%%%%%%%%%%%%%%%%%%%%%%%%%%%%%%%%%%
%%%%%%%%%%%%%%%%%%%%%%%%%%%%%%%%%%%%%%%%%%%%%%%%%%%%%%%%%%%%%%%%%%%%%%%%%%%%%%%%%%%%%%%%%%
\section{4-GDDs for the proof of Lemma \ref{lem:4-GDD 69^8 m^1}}
\label{app:4-GDD 69^8 m^1}
\adfnull{
$ 69^8 222^1 $,
$ 69^8 225^1 $,
$ 69^8 228^1 $,
$ 69^8 231^1 $,
$ 69^8 234^1 $ and
$ 69^8 237^1 $.
}

% Charlotte:GDD4-1-3-5-mod-6-TeX-gen-A:HITS-fun:4.10
\adfDgap
%ADFvfyBlocksStart {69,69,69,69,69,69,69,69,222}
\noindent{\boldmath $ 69^{8} 222^{1} $}~
With the point set $Z_{774}$ partitioned into
 residue classes modulo $8$ for $\{0, 1, \dots, 551\}$, and
 $\{552, 553, \dots, 773\}$,
 the design is generated from

\adfLgap %ADFvfyDesignStart
$(768, 527, 259, 30)$,
$(768, 356, 435, 262)$,
$(552, 461, 539, 172)$,\adfsplit
$(552, 194, 492, 41)$,
$(552, 38, 284, 250)$,
$(552, 320, 237, 330)$,\adfsplit
$(552, 214, 267, 441)$,
$(552, 235, 271, 304)$,
$(552, 399, 37, 342)$,\adfsplit
$(552, 503, 192, 265)$,
$(553, 162, 15, 213)$,
$(553, 99, 73, 55)$,\adfsplit
$(553, 30, 143, 154)$,
$(553, 403, 2, 196)$,
$(553, 176, 513, 461)$,\adfsplit
$(553, 277, 395, 444)$,
$(553, 232, 380, 110)$,
$(553, 257, 382, 48)$,\adfsplit
$(554, 275, 26, 350)$,
$(554, 56, 73, 159)$,
$(554, 53, 174, 40)$,\adfsplit
$(554, 36, 143, 336)$,
$(554, 89, 334, 499)$,
$(554, 525, 82, 260)$,\adfsplit
$(554, 363, 295, 474)$,
$(554, 268, 201, 37)$,
$(555, 389, 449, 88)$,\adfsplit
$(555, 110, 25, 0)$,
$(555, 501, 162, 80)$,
$(555, 178, 451, 367)$,\adfsplit
$(555, 2, 539, 396)$,
$(555, 174, 75, 471)$,
$(555, 393, 190, 212)$,\adfsplit
$(555, 311, 109, 148)$,
$(556, 351, 125, 443)$,
$(556, 306, 247, 102)$,\adfsplit
$(556, 57, 386, 432)$,
$(556, 383, 195, 93)$,
$(556, 109, 329, 310)$,\adfsplit
$(556, 124, 224, 466)$,
$(556, 110, 108, 145)$,
$(556, 307, 544, 68)$,\adfsplit
$(557, 240, 44, 498)$,
$(557, 130, 103, 541)$,
$(557, 141, 371, 14)$,\adfsplit
$(557, 316, 551, 273)$,
$(557, 64, 483, 102)$,
$(557, 495, 488, 235)$,\adfsplit
$(557, 70, 425, 276)$,
$(557, 53, 362, 433)$,
$(558, 371, 526, 114)$,\adfsplit
$(558, 423, 93, 212)$,
$(558, 139, 276, 170)$,
$(558, 456, 414, 479)$,\adfsplit
$(558, 136, 271, 89)$,
$(558, 493, 465, 4)$,
$(558, 302, 392, 173)$,\adfsplit
$(558, 106, 483, 97)$,
$(559, 27, 177, 7)$,
$(559, 432, 278, 146)$,\adfsplit
$(559, 165, 19, 409)$,
$(559, 279, 396, 395)$,
$(559, 143, 20, 17)$,\adfsplit
$(0, 4, 54, 739)$,
$(0, 6, 205, 631)$,
$(0, 12, 271, 721)$,\adfsplit
$(0, 5, 66, 605)$,
$(0, 14, 183, 704)$,
$(0, 21, 108, 587)$,\adfsplit
$(0, 130, 291, 740)$,
$(0, 105, 319, 596)$,
$(0, 89, 269, 366)$,\adfsplit
$(0, 29, 70, 507)$,
$(0, 77, 302, 668)$,
$(0, 62, 157, 393)$,\adfsplit
$(0, 58, 471, 641)$,
$(0, 30, 217, 686)$,
$(0, 138, 276, 414)$

%ADFvfyBlocksEnd
\adfLgap \noindent by the mapping:
$x \mapsto x +  j \adfmod{552}$ for $x < 552$,
$x \mapsto (x - 552 + 9 j \adfmod{216}) + 552$ for $552 \le x < 768$,
$x \mapsto (x +  j \adfmod{6}) + 768$ for $x \ge 768$,
$0 \le j < 552$
 for the first 77 blocks,
$0 \le j < 138$
 for the last block.
\ADFvfyParStart{(774, ((77, 552, ((552, 1), (216, 9), (6, 1))), (1, 138, ((552, 1), (216, 9), (6, 1)))), ((69, 8), (222, 1)))} %ADFvfyParEnd
% End of 69^8 222^1
%%%%%%%%%%%%%%%%%%%%%%%%%%%%%%%%%%%%%%%%%%%%%%%%%%%%%%%%%%%%%%%%%%%%%%%%%%%%%%%%%%%%%%%%%%
%%%%%%%%%%%%%%%%%%%%%%%%%%%%%%%%%%%%%%%%%%%%%%%%%%%%%%%%%%%%%%%%%%%%%%%%%%%%%%%%%%%%%%%%%%

% Charlotte:GDD4-1-3-5-mod-6-TeX-gen-A:HITS-fun:4.10
\adfDgap
%ADFvfyBlocksStart {69,69,69,69,69,69,69,69,225}
\noindent{\boldmath $ 69^{8} 225^{1} $}~
With the point set $Z_{777}$ partitioned into
 residue classes modulo $8$ for $\{0, 1, \dots, 551\}$, and
 $\{552, 553, \dots, 776\}$,
 the design is generated from

\adfLgap %ADFvfyDesignStart
$(768, 178, 174, 551)$,
$(768, 213, 446, 313)$,
$(769, 333, 454, 212)$,\adfsplit
$(769, 187, 510, 503)$,
$(770, 89, 458, 238)$,
$(770, 297, 247, 96)$,\adfsplit
$(552, 436, 272, 443)$,
$(552, 85, 372, 462)$,
$(552, 242, 496, 199)$,\adfsplit
$(552, 140, 453, 214)$,
$(552, 298, 209, 23)$,
$(552, 350, 177, 264)$,\adfsplit
$(552, 265, 375, 363)$,
$(552, 114, 403, 53)$,
$(553, 110, 440, 135)$,\adfsplit
$(553, 473, 102, 213)$,
$(553, 334, 175, 250)$,
$(553, 124, 119, 433)$,\adfsplit
$(553, 499, 465, 29)$,
$(553, 336, 324, 83)$,
$(553, 500, 349, 186)$,\adfsplit
$(553, 400, 410, 3)$,
$(554, 139, 145, 87)$,
$(554, 356, 414, 458)$,\adfsplit
$(554, 384, 413, 531)$,
$(554, 535, 179, 40)$,
$(554, 349, 100, 489)$,\adfsplit
$(554, 281, 82, 446)$,
$(554, 549, 114, 214)$,
$(554, 176, 156, 407)$,\adfsplit
$(555, 505, 335, 448)$,
$(555, 29, 183, 60)$,
$(555, 322, 476, 491)$,\adfsplit
$(555, 248, 213, 297)$,
$(555, 402, 78, 411)$,
$(555, 163, 41, 205)$,\adfsplit
$(555, 38, 199, 144)$,
$(555, 410, 28, 382)$,
$(556, 169, 194, 123)$,\adfsplit
$(556, 525, 246, 442)$,
$(556, 489, 52, 245)$,
$(556, 284, 383, 24)$,\adfsplit
$(556, 391, 517, 310)$,
$(556, 83, 208, 351)$,
$(556, 398, 403, 377)$,\adfsplit
$(556, 320, 18, 252)$,
$(557, 183, 413, 464)$,
$(557, 271, 372, 169)$,\adfsplit
$(557, 454, 268, 163)$,
$(557, 432, 17, 93)$,
$(557, 47, 491, 92)$,\adfsplit
$(557, 282, 304, 75)$,
$(557, 182, 153, 218)$,
$(557, 229, 270, 442)$,\adfsplit
$(558, 440, 421, 124)$,
$(558, 241, 95, 155)$,
$(558, 40, 302, 43)$,\adfsplit
$(558, 351, 21, 89)$,
$(558, 444, 510, 463)$,
$(558, 389, 288, 274)$,\adfsplit
$(558, 459, 46, 458)$,
$(558, 116, 297, 234)$,
$(559, 148, 463, 289)$,\adfsplit
$(559, 159, 438, 168)$,
$(559, 277, 527, 19)$,
$(559, 459, 254, 521)$,\adfsplit
$(559, 404, 443, 466)$,
$(559, 405, 320, 228)$,
$(559, 378, 297, 16)$,\adfsplit
$(559, 317, 74, 478)$,
$(560, 188, 5, 535)$,
$(560, 73, 158, 532)$,\adfsplit
$(560, 107, 80, 213)$,
$(560, 27, 162, 161)$,
$(560, 202, 135, 46)$,\adfsplit
$(560, 228, 469, 496)$,
$(560, 95, 369, 48)$,
$(560, 434, 115, 126)$,\adfsplit
$(561, 22, 154, 171)$,
$(561, 354, 169, 278)$,
$(561, 199, 328, 146)$,\adfsplit
$(561, 375, 65, 547)$,
$(561, 180, 333, 177)$,
$(561, 284, 104, 6)$,\adfsplit
$(561, 268, 311, 341)$,
$(561, 517, 72, 59)$,
$(562, 487, 358, 281)$,\adfsplit
$(562, 528, 133, 394)$,
$(562, 340, 527, 125)$,
$(562, 443, 174, 156)$,\adfsplit
$(562, 488, 374, 548)$,
$(562, 474, 525, 87)$,
$(562, 256, 219, 145)$,\adfsplit
$(562, 513, 482, 139)$,
$(563, 233, 167, 390)$,
$(563, 364, 49, 261)$,\adfsplit
$(563, 378, 504, 166)$,
$(563, 242, 317, 187)$,
$(563, 57, 212, 3)$,\adfsplit
$(563, 540, 58, 247)$,
$(563, 183, 64, 155)$,
$(563, 416, 446, 85)$,\adfsplit
$(564, 38, 257, 84)$,
$(564, 207, 59, 306)$,
$(564, 120, 373, 242)$,\adfsplit
$(564, 265, 99, 20)$,
$(564, 263, 104, 154)$,
$(564, 309, 510, 124)$,\adfsplit
$(564, 376, 283, 297)$,
$(564, 151, 485, 190)$,
$(565, 17, 477, 442)$,\adfsplit
$(565, 11, 138, 340)$,
$(565, 200, 254, 404)$,
$(565, 37, 372, 195)$,\adfsplit
$(565, 111, 438, 360)$,
$(565, 153, 142, 352)$,
$(565, 269, 359, 145)$,\adfsplit
$(565, 55, 259, 146)$,
$(566, 186, 383, 76)$,
$(566, 315, 30, 125)$,\adfsplit
$(566, 103, 85, 385)$,
$(566, 242, 377, 500)$,
$(566, 323, 310, 392)$,\adfsplit
$(566, 136, 501, 177)$,
$(566, 228, 446, 120)$,
$(566, 523, 159, 34)$,\adfsplit
$(567, 505, 204, 434)$,
$(567, 67, 287, 77)$,
$(567, 426, 255, 20)$,\adfsplit
$(567, 381, 343, 201)$,
$(567, 315, 456, 294)$,
$(567, 232, 532, 85)$,\adfsplit
$(567, 371, 449, 502)$,
$(0, 42, 500, 657)$,
$(0, 34, 317, 658)$,\adfsplit
$(0, 59, 479, 569)$,
$(0, 2, 537, 604)$,
$(0, 26, 433, 605)$,\adfsplit
$(0, 23, 257, 641)$,
$(0, 38, 103, 107)$,
$(0, 6, 361, 363)$,\adfsplit
$(0, 67, 116, 457)$,
$(0, 77, 246, 676)$,
$(0, 61, 415, 623)$,\adfsplit
$(0, 117, 394, 730)$,
$(0, 179, 346, 767)$,
$(0, 87, 428, 712)$,\adfsplit
$(0, 142, 469, 731)$,
$(0, 45, 266, 713)$,
$(0, 97, 130, 455)$,\adfsplit
$(0, 195, 493, 659)$,
$(0, 37, 263, 766)$,
$(0, 167, 203, 694)$,\adfsplit
$(0, 33, 215, 358)$,
$(1, 21, 267, 694)$,
$(0, 138, 276, 414)$,\adfsplit
$(1, 139, 277, 415)$

%ADFvfyBlocksEnd
\adfLgap \noindent by the mapping:
$x \mapsto x + 2 j \adfmod{552}$ for $x < 552$,
$x \mapsto (x - 552 + 18 j \adfmod{216}) + 552$ for $552 \le x < 768$,
$x \mapsto (x - 768 + 3 j \adfmod{9}) + 768$ for $x \ge 768$,
$0 \le j < 276$
 for the first 155 blocks,
$0 \le j < 69$
 for the last two blocks.
\ADFvfyParStart{(777, ((155, 276, ((552, 2), (216, 18), (9, 3))), (2, 69, ((552, 2), (216, 18), (9, 3)))), ((69, 8), (225, 1)))} %ADFvfyParEnd
% End of 69^8 225^1
%%%%%%%%%%%%%%%%%%%%%%%%%%%%%%%%%%%%%%%%%%%%%%%%%%%%%%%%%%%%%%%%%%%%%%%%%%%%%%%%%%%%%%%%%%
%%%%%%%%%%%%%%%%%%%%%%%%%%%%%%%%%%%%%%%%%%%%%%%%%%%%%%%%%%%%%%%%%%%%%%%%%%%%%%%%%%%%%%%%%%

% Charlotte:GDD4-1-3-5-mod-6-TeX-gen-A:HITS-fun:4.10
\adfDgap
%ADFvfyBlocksStart {69,69,69,69,69,69,69,69,228}
\noindent{\boldmath $ 69^{8} 228^{1} $}~
With the point set $Z_{780}$ partitioned into
 residue classes modulo $8$ for $\{0, 1, \dots, 551\}$, and
 $\{552, 553, \dots, 779\}$,
 the design is generated from

\adfLgap %ADFvfyDesignStart
$(768, 507, 20, 550)$,
$(768, 403, 462, 173)$,
$(769, 372, 111, 245)$,\adfsplit
$(769, 500, 328, 487)$,
$(552, 277, 489, 379)$,
$(552, 344, 243, 521)$,\adfsplit
$(552, 311, 234, 97)$,
$(552, 269, 410, 94)$,
$(552, 364, 7, 534)$,\adfsplit
$(552, 132, 549, 370)$,
$(552, 326, 231, 504)$,
$(552, 179, 424, 428)$,\adfsplit
$(553, 414, 441, 525)$,
$(553, 517, 324, 282)$,
$(553, 303, 267, 182)$,\adfsplit
$(553, 71, 19, 217)$,
$(553, 340, 509, 449)$,
$(553, 7, 136, 10)$,\adfsplit
$(553, 168, 299, 386)$,
$(553, 190, 476, 272)$,
$(554, 406, 159, 380)$,\adfsplit
$(554, 531, 33, 246)$,
$(554, 481, 527, 180)$,
$(554, 326, 152, 546)$,\adfsplit
$(554, 391, 0, 50)$,
$(554, 413, 347, 425)$,
$(554, 520, 117, 76)$,\adfsplit
$(554, 547, 466, 277)$,
$(555, 478, 536, 210)$,
$(555, 194, 29, 513)$,\adfsplit
$(555, 516, 206, 21)$,
$(555, 427, 277, 532)$,
$(555, 298, 267, 116)$,\adfsplit
$(555, 127, 112, 371)$,
$(555, 305, 24, 375)$,
$(555, 191, 73, 246)$,\adfsplit
$(556, 504, 217, 46)$,
$(556, 165, 23, 210)$,
$(556, 416, 469, 540)$,\adfsplit
$(556, 305, 507, 172)$,
$(556, 514, 542, 467)$,
$(556, 509, 151, 547)$,\adfsplit
$(556, 208, 20, 510)$,
$(556, 434, 513, 327)$,
$(557, 324, 83, 263)$,\adfsplit
$(557, 22, 188, 328)$,
$(557, 193, 346, 54)$,
$(557, 50, 13, 43)$,\adfsplit
$(557, 230, 387, 393)$,
$(557, 511, 234, 28)$,
$(557, 161, 485, 152)$,\adfsplit
$(557, 432, 471, 381)$,
$(558, 389, 274, 155)$,
$(558, 52, 291, 69)$,\adfsplit
$(558, 439, 12, 421)$,
$(558, 104, 446, 441)$,
$(558, 111, 406, 258)$,\adfsplit
$(558, 328, 551, 458)$,
$(558, 126, 20, 449)$,
$(558, 289, 312, 403)$,\adfsplit
$(559, 47, 209, 499)$,
$(559, 241, 4, 87)$,
$(559, 301, 70, 322)$,\adfsplit
$(559, 443, 174, 357)$,
$(0, 122, 325, 578)$,
$(0, 2, 76, 605)$,\adfsplit
$(0, 10, 209, 649)$,
$(0, 113, 371, 694)$,
$(0, 11, 356, 631)$,\adfsplit
$(0, 99, 397, 730)$,
$(0, 14, 103, 722)$,
$(0, 33, 197, 695)$,\adfsplit
$(0, 49, 116, 632)$,
$(0, 1, 35, 362)$,
$(0, 132, 385, 749)$,\adfsplit
$(0, 20, 117, 731)$,
$(0, 98, 407, 704)$,
$(0, 19, 63, 92)$,\adfsplit
$(0, 138, 276, 414)$

%ADFvfyBlocksEnd
\adfLgap \noindent by the mapping:
$x \mapsto x +  j \adfmod{552}$ for $x < 552$,
$x \mapsto (x - 552 + 9 j \adfmod{216}) + 552$ for $552 \le x < 768$,
$x \mapsto (x + 2 j \adfmod{12}) + 768$ for $x \ge 768$,
$0 \le j < 552$
 for the first 78 blocks,
$0 \le j < 138$
 for the last block.
\ADFvfyParStart{(780, ((78, 552, ((552, 1), (216, 9), (12, 2))), (1, 138, ((552, 1), (216, 9), (12, 2)))), ((69, 8), (228, 1)))} %ADFvfyParEnd
% End of 69^8 228^1
%%%%%%%%%%%%%%%%%%%%%%%%%%%%%%%%%%%%%%%%%%%%%%%%%%%%%%%%%%%%%%%%%%%%%%%%%%%%%%%%%%%%%%%%%%
%%%%%%%%%%%%%%%%%%%%%%%%%%%%%%%%%%%%%%%%%%%%%%%%%%%%%%%%%%%%%%%%%%%%%%%%%%%%%%%%%%%%%%%%%%

% Charlotte:GDD4-1-3-5-mod-6-TeX-gen-A:HITS-fun:4.10
\adfDgap
%ADFvfyBlocksStart {69,69,69,69,69,69,69,69,231}
\noindent{\boldmath $ 69^{8} 231^{1} $}~
With the point set $Z_{783}$ partitioned into
 residue classes modulo $7$ for $\{0, 1, \dots, 482\}$,
 $\{483, 484, \dots, 551\}$, and
 $\{552, 553, \dots, 782\}$,
 the design is generated from

\adfLgap %ADFvfyDesignStart
$(91, 201, 536, 693)$,
$(381, 427, 510, 755)$,
$(206, 415, 495, 592)$,\adfsplit
$(320, 447, 486, 760)$,
$(170, 266, 514, 740)$,
$(131, 206, 280, 585)$,\adfsplit
$(320, 396, 450, 555)$,
$(157, 268, 440, 670)$,
$(115, 200, 420, 617)$,\adfsplit
$(22, 227, 287, 638)$,
$(103, 140, 422, 618)$,
$(72, 187, 344, 604)$,\adfsplit
$(145, 482, 494, 646)$,
$(67, 259, 467, 662)$,
$(363, 392, 511, 774)$,\adfsplit
$(398, 455, 512, 748)$,
$(92, 347, 502, 746)$,
$(128, 194, 368, 726)$,\adfsplit
$(394, 406, 530, 578)$,
$(57, 422, 540, 584)$,
$(232, 364, 370, 589)$,\adfsplit
$(116, 222, 372, 664)$,
$(40, 287, 309, 586)$,
$(133, 156, 180, 591)$,\adfsplit
$(117, 126, 144, 769)$,
$(0, 253, 503, 556)$,
$(13, 51, 95, 722)$,\adfsplit
$(246, 278, 546, 776)$,
$(155, 403, 501, 761)$,
$(95, 243, 528, 781)$,\adfsplit
$(81, 122, 542, 708)$,
$(202, 213, 462, 615)$,
$(113, 114, 222, 594)$,\adfsplit
$(46, 234, 236, 564)$,
$(191, 291, 412, 668)$,
$(301, 354, 388, 629)$,\adfsplit
$(63, 121, 467, 766)$,
$(168, 338, 474, 613)$,
$(90, 366, 480, 740)$,\adfsplit
$(12, 296, 541, 586)$,
$(9, 280, 372, 705)$,
$(178, 337, 341, 608)$,\adfsplit
$(133, 198, 518, 560)$,
$(123, 254, 420, 764)$,
$(295, 397, 487, 678)$,\adfsplit
$(253, 263, 486, 560)$,
$(23, 330, 539, 641)$,
$(210, 307, 408, 584)$,\adfsplit
$(149, 407, 529, 617)$,
$(204, 357, 461, 768)$,
$(58, 260, 494, 739)$,\adfsplit
$(273, 337, 380, 621)$,
$(164, 215, 235, 636)$,
$(63, 307, 402, 645)$,\adfsplit
$(26, 280, 461, 738)$,
$(8, 166, 202, 777)$,
$(93, 98, 339, 729)$,\adfsplit
$(221, 237, 486, 745)$,
$(66, 128, 482, 687)$,
$(118, 283, 333, 732)$,\adfsplit
$(84, 264, 535, 554)$,
$(280, 305, 535, 725)$,
$(85, 111, 130, 570)$,\adfsplit
$(176, 304, 363, 779)$,
$(135, 215, 339, 699)$,
$(300, 369, 442, 649)$,\adfsplit
$(78, 166, 427, 774)$,
$(237, 326, 500, 698)$,
$(0, 15, 156, 758)$,\adfsplit
$(0, 39, 125, 770)$,
$(0, 8, 300, 639)$,
$(0, 3, 33, 557)$,\adfsplit
$(0, 81, 498, 744)$,
$(0, 116, 513, 574)$,
$(0, 99, 250, 706)$,\adfsplit
$(0, 145, 508, 724)$,
$(0, 117, 323, 702)$,
$(0, 113, 286, 490)$,\adfsplit
$(0, 179, 517, 676)$,
$(0, 72, 491, 643)$,
$(0, 123, 316, 720)$,\adfsplit
$(0, 31, 298, 686)$,
$(0, 13, 68, 710)$,
$(0, 152, 321, 754)$,\adfsplit
$(0, 90, 184, 664)$,
$(0, 40, 143, 195)$,
$(0, 17, 139, 587)$,\adfsplit
$(0, 61, 232, 558)$,
$(0, 78, 348, 734)$,
$(0, 219, 546, 635)$

%ADFvfyBlocksEnd
\adfLgap \noindent by the mapping:
$x \mapsto x +  j \adfmod{483}$ for $x < 483$,
$x \mapsto (x +  j \adfmod{69}) + 483$ for $483 \le x < 552$,
$x \mapsto (x - 552 + 11 j \adfmod{231}) + 552$ for $x \ge 552$,
$0 \le j < 483$.
\ADFvfyParStart{(783, ((90, 483, ((483, 1), (69, 1), (231, 11)))), ((69, 7), (69, 1), (231, 1)))} %ADFvfyParEnd
% End of 69^8 231^1
%%%%%%%%%%%%%%%%%%%%%%%%%%%%%%%%%%%%%%%%%%%%%%%%%%%%%%%%%%%%%%%%%%%%%%%%%%%%%%%%%%%%%%%%%%
%%%%%%%%%%%%%%%%%%%%%%%%%%%%%%%%%%%%%%%%%%%%%%%%%%%%%%%%%%%%%%%%%%%%%%%%%%%%%%%%%%%%%%%%%%

% Charlotte:GDD4-1-3-5-mod-6-TeX-gen-A:HITS-fun:4.10
\adfDgap
%ADFvfyBlocksStart {69,69,69,69,69,69,69,69,234}
\noindent{\boldmath $ 69^{8} 234^{1} $}~
With the point set $Z_{786}$ partitioned into
 residue classes modulo $8$ for $\{0, 1, \dots, 551\}$, and
 $\{552, 553, \dots, 785\}$,
 the design is generated from

\adfLgap %ADFvfyDesignStart
$(768, 448, 127, 276)$,
$(768, 203, 440, 177)$,
$(769, 160, 349, 147)$,\adfsplit
$(769, 150, 437, 482)$,
$(770, 183, 109, 440)$,
$(770, 166, 527, 252)$,\adfsplit
$(552, 114, 75, 157)$,
$(552, 113, 268, 127)$,
$(552, 26, 285, 503)$,\adfsplit
$(552, 499, 296, 414)$,
$(552, 548, 490, 153)$,
$(552, 515, 173, 254)$,\adfsplit
$(552, 280, 135, 505)$,
$(552, 288, 156, 310)$,
$(553, 348, 462, 392)$,\adfsplit
$(553, 158, 175, 394)$,
$(553, 91, 543, 292)$,
$(553, 217, 184, 77)$,\adfsplit
$(553, 450, 120, 497)$,
$(553, 382, 513, 363)$,
$(553, 23, 404, 133)$,\adfsplit
$(553, 290, 203, 477)$,
$(554, 470, 69, 227)$,
$(554, 106, 468, 102)$,\adfsplit
$(554, 452, 530, 94)$,
$(554, 447, 104, 523)$,
$(554, 196, 304, 201)$,\adfsplit
$(554, 507, 101, 265)$,
$(554, 527, 330, 301)$,
$(554, 168, 535, 281)$,\adfsplit
$(555, 420, 549, 450)$,
$(555, 277, 275, 337)$,
$(555, 0, 46, 53)$,\adfsplit
$(555, 249, 304, 231)$,
$(555, 497, 407, 459)$,
$(555, 110, 244, 56)$,\adfsplit
$(555, 246, 530, 463)$,
$(555, 500, 466, 187)$,
$(556, 277, 185, 212)$,\adfsplit
$(556, 386, 174, 0)$,
$(556, 145, 3, 156)$,
$(556, 320, 447, 34)$,\adfsplit
$(556, 88, 31, 402)$,
$(556, 525, 388, 143)$,
$(556, 158, 57, 179)$,\adfsplit
$(556, 190, 221, 43)$,
$(557, 319, 16, 517)$,
$(557, 120, 245, 239)$,\adfsplit
$(557, 11, 258, 294)$,
$(557, 309, 49, 214)$,
$(557, 523, 185, 194)$,\adfsplit
$(557, 348, 81, 363)$,
$(557, 519, 470, 44)$,
$(557, 316, 154, 248)$,\adfsplit
$(558, 274, 454, 97)$,
$(558, 451, 540, 223)$,
$(558, 196, 519, 173)$,\adfsplit
$(558, 472, 81, 531)$,
$(558, 140, 440, 374)$,
$(558, 246, 453, 450)$,\adfsplit
$(558, 37, 98, 233)$,
$(558, 360, 443, 431)$,
$(559, 495, 532, 73)$,\adfsplit
$(559, 430, 545, 402)$,
$(0, 1, 447, 568)$,
$(0, 41, 502, 560)$,\adfsplit
$(0, 97, 341, 614)$,
$(0, 20, 319, 650)$,
$(0, 123, 373, 703)$,\adfsplit
$(0, 10, 121, 379)$,
$(0, 156, 325, 712)$,
$(0, 148, 347, 730)$,\adfsplit
$(0, 69, 359, 658)$,
$(0, 117, 428, 721)$,
$(0, 42, 339, 695)$,\adfsplit
$(0, 25, 468, 740)$,
$(0, 35, 489, 758)$,
$(0, 79, 246, 767)$,\adfsplit
$(0, 159, 389, 641)$,
$(0, 138, 276, 414)$

%ADFvfyBlocksEnd
\adfLgap \noindent by the mapping:
$x \mapsto x +  j \adfmod{552}$ for $x < 552$,
$x \mapsto (x - 552 + 9 j \adfmod{216}) + 552$ for $552 \le x < 768$,
$x \mapsto (x - 768 + 3 j \adfmod{18}) + 768$ for $x \ge 768$,
$0 \le j < 552$
 for the first 79 blocks,
$0 \le j < 138$
 for the last block.
\ADFvfyParStart{(786, ((79, 552, ((552, 1), (216, 9), (18, 3))), (1, 138, ((552, 1), (216, 9), (18, 3)))), ((69, 8), (234, 1)))} %ADFvfyParEnd
% End of 69^8 234^1
%%%%%%%%%%%%%%%%%%%%%%%%%%%%%%%%%%%%%%%%%%%%%%%%%%%%%%%%%%%%%%%%%%%%%%%%%%%%%%%%%%%%%%%%%%
%%%%%%%%%%%%%%%%%%%%%%%%%%%%%%%%%%%%%%%%%%%%%%%%%%%%%%%%%%%%%%%%%%%%%%%%%%%%%%%%%%%%%%%%%%

% Charlotte:GDD4-1-3-5-mod-6-TeX-gen-A:HITS-fun:4.10
\adfDgap
%ADFvfyBlocksStart {69,69,69,69,69,69,69,69,237}
\noindent{\boldmath $ 69^{8} 237^{1} $}~
With the point set $Z_{789}$ partitioned into
 residue classes modulo $8$ for $\{0, 1, \dots, 551\}$, and
 $\{552, 553, \dots, 788\}$,
 the design is generated from

\adfLgap %ADFvfyDesignStart
$(768, 528, 357, 410)$,
$(768, 430, 503, 349)$,
$(769, 126, 523, 465)$,\adfsplit
$(769, 422, 437, 16)$,
$(770, 402, 267, 61)$,
$(770, 544, 275, 494)$,\adfsplit
$(771, 394, 216, 169)$,
$(771, 255, 260, 17)$,
$(772, 500, 161, 162)$,\adfsplit
$(772, 21, 439, 538)$,
$(773, 79, 90, 440)$,
$(773, 539, 93, 196)$,\adfsplit
$(774, 466, 512, 51)$,
$(774, 169, 54, 479)$,
$(552, 318, 359, 9)$,\adfsplit
$(552, 173, 387, 154)$,
$(552, 209, 288, 333)$,
$(552, 63, 541, 146)$,\adfsplit
$(552, 388, 470, 224)$,
$(552, 452, 427, 433)$,
$(552, 352, 142, 203)$,\adfsplit
$(552, 492, 138, 367)$,
$(553, 513, 523, 245)$,
$(553, 267, 338, 406)$,\adfsplit
$(553, 223, 164, 56)$,
$(553, 239, 426, 326)$,
$(553, 246, 157, 264)$,\adfsplit
$(553, 409, 256, 124)$,
$(553, 425, 180, 15)$,
$(553, 237, 275, 466)$,\adfsplit
$(554, 236, 523, 494)$,
$(554, 517, 471, 27)$,
$(554, 418, 143, 540)$,\adfsplit
$(554, 385, 237, 40)$,
$(554, 52, 104, 342)$,
$(554, 489, 408, 506)$,\adfsplit
$(554, 103, 214, 18)$,
$(554, 161, 419, 149)$,
$(555, 60, 184, 495)$,\adfsplit
$(555, 395, 129, 205)$,
$(555, 68, 529, 138)$,
$(555, 454, 248, 274)$,\adfsplit
$(555, 494, 379, 381)$,
$(555, 247, 528, 291)$,
$(555, 486, 194, 551)$,\adfsplit
$(555, 76, 209, 461)$,
$(556, 120, 446, 84)$,
$(556, 121, 464, 247)$,\adfsplit
$(556, 499, 269, 105)$,
$(556, 274, 308, 239)$,
$(556, 142, 213, 363)$,\adfsplit
$(556, 347, 390, 88)$,
$(556, 111, 436, 498)$,
$(556, 2, 253, 113)$,\adfsplit
$(557, 492, 241, 359)$,
$(557, 199, 539, 472)$,
$(557, 139, 344, 398)$,\adfsplit
$(557, 434, 528, 358)$,
$(557, 365, 315, 345)$,
$(557, 514, 236, 294)$,\adfsplit
$(557, 133, 522, 412)$,
$(557, 281, 21, 183)$,
$(558, 55, 512, 29)$,\adfsplit
$(558, 423, 34, 493)$,
$(558, 434, 459, 278)$,
$(558, 308, 81, 448)$,\adfsplit
$(558, 336, 449, 522)$,
$(558, 454, 36, 59)$,
$(558, 263, 91, 294)$,\adfsplit
$(558, 189, 172, 433)$,
$(559, 58, 160, 379)$,
$(559, 77, 340, 22)$,\adfsplit
$(559, 78, 282, 295)$,
$(559, 71, 123, 449)$,
$(559, 433, 80, 228)$,\adfsplit
$(559, 218, 321, 549)$,
$(559, 188, 528, 491)$,
$(559, 543, 325, 494)$,\adfsplit
$(560, 517, 121, 310)$,
$(560, 339, 367, 476)$,
$(560, 137, 80, 303)$,\adfsplit
$(560, 456, 245, 76)$,
$(560, 69, 446, 386)$,
$(560, 203, 60, 126)$,\adfsplit
$(560, 167, 402, 201)$,
$(560, 280, 547, 442)$,
$(561, 63, 206, 186)$,\adfsplit
$(561, 335, 139, 60)$,
$(561, 245, 75, 104)$,
$(561, 448, 79, 25)$,\adfsplit
$(561, 17, 205, 10)$,
$(561, 54, 21, 44)$,
$(561, 292, 214, 408)$,\adfsplit
$(561, 489, 83, 122)$,
$(562, 193, 15, 162)$,
$(562, 536, 226, 161)$,\adfsplit
$(562, 461, 262, 465)$,
$(562, 405, 432, 108)$,
$(562, 26, 280, 307)$,\adfsplit
$(562, 11, 116, 95)$,
$(562, 133, 175, 534)$,
$(562, 374, 460, 315)$,\adfsplit
$(563, 272, 356, 281)$,
$(563, 158, 303, 4)$,
$(563, 445, 343, 25)$,\adfsplit
$(563, 259, 184, 34)$,
$(563, 501, 366, 362)$,
$(563, 216, 227, 149)$,\adfsplit
$(563, 214, 444, 129)$,
$(563, 42, 167, 75)$,
$(564, 7, 460, 309)$,\adfsplit
$(564, 62, 472, 298)$,
$(564, 416, 355, 108)$,
$(564, 314, 219, 48)$,\adfsplit
$(564, 341, 377, 155)$,
$(564, 20, 471, 94)$,
$(564, 481, 287, 541)$,\adfsplit
$(564, 102, 114, 513)$,
$(565, 297, 229, 474)$,
$(565, 192, 460, 213)$,\adfsplit
$(565, 311, 188, 318)$,
$(565, 38, 25, 338)$,
$(565, 83, 293, 415)$,\adfsplit
$(565, 430, 544, 257)$,
$(565, 435, 444, 80)$,
$(565, 43, 226, 495)$,\adfsplit
$(566, 278, 436, 387)$,
$(566, 502, 284, 207)$,
$(566, 0, 282, 319)$,\adfsplit
$(566, 146, 185, 91)$,
$(566, 1, 224, 58)$,
$(566, 299, 119, 365)$,\adfsplit
$(566, 81, 30, 285)$,
$(566, 108, 16, 205)$,
$(567, 31, 464, 13)$,\adfsplit
$(0, 5, 295, 510)$,
$(0, 6, 507, 621)$,
$(0, 106, 537, 658)$,\adfsplit
$(0, 22, 455, 569)$,
$(0, 30, 423, 640)$,
$(0, 53, 426, 585)$,\adfsplit
$(0, 1, 182, 604)$,
$(0, 43, 222, 657)$,
$(0, 131, 313, 605)$,\adfsplit
$(0, 253, 489, 693)$,
$(0, 93, 107, 747)$,
$(0, 14, 161, 675)$,\adfsplit
$(0, 83, 173, 729)$,
$(0, 2, 89, 641)$,
$(0, 63, 149, 748)$,\adfsplit
$(0, 28, 215, 766)$,
$(0, 249, 359, 622)$,
$(0, 241, 263, 676)$,\adfsplit
$(0, 127, 549, 712)$,
$(0, 38, 159, 677)$,
$(0, 3, 44, 695)$,\adfsplit
$(0, 90, 291, 659)$,
$(0, 117, 231, 749)$,
$(1, 83, 199, 749)$,\adfsplit
$(0, 138, 276, 414)$,
$(1, 139, 277, 415)$

%ADFvfyBlocksEnd
\adfLgap \noindent by the mapping:
$x \mapsto x + 2 j \adfmod{552}$ for $x < 552$,
$x \mapsto (x - 552 + 18 j \adfmod{216}) + 552$ for $552 \le x < 768$,
$x \mapsto (x - 768 + 7 j \adfmod{21}) + 768$ for $x \ge 768$,
$0 \le j < 276$
 for the first 159 blocks,
$0 \le j < 69$
 for the last two blocks.
\ADFvfyParStart{(789, ((159, 276, ((552, 2), (216, 18), (21, 7))), (2, 69, ((552, 2), (216, 18), (21, 7)))), ((69, 8), (237, 1)))} %ADFvfyParEnd
% End of 69^8 237^1
%%%%%%%%%%%%%%%%%%%%%%%%%%%%%%%%%%%%%%%%%%%%%%%%%%%%%%%%%%%%%%%%%%%%%%%%%%%%%%%%%%%%%%%%%%
%%%%%%%%%%%%%%%%%%%%%%%%%%%%%%%%%%%%%%%%%%%%%%%%%%%%%%%%%%%%%%%%%%%%%%%%%%%%%%%%%%%%%%%%%%

%%%%%%%%%%%%%%%%%%%%%%%%%%%%%%%%%%%%%%%%%%%%%%%%%%%%%%%%%%%%%%%%%%%%%%%%%%%%%%%%%%%%%%%%%%
%%%%%%%%%%%%%%%%%%%%%%%%%%%%%%%%%%%%%%%%%%%%%%%%%%%%%%%%%%%%%%%%%%%%%%%%%%%%%%%%%%%%%%%%%%
\section{4-GDDs for the proof of Lemma \ref{lem:4-GDD 87^8 m^1}}
\label{app:4-GDD 87^8 m^1}
\adfnull{
$ 87^8 297^1 $ and
$ 87^8 300^1 $.
}

% Charlotte:GDD4-1-3-5-mod-6-TeX-gen-A:HITS-fun:4.10
\adfDgap
%ADFvfyBlocksStart {87,87,87,87,87,87,87,87,297}
\noindent{\boldmath $ 87^{8} 297^{1} $}~
With the point set $Z_{993}$ partitioned into
 residue classes modulo $8$ for $\{0, 1, \dots, 695\}$, and
 $\{696, 697, \dots, 992\}$,
 the design is generated from

\adfLgap %ADFvfyDesignStart
$(984, 244, 659, 417)$,
$(985, 377, 174, 205)$,
$(986, 420, 398, 421)$,\adfsplit
$(696, 256, 114, 487)$,
$(696, 686, 273, 548)$,
$(696, 141, 670, 490)$,\adfsplit
$(696, 560, 100, 493)$,
$(696, 671, 603, 98)$,
$(696, 491, 612, 17)$,\adfsplit
$(696, 101, 240, 351)$,
$(696, 174, 169, 19)$,
$(697, 217, 119, 636)$,\adfsplit
$(697, 122, 477, 355)$,
$(697, 415, 682, 493)$,
$(697, 368, 537, 423)$,\adfsplit
$(697, 90, 374, 603)$,
$(697, 0, 364, 113)$,
$(697, 582, 83, 112)$,\adfsplit
$(697, 461, 454, 524)$,
$(698, 437, 577, 290)$,
$(698, 152, 479, 21)$,\adfsplit
$(698, 483, 158, 593)$,
$(698, 132, 432, 421)$,
$(698, 574, 634, 35)$,\adfsplit
$(698, 426, 543, 208)$,
$(698, 679, 316, 571)$,
$(698, 476, 654, 177)$,\adfsplit
$(699, 42, 237, 196)$,
$(699, 227, 493, 663)$,
$(699, 235, 288, 193)$,\adfsplit
$(699, 441, 438, 365)$,
$(699, 497, 562, 80)$,
$(699, 136, 79, 156)$,\adfsplit
$(699, 483, 239, 670)$,
$(699, 398, 194, 452)$,
$(700, 193, 368, 22)$,\adfsplit
$(700, 576, 610, 129)$,
$(700, 663, 644, 443)$,
$(700, 482, 88, 619)$,\adfsplit
$(700, 295, 90, 604)$,
$(700, 557, 662, 647)$,
$(700, 30, 693, 99)$,\adfsplit
$(700, 13, 545, 324)$,
$(701, 136, 535, 669)$,
$(701, 627, 310, 434)$,\adfsplit
$(701, 293, 311, 564)$,
$(701, 481, 471, 419)$,
$(701, 206, 449, 212)$,\adfsplit
$(701, 6, 115, 445)$,
$(701, 538, 216, 508)$,
$(701, 594, 417, 392)$,\adfsplit
$(702, 401, 530, 438)$,
$(702, 77, 566, 539)$,
$(702, 213, 168, 249)$,\adfsplit
$(702, 247, 306, 147)$,
$(702, 193, 514, 328)$,
$(702, 685, 228, 142)$,\adfsplit
$(702, 647, 572, 296)$,
$(702, 375, 148, 571)$,
$(703, 278, 404, 325)$,\adfsplit
$(703, 311, 499, 168)$,
$(703, 694, 179, 477)$,
$(703, 529, 147, 485)$,\adfsplit
$(703, 471, 425, 132)$,
$(703, 202, 343, 488)$,
$(703, 364, 40, 570)$,\adfsplit
$(703, 414, 33, 410)$,
$(704, 197, 169, 218)$,
$(704, 351, 402, 277)$,\adfsplit
$(704, 499, 549, 511)$,
$(704, 524, 154, 510)$,
$(704, 280, 396, 387)$,\adfsplit
$(704, 124, 566, 81)$,
$(704, 622, 56, 17)$,
$(704, 335, 587, 648)$,\adfsplit
$(705, 589, 283, 366)$,
$(705, 101, 585, 8)$,
$(705, 383, 155, 498)$,\adfsplit
$(705, 529, 447, 118)$,
$(705, 436, 521, 50)$,
$(705, 470, 444, 175)$,\adfsplit
$(705, 260, 363, 576)$,
$(705, 93, 328, 394)$,
$(706, 437, 3, 454)$,\adfsplit
$(706, 235, 381, 672)$,
$(0, 94, 362, 731)$,
$(0, 151, 433, 802)$,\adfsplit
$(0, 2, 89, 898)$,
$(0, 99, 406, 946)$,
$(0, 149, 486, 886)$,\adfsplit
$(0, 71, 318, 850)$,
$(0, 106, 305, 754)$,
$(0, 158, 388, 887)$,\adfsplit
$(0, 118, 274, 767)$,
$(0, 58, 190, 923)$,
$(0, 209, 450, 779)$,\adfsplit
$(0, 84, 278, 426)$,
$(0, 13, 198, 971)$,
$(0, 35, 569, 875)$,\adfsplit
$(0, 133, 294, 983)$,
$(0, 174, 348, 522)$

%ADFvfyBlocksEnd
\adfLgap \noindent by the mapping:
$x \mapsto x +  j \adfmod{696}$ for $x < 696$,
$x \mapsto (x - 696 + 12 j \adfmod{288}) + 696$ for $696 \le x < 984$,
$x \mapsto (x - 984 + 3 j \adfmod{9}) + 984$ for $x \ge 984$,
$0 \le j < 696$
 for the first 100 blocks,
$0 \le j < 174$
 for the last block.
\ADFvfyParStart{(993, ((100, 696, ((696, 1), (288, 12), (9, 3))), (1, 174, ((696, 1), (288, 12), (9, 3)))), ((87, 8), (297, 1)))} %ADFvfyParEnd
% End of 87^8 297^1
%%%%%%%%%%%%%%%%%%%%%%%%%%%%%%%%%%%%%%%%%%%%%%%%%%%%%%%%%%%%%%%%%%%%%%%%%%%%%%%%%%%%%%%%%%
%%%%%%%%%%%%%%%%%%%%%%%%%%%%%%%%%%%%%%%%%%%%%%%%%%%%%%%%%%%%%%%%%%%%%%%%%%%%%%%%%%%%%%%%%%

% Charlotte:GDD4-1-3-5-mod-6-TeX-gen-A:HITS-fun:4.10
\adfDgap
%ADFvfyBlocksStart {87,87,87,87,87,87,87,87,300}
\noindent{\boldmath $ 87^{8} 300^{1} $}~
With the point set $Z_{996}$ partitioned into
 residue classes modulo $8$ for $\{0, 1, \dots, 695\}$, and
 $\{696, 697, \dots, 995\}$,
 the design is generated from

\adfLgap %ADFvfyDesignStart
$(696, 281, 202, 509)$,
$(696, 527, 565, 302)$,
$(696, 690, 334, 112)$,\adfsplit
$(696, 572, 243, 560)$,
$(696, 146, 81, 576)$,
$(696, 534, 655, 529)$,\adfsplit
$(696, 155, 429, 684)$,
$(696, 76, 567, 259)$,
$(697, 380, 474, 488)$,\adfsplit
$(697, 455, 219, 100)$,
$(697, 69, 391, 510)$,
$(697, 113, 475, 365)$,\adfsplit
$(697, 369, 397, 314)$,
$(697, 302, 444, 424)$,
$(697, 22, 337, 0)$,\adfsplit
$(697, 34, 543, 59)$,
$(698, 207, 238, 505)$,
$(698, 346, 329, 160)$,\adfsplit
$(698, 410, 189, 609)$,
$(698, 342, 258, 157)$,
$(698, 365, 134, 679)$,\adfsplit
$(698, 260, 611, 312)$,
$(698, 52, 119, 139)$,
$(698, 584, 132, 291)$,\adfsplit
$(699, 105, 658, 664)$,
$(699, 95, 235, 169)$,
$(699, 530, 504, 252)$,\adfsplit
$(699, 188, 142, 450)$,
$(699, 493, 175, 462)$,
$(699, 56, 495, 53)$,\adfsplit
$(699, 261, 131, 185)$,
$(699, 676, 110, 3)$,
$(700, 48, 201, 189)$,\adfsplit
$(700, 615, 385, 596)$,
$(700, 286, 143, 154)$,
$(700, 508, 570, 512)$,\adfsplit
$(700, 510, 132, 35)$,
$(700, 31, 77, 257)$,
$(700, 74, 603, 109)$,\adfsplit
$(700, 91, 136, 38)$,
$(701, 589, 366, 235)$,
$(701, 577, 368, 477)$,\adfsplit
$(701, 636, 298, 432)$,
$(701, 513, 591, 90)$,
$(701, 116, 243, 46)$,\adfsplit
$(701, 463, 11, 544)$,
$(701, 422, 532, 383)$,
$(701, 314, 173, 329)$,\adfsplit
$(702, 430, 221, 313)$,
$(702, 88, 549, 177)$,
$(702, 349, 467, 591)$,\adfsplit
$(702, 307, 58, 686)$,
$(702, 156, 654, 528)$,
$(702, 146, 295, 292)$,\adfsplit
$(702, 449, 402, 267)$,
$(702, 680, 20, 47)$,
$(703, 273, 400, 187)$,\adfsplit
$(703, 314, 190, 69)$,
$(703, 385, 603, 263)$,
$(703, 656, 324, 655)$,\adfsplit
$(703, 162, 183, 205)$,
$(703, 563, 269, 144)$,
$(703, 548, 82, 174)$,\adfsplit
$(703, 254, 617, 220)$,
$(704, 489, 651, 194)$,
$(704, 432, 148, 34)$,\adfsplit
$(704, 41, 357, 234)$,
$(704, 390, 229, 660)$,
$(704, 40, 343, 341)$,\adfsplit
$(704, 191, 83, 470)$,
$(704, 32, 382, 307)$,
$(704, 409, 164, 399)$,\adfsplit
$(705, 108, 5, 610)$,
$(705, 613, 254, 548)$,
$(705, 576, 499, 52)$,\adfsplit
$(705, 176, 105, 635)$,
$(705, 103, 137, 334)$,
$(705, 97, 378, 159)$,\adfsplit
$(705, 359, 630, 184)$,
$(705, 291, 381, 98)$,
$(706, 622, 562, 28)$,\adfsplit
$(706, 402, 245, 112)$,
$(706, 686, 93, 631)$,
$(706, 407, 308, 523)$,\adfsplit
$(706, 684, 659, 615)$,
$(706, 685, 627, 601)$,
$(706, 242, 320, 561)$,\adfsplit
$(706, 270, 137, 24)$,
$(707, 444, 186, 683)$,
$(707, 375, 597, 89)$,\adfsplit
$(707, 218, 169, 560)$,
$(707, 239, 94, 341)$,
$(707, 178, 150, 307)$,\adfsplit
$(707, 328, 110, 127)$,
$(707, 99, 244, 249)$,
$(707, 168, 445, 68)$,\adfsplit
$(708, 638, 53, 336)$,
$(708, 381, 424, 339)$,
$(708, 547, 294, 181)$,\adfsplit
$(708, 511, 452, 118)$,
$(708, 82, 89, 599)$,
$(708, 612, 681, 323)$,\adfsplit
$(708, 426, 652, 512)$,
$(708, 650, 591, 145)$,
$(709, 189, 499, 680)$,\adfsplit
$(709, 467, 517, 17)$,
$(709, 96, 338, 339)$,
$(709, 324, 186, 23)$,\adfsplit
$(709, 124, 520, 553)$,
$(709, 206, 81, 679)$,
$(709, 567, 394, 582)$,\adfsplit
$(709, 694, 284, 653)$,
$(710, 342, 512, 173)$,
$(710, 535, 603, 254)$,\adfsplit
$(710, 547, 430, 10)$,
$(710, 231, 213, 472)$,
$(710, 361, 661, 484)$,\adfsplit
$(710, 42, 353, 524)$,
$(710, 122, 465, 228)$,
$(710, 360, 491, 695)$,\adfsplit
$(711, 611, 404, 560)$,
$(711, 158, 111, 75)$,
$(711, 671, 574, 564)$,\adfsplit
$(711, 0, 436, 261)$,
$(711, 79, 497, 114)$,
$(711, 256, 637, 226)$,\adfsplit
$(711, 510, 50, 481)$,
$(711, 197, 609, 91)$,
$(712, 173, 25, 636)$,\adfsplit
$(712, 135, 466, 621)$,
$(712, 74, 620, 72)$,
$(712, 580, 225, 571)$,\adfsplit
$(712, 570, 431, 142)$,
$(712, 632, 590, 123)$,
$(712, 256, 305, 102)$,\adfsplit
$(712, 275, 223, 421)$,
$(713, 557, 122, 396)$,
$(713, 52, 570, 381)$,\adfsplit
$(713, 670, 417, 347)$,
$(713, 577, 13, 150)$,
$(713, 202, 164, 398)$,\adfsplit
$(713, 344, 615, 425)$,
$(713, 336, 451, 623)$,
$(713, 448, 319, 99)$,\adfsplit
$(714, 615, 136, 585)$,
$(714, 566, 512, 475)$,
$(714, 500, 237, 294)$,\adfsplit
$(714, 268, 598, 41)$,
$(714, 34, 361, 247)$,
$(714, 517, 672, 651)$,\adfsplit
$(714, 335, 450, 77)$,
$(714, 122, 660, 35)$,
$(715, 197, 160, 411)$,\adfsplit
$(715, 475, 45, 482)$,
$(715, 260, 449, 455)$,
$(715, 175, 652, 633)$,\adfsplit
$(715, 636, 471, 426)$,
$(715, 625, 440, 678)$,
$(715, 96, 278, 683)$,\adfsplit
$(715, 322, 94, 493)$,
$(716, 195, 470, 154)$,
$(716, 622, 77, 0)$,\adfsplit
$(716, 567, 161, 429)$,
$(716, 347, 652, 253)$,
$(716, 338, 56, 191)$,\adfsplit
$(716, 102, 625, 163)$,
$(716, 420, 177, 330)$,
$(716, 535, 640, 524)$,\adfsplit
$(717, 565, 664, 338)$,
$(717, 56, 366, 586)$,
$(717, 617, 690, 347)$,\adfsplit
$(717, 660, 221, 483)$,
$(717, 672, 289, 20)$,
$(717, 573, 71, 403)$,\adfsplit
$(717, 687, 206, 4)$,
$(0, 39, 66, 819)$,
$(0, 292, 585, 743)$,\adfsplit
$(0, 29, 212, 919)$,
$(0, 57, 493, 843)$,
$(0, 190, 463, 793)$,\adfsplit
$(0, 467, 673, 918)$,
$(0, 179, 321, 742)$,
$(0, 73, 382, 868)$,\adfsplit
$(0, 217, 371, 968)$,
$(0, 18, 663, 943)$,
$(0, 165, 471, 768)$,\adfsplit
$(0, 147, 442, 844)$,
$(0, 211, 306, 720)$,
$(0, 333, 635, 770)$,\adfsplit
$(0, 9, 76, 945)$,
$(0, 291, 617, 794)$,
$(0, 50, 375, 389)$,\adfsplit
$(0, 75, 489, 944)$,
$(0, 93, 385, 744)$,
$(0, 82, 361, 994)$,\adfsplit
$(0, 13, 95, 845)$,
$(0, 105, 516, 820)$,
$(0, 101, 633, 769)$,\adfsplit
$(0, 603, 607, 995)$,
$(0, 164, 345, 870)$,
$(0, 445, 505, 745)$,\adfsplit
$(0, 174, 348, 522)$,
$(1, 175, 349, 523)$

%ADFvfyBlocksEnd
\adfLgap \noindent by the mapping:
$x \mapsto x + 2 j \adfmod{696}$ for $x < 696$,
$x \mapsto (x - 696 + 25 j \adfmod{300}) + 696$ for $x \ge 696$,
$0 \le j < 348$
 for the first 201 blocks,
$0 \le j < 87$
 for the last two blocks.
\ADFvfyParStart{(996, ((201, 348, ((696, 2), (300, 25))), (2, 87, ((696, 2), (300, 25)))), ((87, 8), (300, 1)))} %ADFvfyParEnd
% End of 87^8 300^1
%%%%%%%%%%%%%%%%%%%%%%%%%%%%%%%%%%%%%%%%%%%%%%%%%%%%%%%%%%%%%%%%%%%%%%%%%%%%%%%%%%%%%%%%%%
%%%%%%%%%%%%%%%%%%%%%%%%%%%%%%%%%%%%%%%%%%%%%%%%%%%%%%%%%%%%%%%%%%%%%%%%%%%%%%%%%%%%%%%%%%

%%%%%%%%%%%%%%%%%%%%%%%%%%%%%%%%%%%%%%%%%%%%%%%%%%%%%%%%%%%%%%%%%%%%%%%%%%%%%%%%%%%%%%%%%%
%%%%%%%%%%%%%%%%%%%%%%%%%%%%%%%%%%%%%%%%%%%%%%%%%%%%%%%%%%%%%%%%%%%%%%%%%%%%%%%%%%%%%%%%%%
\section{4-GDDs for the proof of Lemma \ref{lem:4-GDD 93^8 m^1}}
\label{app:4-GDD 93^8 m^1}
\adfnull{
$ 93^8 318^1 $ and
$ 93^8 321^1 $.
}

% Charlotte:GDD4-1-3-5-mod-6-TeX-gen-A:HITS-fun:4.10
\adfDgap
%ADFvfyBlocksStart {93,93,93,93,93,93,93,93,318}
\noindent{\boldmath $ 93^{8} 318^{1} $}~
With the point set $Z_{1062}$ partitioned into
 residue classes modulo $8$ for $\{0, 1, \dots, 743\}$, and
 $\{744, 745, \dots, 1061\}$,
 the design is generated from

\adfLgap %ADFvfyDesignStart
$(1056, 430, 695, 561)$,
$(1057, 733, 176, 420)$,
$(744, 619, 112, 622)$,\adfsplit
$(744, 689, 363, 464)$,
$(744, 117, 44, 630)$,
$(744, 671, 626, 696)$,\adfsplit
$(744, 97, 439, 682)$,
$(744, 66, 725, 708)$,
$(744, 148, 231, 345)$,\adfsplit
$(744, 662, 301, 227)$,
$(745, 716, 280, 217)$,
$(745, 511, 540, 29)$,\adfsplit
$(745, 75, 344, 38)$,
$(745, 619, 410, 637)$,
$(745, 201, 580, 202)$,\adfsplit
$(745, 423, 382, 11)$,
$(745, 593, 270, 216)$,
$(745, 90, 405, 623)$,\adfsplit
$(746, 652, 110, 91)$,
$(746, 26, 164, 103)$,
$(746, 15, 126, 184)$,\adfsplit
$(746, 12, 557, 546)$,
$(746, 704, 59, 370)$,
$(746, 286, 81, 552)$,\adfsplit
$(746, 263, 305, 109)$,
$(746, 621, 385, 555)$,
$(747, 592, 38, 533)$,\adfsplit
$(747, 517, 342, 457)$,
$(747, 183, 273, 322)$,
$(747, 36, 367, 3)$,\adfsplit
$(747, 32, 236, 286)$,
$(747, 443, 50, 405)$,
$(747, 162, 172, 456)$,\adfsplit
$(747, 569, 427, 719)$,
$(748, 287, 8, 34)$,
$(748, 111, 125, 579)$,\adfsplit
$(748, 81, 522, 301)$,
$(748, 74, 592, 430)$,
$(748, 117, 174, 457)$,\adfsplit
$(748, 355, 350, 708)$,
$(748, 539, 41, 436)$,
$(748, 164, 480, 343)$,\adfsplit
$(749, 468, 254, 226)$,
$(749, 553, 47, 392)$,
$(749, 257, 724, 101)$,\adfsplit
$(749, 366, 530, 3)$,
$(749, 19, 40, 495)$,
$(749, 177, 96, 131)$,\adfsplit
$(749, 477, 20, 642)$,
$(749, 478, 493, 487)$,
$(750, 302, 99, 170)$,\adfsplit
$(750, 652, 739, 439)$,
$(750, 81, 61, 342)$,
$(750, 191, 20, 317)$,\adfsplit
$(750, 491, 632, 570)$,
$(750, 108, 449, 208)$,
$(750, 183, 381, 433)$,\adfsplit
$(750, 106, 358, 672)$,
$(751, 545, 158, 643)$,
$(751, 261, 409, 358)$,\adfsplit
$(751, 255, 54, 576)$,
$(751, 506, 132, 57)$,
$(751, 154, 411, 29)$,\adfsplit
$(751, 479, 592, 426)$,
$(751, 421, 151, 268)$,
$(751, 512, 539, 500)$,\adfsplit
$(752, 649, 565, 15)$,
$(752, 551, 125, 680)$,
$(752, 683, 588, 154)$,\adfsplit
$(752, 285, 722, 308)$,
$(752, 705, 376, 675)$,
$(752, 700, 331, 102)$,\adfsplit
$(752, 22, 594, 24)$,
$(752, 607, 182, 377)$,
$(753, 537, 581, 130)$,\adfsplit
$(753, 652, 583, 193)$,
$(753, 644, 735, 568)$,
$(753, 215, 339, 72)$,\adfsplit
$(753, 156, 707, 238)$,
$(753, 85, 606, 258)$,
$(753, 98, 737, 715)$,\adfsplit
$(753, 464, 597, 542)$,
$(754, 122, 373, 257)$,
$(754, 53, 40, 462)$,\adfsplit
$(754, 444, 552, 297)$,
$(754, 505, 299, 207)$,
$(754, 479, 547, 188)$,\adfsplit
$(754, 487, 453, 518)$,
$(754, 435, 316, 34)$,
$(754, 320, 718, 234)$,\adfsplit
$(755, 216, 379, 93)$,
$(755, 219, 223, 404)$,
$(755, 704, 234, 54)$,\adfsplit
$(0, 107, 247, 899)$,
$(0, 151, 333, 859)$,
$(0, 36, 191, 1041)$,\adfsplit
$(0, 43, 599, 937)$,
$(0, 118, 419, 846)$,
$(0, 106, 263, 911)$,\adfsplit
$(0, 207, 509, 834)$,
$(0, 109, 614, 964)$,
$(0, 7, 219, 912)$,\adfsplit
$(0, 149, 427, 990)$,
$(0, 67, 473, 938)$,
$(0, 47, 305, 394)$,\adfsplit
$(0, 93, 417, 977)$,
$(0, 177, 405, 847)$,
$(0, 186, 372, 558)$

%ADFvfyBlocksEnd
\adfLgap \noindent by the mapping:
$x \mapsto x +  j \adfmod{744}$ for $x < 744$,
$x \mapsto (x - 744 + 13 j \adfmod{312}) + 744$ for $744 \le x < 1056$,
$x \mapsto (x + 2 j \adfmod{6}) + 1056$ for $x \ge 1056$,
$0 \le j < 744$
 for the first 107 blocks,
$0 \le j < 186$
 for the last block.
\ADFvfyParStart{(1062, ((107, 744, ((744, 1), (312, 13), (6, 2))), (1, 186, ((744, 1), (312, 13), (6, 2)))), ((93, 8), (318, 1)))} %ADFvfyParEnd
% End of 93^8 318^1
%%%%%%%%%%%%%%%%%%%%%%%%%%%%%%%%%%%%%%%%%%%%%%%%%%%%%%%%%%%%%%%%%%%%%%%%%%%%%%%%%%%%%%%%%%
%%%%%%%%%%%%%%%%%%%%%%%%%%%%%%%%%%%%%%%%%%%%%%%%%%%%%%%%%%%%%%%%%%%%%%%%%%%%%%%%%%%%%%%%%%

% Charlotte:GDD4-1-3-5-mod-6-TeX-gen-A:HITS-fun:4.10
\adfDgap
%ADFvfyBlocksStart {93,93,93,93,93,93,93,93,321}
\noindent{\boldmath $ 93^{8} 321^{1} $}~
With the point set $Z_{1065}$ partitioned into
 residue classes modulo $8$ for $\{0, 1, \dots, 743\}$, and
 $\{744, 745, \dots, 1064\}$,
 the design is generated from

\adfLgap %ADFvfyDesignStart
$(1056, 579, 533, 740)$,
$(1056, 709, 90, 742)$,
$(1057, 683, 18, 424)$,\adfsplit
$(1057, 123, 367, 602)$,
$(1058, 51, 190, 719)$,
$(1058, 218, 703, 342)$,\adfsplit
$(744, 442, 61, 654)$,
$(744, 166, 231, 332)$,
$(744, 27, 57, 642)$,\adfsplit
$(744, 21, 638, 347)$,
$(744, 152, 463, 532)$,
$(744, 300, 648, 601)$,\adfsplit
$(744, 2, 473, 115)$,
$(744, 328, 437, 263)$,
$(745, 326, 272, 419)$,\adfsplit
$(745, 454, 393, 192)$,
$(745, 433, 434, 148)$,
$(745, 719, 380, 733)$,\adfsplit
$(745, 303, 90, 520)$,
$(745, 31, 538, 150)$,
$(745, 19, 329, 29)$,\adfsplit
$(745, 213, 396, 219)$,
$(746, 189, 184, 739)$,
$(746, 528, 20, 135)$,\adfsplit
$(746, 373, 70, 272)$,
$(746, 362, 485, 609)$,
$(746, 446, 71, 707)$,\adfsplit
$(746, 498, 540, 3)$,
$(746, 97, 391, 316)$,
$(746, 610, 713, 126)$,\adfsplit
$(747, 397, 415, 566)$,
$(747, 159, 140, 10)$,
$(747, 238, 609, 71)$,\adfsplit
$(747, 28, 48, 43)$,
$(747, 97, 66, 347)$,
$(747, 416, 492, 437)$,\adfsplit
$(747, 232, 194, 150)$,
$(747, 189, 593, 699)$,
$(748, 471, 529, 568)$,\adfsplit
$(748, 281, 2, 512)$,
$(748, 430, 312, 125)$,
$(748, 55, 630, 402)$,\adfsplit
$(748, 513, 383, 419)$,
$(748, 331, 541, 572)$,
$(748, 244, 302, 51)$,\adfsplit
$(748, 154, 213, 564)$,
$(749, 418, 457, 239)$,
$(749, 502, 474, 224)$,\adfsplit
$(749, 231, 243, 492)$,
$(749, 112, 709, 506)$,
$(749, 153, 55, 100)$,\adfsplit
$(749, 187, 470, 668)$,
$(749, 125, 552, 353)$,
$(749, 275, 477, 294)$,\adfsplit
$(750, 56, 612, 483)$,
$(750, 481, 376, 347)$,
$(750, 327, 254, 165)$,\adfsplit
$(750, 412, 626, 125)$,
$(750, 451, 394, 609)$,
$(750, 407, 718, 37)$,\adfsplit
$(750, 620, 666, 630)$,
$(750, 367, 89, 96)$,
$(751, 119, 414, 132)$,\adfsplit
$(751, 394, 57, 220)$,
$(751, 491, 645, 16)$,
$(751, 704, 548, 471)$,\adfsplit
$(751, 665, 310, 157)$,
$(751, 602, 168, 350)$,
$(751, 247, 629, 219)$,\adfsplit
$(751, 235, 234, 289)$,
$(752, 735, 33, 83)$,
$(752, 723, 430, 546)$,\adfsplit
$(752, 554, 278, 192)$,
$(752, 64, 58, 407)$,
$(752, 222, 116, 541)$,\adfsplit
$(752, 689, 156, 19)$,
$(752, 632, 388, 597)$,
$(752, 511, 197, 25)$,\adfsplit
$(753, 666, 568, 252)$,
$(753, 4, 277, 441)$,
$(753, 645, 47, 99)$,\adfsplit
$(753, 31, 434, 545)$,
$(753, 155, 461, 15)$,
$(753, 80, 553, 230)$,\adfsplit
$(753, 548, 358, 283)$,
$(753, 510, 730, 288)$,
$(754, 379, 262, 345)$,\adfsplit
$(754, 125, 18, 15)$,
$(754, 575, 219, 706)$,
$(754, 88, 14, 329)$,\adfsplit
$(754, 141, 720, 25)$,
$(754, 246, 636, 176)$,
$(754, 668, 175, 85)$,\adfsplit
$(754, 100, 722, 587)$,
$(755, 642, 719, 19)$,
$(755, 96, 289, 597)$,\adfsplit
$(755, 438, 665, 461)$,
$(755, 495, 205, 724)$,
$(755, 411, 732, 206)$,\adfsplit
$(755, 489, 347, 310)$,
$(755, 392, 7, 394)$,
$(755, 44, 578, 208)$,\adfsplit
$(756, 620, 448, 590)$,
$(756, 442, 120, 737)$,
$(756, 43, 726, 311)$,\adfsplit
$(756, 114, 557, 516)$,
$(756, 59, 439, 97)$,
$(756, 310, 225, 549)$,\adfsplit
$(756, 584, 613, 171)$,
$(756, 674, 615, 76)$,
$(757, 30, 671, 609)$,\adfsplit
$(757, 392, 437, 259)$,
$(757, 38, 648, 179)$,
$(757, 306, 687, 333)$,\adfsplit
$(757, 250, 151, 517)$,
$(757, 404, 737, 387)$,
$(757, 612, 472, 313)$,\adfsplit
$(757, 554, 100, 214)$,
$(758, 329, 447, 400)$,
$(758, 150, 523, 249)$,\adfsplit
$(758, 422, 372, 169)$,
$(758, 338, 93, 7)$,
$(758, 570, 395, 728)$,\adfsplit
$(758, 349, 22, 460)$,
$(758, 384, 603, 629)$,
$(758, 263, 130, 476)$,\adfsplit
$(759, 546, 534, 237)$,
$(759, 389, 537, 651)$,
$(759, 316, 715, 312)$,\adfsplit
$(759, 109, 639, 641)$,
$(759, 695, 320, 404)$,
$(759, 97, 386, 86)$,\adfsplit
$(759, 660, 64, 199)$,
$(759, 130, 670, 107)$,
$(760, 36, 648, 555)$,\adfsplit
$(760, 193, 668, 477)$,
$(760, 280, 229, 638)$,
$(760, 153, 659, 334)$,\adfsplit
$(760, 535, 126, 355)$,
$(760, 330, 4, 239)$,
$(760, 377, 629, 178)$,\adfsplit
$(760, 506, 80, 495)$,
$(761, 356, 150, 545)$,
$(761, 436, 341, 607)$,\adfsplit
$(761, 501, 528, 334)$,
$(761, 64, 708, 159)$,
$(761, 290, 603, 265)$,\adfsplit
$(761, 346, 503, 499)$,
$(761, 666, 61, 632)$,
$(761, 561, 518, 395)$,\adfsplit
$(762, 225, 288, 653)$,
$(762, 364, 651, 274)$,
$(762, 228, 439, 406)$,\adfsplit
$(762, 737, 650, 541)$,
$(762, 443, 318, 95)$,
$(762, 640, 721, 619)$,\adfsplit
$(762, 429, 495, 536)$,
$(762, 230, 332, 162)$,
$(763, 73, 173, 695)$,\adfsplit
$(763, 86, 471, 8)$,
$(763, 7, 322, 465)$,
$(763, 493, 636, 184)$,\adfsplit
$(763, 220, 314, 451)$,
$(763, 366, 620, 563)$,
$(763, 450, 0, 401)$,\adfsplit
$(763, 165, 411, 478)$,
$(764, 499, 25, 684)$,
$(764, 657, 310, 335)$,\adfsplit
$(764, 226, 733, 72)$,
$(764, 628, 29, 110)$,
$(764, 131, 199, 690)$,\adfsplit
$(764, 69, 296, 387)$,
$(764, 351, 654, 424)$,
$(764, 521, 692, 626)$,\adfsplit
$(765, 209, 608, 196)$,
$(765, 273, 355, 333)$,
$(765, 311, 90, 389)$,\adfsplit
$(765, 458, 472, 133)$,
$(765, 291, 583, 121)$,
$(765, 518, 250, 587)$,\adfsplit
$(765, 558, 420, 312)$,
$(765, 188, 430, 735)$,
$(766, 156, 642, 319)$,\adfsplit
$(766, 538, 17, 173)$,
$(766, 659, 650, 97)$,
$(766, 428, 624, 615)$,\adfsplit
$(766, 430, 501, 585)$,
$(766, 739, 688, 479)$,
$(766, 470, 196, 315)$,\adfsplit
$(0, 22, 389, 974)$,
$(0, 131, 543, 767)$,
$(0, 79, 436, 691)$,\adfsplit
$(0, 89, 366, 846)$,
$(0, 17, 571, 769)$,
$(0, 175, 324, 949)$,\adfsplit
$(0, 35, 682, 924)$,
$(0, 110, 599, 795)$,
$(0, 52, 469, 871)$,\adfsplit
$(0, 60, 506, 976)$,
$(0, 266, 657, 1029)$,
$(0, 277, 297, 847)$,\adfsplit
$(0, 18, 307, 873)$,
$(0, 162, 627, 821)$,
$(0, 126, 623, 1055)$,\adfsplit
$(1, 189, 415, 977)$,
$(0, 7, 270, 1053)$,
$(0, 26, 549, 1027)$,\adfsplit
$(0, 67, 180, 1054)$,
$(0, 217, 707, 923)$,
$(0, 233, 701, 897)$,\adfsplit
$(0, 321, 459, 845)$,
$(0, 477, 603, 1002)$,
$(0, 383, 729, 872)$,\adfsplit
$(0, 3, 505, 1028)$,
$(1, 71, 595, 846)$,
$(0, 186, 372, 558)$,\adfsplit
$(1, 187, 373, 559)$

%ADFvfyBlocksEnd
\adfLgap \noindent by the mapping:
$x \mapsto x + 2 j \adfmod{744}$ for $x < 744$,
$x \mapsto (x - 744 + 26 j \adfmod{312}) + 744$ for $744 \le x < 1056$,
$x \mapsto (x - 1056 + 3 j \adfmod{9}) + 1056$ for $x \ge 1056$,
$0 \le j < 372$
 for the first 215 blocks,
$0 \le j < 93$
 for the last two blocks.
\ADFvfyParStart{(1065, ((215, 372, ((744, 2), (312, 26), (9, 3))), (2, 93, ((744, 2), (312, 26), (9, 3)))), ((93, 8), (321, 1)))} %ADFvfyParEnd
% End of 93^8 321^1
%%%%%%%%%%%%%%%%%%%%%%%%%%%%%%%%%%%%%%%%%%%%%%%%%%%%%%%%%%%%%%%%%%%%%%%%%%%%%%%%%%%%%%%%%%
%%%%%%%%%%%%%%%%%%%%%%%%%%%%%%%%%%%%%%%%%%%%%%%%%%%%%%%%%%%%%%%%%%%%%%%%%%%%%%%%%%%%%%%%%%

%%%%%%%%%%%%%%%%%%%%%%%%%%%%%%%%%%%%%%%%%%%%%%%%%%%%%%%%%%%%%%%%%%%%%%%%%%%%%%%%%%%%%%%%%%
%%%%%%%%%%%%%%%%%%%%%%%%%%%%%%%%%%%%%%%%%%%%%%%%%%%%%%%%%%%%%%%%%%%%%%%%%%%%%%%%%%%%%%%%%%
\section{4-GDDs for the proof of Lemma \ref{lem:4-GDD 13^u m^1}}
\label{app:4-GDD 13^u m^1}
\adfnull{
$ 13^{12} 7^1 $,
$ 13^{12} 10^1 $,
$ 13^9 10^1 $,
$ 13^9 16^1 $,
$ 13^9 22^1 $,
$ 13^9 28^1 $,
$ 13^9 34^1 $,
$ 13^9 40^1 $ and
$ 13^9 46^1 $.
}

% Charlotte:GDD4-1-3-5-mod-6-TeX-gen-A:HITS-fun:4.10
\adfDgap
%ADFvfyBlocksStart {13,13,13,13,13,13,13,13,13,13,13,13,7}
\noindent{\boldmath $ 13^{12} 7^{1} $}~
With the point set $Z_{163}$ partitioned into
 residue classes modulo $12$ for $\{0, 1, \dots, 155\}$, and
 $\{156, 157, \dots, 162\}$,
 the design is generated from

\adfLgap %ADFvfyDesignStart
$(156, 47, 154, 123)$,
$(156, 120, 139, 32)$,
$(157, 139, 153, 136)$,
$(157, 42, 71, 98)$,\adfsplit
$(41, 146, 132, 75)$,
$(0, 1, 2, 94)$,
$(0, 4, 9, 150)$,
$(0, 7, 118, 126)$,\adfsplit
$(0, 11, 13, 74)$,
$(0, 16, 128, 149)$,
$(0, 20, 42, 67)$,
$(0, 23, 76, 102)$,\adfsplit
$(0, 18, 58, 89)$,
$(0, 33, 41, 106)$,
$(0, 27, 110, 145)$,
$(0, 43, 86, 141)$,\adfsplit
$(0, 32, 66, 147)$,
$(0, 131, 135, 151)$,
$(0, 53, 109, 139)$,
$(0, 59, 105, 123)$,\adfsplit
$(0, 57, 101, 111)$,
$(0, 69, 97, 119)$,
$(0, 61, 87, 127)$,
$(0, 79, 121, 153)$,\adfsplit
$(0, 75, 137, 143)$,
$(0, 39, 78, 117)$,
$(162, 0, 52, 104)$,
$(162, 1, 53, 105)$

%ADFvfyBlocksEnd
\adfLgap \noindent by the mapping:
$x \mapsto x + 2 j \adfmod{156}$ for $x < 156$,
$x \mapsto (x + 2 j \adfmod{6}) + 156$ for $156 \le x < 162$,
$162 \mapsto 162$,
$0 \le j < 78$
 for the first 25 blocks,
$0 \le j < 39$
 for the next block,
$0 \le j < 26$
 for the last two blocks.
\ADFvfyParStart{(163, ((25, 78, ((156, 2), (6, 2), (1, 1))), (1, 39, ((156, 2), (6, 2), (1, 1))), (2, 26, ((156, 2), (6, 2), (1, 1)))), ((13, 12), (7, 1)))} %ADFvfyParEnd
% End of 13^12 7^1
%%%%%%%%%%%%%%%%%%%%%%%%%%%%%%%%%%%%%%%%%%%%%%%%%%%%%%%%%%%%%%%%%%%%%%%%%%%%%%%%%%%%%%%%%%
%%%%%%%%%%%%%%%%%%%%%%%%%%%%%%%%%%%%%%%%%%%%%%%%%%%%%%%%%%%%%%%%%%%%%%%%%%%%%%%%%%%%%%%%%%

% Charlotte:GDD4-1-3-5-mod-6-TeX-gen-A:HITS-fun:4.10
\adfDgap
%ADFvfyBlocksStart {13,13,13,13,13,13,13,13,13,13,13,13,10}
\noindent{\boldmath $ 13^{12} 10^{1} $}~
With the point set $Z_{166}$ partitioned into
 residue classes modulo $12$ for $\{0, 1, \dots, 155\}$, and
 $\{156, 157, \dots, 165\}$,
 the design is generated from

\adfLgap %ADFvfyDesignStart
$(156, 145, 47, 135)$,
$(157, 1, 98, 42)$,
$(158, 86, 64, 60)$,
$(0, 1, 3, 136)$,\adfsplit
$(0, 5, 11, 114)$,
$(0, 7, 15, 92)$,
$(0, 9, 25, 43)$,
$(0, 13, 40, 95)$,\adfsplit
$(0, 28, 57, 119)$,
$(0, 14, 33, 87)$,
$(0, 32, 67, 118)$,
$(0, 30, 75, 106)$,\adfsplit
$(0, 17, 63, 107)$,
$(0, 39, 78, 117)$,
$(165, 0, 52, 104)$

%ADFvfyBlocksEnd
\adfLgap \noindent by the mapping:
$x \mapsto x +  j \adfmod{156}$ for $x < 156$,
$x \mapsto (x - 156 + 3 j \adfmod{9}) + 156$ for $156 \le x < 165$,
$165 \mapsto 165$,
$0 \le j < 156$
 for the first 13 blocks,
$0 \le j < 39$
 for the next block,
$0 \le j < 52$
 for the last block.
\ADFvfyParStart{(166, ((13, 156, ((156, 1), (9, 3), (1, 1))), (1, 39, ((156, 1), (9, 3), (1, 1))), (1, 52, ((156, 1), (9, 3), (1, 1)))), ((13, 12), (10, 1)))} %ADFvfyParEnd
% End of 13^12 10^1
%%%%%%%%%%%%%%%%%%%%%%%%%%%%%%%%%%%%%%%%%%%%%%%%%%%%%%%%%%%%%%%%%%%%%%%%%%%%%%%%%%%%%%%%%%
%%%%%%%%%%%%%%%%%%%%%%%%%%%%%%%%%%%%%%%%%%%%%%%%%%%%%%%%%%%%%%%%%%%%%%%%%%%%%%%%%%%%%%%%%%

% Charlotte:GDD4-1-3-5-mod-6-TeX-gen-A:HITS-fun:4.10
\adfDgap
%ADFvfyBlocksStart {13,13,13,13,13,13,13,13,13,10}
\noindent{\boldmath $ 13^{9} 10^{1} $}~
With the point set $Z_{127}$ partitioned into
 residue classes modulo $9$ for $\{0, 1, \dots, 116\}$, and
 $\{117, 118, \dots, 126\}$,
 the design is generated from

\adfLgap %ADFvfyDesignStart
$(117, 77, 94, 63)$,
$(118, 96, 109, 17)$,
$(119, 110, 45, 16)$,
$(0, 1, 3, 33)$,\adfsplit
$(0, 4, 10, 80)$,
$(0, 8, 20, 64)$,
$(0, 15, 49, 75)$,
$(0, 21, 43, 71)$,\adfsplit
$(0, 7, 55, 66)$,
$(0, 5, 24, 40)$,
$(126, 0, 39, 78)$

%ADFvfyBlocksEnd
\adfLgap \noindent by the mapping:
$x \mapsto x +  j \adfmod{117}$ for $x < 117$,
$x \mapsto (x + 3 j \adfmod{9}) + 117$ for $117 \le x < 126$,
$126 \mapsto 126$,
$0 \le j < 117$
 for the first ten blocks,
$0 \le j < 39$
 for the last block.
\ADFvfyParStart{(127, ((10, 117, ((117, 1), (9, 3), (1, 1))), (1, 39, ((117, 1), (9, 3), (1, 1)))), ((13, 9), (10, 1)))} %ADFvfyParEnd
% End of 13^9 10^1
%%%%%%%%%%%%%%%%%%%%%%%%%%%%%%%%%%%%%%%%%%%%%%%%%%%%%%%%%%%%%%%%%%%%%%%%%%%%%%%%%%%%%%%%%%
%%%%%%%%%%%%%%%%%%%%%%%%%%%%%%%%%%%%%%%%%%%%%%%%%%%%%%%%%%%%%%%%%%%%%%%%%%%%%%%%%%%%%%%%%%

% Charlotte:GDD4-1-3-5-mod-6-TeX-gen-A:HITS-fun:4.10
\adfDgap
%ADFvfyBlocksStart {13,13,13,13,13,13,13,13,13,16}
\noindent{\boldmath $ 13^{9} 16^{1} $}~
With the point set $Z_{133}$ partitioned into
 residue classes modulo $9$ for $\{0, 1, \dots, 116\}$, and
 $\{117, 118, \dots, 132\}$,
 the design is generated from

\adfLgap %ADFvfyDesignStart
$(117, 68, 78, 22)$,
$(118, 70, 102, 71)$,
$(119, 4, 116, 27)$,
$(120, 4, 24, 2)$,\adfsplit
$(121, 38, 54, 4)$,
$(0, 3, 7, 82)$,
$(0, 6, 49, 64)$,
$(0, 13, 37, 70)$,\adfsplit
$(0, 8, 25, 73)$,
$(0, 12, 26, 88)$,
$(0, 11, 30, 51)$,
$(132, 0, 39, 78)$

%ADFvfyBlocksEnd
\adfLgap \noindent by the mapping:
$x \mapsto x +  j \adfmod{117}$ for $x < 117$,
$x \mapsto (x - 117 + 5 j \adfmod{15}) + 117$ for $117 \le x < 132$,
$132 \mapsto 132$,
$0 \le j < 117$
 for the first 11 blocks,
$0 \le j < 39$
 for the last block.
\ADFvfyParStart{(133, ((11, 117, ((117, 1), (15, 5), (1, 1))), (1, 39, ((117, 1), (15, 5), (1, 1)))), ((13, 9), (16, 1)))} %ADFvfyParEnd
% End of 13^9 16^1
%%%%%%%%%%%%%%%%%%%%%%%%%%%%%%%%%%%%%%%%%%%%%%%%%%%%%%%%%%%%%%%%%%%%%%%%%%%%%%%%%%%%%%%%%%
%%%%%%%%%%%%%%%%%%%%%%%%%%%%%%%%%%%%%%%%%%%%%%%%%%%%%%%%%%%%%%%%%%%%%%%%%%%%%%%%%%%%%%%%%%

% Charlotte:GDD4-1-3-5-mod-6-TeX-gen-A:HITS-fun:4.10
\adfDgap
%ADFvfyBlocksStart {13,13,13,13,13,13,13,13,13,22}
\noindent{\boldmath $ 13^{9} 22^{1} $}~
With the point set $Z_{139}$ partitioned into
 residue classes modulo $9$ for $\{0, 1, \dots, 116\}$, and
 $\{117, 118, \dots, 138\}$,
 the design is generated from

\adfLgap %ADFvfyDesignStart
$(117, 36, 80, 1)$,
$(118, 98, 100, 15)$,
$(119, 11, 33, 88)$,
$(120, 5, 18, 34)$,\adfsplit
$(121, 109, 108, 59)$,
$(122, 79, 93, 23)$,
$(123, 87, 77, 46)$,
$(0, 3, 7, 96)$,\adfsplit
$(0, 15, 48, 74)$,
$(0, 5, 51, 57)$,
$(0, 8, 20, 100)$,
$(0, 11, 30, 53)$,\adfsplit
$(138, 0, 39, 78)$

%ADFvfyBlocksEnd
\adfLgap \noindent by the mapping:
$x \mapsto x +  j \adfmod{117}$ for $x < 117$,
$x \mapsto (x - 117 + 7 j \adfmod{21}) + 117$ for $117 \le x < 138$,
$138 \mapsto 138$,
$0 \le j < 117$
 for the first 12 blocks,
$0 \le j < 39$
 for the last block.
\ADFvfyParStart{(139, ((12, 117, ((117, 1), (21, 7), (1, 1))), (1, 39, ((117, 1), (21, 7), (1, 1)))), ((13, 9), (22, 1)))} %ADFvfyParEnd
% End of 13^9 22^1
%%%%%%%%%%%%%%%%%%%%%%%%%%%%%%%%%%%%%%%%%%%%%%%%%%%%%%%%%%%%%%%%%%%%%%%%%%%%%%%%%%%%%%%%%%
%%%%%%%%%%%%%%%%%%%%%%%%%%%%%%%%%%%%%%%%%%%%%%%%%%%%%%%%%%%%%%%%%%%%%%%%%%%%%%%%%%%%%%%%%%

% Charlotte:GDD4-1-3-5-mod-6-TeX-gen-A:HITS-fun:4.10
\adfDgap
%ADFvfyBlocksStart {13,13,13,13,13,13,13,13,13,28}
\noindent{\boldmath $ 13^{9} 28^{1} $}~
With the point set $Z_{145}$ partitioned into
 residue classes modulo $9$ for $\{0, 1, \dots, 116\}$, and
 $\{117, 118, \dots, 144\}$,
 the design is generated from

\adfLgap %ADFvfyDesignStart
$(117, 59, 25, 57)$,
$(118, 74, 12, 55)$,
$(119, 15, 92, 85)$,
$(120, 34, 56, 42)$,\adfsplit
$(121, 29, 3, 13)$,
$(122, 31, 90, 77)$,
$(123, 59, 10, 30)$,
$(0, 1, 5, 124)$,\adfsplit
$(0, 23, 56, 80)$,
$(0, 21, 51, 86)$,
$(0, 6, 48, 73)$,
$(0, 11, 28, 125)$,\adfsplit
$(0, 3, 15, 79)$,
$(144, 0, 39, 78)$

%ADFvfyBlocksEnd
\adfLgap \noindent by the mapping:
$x \mapsto x +  j \adfmod{117}$ for $x < 117$,
$x \mapsto (x - 117 + 9 j \adfmod{27}) + 117$ for $117 \le x < 144$,
$144 \mapsto 144$,
$0 \le j < 117$
 for the first 13 blocks,
$0 \le j < 39$
 for the last block.
\ADFvfyParStart{(145, ((13, 117, ((117, 1), (27, 9), (1, 1))), (1, 39, ((117, 1), (27, 9), (1, 1)))), ((13, 9), (28, 1)))} %ADFvfyParEnd
% End of 13^9 28^1
%%%%%%%%%%%%%%%%%%%%%%%%%%%%%%%%%%%%%%%%%%%%%%%%%%%%%%%%%%%%%%%%%%%%%%%%%%%%%%%%%%%%%%%%%%
%%%%%%%%%%%%%%%%%%%%%%%%%%%%%%%%%%%%%%%%%%%%%%%%%%%%%%%%%%%%%%%%%%%%%%%%%%%%%%%%%%%%%%%%%%

% Charlotte:GDD4-1-3-5-mod-6-TeX-gen-A:HITS-fun:4.10
\adfDgap
%ADFvfyBlocksStart {13,13,13,13,13,13,13,13,13,34}
\noindent{\boldmath $ 13^{9} 34^{1} $}~
With the point set $Z_{151}$ partitioned into
 residue classes modulo $9$ for $\{0, 1, \dots, 116\}$, and
 $\{117, 118, \dots, 150\}$,
 the design is generated from

\adfLgap %ADFvfyDesignStart
$(147, 2, 73, 21)$,
$(144, 48, 95, 52)$,
$(135, 24, 39, 10)$,
$(135, 35, 52, 77)$,\adfsplit
$(135, 13, 47, 63)$,
$(0, 1, 3, 23)$,
$(0, 5, 11, 60)$,
$(0, 12, 33, 73)$,\adfsplit
$(0, 26, 64, 118)$,
$(0, 7, 87, 120)$,
$(0, 8, 93, 124)$,
$(0, 13, 41, 131)$,\adfsplit
$(0, 10, 58, 134)$,
$(0, 31, 66, 132)$,
$(150, 0, 39, 78)$

%ADFvfyBlocksEnd
\adfLgap \noindent by the mapping:
$x \mapsto x +  j \adfmod{117}$ for $x < 117$,
$x \mapsto (x +  j \adfmod{9}) + 117$ for $117 \le x < 126$,
$x \mapsto (x +  j \adfmod{9}) + 126$ for $126 \le x < 135$,
$x \mapsto (x +  j \adfmod{9}) + 135$ for $135 \le x < 144$,
$x \mapsto (x +  j \adfmod{3}) + 144$ for $144 \le x < 147$,
$x \mapsto (x +  j \adfmod{3}) + 147$ for $147 \le x < 150$,
$150 \mapsto 150$,
$0 \le j < 117$
 for the first 14 blocks,
$0 \le j < 39$
 for the last block.
\ADFvfyParStart{(151, ((14, 117, ((117, 1), (9, 1), (9, 1), (9, 1), (3, 1), (3, 1), (1, 1))), (1, 39, ((117, 1), (9, 1), (9, 1), (9, 1), (3, 1), (3, 1), (1, 1)))), ((13, 9), (34, 1)))} %ADFvfyParEnd
% End of 13^9 34^1
%%%%%%%%%%%%%%%%%%%%%%%%%%%%%%%%%%%%%%%%%%%%%%%%%%%%%%%%%%%%%%%%%%%%%%%%%%%%%%%%%%%%%%%%%%
%%%%%%%%%%%%%%%%%%%%%%%%%%%%%%%%%%%%%%%%%%%%%%%%%%%%%%%%%%%%%%%%%%%%%%%%%%%%%%%%%%%%%%%%%%

% Charlotte:GDD4-1-3-5-mod-6-TeX-gen-A:HITS-fun:4.10
\adfDgap
%ADFvfyBlocksStart {13,13,13,13,13,13,13,13,13,40}
\noindent{\boldmath $ 13^{9} 40^{1} $}~
With the point set $Z_{157}$ partitioned into
 residue classes modulo $9$ for $\{0, 1, \dots, 116\}$, and
 $\{117, 118, \dots, 156\}$,
 the design is generated from

\adfLgap %ADFvfyDesignStart
$(117, 78, 52, 8)$,
$(117, 64, 99, 80)$,
$(117, 22, 51, 111)$,
$(117, 18, 11, 114)$,\adfsplit
$(117, 104, 67, 73)$,
$(117, 109, 44, 113)$,
$(0, 1, 3, 25)$,
$(0, 8, 40, 51)$,\adfsplit
$(0, 5, 55, 155)$,
$(0, 17, 59, 149)$,
$(0, 12, 46, 153)$,
$(0, 15, 56, 118)$,\adfsplit
$(0, 20, 53, 133)$,
$(0, 30, 79, 127)$,
$(0, 10, 23, 147)$,
$(156, 0, 39, 78)$

%ADFvfyBlocksEnd
\adfLgap \noindent by the mapping:
$x \mapsto x +  j \adfmod{117}$ for $x < 117$,
$x \mapsto (x +  j \adfmod{39}) + 117$ for $117 \le x < 156$,
$156 \mapsto 156$,
$0 \le j < 117$
 for the first 15 blocks,
$0 \le j < 39$
 for the last block.
\ADFvfyParStart{(157, ((15, 117, ((117, 1), (39, 1), (1, 1))), (1, 39, ((117, 1), (39, 1), (1, 1)))), ((13, 9), (40, 1)))} %ADFvfyParEnd
% End of 13^9 40^1
%%%%%%%%%%%%%%%%%%%%%%%%%%%%%%%%%%%%%%%%%%%%%%%%%%%%%%%%%%%%%%%%%%%%%%%%%%%%%%%%%%%%%%%%%%
%%%%%%%%%%%%%%%%%%%%%%%%%%%%%%%%%%%%%%%%%%%%%%%%%%%%%%%%%%%%%%%%%%%%%%%%%%%%%%%%%%%%%%%%%%

% Charlotte:GDD4-1-3-5-mod-6-TeX-gen-A:HITS-fun:4.10
\adfDgap
%ADFvfyBlocksStart {13,13,13,13,13,13,13,13,13,46}
\noindent{\boldmath $ 13^{9} 46^{1} $}~
With the point set $Z_{163}$ partitioned into
 residue classes modulo $9$ for $\{0, 1, \dots, 116\}$, and
 $\{117, 118, \dots, 162\}$,
 the design is generated from

\adfLgap %ADFvfyDesignStart
$(159, 20, 73, 72)$,
$(156, 33, 109, 20)$,
$(117, 24, 44, 59)$,
$(117, 107, 97, 74)$,\adfsplit
$(117, 78, 7, 80)$,
$(117, 76, 50, 28)$,
$(117, 31, 25, 93)$,
$(0, 3, 17, 77)$,\adfsplit
$(0, 4, 12, 120)$,
$(0, 7, 93, 153)$,
$(0, 16, 37, 142)$,
$(0, 19, 61, 152)$,\adfsplit
$(0, 32, 66, 122)$,
$(0, 11, 70, 135)$,
$(0, 5, 92, 123)$,
$(0, 29, 79, 140)$,\adfsplit
$(162, 0, 39, 78)$

%ADFvfyBlocksEnd
\adfLgap \noindent by the mapping:
$x \mapsto x +  j \adfmod{117}$ for $x < 117$,
$x \mapsto (x +  j \adfmod{39}) + 117$ for $117 \le x < 156$,
$x \mapsto (x +  j \adfmod{3}) + 156$ for $156 \le x < 159$,
$x \mapsto (x +  j \adfmod{3}) + 159$ for $159 \le x < 162$,
$162 \mapsto 162$,
$0 \le j < 117$
 for the first 16 blocks,
$0 \le j < 39$
 for the last block.
\ADFvfyParStart{(163, ((16, 117, ((117, 1), (39, 1), (3, 1), (3, 1), (1, 1))), (1, 39, ((117, 1), (39, 1), (3, 1), (3, 1), (1, 1)))), ((13, 9), (46, 1)))} %ADFvfyParEnd
% End of 13^9 46^1
%%%%%%%%%%%%%%%%%%%%%%%%%%%%%%%%%%%%%%%%%%%%%%%%%%%%%%%%%%%%%%%%%%%%%%%%%%%%%%%%%%%%%%%%%%
%%%%%%%%%%%%%%%%%%%%%%%%%%%%%%%%%%%%%%%%%%%%%%%%%%%%%%%%%%%%%%%%%%%%%%%%%%%%%%%%%%%%%%%%%%

%%%%%%%%%%%%%%%%%%%%%%%%%%%%%%%%%%%%%%%%%%%%%%%%%%%%%%%%%%%%%%%%%%%%%%%%%%%%%%%%%%%%%%%%%%
%%%%%%%%%%%%%%%%%%%%%%%%%%%%%%%%%%%%%%%%%%%%%%%%%%%%%%%%%%%%%%%%%%%%%%%%%%%%%%%%%%%%%%%%%%
\section{4-GDDs for the proof of Lemma \ref{lem:4-GDD 17^u m^1}}
\label{app:4-GDD 17^u m^1}
\adfnull{
$ 17^{12} 2^1 $,
$ 17^{12} 5^1 $,
$ 17^{12} 8^1 $,
$ 17^{12} 11^1 $,
$ 17^{12} 14^1 $,
$ 17^{15} 11^1 $,
$ 17^9 8^1 $,
$ 17^9 14^1 $,
$ 17^9 20^1 $,
$ 17^9 26^1 $,
$ 17^9 32^1 $,
$ 17^9 38^1 $,
$ 17^9 44^1 $,
$ 17^9 50^1 $,
$ 17^9 56^1 $ and
$ 17^9 62^1 $.
}

% Charlotte:GDD4-1-3-5-mod-6-TeX-gen-A:HITS-fun:4.10
\adfDgap
%ADFvfyBlocksStart {17,17,17,17,17,17,17,17,17,17,17,17,2}
\noindent{\boldmath $ 17^{12} 2^{1} $}~
With the point set $Z_{206}$ partitioned into
 residue classes modulo $12$ for $\{0, 1, \dots, 203\}$, and
 $\{204, 205\}$,
 the design is generated from

\adfLgap %ADFvfyDesignStart
$(204, 91, 119, 132)$,
$(205, 118, 194, 15)$,
$(11, 8, 15, 186)$,
$(116, 38, 102, 96)$,\adfsplit
$(92, 126, 169, 87)$,
$(35, 22, 132, 102)$,
$(27, 127, 124, 191)$,
$(58, 59, 163, 1)$,\adfsplit
$(9, 92, 51, 85)$,
$(58, 131, 148, 147)$,
$(42, 130, 69, 147)$,
$(165, 147, 78, 95)$,\adfsplit
$(9, 182, 136, 83)$,
$(54, 65, 86, 91)$,
$(174, 23, 91, 62)$,
$(24, 10, 43, 160)$,\adfsplit
$(53, 162, 83, 49)$,
$(38, 66, 124, 155)$,
$(106, 141, 122, 187)$,
$(114, 53, 68, 163)$,\adfsplit
$(40, 185, 135, 133)$,
$(117, 25, 190, 188)$,
$(166, 184, 177, 176)$,
$(139, 61, 201, 116)$,\adfsplit
$(49, 161, 19, 57)$,
$(1, 120, 135, 160)$,
$(109, 128, 46, 194)$,
$(0, 2, 129, 184)$,\adfsplit
$(0, 3, 16, 31)$,
$(0, 4, 9, 90)$,
$(0, 23, 68, 141)$,
$(0, 8, 106, 175)$,\adfsplit
$(0, 10, 57, 113)$,
$(0, 44, 79, 150)$,
$(0, 62, 111, 178)$,
$(1, 50, 92, 139)$,\adfsplit
$(0, 52, 149, 181)$,
$(0, 145, 151, 172)$,
$(0, 80, 154, 194)$,
$(1, 134, 161, 167)$,\adfsplit
$(0, 21, 59, 182)$,
$(0, 47, 116, 130)$,
$(0, 22, 66, 101)$,
$(0, 45, 110, 167)$,\adfsplit
$(1, 26, 44, 155)$,
$(1, 62, 71, 125)$,
$(0, 128, 139, 148)$,
$(0, 51, 102, 153)$,\adfsplit
$(1, 52, 103, 154)$,
$(2, 53, 104, 155)$

%ADFvfyBlocksEnd
\adfLgap \noindent by the mapping:
$x \mapsto x + 3 j \adfmod{204}$ for $x < 204$,
$x \mapsto x$ for $x \ge 204$,
$0 \le j < 68$
 for the first 47 blocks,
$0 \le j < 17$
 for the last three blocks.
\ADFvfyParStart{(206, ((47, 68, ((204, 3), (2, 2))), (3, 17, ((204, 3), (2, 2)))), ((17, 12), (2, 1)))} %ADFvfyParEnd
% End of 17^12 2^1
%%%%%%%%%%%%%%%%%%%%%%%%%%%%%%%%%%%%%%%%%%%%%%%%%%%%%%%%%%%%%%%%%%%%%%%%%%%%%%%%%%%%%%%%%%
%%%%%%%%%%%%%%%%%%%%%%%%%%%%%%%%%%%%%%%%%%%%%%%%%%%%%%%%%%%%%%%%%%%%%%%%%%%%%%%%%%%%%%%%%%

% Charlotte:GDD4-1-3-5-mod-6-TeX-gen-A:HITS-fun:4.10
\adfDgap
%ADFvfyBlocksStart {17,17,17,17,17,17,17,17,17,17,17,17,5}
\noindent{\boldmath $ 17^{12} 5^{1} $}~
With the point set $Z_{209}$ partitioned into
 residue classes modulo $12$ for $\{0, 1, \dots, 203\}$, and
 $\{204, 205, 206, 207, 208\}$,
 the design is generated from

\adfLgap %ADFvfyDesignStart
$(204, 0, 1, 2)$,
$(205, 0, 67, 137)$,
$(206, 0, 70, 68)$,
$(207, 0, 136, 203)$,\adfsplit
$(208, 0, 202, 134)$,
$(196, 51, 45, 65)$,
$(84, 125, 93, 138)$,
$(104, 196, 109, 149)$,\adfsplit
$(52, 66, 182, 142)$,
$(57, 72, 190, 67)$,
$(39, 168, 116, 143)$,
$(126, 129, 100, 116)$,\adfsplit
$(138, 172, 92, 130)$,
$(75, 36, 56, 155)$,
$(99, 8, 117, 191)$,
$(191, 37, 66, 44)$,\adfsplit
$(100, 82, 189, 11)$,
$(126, 5, 122, 57)$,
$(86, 37, 165, 80)$,
$(0, 7, 100, 121)$,\adfsplit
$(0, 11, 66, 122)$,
$(0, 17, 28, 174)$,
$(0, 15, 154, 185)$,
$(0, 43, 78, 201)$,\adfsplit
$(0, 32, 173, 177)$,
$(0, 37, 95, 103)$,
$(0, 33, 63, 181)$,
$(0, 23, 64, 165)$,\adfsplit
$(0, 55, 97, 195)$,
$(0, 25, 98, 160)$,
$(0, 35, 57, 151)$,
$(0, 113, 129, 167)$,\adfsplit
$(0, 61, 105, 187)$,
$(0, 49, 77, 110)$,
$(0, 51, 102, 153)$

%ADFvfyBlocksEnd
\adfLgap \noindent by the mapping:
$x \mapsto x + 3 j \adfmod{204}$ for $x < 204$,
$x \mapsto x$ for $x \ge 204$,
$0 \le j < 68$
 for the first five blocks;
$x \mapsto x + 2 j \adfmod{204}$ for $x < 204$,
$x \mapsto x$ for $x \ge 204$,
$0 \le j < 102$
 for the next 29 blocks,
$0 \le j < 51$
 for the last block.
\ADFvfyParStart{(209, ((5, 68, ((204, 3), (5, 5))), (29, 102, ((204, 2), (5, 5))), (1, 51, ((204, 2), (5, 5)))), ((17, 12), (5, 1)))} %ADFvfyParEnd
% End of 17^12 5^1
%%%%%%%%%%%%%%%%%%%%%%%%%%%%%%%%%%%%%%%%%%%%%%%%%%%%%%%%%%%%%%%%%%%%%%%%%%%%%%%%%%%%%%%%%%
%%%%%%%%%%%%%%%%%%%%%%%%%%%%%%%%%%%%%%%%%%%%%%%%%%%%%%%%%%%%%%%%%%%%%%%%%%%%%%%%%%%%%%%%%%

% Charlotte:GDD4-1-3-5-mod-6-TeX-gen-A:HITS-fun:4.10
\adfDgap
%ADFvfyBlocksStart {17,17,17,17,17,17,17,17,17,17,17,17,8}
\noindent{\boldmath $ 17^{12} 8^{1} $}~
With the point set $Z_{212}$ partitioned into
 residue classes modulo $12$ for $\{0, 1, \dots, 203\}$, and
 $\{204, 205, \dots, 211\}$,
 the design is generated from

\adfLgap %ADFvfyDesignStart
$(204, 136, 183, 38)$,
$(205, 60, 32, 109)$,
$(206, 125, 103, 33)$,
$(207, 23, 88, 72)$,\adfsplit
$(208, 61, 164, 27)$,
$(209, 195, 179, 106)$,
$(210, 98, 16, 159)$,
$(211, 25, 65, 69)$,\adfsplit
$(73, 180, 40, 31)$,
$(117, 75, 138, 11)$,
$(126, 11, 170, 163)$,
$(152, 182, 102, 202)$,\adfsplit
$(179, 80, 12, 6)$,
$(188, 43, 30, 112)$,
$(74, 4, 167, 45)$,
$(110, 154, 124, 3)$,\adfsplit
$(151, 141, 47, 15)$,
$(183, 157, 10, 161)$,
$(65, 1, 99, 186)$,
$(202, 97, 128, 79)$,\adfsplit
$(32, 94, 72, 139)$,
$(88, 25, 102, 36)$,
$(81, 36, 175, 92)$,
$(38, 179, 126, 108)$,\adfsplit
$(1, 130, 201, 92)$,
$(67, 102, 75, 50)$,
$(180, 51, 66, 89)$,
$(43, 36, 201, 122)$,\adfsplit
$(96, 79, 104, 39)$,
$(49, 78, 65, 55)$,
$(0, 2, 9, 35)$,
$(0, 14, 69, 99)$,\adfsplit
$(0, 54, 184, 185)$,
$(0, 25, 111, 202)$,
$(0, 3, 20, 88)$,
$(0, 103, 142, 161)$,\adfsplit
$(0, 58, 124, 145)$,
$(0, 1, 33, 112)$,
$(0, 43, 81, 128)$,
$(0, 41, 73, 76)$,\adfsplit
$(0, 62, 101, 148)$,
$(0, 5, 28, 154)$,
$(0, 109, 146, 203)$,
$(0, 104, 194, 199)$,\adfsplit
$(1, 16, 44, 194)$,
$(1, 35, 113, 115)$,
$(1, 47, 55, 116)$,
$(2, 5, 11, 188)$,\adfsplit
$(1, 14, 89, 176)$,
$(1, 68, 134, 149)$,
$(0, 51, 102, 153)$,
$(1, 52, 103, 154)$,\adfsplit
$(2, 53, 104, 155)$

%ADFvfyBlocksEnd
\adfLgap \noindent by the mapping:
$x \mapsto x + 3 j \adfmod{204}$ for $x < 204$,
$x \mapsto x$ for $x \ge 204$,
$0 \le j < 68$
 for the first 50 blocks,
$0 \le j < 17$
 for the last three blocks.
\ADFvfyParStart{(212, ((50, 68, ((204, 3), (8, 8))), (3, 17, ((204, 3), (8, 8)))), ((17, 12), (8, 1)))} %ADFvfyParEnd
% End of 17^12 8^1
%%%%%%%%%%%%%%%%%%%%%%%%%%%%%%%%%%%%%%%%%%%%%%%%%%%%%%%%%%%%%%%%%%%%%%%%%%%%%%%%%%%%%%%%%%
%%%%%%%%%%%%%%%%%%%%%%%%%%%%%%%%%%%%%%%%%%%%%%%%%%%%%%%%%%%%%%%%%%%%%%%%%%%%%%%%%%%%%%%%%%

% Charlotte:GDD4-1-3-5-mod-6-TeX-gen-A:HITS-fun:4.10
\adfDgap
%ADFvfyBlocksStart {17,17,17,17,17,17,17,17,17,17,17,17,11}
\noindent{\boldmath $ 17^{12} 11^{1} $}~
With the point set $Z_{215}$ partitioned into
 residue classes modulo $11$ for $\{0, 1, \dots, 186\}$,
 $\{187, 188, \dots, 203\}$, and
 $\{204, 205, \dots, 214\}$,
 the design is generated from

\adfLgap %ADFvfyDesignStart
$(204, 188, 108, 62)$,
$(204, 136, 126, 175)$,
$(204, 66, 149, 179)$,
$(204, 34, 145, 140)$,\adfsplit
$(12, 24, 11, 80)$,
$(133, 162, 75, 137)$,
$(160, 175, 144, 123)$,
$(39, 202, 80, 134)$,\adfsplit
$(202, 17, 125, 62)$,
$(176, 42, 194, 169)$,
$(0, 2, 75, 200)$,
$(0, 3, 93, 192)$,\adfsplit
$(0, 14, 38, 85)$,
$(0, 20, 70, 98)$,
$(0, 19, 51, 115)$,
$(0, 8, 26, 152)$,\adfsplit
$(0, 9, 36, 139)$,
$(0, 6, 40, 107)$,
$(0, 17, 59, 82)$

%ADFvfyBlocksEnd
\adfLgap \noindent by the mapping:
$x \mapsto x +  j \adfmod{187}$ for $x < 187$,
$x \mapsto (x +  j \adfmod{17}) + 187$ for $187 \le x < 204$,
$x \mapsto (x - 204 +  j \adfmod{11}) + 204$ for $x \ge 204$,
$0 \le j < 187$.
\ADFvfyParStart{(215, ((19, 187, ((187, 1), (17, 1), (11, 1)))), ((17, 11), (17, 1), (11, 1)))} %ADFvfyParEnd
% End of 17^12 11^1
%%%%%%%%%%%%%%%%%%%%%%%%%%%%%%%%%%%%%%%%%%%%%%%%%%%%%%%%%%%%%%%%%%%%%%%%%%%%%%%%%%%%%%%%%%
%%%%%%%%%%%%%%%%%%%%%%%%%%%%%%%%%%%%%%%%%%%%%%%%%%%%%%%%%%%%%%%%%%%%%%%%%%%%%%%%%%%%%%%%%%

% Charlotte:GDD4-1-3-5-mod-6-TeX-gen-A:HITS-fun:4.10
\adfDgap
%ADFvfyBlocksStart {17,17,17,17,17,17,17,17,17,17,17,17,14}
\noindent{\boldmath $ 17^{12} 14^{1} $}~
With the point set $Z_{218}$ partitioned into
 residue classes modulo $12$ for $\{0, 1, \dots, 203\}$, and
 $\{204, 205, \dots, 217\}$,
 the design is generated from

\adfLgap %ADFvfyDesignStart
$(204, 187, 128, 180)$,
$(205, 154, 185, 27)$,
$(206, 44, 202, 15)$,
$(207, 34, 186, 74)$,\adfsplit
$(208, 67, 48, 41)$,
$(209, 19, 39, 185)$,
$(210, 13, 48, 14)$,
$(211, 195, 199, 2)$,\adfsplit
$(212, 74, 96, 127)$,
$(213, 145, 128, 168)$,
$(214, 96, 134, 4)$,
$(215, 118, 21, 110)$,\adfsplit
$(216, 122, 42, 106)$,
$(217, 157, 95, 63)$,
$(81, 99, 91, 26)$,
$(77, 22, 138, 119)$,\adfsplit
$(28, 55, 111, 102)$,
$(48, 134, 114, 82)$,
$(15, 58, 73, 12)$,
$(23, 90, 166, 160)$,\adfsplit
$(87, 6, 104, 48)$,
$(192, 117, 133, 147)$,
$(47, 72, 62, 137)$,
$(73, 48, 117, 161)$,\adfsplit
$(141, 174, 61, 137)$,
$(179, 136, 3, 44)$,
$(149, 39, 33, 167)$,
$(30, 180, 31, 101)$,\adfsplit
$(202, 96, 200, 91)$,
$(71, 168, 139, 181)$,
$(21, 29, 100, 19)$,
$(138, 195, 164, 73)$,\adfsplit
$(40, 199, 68, 73)$,
$(0, 2, 22, 111)$,
$(0, 5, 21, 35)$,
$(0, 27, 74, 117)$,\adfsplit
$(0, 37, 105, 154)$,
$(0, 23, 126, 193)$,
$(0, 59, 140, 141)$,
$(0, 15, 83, 166)$,\adfsplit
$(0, 28, 118, 155)$,
$(0, 50, 73, 77)$,
$(0, 62, 119, 163)$,
$(0, 53, 100, 109)$,\adfsplit
$(0, 40, 103, 178)$,
$(0, 85, 142, 191)$,
$(1, 55, 170, 173)$,
$(1, 4, 83, 187)$,\adfsplit
$(1, 14, 23, 164)$,
$(1, 59, 65, 176)$,
$(2, 23, 68, 101)$,
$(1, 26, 106, 191)$,\adfsplit
$(1, 20, 31, 70)$,
$(0, 51, 102, 153)$,
$(1, 52, 103, 154)$,
$(2, 53, 104, 155)$

%ADFvfyBlocksEnd
\adfLgap \noindent by the mapping:
$x \mapsto x + 3 j \adfmod{204}$ for $x < 204$,
$x \mapsto x$ for $x \ge 204$,
$0 \le j < 68$
 for the first 53 blocks,
$0 \le j < 17$
 for the last three blocks.
\ADFvfyParStart{(218, ((53, 68, ((204, 3), (14, 14))), (3, 17, ((204, 3), (14, 14)))), ((17, 12), (14, 1)))} %ADFvfyParEnd
% End of 17^12 14^1
%%%%%%%%%%%%%%%%%%%%%%%%%%%%%%%%%%%%%%%%%%%%%%%%%%%%%%%%%%%%%%%%%%%%%%%%%%%%%%%%%%%%%%%%%%
%%%%%%%%%%%%%%%%%%%%%%%%%%%%%%%%%%%%%%%%%%%%%%%%%%%%%%%%%%%%%%%%%%%%%%%%%%%%%%%%%%%%%%%%%%

% Charlotte:GDD4-1-3-5-mod-6-TeX-gen-A:HITS-fun:4.10
\adfDgap
%ADFvfyBlocksStart {17,17,17,17,17,17,17,17,17,17,17,17,17,17,17,11}
\noindent{\boldmath $ 17^{15} 11^{1} $}~
With the point set $Z_{266}$ partitioned into
 residue classes modulo $15$ for $\{0, 1, \dots, 254\}$, and
 $\{255, 256, \dots, 265\}$,
 the design is generated from

\adfLgap %ADFvfyDesignStart
$(255, 0, 1, 2)$,
$(256, 0, 85, 83)$,
$(257, 0, 169, 254)$,
$(258, 0, 172, 86)$,\adfsplit
$(259, 0, 253, 170)$,
$(260, 0, 4, 11)$,
$(261, 0, 7, 251)$,
$(262, 0, 244, 248)$,\adfsplit
$(263, 0, 10, 23)$,
$(264, 0, 13, 245)$,
$(265, 0, 232, 242)$,
$(159, 215, 55, 246)$,\adfsplit
$(99, 72, 51, 48)$,
$(33, 158, 164, 146)$,
$(184, 86, 234, 95)$,
$(252, 88, 41, 159)$,\adfsplit
$(42, 170, 144, 198)$,
$(35, 111, 217, 152)$,
$(139, 158, 96, 197)$,
$(22, 239, 94, 143)$,\adfsplit
$(9, 51, 91, 120)$,
$(0, 5, 68, 84)$,
$(0, 17, 53, 114)$,
$(0, 33, 70, 136)$,\adfsplit
$(0, 8, 22, 100)$,
$(0, 55, 122, 181)$,
$(0, 34, 80, 174)$,
$(0, 20, 108, 143)$,\adfsplit
$(0, 25, 57, 203)$

%ADFvfyBlocksEnd
\adfLgap \noindent by the mapping:
$x \mapsto x + 3 j \adfmod{255}$ for $x < 255$,
$x \mapsto x$ for $x \ge 255$,
$0 \le j < 85$
 for the first 11 blocks;
$x \mapsto x +  j \adfmod{255}$ for $x < 255$,
$x \mapsto x$ for $x \ge 255$,
$0 \le j < 255$
 for the last 18 blocks.
\ADFvfyParStart{(266, ((11, 85, ((255, 3), (11, 11))), (18, 255, ((255, 1), (11, 11)))), ((17, 15), (11, 1)))} %ADFvfyParEnd
% End of 17^15 11^1
%%%%%%%%%%%%%%%%%%%%%%%%%%%%%%%%%%%%%%%%%%%%%%%%%%%%%%%%%%%%%%%%%%%%%%%%%%%%%%%%%%%%%%%%%%
%%%%%%%%%%%%%%%%%%%%%%%%%%%%%%%%%%%%%%%%%%%%%%%%%%%%%%%%%%%%%%%%%%%%%%%%%%%%%%%%%%%%%%%%%%

% Charlotte:GDD4-1-3-5-mod-6-TeX-gen-A:HITS-fun:4.10
\adfDgap
%ADFvfyBlocksStart {17,17,17,17,17,17,17,17,17,8}
\noindent{\boldmath $ 17^{9} 8^{1} $}~
With the point set $Z_{161}$ partitioned into
 residue classes modulo $9$ for $\{0, 1, \dots, 152\}$, and
 $\{153, 154, \dots, 160\}$,
 the design is generated from

\adfLgap %ADFvfyDesignStart
$(159, 68, 34, 69)$,
$(160, 130, 20, 90)$,
$(153, 77, 115, 51)$,
$(153, 76, 116, 48)$,\adfsplit
$(153, 65, 82, 108)$,
$(154, 107, 117, 85)$,
$(154, 101, 66, 6)$,
$(154, 5, 109, 106)$,\adfsplit
$(146, 72, 129, 5)$,
$(122, 9, 21, 29)$,
$(4, 57, 68, 146)$,
$(50, 111, 56, 19)$,\adfsplit
$(137, 150, 10, 121)$,
$(27, 33, 89, 16)$,
$(34, 91, 138, 99)$,
$(84, 130, 107, 77)$,\adfsplit
$(55, 60, 113, 85)$,
$(76, 145, 47, 66)$,
$(61, 112, 18, 42)$,
$(106, 129, 150, 73)$,\adfsplit
$(135, 14, 73, 29)$,
$(40, 119, 3, 141)$,
$(0, 14, 42, 119)$,
$(0, 1, 3, 123)$,\adfsplit
$(0, 22, 34, 102)$,
$(0, 38, 69, 149)$,
$(0, 4, 75, 133)$,
$(0, 16, 66, 137)$,\adfsplit
$(0, 7, 104, 139)$,
$(0, 41, 55, 65)$,
$(0, 5, 44, 112)$,
$(0, 2, 59, 105)$,\adfsplit
$(0, 67, 115, 128)$,
$(1, 20, 40, 152)$,
$(1, 16, 104, 107)$,
$(1, 5, 71, 76)$,\adfsplit
$(1, 2, 88, 148)$,
$(1, 26, 77, 146)$

%ADFvfyBlocksEnd
\adfLgap \noindent by the mapping:
$x \mapsto x + 3 j \adfmod{153}$ for $x < 153$,
$x \mapsto (x - 153 + 2 j \adfmod{6}) + 153$ for $153 \le x < 159$,
$x \mapsto x$ for $x \ge 159$,
$0 \le j < 51$.
\ADFvfyParStart{(161, ((38, 51, ((153, 3), (6, 2), (2, 2)))), ((17, 9), (8, 1)))} %ADFvfyParEnd
% End of 17^9 8^1
%%%%%%%%%%%%%%%%%%%%%%%%%%%%%%%%%%%%%%%%%%%%%%%%%%%%%%%%%%%%%%%%%%%%%%%%%%%%%%%%%%%%%%%%%%
%%%%%%%%%%%%%%%%%%%%%%%%%%%%%%%%%%%%%%%%%%%%%%%%%%%%%%%%%%%%%%%%%%%%%%%%%%%%%%%%%%%%%%%%%%

% Charlotte:GDD4-1-3-5-mod-6-TeX-gen-A:HITS-fun:4.10
\adfDgap
%ADFvfyBlocksStart {17,17,17,17,17,17,17,17,17,14}
\noindent{\boldmath $ 17^{9} 14^{1} $}~
With the point set $Z_{167}$ partitioned into
 residue classes modulo $9$ for $\{0, 1, \dots, 152\}$, and
 $\{153, 154, \dots, 166\}$,
 the design is generated from

\adfLgap %ADFvfyDesignStart
$(165, 81, 134, 37)$,
$(166, 104, 124, 102)$,
$(153, 35, 118, 47)$,
$(153, 57, 54, 52)$,\adfsplit
$(153, 78, 49, 149)$,
$(154, 116, 122, 33)$,
$(154, 16, 129, 13)$,
$(154, 65, 36, 100)$,\adfsplit
$(155, 113, 56, 1)$,
$(155, 78, 112, 143)$,
$(155, 12, 61, 36)$,
$(156, 122, 102, 146)$,\adfsplit
$(156, 115, 141, 58)$,
$(156, 144, 125, 82)$,
$(63, 78, 95, 57)$,
$(9, 8, 93, 94)$,\adfsplit
$(38, 145, 106, 140)$,
$(76, 123, 93, 136)$,
$(35, 72, 149, 128)$,
$(137, 71, 102, 79)$,\adfsplit
$(17, 79, 51, 13)$,
$(45, 32, 138, 145)$,
$(69, 74, 124, 49)$,
$(24, 50, 111, 20)$,\adfsplit
$(0, 4, 78, 88)$,
$(0, 11, 39, 141)$,
$(0, 48, 121, 143)$,
$(0, 8, 41, 57)$,\adfsplit
$(0, 42, 118, 128)$,
$(1, 8, 25, 77)$,
$(0, 14, 52, 58)$,
$(0, 82, 101, 103)$,\adfsplit
$(0, 50, 97, 98)$,
$(0, 59, 74, 139)$,
$(0, 33, 113, 145)$,
$(0, 94, 131, 142)$,\adfsplit
$(0, 16, 67, 146)$,
$(0, 19, 31, 61)$,
$(1, 17, 95, 121)$,
$(1, 14, 125, 139)$,\adfsplit
$(0, 46, 107, 110)$

%ADFvfyBlocksEnd
\adfLgap \noindent by the mapping:
$x \mapsto x + 3 j \adfmod{153}$ for $x < 153$,
$x \mapsto (x - 153 + 4 j \adfmod{12}) + 153$ for $153 \le x < 165$,
$x \mapsto x$ for $x \ge 165$,
$0 \le j < 51$.
\ADFvfyParStart{(167, ((41, 51, ((153, 3), (12, 4), (2, 2)))), ((17, 9), (14, 1)))} %ADFvfyParEnd
% End of 17^9 14^1
%%%%%%%%%%%%%%%%%%%%%%%%%%%%%%%%%%%%%%%%%%%%%%%%%%%%%%%%%%%%%%%%%%%%%%%%%%%%%%%%%%%%%%%%%%
%%%%%%%%%%%%%%%%%%%%%%%%%%%%%%%%%%%%%%%%%%%%%%%%%%%%%%%%%%%%%%%%%%%%%%%%%%%%%%%%%%%%%%%%%%

% Charlotte:GDD4-1-3-5-mod-6-TeX-gen-A:HITS-fun:4.10
\adfDgap
%ADFvfyBlocksStart {17,17,17,17,17,17,17,17,17,20}
\noindent{\boldmath $ 17^{9} 20^{1} $}~
With the point set $Z_{173}$ partitioned into
 residue classes modulo $9$ for $\{0, 1, \dots, 152\}$, and
 $\{153, 154, \dots, 172\}$,
 the design is generated from

\adfLgap %ADFvfyDesignStart
$(170, 24, 85, 53)$,
$(170, 48, 70, 100)$,
$(170, 119, 99, 95)$,
$(153, 48, 116, 128)$,\adfsplit
$(153, 126, 40, 125)$,
$(153, 6, 85, 71)$,
$(153, 53, 110, 73)$,
$(153, 16, 130, 39)$,\adfsplit
$(153, 137, 135, 152)$,
$(153, 98, 117, 147)$,
$(153, 127, 87, 106)$,
$(153, 113, 42, 63)$,\adfsplit
$(153, 32, 60, 109)$,
$(153, 78, 121, 37)$,
$(153, 43, 120, 100)$,
$(153, 3, 103, 146)$,\adfsplit
$(153, 80, 119, 72)$,
$(153, 38, 51, 133)$,
$(153, 30, 148, 5)$,
$(153, 64, 112, 143)$,\adfsplit
$(48, 105, 59, 81)$,
$(48, 6, 92, 9)$,
$(131, 72, 103, 70)$,
$(67, 114, 63, 23)$,\adfsplit
$(122, 115, 0, 127)$,
$(103, 80, 48, 149)$,
$(117, 1, 56, 143)$,
$(0, 1, 60, 147)$,\adfsplit
$(0, 5, 48, 58)$,
$(0, 12, 46, 121)$,
$(0, 15, 38, 139)$,
$(0, 13, 73, 78)$,\adfsplit
$(0, 16, 84, 119)$,
$(0, 64, 70, 89)$,
$(0, 53, 142, 146)$,
$(0, 56, 77, 103)$,\adfsplit
$(0, 14, 97, 137)$,
$(0, 95, 98, 136)$,
$(1, 2, 88, 139)$,
$(1, 35, 43, 137)$,\adfsplit
$(0, 25, 28, 116)$,
$(1, 23, 25, 98)$,
$(0, 41, 74, 145)$,
$(1, 14, 119, 125)$

%ADFvfyBlocksEnd
\adfLgap \noindent by the mapping:
$x \mapsto x + 3 j \adfmod{153}$ for $x < 153$,
$x \mapsto (x +  j \adfmod{17}) + 153$ for $153 \le x < 170$,
$x \mapsto (x - 170 +  j \adfmod{3}) + 170$ for $x \ge 170$,
$0 \le j < 51$.
\ADFvfyParStart{(173, ((44, 51, ((153, 3), (17, 1), (3, 1)))), ((17, 9), (20, 1)))} %ADFvfyParEnd
% End of 17^9 20^1
%%%%%%%%%%%%%%%%%%%%%%%%%%%%%%%%%%%%%%%%%%%%%%%%%%%%%%%%%%%%%%%%%%%%%%%%%%%%%%%%%%%%%%%%%%
%%%%%%%%%%%%%%%%%%%%%%%%%%%%%%%%%%%%%%%%%%%%%%%%%%%%%%%%%%%%%%%%%%%%%%%%%%%%%%%%%%%%%%%%%%

% Charlotte:GDD4-1-3-5-mod-6-TeX-gen-A:HITS-fun:4.10
\adfDgap
%ADFvfyBlocksStart {17,17,17,17,17,17,17,17,17,26}
\noindent{\boldmath $ 17^{9} 26^{1} $}~
With the point set $Z_{179}$ partitioned into
 residue classes modulo $9$ for $\{0, 1, \dots, 152\}$, and
 $\{153, 154, \dots, 178\}$,
 the design is generated from

\adfLgap %ADFvfyDesignStart
$(170, 88, 54, 11)$,
$(170, 145, 85, 113)$,
$(170, 53, 75, 78)$,
$(171, 95, 33, 16)$,\adfsplit
$(171, 49, 71, 82)$,
$(171, 54, 38, 84)$,
$(172, 149, 138, 35)$,
$(172, 34, 126, 28)$,\adfsplit
$(172, 121, 78, 137)$,
$(153, 110, 37, 105)$,
$(153, 47, 136, 26)$,
$(153, 16, 104, 71)$,\adfsplit
$(153, 89, 65, 106)$,
$(153, 133, 103, 95)$,
$(153, 152, 5, 151)$,
$(153, 141, 25, 81)$,\adfsplit
$(153, 23, 80, 135)$,
$(153, 130, 102, 61)$,
$(153, 94, 42, 119)$,
$(153, 62, 109, 63)$,\adfsplit
$(153, 75, 117, 13)$,
$(153, 46, 124, 41)$,
$(153, 138, 18, 91)$,
$(153, 45, 129, 137)$,\adfsplit
$(153, 57, 9, 19)$,
$(153, 32, 72, 150)$,
$(91, 23, 89, 9)$,
$(84, 74, 5, 124)$,\adfsplit
$(132, 26, 97, 74)$,
$(0, 1, 20, 32)$,
$(0, 2, 39, 140)$,
$(0, 7, 56, 134)$,\adfsplit
$(0, 76, 86, 89)$,
$(1, 43, 110, 140)$,
$(0, 23, 65, 125)$,
$(0, 29, 58, 124)$,\adfsplit
$(0, 44, 103, 151)$,
$(0, 24, 94, 146)$,
$(0, 12, 83, 145)$,
$(0, 22, 26, 79)$,\adfsplit
$(0, 17, 109, 132)$,
$(0, 15, 68, 102)$,
$(0, 4, 25, 57)$,
$(0, 64, 67, 104)$,\adfsplit
$(0, 88, 127, 139)$,
$(0, 16, 31, 149)$,
$(0, 6, 19, 148)$

%ADFvfyBlocksEnd
\adfLgap \noindent by the mapping:
$x \mapsto x + 3 j \adfmod{153}$ for $x < 153$,
$x \mapsto (x +  j \adfmod{17}) + 153$ for $153 \le x < 170$,
$x \mapsto (x - 170 + 3 j \adfmod{9}) + 170$ for $x \ge 170$,
$0 \le j < 51$.
\ADFvfyParStart{(179, ((47, 51, ((153, 3), (17, 1), (9, 3)))), ((17, 9), (26, 1)))} %ADFvfyParEnd
% End of 17^9 26^1
%%%%%%%%%%%%%%%%%%%%%%%%%%%%%%%%%%%%%%%%%%%%%%%%%%%%%%%%%%%%%%%%%%%%%%%%%%%%%%%%%%%%%%%%%%
%%%%%%%%%%%%%%%%%%%%%%%%%%%%%%%%%%%%%%%%%%%%%%%%%%%%%%%%%%%%%%%%%%%%%%%%%%%%%%%%%%%%%%%%%%

% Charlotte:GDD4-1-3-5-mod-6-TeX-gen-A:HITS-fun:4.10
\adfDgap
%ADFvfyBlocksStart {17,17,17,17,17,17,17,17,17,32}
\noindent{\boldmath $ 17^{9} 32^{1} $}~
With the point set $Z_{185}$ partitioned into
 residue classes modulo $9$ for $\{0, 1, \dots, 152\}$, and
 $\{153, 154, \dots, 184\}$,
 the design is generated from

\adfLgap %ADFvfyDesignStart
$(170, 107, 148, 126)$,
$(170, 97, 42, 95)$,
$(170, 111, 101, 118)$,
$(171, 133, 96, 144)$,\adfsplit
$(171, 107, 139, 32)$,
$(171, 100, 102, 29)$,
$(172, 122, 120, 98)$,
$(172, 88, 65, 130)$,\adfsplit
$(172, 28, 123, 81)$,
$(173, 33, 99, 43)$,
$(173, 122, 20, 100)$,
$(173, 84, 85, 26)$,\adfsplit
$(174, 44, 74, 49)$,
$(174, 81, 86, 124)$,
$(174, 129, 145, 60)$,
$(153, 150, 116, 16)$,\adfsplit
$(153, 80, 52, 40)$,
$(153, 93, 137, 41)$,
$(153, 102, 88, 149)$,
$(153, 9, 119, 26)$,\adfsplit
$(153, 13, 126, 62)$,
$(153, 56, 124, 140)$,
$(153, 21, 44, 15)$,
$(153, 100, 39, 18)$,\adfsplit
$(153, 71, 19, 78)$,
$(153, 33, 36, 101)$,
$(153, 81, 57, 133)$,
$(153, 127, 130, 134)$,\adfsplit
$(153, 63, 109, 23)$,
$(153, 43, 148, 53)$,
$(153, 112, 136, 3)$,
$(153, 147, 8, 106)$,\adfsplit
$(0, 4, 73, 123)$,
$(0, 8, 120, 140)$,
$(0, 11, 25, 145)$,
$(0, 51, 115, 149)$,\adfsplit
$(0, 56, 91, 130)$,
$(0, 13, 28, 88)$,
$(0, 12, 79, 116)$,
$(0, 39, 122, 125)$,\adfsplit
$(0, 35, 118, 124)$,
$(1, 14, 22, 128)$,
$(0, 38, 71, 148)$,
$(0, 26, 32, 74)$,\adfsplit
$(0, 50, 70, 96)$,
$(1, 2, 52, 143)$,
$(0, 49, 106, 136)$,
$(0, 31, 62, 128)$,\adfsplit
$(0, 15, 93, 152)$,
$(0, 77, 92, 121)$

%ADFvfyBlocksEnd
\adfLgap \noindent by the mapping:
$x \mapsto x + 3 j \adfmod{153}$ for $x < 153$,
$x \mapsto (x +  j \adfmod{17}) + 153$ for $153 \le x < 170$,
$x \mapsto (x - 170 + 5 j \adfmod{15}) + 170$ for $x \ge 170$,
$0 \le j < 51$.
\ADFvfyParStart{(185, ((50, 51, ((153, 3), (17, 1), (15, 5)))), ((17, 9), (32, 1)))} %ADFvfyParEnd
% End of 17^9 32^1
%%%%%%%%%%%%%%%%%%%%%%%%%%%%%%%%%%%%%%%%%%%%%%%%%%%%%%%%%%%%%%%%%%%%%%%%%%%%%%%%%%%%%%%%%%
%%%%%%%%%%%%%%%%%%%%%%%%%%%%%%%%%%%%%%%%%%%%%%%%%%%%%%%%%%%%%%%%%%%%%%%%%%%%%%%%%%%%%%%%%%

% Charlotte:GDD4-1-3-5-mod-6-TeX-gen-A:HITS-fun:4.10
\adfDgap
%ADFvfyBlocksStart {17,17,17,17,17,17,17,17,17,38}
\noindent{\boldmath $ 17^{9} 38^{1} $}~
With the point set $Z_{191}$ partitioned into
 residue classes modulo $9$ for $\{0, 1, \dots, 152\}$, and
 $\{153, 154, \dots, 190\}$,
 the design is generated from

\adfLgap %ADFvfyDesignStart
$(190, 146, 91, 105)$,
$(187, 130, 11, 134)$,
$(187, 1, 108, 32)$,
$(187, 6, 88, 84)$,\adfsplit
$(153, 15, 2, 49)$,
$(153, 13, 133, 77)$,
$(153, 95, 118, 98)$,
$(153, 20, 57, 124)$,\adfsplit
$(153, 30, 38, 63)$,
$(153, 138, 137, 8)$,
$(153, 103, 144, 74)$,
$(153, 84, 61, 150)$,\adfsplit
$(153, 56, 141, 134)$,
$(153, 19, 79, 54)$,
$(153, 72, 75, 119)$,
$(153, 145, 11, 97)$,\adfsplit
$(153, 101, 109, 120)$,
$(153, 92, 147, 55)$,
$(153, 40, 131, 60)$,
$(153, 37, 0, 65)$,\adfsplit
$(153, 85, 27, 127)$,
$(154, 49, 81, 88)$,
$(154, 40, 113, 53)$,
$(154, 115, 103, 12)$,\adfsplit
$(154, 31, 7, 143)$,
$(154, 144, 19, 71)$,
$(154, 138, 47, 9)$,
$(154, 83, 131, 126)$,\adfsplit
$(154, 108, 150, 50)$,
$(154, 3, 112, 35)$,
$(154, 68, 106, 33)$,
$(154, 85, 123, 74)$,\adfsplit
$(154, 102, 116, 128)$,
$(154, 97, 18, 140)$,
$(0, 1, 132, 178)$,
$(0, 17, 39, 162)$,\adfsplit
$(1, 4, 122, 172)$,
$(1, 17, 31, 160)$,
$(0, 11, 15, 55)$,
$(1, 71, 86, 128)$,\adfsplit
$(0, 13, 92, 113)$,
$(1, 11, 62, 95)$,
$(0, 43, 125, 127)$,
$(0, 12, 86, 136)$,\adfsplit
$(0, 31, 97, 119)$,
$(1, 7, 22, 47)$,
$(0, 101, 106, 107)$,
$(0, 2, 89, 105)$,\adfsplit
$(0, 52, 60, 148)$,
$(0, 10, 69, 85)$,
$(0, 49, 56, 151)$,
$(0, 20, 30, 59)$,\adfsplit
$(0, 6, 57, 76)$

%ADFvfyBlocksEnd
\adfLgap \noindent by the mapping:
$x \mapsto x + 3 j \adfmod{153}$ for $x < 153$,
$x \mapsto (x - 153 + 2 j \adfmod{34}) + 153$ for $153 \le x < 187$,
$x \mapsto (x - 187 +  j \adfmod{3}) + 187$ for $187 \le x < 190$,
$190 \mapsto 190$,
$0 \le j < 51$.
\ADFvfyParStart{(191, ((53, 51, ((153, 3), (34, 2), (3, 1), (1, 1)))), ((17, 9), (38, 1)))} %ADFvfyParEnd
% End of 17^9 38^1
%%%%%%%%%%%%%%%%%%%%%%%%%%%%%%%%%%%%%%%%%%%%%%%%%%%%%%%%%%%%%%%%%%%%%%%%%%%%%%%%%%%%%%%%%%
%%%%%%%%%%%%%%%%%%%%%%%%%%%%%%%%%%%%%%%%%%%%%%%%%%%%%%%%%%%%%%%%%%%%%%%%%%%%%%%%%%%%%%%%%%

% Charlotte:GDD4-1-3-5-mod-6-TeX-gen-A:HITS-fun:4.10
\adfDgap
%ADFvfyBlocksStart {17,17,17,17,17,17,17,17,17,44}
\noindent{\boldmath $ 17^{9} 44^{1} $}~
With the point set $Z_{197}$ partitioned into
 residue classes modulo $9$ for $\{0, 1, \dots, 152\}$, and
 $\{153, 154, \dots, 196\}$,
 the design is generated from

\adfLgap %ADFvfyDesignStart
$(196, 32, 21, 148)$,
$(187, 99, 102, 137)$,
$(187, 58, 23, 44)$,
$(187, 34, 96, 1)$,\adfsplit
$(188, 123, 88, 140)$,
$(188, 38, 144, 111)$,
$(188, 127, 53, 139)$,
$(189, 65, 97, 68)$,\adfsplit
$(189, 21, 126, 118)$,
$(189, 78, 26, 76)$,
$(153, 54, 119, 94)$,
$(153, 118, 74, 35)$,\adfsplit
$(153, 127, 34, 38)$,
$(153, 58, 79, 41)$,
$(153, 51, 111, 31)$,
$(153, 39, 78, 40)$,\adfsplit
$(153, 84, 71, 13)$,
$(153, 36, 42, 14)$,
$(153, 128, 126, 37)$,
$(153, 120, 95, 62)$,\adfsplit
$(153, 132, 49, 83)$,
$(153, 6, 52, 29)$,
$(153, 10, 124, 104)$,
$(153, 96, 148, 56)$,\adfsplit
$(153, 110, 66, 98)$,
$(153, 99, 123, 55)$,
$(153, 19, 101, 12)$,
$(154, 117, 64, 110)$,\adfsplit
$(154, 39, 81, 49)$,
$(154, 78, 67, 152)$,
$(154, 99, 107, 1)$,
$(154, 54, 70, 80)$,\adfsplit
$(154, 128, 85, 114)$,
$(154, 149, 130, 73)$,
$(154, 76, 65, 135)$,
$(154, 58, 146, 71)$,\adfsplit
$(154, 91, 108, 138)$,
$(0, 12, 49, 79)$,
$(1, 25, 98, 103)$,
$(0, 61, 148, 154)$,\adfsplit
$(0, 31, 34, 62)$,
$(0, 43, 112, 158)$,
$(0, 86, 137, 139)$,
$(0, 19, 25, 130)$,\adfsplit
$(0, 13, 28, 78)$,
$(1, 17, 41, 128)$,
$(0, 68, 76, 176)$,
$(1, 8, 23, 168)$,\adfsplit
$(0, 4, 5, 53)$,
$(0, 22, 77, 134)$,
$(0, 56, 57, 149)$,
$(0, 110, 116, 182)$,\adfsplit
$(0, 15, 102, 122)$,
$(0, 59, 69, 172)$,
$(0, 41, 71, 174)$,
$(0, 21, 50, 119)$

%ADFvfyBlocksEnd
\adfLgap \noindent by the mapping:
$x \mapsto x + 3 j \adfmod{153}$ for $x < 153$,
$x \mapsto (x - 153 + 2 j \adfmod{34}) + 153$ for $153 \le x < 187$,
$x \mapsto (x - 187 + 3 j \adfmod{9}) + 187$ for $187 \le x < 196$,
$196 \mapsto 196$,
$0 \le j < 51$.
\ADFvfyParStart{(197, ((56, 51, ((153, 3), (34, 2), (9, 3), (1, 1)))), ((17, 9), (44, 1)))} %ADFvfyParEnd
% End of 17^9 44^1
%%%%%%%%%%%%%%%%%%%%%%%%%%%%%%%%%%%%%%%%%%%%%%%%%%%%%%%%%%%%%%%%%%%%%%%%%%%%%%%%%%%%%%%%%%
%%%%%%%%%%%%%%%%%%%%%%%%%%%%%%%%%%%%%%%%%%%%%%%%%%%%%%%%%%%%%%%%%%%%%%%%%%%%%%%%%%%%%%%%%%

% Charlotte:GDD4-1-3-5-mod-6-TeX-gen-A:HITS-fun:4.10
\adfDgap
%ADFvfyBlocksStart {17,17,17,17,17,17,17,17,17,50}
\noindent{\boldmath $ 17^{9} 50^{1} $}~
With the point set $Z_{203}$ partitioned into
 residue classes modulo $9$ for $\{0, 1, \dots, 152\}$, and
 $\{153, 154, \dots, 202\}$,
 the design is generated from

\adfLgap %ADFvfyDesignStart
$(202, 123, 5, 76)$,
$(187, 62, 11, 21)$,
$(187, 112, 133, 99)$,
$(187, 14, 51, 37)$,\adfsplit
$(188, 86, 115, 67)$,
$(188, 30, 60, 83)$,
$(188, 46, 0, 71)$,
$(189, 48, 79, 110)$,\adfsplit
$(189, 78, 40, 125)$,
$(189, 1, 122, 117)$,
$(190, 120, 81, 61)$,
$(190, 20, 127, 51)$,\adfsplit
$(190, 134, 148, 50)$,
$(191, 4, 15, 34)$,
$(191, 144, 129, 20)$,
$(191, 28, 26, 104)$,\adfsplit
$(153, 9, 23, 67)$,
$(153, 73, 20, 115)$,
$(153, 94, 12, 86)$,
$(153, 88, 49, 92)$,\adfsplit
$(153, 3, 132, 59)$,
$(153, 29, 53, 32)$,
$(153, 112, 140, 45)$,
$(153, 123, 39, 149)$,\adfsplit
$(153, 57, 56, 17)$,
$(153, 152, 102, 69)$,
$(153, 109, 34, 99)$,
$(153, 121, 25, 116)$,\adfsplit
$(153, 52, 26, 24)$,
$(153, 66, 40, 62)$,
$(153, 130, 138, 135)$,
$(153, 144, 55, 148)$,\adfsplit
$(153, 95, 31, 129)$,
$(154, 31, 134, 92)$,
$(154, 89, 1, 24)$,
$(154, 152, 136, 56)$,\adfsplit
$(154, 145, 128, 79)$,
$(154, 54, 66, 25)$,
$(154, 149, 8, 57)$,
$(154, 97, 132, 36)$,\adfsplit
$(0, 6, 131, 146)$,
$(0, 8, 49, 128)$,
$(0, 101, 121, 176)$,
$(0, 7, 40, 91)$,\adfsplit
$(0, 42, 85, 102)$,
$(1, 4, 142, 182)$,
$(0, 25, 77, 156)$,
$(0, 16, 22, 66)$,\adfsplit
$(0, 38, 68, 151)$,
$(0, 17, 73, 158)$,
$(0, 79, 89, 170)$,
$(0, 1, 107, 162)$,\adfsplit
$(0, 59, 97, 180)$,
$(0, 11, 132, 160)$,
$(0, 20, 75, 182)$,
$(0, 70, 137, 154)$,\adfsplit
$(1, 2, 8, 95)$,
$(0, 48, 100, 134)$,
$(1, 14, 119, 130)$

%ADFvfyBlocksEnd
\adfLgap \noindent by the mapping:
$x \mapsto x + 3 j \adfmod{153}$ for $x < 153$,
$x \mapsto (x - 153 + 2 j \adfmod{34}) + 153$ for $153 \le x < 187$,
$x \mapsto (x - 187 + 5 j \adfmod{15}) + 187$ for $187 \le x < 202$,
$202 \mapsto 202$,
$0 \le j < 51$.
\ADFvfyParStart{(203, ((59, 51, ((153, 3), (34, 2), (15, 5), (1, 1)))), ((17, 9), (50, 1)))} %ADFvfyParEnd
% End of 17^9 50^1
%%%%%%%%%%%%%%%%%%%%%%%%%%%%%%%%%%%%%%%%%%%%%%%%%%%%%%%%%%%%%%%%%%%%%%%%%%%%%%%%%%%%%%%%%%
%%%%%%%%%%%%%%%%%%%%%%%%%%%%%%%%%%%%%%%%%%%%%%%%%%%%%%%%%%%%%%%%%%%%%%%%%%%%%%%%%%%%%%%%%%

% Charlotte:GDD4-1-3-5-mod-6-TeX-gen-A:HITS-fun:4.10
\adfDgap
%ADFvfyBlocksStart {17,17,17,17,17,17,17,17,17,56}
\noindent{\boldmath $ 17^{9} 56^{1} $}~
With the point set $Z_{209}$ partitioned into
 residue classes modulo $9$ for $\{0, 1, \dots, 152\}$, and
 $\{153, 154, \dots, 208\}$,
 the design is generated from

\adfLgap %ADFvfyDesignStart
$(207, 98, 124, 45)$,
$(208, 93, 97, 41)$,
$(204, 74, 107, 115)$,
$(204, 46, 15, 22)$,\adfsplit
$(204, 122, 0, 111)$,
$(153, 74, 25, 96)$,
$(153, 59, 146, 21)$,
$(153, 133, 105, 18)$,\adfsplit
$(153, 134, 46, 67)$,
$(153, 71, 119, 57)$,
$(153, 92, 28, 86)$,
$(153, 87, 22, 88)$,\adfsplit
$(153, 7, 9, 15)$,
$(153, 55, 52, 81)$,
$(153, 61, 12, 145)$,
$(153, 11, 70, 84)$,\adfsplit
$(153, 49, 93, 29)$,
$(153, 142, 39, 14)$,
$(153, 26, 101, 27)$,
$(153, 107, 126, 136)$,\adfsplit
$(153, 140, 47, 0)$,
$(153, 53, 13, 150)$,
$(154, 145, 51, 143)$,
$(154, 58, 135, 120)$,\adfsplit
$(154, 8, 85, 117)$,
$(154, 113, 98, 74)$,
$(154, 50, 108, 130)$,
$(154, 126, 89, 5)$,\adfsplit
$(154, 139, 100, 3)$,
$(154, 17, 20, 127)$,
$(154, 82, 32, 48)$,
$(154, 97, 12, 55)$,\adfsplit
$(154, 29, 112, 123)$,
$(154, 144, 77, 132)$,
$(154, 14, 67, 19)$,
$(154, 103, 35, 73)$,\adfsplit
$(154, 64, 87, 53)$,
$(154, 60, 95, 90)$,
$(154, 96, 91, 27)$,
$(155, 90, 33, 57)$,\adfsplit
$(155, 149, 81, 142)$,
$(155, 150, 109, 52)$,
$(155, 94, 128, 24)$,
$(0, 3, 105, 113)$,\adfsplit
$(0, 20, 67, 118)$,
$(0, 25, 26, 77)$,
$(0, 71, 75, 161)$,
$(0, 46, 65, 107)$,\adfsplit
$(0, 23, 146, 155)$,
$(0, 40, 93, 170)$,
$(0, 73, 83, 106)$,
$(0, 21, 58, 164)$,\adfsplit
$(0, 2, 114, 197)$,
$(0, 52, 143, 191)$,
$(0, 13, 19, 50)$,
$(0, 17, 29, 203)$,\adfsplit
$(1, 23, 94, 161)$,
$(1, 16, 137, 188)$,
$(1, 5, 79, 203)$,
$(1, 13, 29, 197)$,\adfsplit
$(1, 14, 110, 179)$,
$(1, 119, 140, 194)$

%ADFvfyBlocksEnd
\adfLgap \noindent by the mapping:
$x \mapsto x + 3 j \adfmod{153}$ for $x < 153$,
$x \mapsto (x + 3 j \adfmod{51}) + 153$ for $153 \le x < 204$,
$x \mapsto (x +  j \adfmod{3}) + 204$ for $204 \le x < 207$,
$x \mapsto x$ for $x \ge 207$,
$0 \le j < 51$.
\ADFvfyParStart{(209, ((62, 51, ((153, 3), (51, 3), (3, 1), (2, 2)))), ((17, 9), (56, 1)))} %ADFvfyParEnd
% End of 17^9 56^1
%%%%%%%%%%%%%%%%%%%%%%%%%%%%%%%%%%%%%%%%%%%%%%%%%%%%%%%%%%%%%%%%%%%%%%%%%%%%%%%%%%%%%%%%%%
%%%%%%%%%%%%%%%%%%%%%%%%%%%%%%%%%%%%%%%%%%%%%%%%%%%%%%%%%%%%%%%%%%%%%%%%%%%%%%%%%%%%%%%%%%

% Charlotte:GDD4-1-3-5-mod-6-TeX-gen-A:HITS-fun:4.10
\adfDgap
%ADFvfyBlocksStart {17,17,17,17,17,17,17,17,17,62}
\noindent{\boldmath $ 17^{9} 62^{1} $}~
With the point set $Z_{215}$ partitioned into
 residue classes modulo $9$ for $\{0, 1, \dots, 152\}$, and
 $\{153, 154, \dots, 214\}$,
 the design is generated from

\adfLgap %ADFvfyDesignStart
$(213, 56, 112, 120)$,
$(214, 97, 122, 18)$,
$(204, 125, 85, 34)$,
$(204, 117, 23, 38)$,\adfsplit
$(204, 114, 136, 30)$,
$(205, 44, 137, 87)$,
$(205, 85, 36, 73)$,
$(205, 142, 113, 12)$,\adfsplit
$(206, 147, 89, 90)$,
$(206, 124, 94, 105)$,
$(206, 127, 41, 146)$,
$(153, 31, 55, 59)$,\adfsplit
$(153, 42, 45, 100)$,
$(153, 124, 20, 44)$,
$(153, 117, 103, 11)$,
$(153, 141, 16, 26)$,\adfsplit
$(153, 57, 150, 135)$,
$(153, 56, 131, 75)$,
$(153, 123, 12, 125)$,
$(153, 65, 3, 32)$,\adfsplit
$(153, 13, 46, 149)$,
$(153, 142, 51, 41)$,
$(153, 104, 69, 38)$,
$(153, 86, 17, 85)$,\adfsplit
$(153, 28, 60, 27)$,
$(153, 101, 58, 37)$,
$(153, 30, 127, 112)$,
$(153, 43, 121, 87)$,\adfsplit
$(154, 92, 114, 85)$,
$(154, 124, 129, 116)$,
$(154, 122, 67, 42)$,
$(154, 2, 132, 23)$,\adfsplit
$(154, 100, 58, 80)$,
$(154, 87, 103, 99)$,
$(154, 79, 8, 120)$,
$(154, 55, 113, 45)$,\adfsplit
$(154, 21, 98, 64)$,
$(154, 50, 33, 133)$,
$(154, 70, 38, 108)$,
$(154, 9, 39, 146)$,\adfsplit
$(154, 25, 105, 137)$,
$(154, 68, 94, 56)$,
$(154, 61, 126, 26)$,
$(154, 51, 97, 117)$,\adfsplit
$(154, 91, 139, 32)$,
$(155, 46, 125, 105)$,
$(0, 13, 52, 102)$,
$(0, 8, 14, 65)$,\adfsplit
$(0, 5, 39, 155)$,
$(0, 6, 92, 176)$,
$(0, 21, 146, 161)$,
$(0, 11, 48, 182)$,\adfsplit
$(0, 24, 151, 197)$,
$(0, 40, 128, 179)$,
$(0, 26, 70, 164)$,
$(0, 71, 85, 200)$,\adfsplit
$(0, 31, 118, 158)$,
$(0, 61, 67, 98)$,
$(1, 4, 61, 185)$,
$(0, 7, 76, 194)$,\adfsplit
$(1, 77, 107, 203)$,
$(0, 136, 149, 188)$,
$(1, 17, 131, 173)$,
$(1, 101, 143, 164)$,\adfsplit
$(1, 149, 152, 158)$

%ADFvfyBlocksEnd
\adfLgap \noindent by the mapping:
$x \mapsto x + 3 j \adfmod{153}$ for $x < 153$,
$x \mapsto (x + 3 j \adfmod{51}) + 153$ for $153 \le x < 204$,
$x \mapsto (x - 204 + 3 j \adfmod{9}) + 204$ for $204 \le x < 213$,
$x \mapsto x$ for $x \ge 213$,
$0 \le j < 51$.
\ADFvfyParStart{(215, ((65, 51, ((153, 3), (51, 3), (9, 3), (2, 2)))), ((17, 9), (62, 1)))} %ADFvfyParEnd
% End of 17^9 62^1
%%%%%%%%%%%%%%%%%%%%%%%%%%%%%%%%%%%%%%%%%%%%%%%%%%%%%%%%%%%%%%%%%%%%%%%%%%%%%%%%%%%%%%%%%%
%%%%%%%%%%%%%%%%%%%%%%%%%%%%%%%%%%%%%%%%%%%%%%%%%%%%%%%%%%%%%%%%%%%%%%%%%%%%%%%%%%%%%%%%%%

%%%%%%%%%%%%%%%%%%%%%%%%%%%%%%%%%%%%%%%%%%%%%%%%%%%%%%%%%%%%%%%%%%%%%%%%%%%%%%%%%%%%%%%%%%
%%%%%%%%%%%%%%%%%%%%%%%%%%%%%%%%%%%%%%%%%%%%%%%%%%%%%%%%%%%%%%%%%%%%%%%%%%%%%%%%%%%%%%%%%%
\section{4-GDDs for the proof of Lemma \ref{lem:4-GDD 19^u m^1}}
\label{app:4-GDD 19^u m^1}
\adfnull{
$ 19^{12} 7^1 $,
$ 19^{12} 10^1 $,
$ 19^{12} 13^1 $,
$ 19^{12} 16^1 $,
$ 19^{24} 13^1 $,
$ 19^{24} 16^1 $,
$ 19^{15} 13^1 $,
$ 19^9 10^1 $,
$ 19^9 16^1 $,
$ 19^9 22^1 $,
$ 19^9 28^1 $,
$ 19^9 34^1 $,
$ 19^9 40^1 $,
$ 19^9 46^1 $,
$ 19^9 52^1 $,
$ 19^9 58^1 $,
$ 19^9 64^1 $,
$ 19^9 70^1 $ and
$ 19^{21} 16^1 $.
}

% Charlotte:GDD4-1-3-5-mod-6-TeX-gen-A:HITS-fun:4.10
\adfDgap
%ADFvfyBlocksStart {19,19,19,19,19,19,19,19,19,19,19,19,7}
\noindent{\boldmath $ 19^{12} 7^{1} $}~
With the point set $Z_{235}$ partitioned into
 residue classes modulo $12$ for $\{0, 1, \dots, 227\}$, and
 $\{228, 229, \dots, 234\}$,
 the design is generated from

\adfLgap %ADFvfyDesignStart
$(228, 103, 198, 110)$,
$(229, 35, 192, 226)$,
$(15, 204, 218, 224)$,
$(104, 158, 173, 83)$,\adfsplit
$(78, 28, 195, 223)$,
$(10, 48, 114, 88)$,
$(202, 16, 110, 197)$,
$(45, 175, 62, 107)$,\adfsplit
$(0, 1, 3, 164)$,
$(0, 4, 13, 210)$,
$(0, 23, 79, 122)$,
$(0, 16, 74, 125)$,\adfsplit
$(0, 35, 81, 158)$,
$(0, 10, 73, 128)$,
$(0, 29, 82, 131)$,
$(0, 11, 41, 127)$,\adfsplit
$(0, 8, 52, 143)$,
$(0, 27, 59, 148)$,
$(0, 57, 114, 171)$,
$(234, 0, 76, 152)$

%ADFvfyBlocksEnd
\adfLgap \noindent by the mapping:
$x \mapsto x +  j \adfmod{228}$ for $x < 228$,
$x \mapsto (x + 2 j \adfmod{6}) + 228$ for $228 \le x < 234$,
$234 \mapsto 234$,
$0 \le j < 228$
 for the first 18 blocks,
$0 \le j < 57$
 for the next block,
$0 \le j < 76$
 for the last block.
\ADFvfyParStart{(235, ((18, 228, ((228, 1), (6, 2), (1, 1))), (1, 57, ((228, 1), (6, 2), (1, 1))), (1, 76, ((228, 1), (6, 2), (1, 1)))), ((19, 12), (7, 1)))} %ADFvfyParEnd
% End of 19^12 7^1
%%%%%%%%%%%%%%%%%%%%%%%%%%%%%%%%%%%%%%%%%%%%%%%%%%%%%%%%%%%%%%%%%%%%%%%%%%%%%%%%%%%%%%%%%%
%%%%%%%%%%%%%%%%%%%%%%%%%%%%%%%%%%%%%%%%%%%%%%%%%%%%%%%%%%%%%%%%%%%%%%%%%%%%%%%%%%%%%%%%%%

% Charlotte:GDD4-1-3-5-mod-6-TeX-gen-A:HITS-fun:4.10
\adfDgap
%ADFvfyBlocksStart {19,19,19,19,19,19,19,19,19,19,19,19,10}
\noindent{\boldmath $ 19^{12} 10^{1} $}~
With the point set $Z_{238}$ partitioned into
 residue classes modulo $12$ for $\{0, 1, \dots, 227\}$, and
 $\{228, 229, \dots, 237\}$,
 the design is generated from

\adfLgap %ADFvfyDesignStart
$(228, 5, 225, 138)$,
$(229, 26, 76, 117)$,
$(230, 167, 97, 64)$,
$(228, 224, 64, 157)$,\adfsplit
$(229, 83, 24, 181)$,
$(230, 218, 204, 123)$,
$(152, 21, 187, 82)$,
$(222, 144, 143, 55)$,\adfsplit
$(23, 200, 208, 54)$,
$(189, 82, 83, 63)$,
$(134, 105, 157, 130)$,
$(203, 212, 111, 40)$,\adfsplit
$(189, 134, 211, 46)$,
$(197, 147, 179, 98)$,
$(131, 13, 2, 176)$,
$(200, 168, 141, 41)$,\adfsplit
$(215, 13, 60, 77)$,
$(208, 89, 131, 42)$,
$(23, 217, 222, 220)$,
$(56, 204, 93, 100)$,\adfsplit
$(24, 162, 87, 124)$,
$(51, 65, 104, 21)$,
$(65, 211, 131, 37)$,
$(78, 44, 9, 65)$,\adfsplit
$(0, 15, 119, 205)$,
$(1, 5, 11, 79)$,
$(0, 5, 22, 217)$,
$(0, 64, 131, 177)$,\adfsplit
$(0, 73, 75, 116)$,
$(0, 7, 98, 123)$,
$(0, 30, 83, 141)$,
$(0, 9, 102, 122)$,\adfsplit
$(0, 3, 18, 218)$,
$(0, 6, 58, 179)$,
$(0, 16, 42, 134)$,
$(0, 19, 162, 207)$,\adfsplit
$(0, 39, 142, 188)$,
$(0, 57, 114, 171)$,
$(237, 0, 76, 152)$,
$(237, 1, 77, 153)$

%ADFvfyBlocksEnd
\adfLgap \noindent by the mapping:
$x \mapsto x + 2 j \adfmod{228}$ for $x < 228$,
$x \mapsto (x - 228 + 3 j \adfmod{9}) + 228$ for $228 \le x < 237$,
$237 \mapsto 237$,
$0 \le j < 114$
 for the first 37 blocks,
$0 \le j < 57$
 for the next block,
$0 \le j < 38$
 for the last two blocks.
\ADFvfyParStart{(238, ((37, 114, ((228, 2), (9, 3), (1, 1))), (1, 57, ((228, 2), (9, 3), (1, 1))), (2, 38, ((228, 2), (9, 3), (1, 1)))), ((19, 12), (10, 1)))} %ADFvfyParEnd
% End of 19^12 10^1
%%%%%%%%%%%%%%%%%%%%%%%%%%%%%%%%%%%%%%%%%%%%%%%%%%%%%%%%%%%%%%%%%%%%%%%%%%%%%%%%%%%%%%%%%%
%%%%%%%%%%%%%%%%%%%%%%%%%%%%%%%%%%%%%%%%%%%%%%%%%%%%%%%%%%%%%%%%%%%%%%%%%%%%%%%%%%%%%%%%%%

% Charlotte:GDD4-1-3-5-mod-6-TeX-gen-A:HITS-fun:4.10
\adfDgap
%ADFvfyBlocksStart {19,19,19,19,19,19,19,19,19,19,19,19,13}
\noindent{\boldmath $ 19^{12} 13^{1} $}~
With the point set $Z_{241}$ partitioned into
 residue classes modulo $12$ for $\{0, 1, \dots, 227\}$, and
 $\{228, 229, \dots, 240\}$,
 the design is generated from

\adfLgap %ADFvfyDesignStart
$(228, 155, 202, 174)$,
$(229, 14, 132, 115)$,
$(230, 40, 80, 114)$,
$(231, 90, 145, 176)$,\adfsplit
$(95, 149, 160, 0)$,
$(35, 152, 29, 33)$,
$(124, 26, 139, 119)$,
$(31, 68, 98, 76)$,\adfsplit
$(168, 47, 186, 22)$,
$(0, 1, 63, 70)$,
$(0, 3, 32, 53)$,
$(0, 13, 56, 153)$,\adfsplit
$(0, 9, 44, 90)$,
$(0, 14, 99, 126)$,
$(0, 33, 91, 157)$,
$(0, 41, 92, 169)$,\adfsplit
$(0, 16, 42, 94)$,
$(0, 23, 61, 148)$,
$(0, 10, 49, 155)$,
$(0, 57, 114, 171)$,\adfsplit
$(240, 0, 76, 152)$

%ADFvfyBlocksEnd
\adfLgap \noindent by the mapping:
$x \mapsto x +  j \adfmod{228}$ for $x < 228$,
$x \mapsto (x + 4 j \adfmod{12}) + 228$ for $228 \le x < 240$,
$240 \mapsto 240$,
$0 \le j < 228$
 for the first 19 blocks,
$0 \le j < 57$
 for the next block,
$0 \le j < 76$
 for the last block.
\ADFvfyParStart{(241, ((19, 228, ((228, 1), (12, 4), (1, 1))), (1, 57, ((228, 1), (12, 4), (1, 1))), (1, 76, ((228, 1), (12, 4), (1, 1)))), ((19, 12), (13, 1)))} %ADFvfyParEnd
% End of 19^12 13^1
%%%%%%%%%%%%%%%%%%%%%%%%%%%%%%%%%%%%%%%%%%%%%%%%%%%%%%%%%%%%%%%%%%%%%%%%%%%%%%%%%%%%%%%%%%
%%%%%%%%%%%%%%%%%%%%%%%%%%%%%%%%%%%%%%%%%%%%%%%%%%%%%%%%%%%%%%%%%%%%%%%%%%%%%%%%%%%%%%%%%%

% Charlotte:GDD4-1-3-5-mod-6-TeX-gen-A:HITS-fun:4.10
\adfDgap
%ADFvfyBlocksStart {19,19,19,19,19,19,19,19,19,19,19,19,16}
\noindent{\boldmath $ 19^{12} 16^{1} $}~
With the point set $Z_{244}$ partitioned into
 residue classes modulo $12$ for $\{0, 1, \dots, 227\}$, and
 $\{228, 229, \dots, 243\}$,
 the design is generated from

\adfLgap %ADFvfyDesignStart
$(228, 47, 56, 114)$,
$(229, 157, 170, 106)$,
$(230, 173, 164, 118)$,
$(231, 202, 199, 14)$,\adfsplit
$(232, 28, 227, 135)$,
$(228, 15, 142, 55)$,
$(229, 75, 11, 156)$,
$(230, 199, 39, 150)$,\adfsplit
$(231, 189, 119, 156)$,
$(232, 7, 174, 26)$,
$(193, 146, 53, 112)$,
$(25, 4, 144, 163)$,\adfsplit
$(80, 171, 202, 209)$,
$(121, 212, 55, 196)$,
$(30, 56, 15, 112)$,
$(182, 87, 24, 44)$,\adfsplit
$(55, 13, 40, 161)$,
$(37, 146, 92, 41)$,
$(87, 200, 218, 82)$,
$(140, 185, 217, 72)$,\adfsplit
$(198, 56, 49, 95)$,
$(122, 83, 133, 211)$,
$(226, 153, 218, 124)$,
$(2, 15, 145, 178)$,\adfsplit
$(42, 44, 43, 147)$,
$(0, 4, 78, 128)$,
$(0, 23, 166, 198)$,
$(0, 22, 95, 184)$,\adfsplit
$(0, 14, 130, 165)$,
$(0, 6, 181, 211)$,
$(0, 25, 27, 186)$,
$(0, 10, 38, 217)$,\adfsplit
$(0, 149, 203, 223)$,
$(0, 3, 37, 59)$,
$(0, 17, 111, 193)$,
$(0, 31, 75, 93)$,\adfsplit
$(0, 41, 99, 127)$,
$(0, 157, 167, 183)$,
$(1, 7, 15, 117)$,
$(0, 57, 114, 171)$,\adfsplit
$(243, 0, 76, 152)$,
$(243, 1, 77, 153)$

%ADFvfyBlocksEnd
\adfLgap \noindent by the mapping:
$x \mapsto x + 2 j \adfmod{228}$ for $x < 228$,
$x \mapsto (x - 228 + 5 j \adfmod{15}) + 228$ for $228 \le x < 243$,
$243 \mapsto 243$,
$0 \le j < 114$
 for the first 39 blocks,
$0 \le j < 57$
 for the next block,
$0 \le j < 38$
 for the last two blocks.
\ADFvfyParStart{(244, ((39, 114, ((228, 2), (15, 5), (1, 1))), (1, 57, ((228, 2), (15, 5), (1, 1))), (2, 38, ((228, 2), (15, 5), (1, 1)))), ((19, 12), (16, 1)))} %ADFvfyParEnd
% End of 19^12 16^1
%%%%%%%%%%%%%%%%%%%%%%%%%%%%%%%%%%%%%%%%%%%%%%%%%%%%%%%%%%%%%%%%%%%%%%%%%%%%%%%%%%%%%%%%%%
%%%%%%%%%%%%%%%%%%%%%%%%%%%%%%%%%%%%%%%%%%%%%%%%%%%%%%%%%%%%%%%%%%%%%%%%%%%%%%%%%%%%%%%%%%

% Charlotte:GDD4-1-3-5-mod-6-TeX-gen-A:HITS-fun:4.10
\adfDgap
%ADFvfyBlocksStart {19,19,19,19,19,19,19,19,19,19,19,19,19,19,19,19,19,19,19,19,19,19,19,19,13}
\noindent{\boldmath $ 19^{24} 13^{1} $}~
With the point set $Z_{469}$ partitioned into
 residue classes modulo $24$ for $\{0, 1, \dots, 455\}$, and
 $\{456, 457, \dots, 468\}$,
 the design is generated from

\adfLgap %ADFvfyDesignStart
$(456, 213, 347, 337)$,
$(457, 73, 152, 321)$,
$(458, 269, 64, 426)$,
$(459, 111, 284, 337)$,\adfsplit
$(178, 387, 299, 146)$,
$(76, 267, 247, 67)$,
$(372, 211, 265, 360)$,
$(216, 261, 119, 378)$,\adfsplit
$(26, 229, 44, 84)$,
$(123, 235, 393, 212)$,
$(353, 326, 175, 394)$,
$(109, 375, 62, 320)$,\adfsplit
$(96, 298, 25, 331)$,
$(42, 249, 376, 415)$,
$(261, 292, 311, 392)$,
$(11, 268, 45, 7)$,\adfsplit
$(375, 49, 449, 95)$,
$(103, 40, 37, 309)$,
$(324, 354, 432, 427)$,
$(88, 102, 94, 242)$,\adfsplit
$(175, 251, 39, 116)$,
$(285, 257, 403, 270)$,
$(226, 420, 295, 330)$,
$(229, 91, 128, 5)$,\adfsplit
$(70, 155, 0, 330)$,
$(0, 1, 52, 395)$,
$(0, 11, 60, 369)$,
$(0, 7, 91, 116)$,\adfsplit
$(0, 17, 132, 291)$,
$(0, 65, 225, 292)$,
$(0, 29, 128, 346)$,
$(0, 2, 44, 286)$,\adfsplit
$(0, 82, 174, 337)$,
$(0, 36, 129, 351)$,
$(0, 57, 213, 277)$,
$(0, 16, 193, 268)$,\adfsplit
$(0, 26, 106, 345)$,
$(0, 21, 43, 210)$,
$(0, 114, 228, 342)$,
$(468, 0, 152, 304)$

%ADFvfyBlocksEnd
\adfLgap \noindent by the mapping:
$x \mapsto x +  j \adfmod{456}$ for $x < 456$,
$x \mapsto (x + 4 j \adfmod{12}) + 456$ for $456 \le x < 468$,
$468 \mapsto 468$,
$0 \le j < 456$
 for the first 38 blocks,
$0 \le j < 114$
 for the next block,
$0 \le j < 152$
 for the last block.
\ADFvfyParStart{(469, ((38, 456, ((456, 1), (12, 4), (1, 1))), (1, 114, ((456, 1), (12, 4), (1, 1))), (1, 152, ((456, 1), (12, 4), (1, 1)))), ((19, 24), (13, 1)))} %ADFvfyParEnd
% End of 19^24 13^1
%%%%%%%%%%%%%%%%%%%%%%%%%%%%%%%%%%%%%%%%%%%%%%%%%%%%%%%%%%%%%%%%%%%%%%%%%%%%%%%%%%%%%%%%%%
%%%%%%%%%%%%%%%%%%%%%%%%%%%%%%%%%%%%%%%%%%%%%%%%%%%%%%%%%%%%%%%%%%%%%%%%%%%%%%%%%%%%%%%%%%

% Charlotte:GDD4-1-3-5-mod-6-TeX-gen-A:HITS-fun:4.10
\adfDgap
%ADFvfyBlocksStart {19,19,19,19,19,19,19,19,19,19,19,19,19,19,19,19,19,19,19,19,19,19,19,19,16}
\noindent{\boldmath $ 19^{24} 16^{1} $}~
With the point set $Z_{472}$ partitioned into
 residue classes modulo $24$ for $\{0, 1, \dots, 455\}$, and
 $\{456, 457, \dots, 471\}$,
 the design is generated from

\adfLgap %ADFvfyDesignStart
$(456, 258, 23, 130)$,
$(456, 97, 345, 68)$,
$(457, 268, 120, 189)$,
$(457, 67, 155, 332)$,\adfsplit
$(458, 340, 363, 114)$,
$(458, 395, 452, 367)$,
$(459, 358, 145, 288)$,
$(459, 98, 141, 143)$,\adfsplit
$(460, 161, 26, 414)$,
$(460, 52, 21, 397)$,
$(101, 47, 89, 145)$,
$(251, 120, 139, 412)$,\adfsplit
$(28, 201, 342, 247)$,
$(153, 101, 438, 188)$,
$(280, 440, 291, 241)$,
$(208, 344, 77, 250)$,\adfsplit
$(191, 142, 354, 376)$,
$(102, 368, 117, 265)$,
$(49, 232, 41, 282)$,
$(397, 110, 141, 163)$,\adfsplit
$(129, 71, 195, 401)$,
$(421, 332, 127, 345)$,
$(308, 79, 379, 28)$,
$(441, 350, 179, 5)$,\adfsplit
$(132, 142, 431, 404)$,
$(222, 454, 216, 298)$,
$(365, 369, 216, 451)$,
$(352, 414, 50, 341)$,\adfsplit
$(381, 390, 188, 448)$,
$(243, 11, 238, 30)$,
$(131, 37, 246, 21)$,
$(108, 71, 247, 96)$,\adfsplit
$(328, 244, 414, 423)$,
$(10, 214, 28, 331)$,
$(299, 264, 407, 98)$,
$(107, 304, 188, 296)$,\adfsplit
$(214, 381, 105, 27)$,
$(106, 13, 183, 102)$,
$(225, 47, 343, 211)$,
$(62, 385, 101, 454)$,\adfsplit
$(277, 79, 260, 422)$,
$(185, 402, 301, 302)$,
$(376, 335, 230, 105)$,
$(63, 229, 84, 452)$,\adfsplit
$(335, 18, 409, 228)$,
$(224, 346, 101, 347)$,
$(175, 272, 238, 297)$,
$(24, 403, 452, 437)$,\adfsplit
$(225, 231, 89, 116)$,
$(111, 439, 129, 28)$,
$(169, 199, 252, 430)$,
$(184, 257, 221, 151)$,\adfsplit
$(23, 28, 325, 386)$,
$(29, 273, 423, 324)$,
$(443, 318, 351, 450)$,
$(1, 76, 402, 128)$,\adfsplit
$(0, 32, 207, 220)$,
$(0, 126, 255, 284)$,
$(0, 133, 159, 263)$,
$(0, 199, 385, 453)$,\adfsplit
$(0, 41, 79, 397)$,
$(0, 21, 241, 281)$,
$(1, 33, 191, 275)$,
$(0, 14, 169, 416)$,\adfsplit
$(0, 65, 75, 344)$,
$(0, 26, 87, 147)$,
$(0, 3, 137, 332)$,
$(0, 102, 282, 439)$,\adfsplit
$(0, 97, 237, 301)$,
$(0, 57, 80, 245)$,
$(0, 66, 319, 409)$,
$(0, 20, 217, 338)$,\adfsplit
$(0, 16, 63, 214)$,
$(0, 7, 322, 378)$,
$(0, 30, 90, 346)$,
$(0, 36, 74, 401)$,\adfsplit
$(0, 2, 46, 352)$,
$(0, 114, 228, 342)$,
$(1, 115, 229, 343)$,
$(471, 0, 152, 304)$,\adfsplit
$(471, 1, 153, 305)$

%ADFvfyBlocksEnd
\adfLgap \noindent by the mapping:
$x \mapsto x + 2 j \adfmod{456}$ for $x < 456$,
$x \mapsto (x - 456 + 5 j \adfmod{15}) + 456$ for $456 \le x < 471$,
$471 \mapsto 471$,
$0 \le j < 228$
 for the first 77 blocks,
$0 \le j < 57$
 for the next two blocks,
$0 \le j < 76$
 for the last two blocks.
\ADFvfyParStart{(472, ((77, 228, ((456, 2), (15, 5), (1, 1))), (2, 57, ((456, 2), (15, 5), (1, 1))), (2, 76, ((456, 2), (15, 5), (1, 1)))), ((19, 24), (16, 1)))} %ADFvfyParEnd
% End of 19^24 16^1
%%%%%%%%%%%%%%%%%%%%%%%%%%%%%%%%%%%%%%%%%%%%%%%%%%%%%%%%%%%%%%%%%%%%%%%%%%%%%%%%%%%%%%%%%%
%%%%%%%%%%%%%%%%%%%%%%%%%%%%%%%%%%%%%%%%%%%%%%%%%%%%%%%%%%%%%%%%%%%%%%%%%%%%%%%%%%%%%%%%%%

% Charlotte:GDD4-1-3-5-mod-6-TeX-gen-A:HITS-fun:4.10
\adfDgap
%ADFvfyBlocksStart {19,19,19,19,19,19,19,19,19,19,19,19,19,19,19,13}
\noindent{\boldmath $ 19^{15} 13^{1} $}~
With the point set $Z_{298}$ partitioned into
 residue classes modulo $15$ for $\{0, 1, \dots, 284\}$, and
 $\{285, 286, \dots, 297\}$,
 the design is generated from

\adfLgap %ADFvfyDesignStart
$(285, 220, 219, 23)$,
$(286, 199, 146, 66)$,
$(287, 100, 171, 185)$,
$(288, 258, 91, 236)$,\adfsplit
$(222, 242, 31, 81)$,
$(31, 189, 161, 57)$,
$(204, 212, 66, 145)$,
$(226, 261, 119, 10)$,\adfsplit
$(257, 239, 215, 27)$,
$(140, 208, 130, 252)$,
$(137, 186, 146, 193)$,
$(179, 166, 249, 272)$,\adfsplit
$(0, 2, 5, 86)$,
$(0, 4, 16, 117)$,
$(0, 6, 25, 234)$,
$(0, 11, 54, 224)$,\adfsplit
$(0, 27, 126, 164)$,
$(0, 39, 103, 149)$,
$(0, 17, 58, 174)$,
$(0, 32, 114, 151)$,\adfsplit
$(0, 21, 87, 123)$,
$(0, 31, 96, 129)$,
$(0, 48, 100, 208)$,
$(0, 29, 91, 222)$,\adfsplit
$(297, 0, 95, 190)$

%ADFvfyBlocksEnd
\adfLgap \noindent by the mapping:
$x \mapsto x +  j \adfmod{285}$ for $x < 285$,
$x \mapsto (x - 285 + 4 j \adfmod{12}) + 285$ for $285 \le x < 297$,
$297 \mapsto 297$,
$0 \le j < 285$
 for the first 24 blocks,
$0 \le j < 95$
 for the last block.
\ADFvfyParStart{(298, ((24, 285, ((285, 1), (12, 4), (1, 1))), (1, 95, ((285, 1), (12, 4), (1, 1)))), ((19, 15), (13, 1)))} %ADFvfyParEnd
% End of 19^15 13^1
%%%%%%%%%%%%%%%%%%%%%%%%%%%%%%%%%%%%%%%%%%%%%%%%%%%%%%%%%%%%%%%%%%%%%%%%%%%%%%%%%%%%%%%%%%
%%%%%%%%%%%%%%%%%%%%%%%%%%%%%%%%%%%%%%%%%%%%%%%%%%%%%%%%%%%%%%%%%%%%%%%%%%%%%%%%%%%%%%%%%%

% Charlotte:GDD4-1-3-5-mod-6-TeX-gen-A:HITS-fun:4.10
\adfDgap
%ADFvfyBlocksStart {19,19,19,19,19,19,19,19,19,10}
\noindent{\boldmath $ 19^{9} 10^{1} $}~
With the point set $Z_{181}$ partitioned into
 residue classes modulo $9$ for $\{0, 1, \dots, 170\}$, and
 $\{171, 172, \dots, 180\}$,
 the design is generated from

\adfLgap %ADFvfyDesignStart
$(171, 80, 48, 32)$,
$(171, 81, 168, 46)$,
$(171, 142, 2, 31)$,
$(95, 92, 136, 117)$,\adfsplit
$(61, 65, 153, 159)$,
$(39, 170, 51, 160)$,
$(0, 1, 56, 67)$,
$(0, 2, 39, 53)$,\adfsplit
$(0, 7, 24, 71)$,
$(0, 5, 26, 96)$,
$(0, 20, 78, 106)$,
$(0, 13, 43, 138)$,\adfsplit
$(0, 8, 42, 110)$,
$(0, 15, 38, 97)$,
$(180, 0, 57, 114)$

%ADFvfyBlocksEnd
\adfLgap \noindent by the mapping:
$x \mapsto x +  j \adfmod{171}$ for $x < 171$,
$x \mapsto (x +  j \adfmod{9}) + 171$ for $171 \le x < 180$,
$180 \mapsto 180$,
$0 \le j < 171$
 for the first 14 blocks,
$0 \le j < 57$
 for the last block.
\ADFvfyParStart{(181, ((14, 171, ((171, 1), (9, 1), (1, 1))), (1, 57, ((171, 1), (9, 1), (1, 1)))), ((19, 9), (10, 1)))} %ADFvfyParEnd
% End of 19^9 10^1
%%%%%%%%%%%%%%%%%%%%%%%%%%%%%%%%%%%%%%%%%%%%%%%%%%%%%%%%%%%%%%%%%%%%%%%%%%%%%%%%%%%%%%%%%%
%%%%%%%%%%%%%%%%%%%%%%%%%%%%%%%%%%%%%%%%%%%%%%%%%%%%%%%%%%%%%%%%%%%%%%%%%%%%%%%%%%%%%%%%%%

% Charlotte:GDD4-1-3-5-mod-6-TeX-gen-A:HITS-fun:4.10
\adfDgap
%ADFvfyBlocksStart {19,19,19,19,19,19,19,19,19,16}
\noindent{\boldmath $ 19^{9} 16^{1} $}~
With the point set $Z_{187}$ partitioned into
 residue classes modulo $9$ for $\{0, 1, \dots, 170\}$, and
 $\{171, 172, \dots, 186\}$,
 the design is generated from

\adfLgap %ADFvfyDesignStart
$(180, 73, 11, 90)$,
$(181, 145, 75, 116)$,
$(171, 99, 113, 17)$,
$(171, 101, 37, 93)$,\adfsplit
$(171, 141, 142, 31)$,
$(110, 99, 138, 136)$,
$(71, 149, 29, 16)$,
$(0, 3, 10, 15)$,\adfsplit
$(0, 20, 52, 85)$,
$(0, 23, 47, 123)$,
$(0, 19, 40, 113)$,
$(0, 16, 69, 104)$,\adfsplit
$(0, 22, 66, 125)$,
$(0, 4, 34, 84)$,
$(0, 6, 31, 128)$,
$(186, 0, 57, 114)$

%ADFvfyBlocksEnd
\adfLgap \noindent by the mapping:
$x \mapsto x +  j \adfmod{171}$ for $x < 171$,
$x \mapsto (x +  j \adfmod{9}) + 171$ for $171 \le x < 180$,
$x \mapsto (x + 2 j \adfmod{6}) + 180$ for $180 \le x < 186$,
$186 \mapsto 186$,
$0 \le j < 171$
 for the first 15 blocks,
$0 \le j < 57$
 for the last block.
\ADFvfyParStart{(187, ((15, 171, ((171, 1), (9, 1), (6, 2), (1, 1))), (1, 57, ((171, 1), (9, 1), (6, 2), (1, 1)))), ((19, 9), (16, 1)))} %ADFvfyParEnd
% End of 19^9 16^1
%%%%%%%%%%%%%%%%%%%%%%%%%%%%%%%%%%%%%%%%%%%%%%%%%%%%%%%%%%%%%%%%%%%%%%%%%%%%%%%%%%%%%%%%%%
%%%%%%%%%%%%%%%%%%%%%%%%%%%%%%%%%%%%%%%%%%%%%%%%%%%%%%%%%%%%%%%%%%%%%%%%%%%%%%%%%%%%%%%%%%

% Charlotte:GDD4-1-3-5-mod-6-TeX-gen-A:HITS-fun:4.10
\adfDgap
%ADFvfyBlocksStart {19,19,19,19,19,19,19,19,19,22}
\noindent{\boldmath $ 19^{9} 22^{1} $}~
With the point set $Z_{193}$ partitioned into
 residue classes modulo $9$ for $\{0, 1, \dots, 170\}$, and
 $\{171, 172, \dots, 192\}$,
 the design is generated from

\adfLgap %ADFvfyDesignStart
$(189, 156, 31, 146)$,
$(171, 109, 47, 132)$,
$(171, 133, 144, 3)$,
$(171, 130, 125, 50)$,\adfsplit
$(172, 67, 141, 119)$,
$(172, 161, 113, 162)$,
$(172, 66, 136, 160)$,
$(55, 42, 84, 131)$,\adfsplit
$(0, 2, 14, 17)$,
$(0, 4, 20, 55)$,
$(0, 19, 50, 83)$,
$(0, 28, 67, 128)$,\adfsplit
$(0, 7, 32, 134)$,
$(0, 8, 68, 106)$,
$(0, 21, 79, 105)$,
$(0, 6, 40, 118)$,\adfsplit
$(192, 0, 57, 114)$

%ADFvfyBlocksEnd
\adfLgap \noindent by the mapping:
$x \mapsto x +  j \adfmod{171}$ for $x < 171$,
$x \mapsto (x - 171 + 2 j \adfmod{18}) + 171$ for $171 \le x < 189$,
$x \mapsto (x +  j \adfmod{3}) + 189$ for $189 \le x < 192$,
$192 \mapsto 192$,
$0 \le j < 171$
 for the first 16 blocks,
$0 \le j < 57$
 for the last block.
\ADFvfyParStart{(193, ((16, 171, ((171, 1), (18, 2), (3, 1), (1, 1))), (1, 57, ((171, 1), (18, 2), (3, 1), (1, 1)))), ((19, 9), (22, 1)))} %ADFvfyParEnd
% End of 19^9 22^1
%%%%%%%%%%%%%%%%%%%%%%%%%%%%%%%%%%%%%%%%%%%%%%%%%%%%%%%%%%%%%%%%%%%%%%%%%%%%%%%%%%%%%%%%%%
%%%%%%%%%%%%%%%%%%%%%%%%%%%%%%%%%%%%%%%%%%%%%%%%%%%%%%%%%%%%%%%%%%%%%%%%%%%%%%%%%%%%%%%%%%

% Charlotte:GDD4-1-3-5-mod-6-TeX-gen-A:HITS-fun:4.10
\adfDgap
%ADFvfyBlocksStart {19,19,19,19,19,19,19,19,19,28}
\noindent{\boldmath $ 19^{9} 28^{1} $}~
With the point set $Z_{199}$ partitioned into
 residue classes modulo $9$ for $\{0, 1, \dots, 170\}$, and
 $\{171, 172, \dots, 198\}$,
 the design is generated from

\adfLgap %ADFvfyDesignStart
$(171, 16, 113, 57)$,
$(171, 51, 19, 85)$,
$(171, 135, 17, 11)$,
$(172, 95, 49, 30)$,\adfsplit
$(172, 151, 108, 116)$,
$(172, 127, 20, 69)$,
$(173, 78, 9, 57)$,
$(173, 88, 143, 11)$,\adfsplit
$(173, 1, 157, 23)$,
$(0, 1, 11, 120)$,
$(0, 4, 29, 79)$,
$(0, 7, 20, 111)$,\adfsplit
$(0, 17, 40, 78)$,
$(0, 12, 71, 95)$,
$(0, 16, 42, 86)$,
$(0, 5, 73, 87)$,\adfsplit
$(0, 2, 30, 33)$,
$(198, 0, 57, 114)$

%ADFvfyBlocksEnd
\adfLgap \noindent by the mapping:
$x \mapsto x +  j \adfmod{171}$ for $x < 171$,
$x \mapsto (x - 171 + 3 j \adfmod{27}) + 171$ for $171 \le x < 198$,
$198 \mapsto 198$,
$0 \le j < 171$
 for the first 17 blocks,
$0 \le j < 57$
 for the last block.
\ADFvfyParStart{(199, ((17, 171, ((171, 1), (27, 3), (1, 1))), (1, 57, ((171, 1), (27, 3), (1, 1)))), ((19, 9), (28, 1)))} %ADFvfyParEnd
% End of 19^9 28^1
%%%%%%%%%%%%%%%%%%%%%%%%%%%%%%%%%%%%%%%%%%%%%%%%%%%%%%%%%%%%%%%%%%%%%%%%%%%%%%%%%%%%%%%%%%
%%%%%%%%%%%%%%%%%%%%%%%%%%%%%%%%%%%%%%%%%%%%%%%%%%%%%%%%%%%%%%%%%%%%%%%%%%%%%%%%%%%%%%%%%%

% Charlotte:GDD4-1-3-5-mod-6-TeX-gen-A:HITS-fun:4.10
\adfDgap
%ADFvfyBlocksStart {19,19,19,19,19,19,19,19,19,34}
\noindent{\boldmath $ 19^{9} 34^{1} $}~
With the point set $Z_{205}$ partitioned into
 residue classes modulo $9$ for $\{0, 1, \dots, 170\}$, and
 $\{171, 172, \dots, 204\}$,
 the design is generated from

\adfLgap %ADFvfyDesignStart
$(198, 16, 30, 128)$,
$(199, 0, 29, 10)$,
$(171, 93, 109, 35)$,
$(171, 13, 2, 169)$,\adfsplit
$(171, 6, 41, 9)$,
$(172, 50, 129, 101)$,
$(172, 63, 150, 26)$,
$(172, 43, 73, 4)$,\adfsplit
$(173, 69, 112, 9)$,
$(173, 74, 95, 7)$,
$(0, 2, 7, 197)$,
$(0, 1, 42, 122)$,\adfsplit
$(0, 22, 48, 137)$,
$(0, 20, 64, 95)$,
$(0, 8, 70, 94)$,
$(0, 13, 38, 138)$,\adfsplit
$(0, 6, 23, 116)$,
$(0, 12, 52, 118)$,
$(204, 0, 57, 114)$

%ADFvfyBlocksEnd
\adfLgap \noindent by the mapping:
$x \mapsto x +  j \adfmod{171}$ for $x < 171$,
$x \mapsto (x - 171 + 3 j \adfmod{27}) + 171$ for $171 \le x < 198$,
$x \mapsto (x + 2 j \adfmod{6}) + 198$ for $198 \le x < 204$,
$204 \mapsto 204$,
$0 \le j < 171$
 for the first 18 blocks,
$0 \le j < 57$
 for the last block.
\ADFvfyParStart{(205, ((18, 171, ((171, 1), (27, 3), (6, 2), (1, 1))), (1, 57, ((171, 1), (27, 3), (6, 2), (1, 1)))), ((19, 9), (34, 1)))} %ADFvfyParEnd
% End of 19^9 34^1
%%%%%%%%%%%%%%%%%%%%%%%%%%%%%%%%%%%%%%%%%%%%%%%%%%%%%%%%%%%%%%%%%%%%%%%%%%%%%%%%%%%%%%%%%%
%%%%%%%%%%%%%%%%%%%%%%%%%%%%%%%%%%%%%%%%%%%%%%%%%%%%%%%%%%%%%%%%%%%%%%%%%%%%%%%%%%%%%%%%%%

% Charlotte:GDD4-1-3-5-mod-6-TeX-gen-A:HITS-fun:4.10
\adfDgap
%ADFvfyBlocksStart {19,19,19,19,19,19,19,19,19,40}
\noindent{\boldmath $ 19^{9} 40^{1} $}~
With the point set $Z_{211}$ partitioned into
 residue classes modulo $9$ for $\{0, 1, \dots, 170\}$, and
 $\{171, 172, \dots, 210\}$,
 the design is generated from

\adfLgap %ADFvfyDesignStart
$(207, 11, 112, 99)$,
$(171, 93, 86, 126)$,
$(171, 79, 15, 164)$,
$(171, 134, 58, 100)$,\adfsplit
$(172, 125, 154, 43)$,
$(172, 22, 23, 47)$,
$(172, 159, 72, 165)$,
$(173, 58, 61, 15)$,\adfsplit
$(173, 144, 91, 107)$,
$(173, 14, 138, 2)$,
$(174, 62, 30, 121)$,
$(0, 4, 14, 19)$,\adfsplit
$(0, 11, 39, 116)$,
$(0, 26, 56, 130)$,
$(0, 2, 51, 71)$,
$(0, 38, 103, 190)$,\adfsplit
$(0, 17, 48, 109)$,
$(0, 21, 44, 119)$,
$(0, 8, 58, 202)$,
$(210, 0, 57, 114)$

%ADFvfyBlocksEnd
\adfLgap \noindent by the mapping:
$x \mapsto x +  j \adfmod{171}$ for $x < 171$,
$x \mapsto (x - 171 + 4 j \adfmod{36}) + 171$ for $171 \le x < 207$,
$x \mapsto (x +  j \adfmod{3}) + 207$ for $207 \le x < 210$,
$210 \mapsto 210$,
$0 \le j < 171$
 for the first 19 blocks,
$0 \le j < 57$
 for the last block.
\ADFvfyParStart{(211, ((19, 171, ((171, 1), (36, 4), (3, 1), (1, 1))), (1, 57, ((171, 1), (36, 4), (3, 1), (1, 1)))), ((19, 9), (40, 1)))} %ADFvfyParEnd
% End of 19^9 40^1
%%%%%%%%%%%%%%%%%%%%%%%%%%%%%%%%%%%%%%%%%%%%%%%%%%%%%%%%%%%%%%%%%%%%%%%%%%%%%%%%%%%%%%%%%%
%%%%%%%%%%%%%%%%%%%%%%%%%%%%%%%%%%%%%%%%%%%%%%%%%%%%%%%%%%%%%%%%%%%%%%%%%%%%%%%%%%%%%%%%%%

% Charlotte:GDD4-1-3-5-mod-6-TeX-gen-A:HITS-fun:4.10
\adfDgap
%ADFvfyBlocksStart {19,19,19,19,19,19,19,19,19,46}
\noindent{\boldmath $ 19^{9} 46^{1} $}~
With the point set $Z_{217}$ partitioned into
 residue classes modulo $9$ for $\{0, 1, \dots, 170\}$, and
 $\{171, 172, \dots, 216\}$,
 the design is generated from

\adfLgap %ADFvfyDesignStart
$(171, 156, 157, 137)$,
$(171, 54, 23, 123)$,
$(171, 28, 170, 7)$,
$(172, 167, 47, 30)$,\adfsplit
$(172, 51, 145, 107)$,
$(172, 135, 70, 76)$,
$(173, 169, 48, 94)$,
$(173, 95, 128, 80)$,\adfsplit
$(173, 144, 154, 78)$,
$(174, 132, 40, 95)$,
$(174, 89, 66, 115)$,
$(0, 2, 64, 174)$,\adfsplit
$(0, 3, 25, 113)$,
$(0, 4, 16, 101)$,
$(0, 24, 68, 111)$,
$(0, 7, 39, 80)$,\adfsplit
$(0, 13, 53, 175)$,
$(0, 11, 104, 195)$,
$(0, 5, 35, 124)$,
$(0, 14, 42, 190)$,\adfsplit
$(216, 0, 57, 114)$

%ADFvfyBlocksEnd
\adfLgap \noindent by the mapping:
$x \mapsto x +  j \adfmod{171}$ for $x < 171$,
$x \mapsto (x - 171 + 5 j \adfmod{45}) + 171$ for $171 \le x < 216$,
$216 \mapsto 216$,
$0 \le j < 171$
 for the first 20 blocks,
$0 \le j < 57$
 for the last block.
\ADFvfyParStart{(217, ((20, 171, ((171, 1), (45, 5), (1, 1))), (1, 57, ((171, 1), (45, 5), (1, 1)))), ((19, 9), (46, 1)))} %ADFvfyParEnd
% End of 19^9 46^1
%%%%%%%%%%%%%%%%%%%%%%%%%%%%%%%%%%%%%%%%%%%%%%%%%%%%%%%%%%%%%%%%%%%%%%%%%%%%%%%%%%%%%%%%%%
%%%%%%%%%%%%%%%%%%%%%%%%%%%%%%%%%%%%%%%%%%%%%%%%%%%%%%%%%%%%%%%%%%%%%%%%%%%%%%%%%%%%%%%%%%

% Charlotte:GDD4-1-3-5-mod-6-TeX-gen-A:HITS-fun:4.10
\adfDgap
%ADFvfyBlocksStart {19,19,19,19,19,19,19,19,19,52}
\noindent{\boldmath $ 19^{9} 52^{1} $}~
With the point set $Z_{223}$ partitioned into
 residue classes modulo $9$ for $\{0, 1, \dots, 170\}$, and
 $\{171, 172, \dots, 222\}$,
 the design is generated from

\adfLgap %ADFvfyDesignStart
$(216, 126, 152, 115)$,
$(217, 23, 15, 79)$,
$(171, 156, 90, 59)$,
$(171, 51, 97, 98)$,\adfsplit
$(171, 100, 11, 112)$,
$(172, 40, 68, 30)$,
$(172, 116, 146, 91)$,
$(172, 27, 51, 7)$,\adfsplit
$(173, 38, 99, 86)$,
$(173, 43, 165, 28)$,
$(173, 132, 40, 44)$,
$(174, 117, 74, 167)$,\adfsplit
$(0, 2, 35, 41)$,
$(0, 3, 17, 103)$,
$(0, 19, 51, 180)$,
$(0, 40, 98, 179)$,\adfsplit
$(0, 21, 80, 200)$,
$(0, 5, 67, 190)$,
$(0, 23, 75, 214)$,
$(0, 7, 29, 94)$,\adfsplit
$(0, 16, 69, 111)$,
$(222, 0, 57, 114)$

%ADFvfyBlocksEnd
\adfLgap \noindent by the mapping:
$x \mapsto x +  j \adfmod{171}$ for $x < 171$,
$x \mapsto (x - 171 + 5 j \adfmod{45}) + 171$ for $171 \le x < 216$,
$x \mapsto (x + 2 j \adfmod{6}) + 216$ for $216 \le x < 222$,
$222 \mapsto 222$,
$0 \le j < 171$
 for the first 21 blocks,
$0 \le j < 57$
 for the last block.
\ADFvfyParStart{(223, ((21, 171, ((171, 1), (45, 5), (6, 2), (1, 1))), (1, 57, ((171, 1), (45, 5), (6, 2), (1, 1)))), ((19, 9), (52, 1)))} %ADFvfyParEnd
% End of 19^9 52^1
%%%%%%%%%%%%%%%%%%%%%%%%%%%%%%%%%%%%%%%%%%%%%%%%%%%%%%%%%%%%%%%%%%%%%%%%%%%%%%%%%%%%%%%%%%
%%%%%%%%%%%%%%%%%%%%%%%%%%%%%%%%%%%%%%%%%%%%%%%%%%%%%%%%%%%%%%%%%%%%%%%%%%%%%%%%%%%%%%%%%%

% Charlotte:GDD4-1-3-5-mod-6-TeX-gen-A:HITS-fun:4.10
\adfDgap
%ADFvfyBlocksStart {19,19,19,19,19,19,19,19,19,58}
\noindent{\boldmath $ 19^{9} 58^{1} $}~
With the point set $Z_{229}$ partitioned into
 residue classes modulo $9$ for $\{0, 1, \dots, 170\}$, and
 $\{171, 172, \dots, 228\}$,
 the design is generated from

\adfLgap %ADFvfyDesignStart
$(225, 50, 123, 91)$,
$(171, 61, 62, 149)$,
$(171, 135, 20, 141)$,
$(171, 22, 100, 138)$,\adfsplit
$(172, 134, 45, 23)$,
$(172, 109, 40, 114)$,
$(172, 111, 128, 52)$,
$(173, 40, 59, 79)$,\adfsplit
$(173, 110, 89, 19)$,
$(173, 42, 90, 57)$,
$(174, 18, 158, 11)$,
$(174, 82, 116, 69)$,\adfsplit
$(174, 112, 3, 70)$,
$(0, 2, 14, 127)$,
$(0, 10, 26, 61)$,
$(0, 8, 85, 175)$,\adfsplit
$(0, 30, 122, 206)$,
$(0, 11, 64, 199)$,
$(0, 28, 71, 217)$,
$(0, 4, 29, 224)$,\adfsplit
$(0, 23, 119, 194)$,
$(0, 3, 40, 106)$,
$(228, 0, 57, 114)$

%ADFvfyBlocksEnd
\adfLgap \noindent by the mapping:
$x \mapsto x +  j \adfmod{171}$ for $x < 171$,
$x \mapsto (x - 171 + 6 j \adfmod{54}) + 171$ for $171 \le x < 225$,
$x \mapsto (x +  j \adfmod{3}) + 225$ for $225 \le x < 228$,
$228 \mapsto 228$,
$0 \le j < 171$
 for the first 22 blocks,
$0 \le j < 57$
 for the last block.
\ADFvfyParStart{(229, ((22, 171, ((171, 1), (54, 6), (3, 1), (1, 1))), (1, 57, ((171, 1), (54, 6), (3, 1), (1, 1)))), ((19, 9), (58, 1)))} %ADFvfyParEnd
% End of 19^9 58^1
%%%%%%%%%%%%%%%%%%%%%%%%%%%%%%%%%%%%%%%%%%%%%%%%%%%%%%%%%%%%%%%%%%%%%%%%%%%%%%%%%%%%%%%%%%
%%%%%%%%%%%%%%%%%%%%%%%%%%%%%%%%%%%%%%%%%%%%%%%%%%%%%%%%%%%%%%%%%%%%%%%%%%%%%%%%%%%%%%%%%%

% Charlotte:GDD4-1-3-5-mod-6-TeX-gen-A:HITS-fun:4.10
\adfDgap
%ADFvfyBlocksStart {19,19,19,19,19,19,19,19,19,64}
\noindent{\boldmath $ 19^{9} 64^{1} $}~
With the point set $Z_{235}$ partitioned into
 residue classes modulo $9$ for $\{0, 1, \dots, 170\}$, and
 $\{171, 172, \dots, 234\}$,
 the design is generated from

\adfLgap %ADFvfyDesignStart
$(171, 24, 165, 101)$,
$(171, 115, 59, 108)$,
$(171, 35, 94, 73)$,
$(172, 39, 164, 81)$,\adfsplit
$(172, 132, 52, 77)$,
$(172, 37, 112, 134)$,
$(173, 165, 13, 151)$,
$(173, 159, 23, 29)$,\adfsplit
$(173, 62, 10, 144)$,
$(174, 137, 125, 145)$,
$(174, 81, 112, 78)$,
$(174, 43, 75, 32)$,\adfsplit
$(175, 56, 156, 33)$,
$(0, 1, 5, 182)$,
$(0, 2, 158, 203)$,
$(0, 26, 73, 102)$,\adfsplit
$(0, 16, 60, 121)$,
$(0, 10, 68, 183)$,
$(0, 24, 86, 204)$,
$(0, 40, 93, 218)$,\adfsplit
$(0, 28, 79, 184)$,
$(0, 39, 104, 226)$,
$(0, 17, 87, 212)$,
$(234, 0, 57, 114)$

%ADFvfyBlocksEnd
\adfLgap \noindent by the mapping:
$x \mapsto x +  j \adfmod{171}$ for $x < 171$,
$x \mapsto (x - 171 + 7 j \adfmod{63}) + 171$ for $171 \le x < 234$,
$234 \mapsto 234$,
$0 \le j < 171$
 for the first 23 blocks,
$0 \le j < 57$
 for the last block.
\ADFvfyParStart{(235, ((23, 171, ((171, 1), (63, 7), (1, 1))), (1, 57, ((171, 1), (63, 7), (1, 1)))), ((19, 9), (64, 1)))} %ADFvfyParEnd
% End of 19^9 64^1
%%%%%%%%%%%%%%%%%%%%%%%%%%%%%%%%%%%%%%%%%%%%%%%%%%%%%%%%%%%%%%%%%%%%%%%%%%%%%%%%%%%%%%%%%%
%%%%%%%%%%%%%%%%%%%%%%%%%%%%%%%%%%%%%%%%%%%%%%%%%%%%%%%%%%%%%%%%%%%%%%%%%%%%%%%%%%%%%%%%%%

% Charlotte:GDD4-1-3-5-mod-6-TeX-gen-A:HITS-fun:4.10
\adfDgap
%ADFvfyBlocksStart {19,19,19,19,19,19,19,19,19,70}
\noindent{\boldmath $ 19^{9} 70^{1} $}~
With the point set $Z_{241}$ partitioned into
 residue classes modulo $9$ for $\{0, 1, \dots, 170\}$, and
 $\{171, 172, \dots, 240\}$,
 the design is generated from

\adfLgap %ADFvfyDesignStart
$(234, 60, 43, 41)$,
$(235, 12, 82, 11)$,
$(171, 113, 109, 146)$,
$(171, 78, 106, 0)$,\adfsplit
$(171, 157, 143, 93)$,
$(172, 99, 133, 14)$,
$(172, 136, 40, 33)$,
$(172, 146, 152, 93)$,\adfsplit
$(173, 136, 47, 124)$,
$(173, 15, 107, 58)$,
$(173, 99, 147, 23)$,
$(174, 12, 63, 25)$,\adfsplit
$(174, 1, 137, 140)$,
$(174, 98, 103, 159)$,
$(0, 8, 23, 155)$,
$(0, 20, 46, 189)$,\adfsplit
$(0, 10, 141, 210)$,
$(0, 22, 109, 224)$,
$(0, 21, 104, 211)$,
$(0, 11, 69, 183)$,\adfsplit
$(0, 31, 91, 225)$,
$(0, 29, 73, 191)$,
$(0, 25, 130, 198)$,
$(0, 42, 97, 205)$,\adfsplit
$(240, 0, 57, 114)$

%ADFvfyBlocksEnd
\adfLgap \noindent by the mapping:
$x \mapsto x +  j \adfmod{171}$ for $x < 171$,
$x \mapsto (x - 171 + 7 j \adfmod{63}) + 171$ for $171 \le x < 234$,
$x \mapsto (x + 2 j \adfmod{6}) + 234$ for $234 \le x < 240$,
$240 \mapsto 240$,
$0 \le j < 171$
 for the first 24 blocks,
$0 \le j < 57$
 for the last block.
\ADFvfyParStart{(241, ((24, 171, ((171, 1), (63, 7), (6, 2), (1, 1))), (1, 57, ((171, 1), (63, 7), (6, 2), (1, 1)))), ((19, 9), (70, 1)))} %ADFvfyParEnd
% End of 19^9 70^1
%%%%%%%%%%%%%%%%%%%%%%%%%%%%%%%%%%%%%%%%%%%%%%%%%%%%%%%%%%%%%%%%%%%%%%%%%%%%%%%%%%%%%%%%%%
%%%%%%%%%%%%%%%%%%%%%%%%%%%%%%%%%%%%%%%%%%%%%%%%%%%%%%%%%%%%%%%%%%%%%%%%%%%%%%%%%%%%%%%%%%

% Charlotte:GDD4-1-3-5-mod-6-TeX-gen-A:HITS-fun:4.10
\adfDgap
%ADFvfyBlocksStart {19,19,19,19,19,19,19,19,19,19,19,19,19,19,19,19,19,19,19,19,19,16}
\noindent{\boldmath $ 19^{21} 16^{1} $}~
With the point set $Z_{415}$ partitioned into
 residue classes modulo $21$ for $\{0, 1, \dots, 398\}$, and
 $\{399, 400, \dots, 414\}$,
 the design is generated from

\adfLgap %ADFvfyDesignStart
$(399, 153, 241, 50)$,
$(400, 142, 204, 290)$,
$(401, 365, 324, 376)$,
$(402, 121, 189, 287)$,\adfsplit
$(403, 84, 134, 364)$,
$(279, 202, 92, 39)$,
$(177, 132, 309, 266)$,
$(192, 353, 356, 388)$,\adfsplit
$(329, 2, 101, 95)$,
$(278, 80, 247, 263)$,
$(184, 78, 221, 202)$,
$(231, 78, 331, 386)$,\adfsplit
$(342, 381, 97, 367)$,
$(4, 320, 261, 377)$,
$(273, 198, 360, 222)$,
$(375, 327, 120, 6)$,\adfsplit
$(7, 254, 35, 249)$,
$(255, 293, 165, 44)$,
$(204, 128, 20, 155)$,
$(19, 349, 371, 248)$,\adfsplit
$(387, 238, 379, 333)$,
$(22, 9, 398, 129)$,
$(262, 8, 88, 44)$,
$(0, 1, 67, 326)$,\adfsplit
$(0, 20, 112, 263)$,
$(0, 34, 113, 224)$,
$(0, 12, 97, 206)$,
$(0, 9, 70, 195)$,\adfsplit
$(0, 17, 82, 199)$,
$(0, 29, 131, 202)$,
$(0, 7, 40, 179)$,
$(0, 56, 137, 295)$,\adfsplit
$(0, 4, 64, 122)$,
$(0, 2, 96, 223)$,
$(414, 0, 133, 266)$

%ADFvfyBlocksEnd
\adfLgap \noindent by the mapping:
$x \mapsto x +  j \adfmod{399}$ for $x < 399$,
$x \mapsto (x - 399 + 5 j \adfmod{15}) + 399$ for $399 \le x < 414$,
$414 \mapsto 414$,
$0 \le j < 399$
 for the first 34 blocks,
$0 \le j < 133$
 for the last block.
\ADFvfyParStart{(415, ((34, 399, ((399, 1), (15, 5), (1, 1))), (1, 133, ((399, 1), (15, 5), (1, 1)))), ((19, 21), (16, 1)))} %ADFvfyParEnd
% End of 19^21 16^1
%%%%%%%%%%%%%%%%%%%%%%%%%%%%%%%%%%%%%%%%%%%%%%%%%%%%%%%%%%%%%%%%%%%%%%%%%%%%%%%%%%%%%%%%%%
%%%%%%%%%%%%%%%%%%%%%%%%%%%%%%%%%%%%%%%%%%%%%%%%%%%%%%%%%%%%%%%%%%%%%%%%%%%%%%%%%%%%%%%%%%

%%%%%%%%%%%%%%%%%%%%%%%%%%%%%%%%%%%%%%%%%%%%%%%%%%%%%%%%%%%%%%%%%%%%%%%%%%%%%%%%%%%%%%%%%%
%%%%%%%%%%%%%%%%%%%%%%%%%%%%%%%%%%%%%%%%%%%%%%%%%%%%%%%%%%%%%%%%%%%%%%%%%%%%%%%%%%%%%%%%%%
\section{4-GDDs for the proof of Lemma \ref{lem:4-GDD 23^u m^1}}
\label{app:4-GDD 23^u m^1}
\adfnull{
$ 23^{12} 14^1 $,
$ 23^{12} 17^1 $,
$ 23^{12} 20^1 $,
$ 23^{15} 17^1 $,
$ 23^9 14^1 $,
$ 23^9 20^1 $,
$ 23^9 26^1 $,
$ 23^9 32^1 $,
$ 23^9 38^1 $,
$ 23^9 44^1 $,
$ 23^9 50^1 $,
$ 23^9 56^1 $,
$ 23^9 62^1 $,
$ 23^9 68^1 $,
$ 23^9 74^1 $,
$ 23^9 80^1 $ and
$ 23^9 86^1 $.
}

% Charlotte:GDD4-1-3-5-mod-6-TeX-gen-A:HITS-fun:4.10
\adfDgap
%ADFvfyBlocksStart {23,23,23,23,23,23,23,23,23,23,23,23,14}
\noindent{\boldmath $ 23^{12} 14^{1} $}~
With the point set $Z_{290}$ partitioned into
 residue classes modulo $12$ for $\{0, 1, \dots, 275\}$, and
 $\{276, 277, \dots, 289\}$,
 the design is generated from

\adfLgap %ADFvfyDesignStart
$(276, 0, 1, 2)$,
$(277, 0, 91, 185)$,
$(278, 0, 94, 92)$,
$(279, 0, 184, 275)$,\adfsplit
$(280, 0, 274, 182)$,
$(281, 0, 4, 11)$,
$(282, 0, 7, 272)$,
$(283, 0, 265, 269)$,\adfsplit
$(284, 0, 10, 23)$,
$(285, 0, 13, 266)$,
$(286, 0, 253, 263)$,
$(287, 0, 16, 35)$,\adfsplit
$(288, 0, 19, 260)$,
$(289, 0, 241, 257)$,
$(88, 204, 41, 49)$,
$(242, 192, 129, 217)$,\adfsplit
$(78, 95, 149, 100)$,
$(214, 107, 150, 272)$,
$(201, 28, 207, 143)$,
$(21, 82, 219, 191)$,\adfsplit
$(270, 190, 61, 204)$,
$(134, 89, 3, 204)$,
$(134, 92, 251, 137)$,
$(105, 126, 167, 135)$,\adfsplit
$(272, 163, 130, 216)$,
$(55, 248, 191, 96)$,
$(244, 177, 140, 107)$,
$(254, 79, 183, 59)$,\adfsplit
$(165, 266, 263, 187)$,
$(87, 129, 251, 40)$,
$(33, 258, 28, 159)$,
$(215, 32, 86, 5)$,\adfsplit
$(0, 6, 93, 193)$,
$(0, 8, 39, 149)$,
$(0, 15, 59, 217)$,
$(0, 21, 151, 236)$,\adfsplit
$(0, 26, 53, 105)$,
$(0, 85, 103, 153)$,
$(0, 125, 165, 267)$,
$(0, 119, 157, 247)$,\adfsplit
$(0, 63, 77, 150)$,
$(0, 181, 227, 261)$,
$(0, 143, 199, 225)$,
$(0, 99, 136, 259)$,\adfsplit
$(0, 30, 62, 127)$,
$(0, 57, 88, 162)$,
$(0, 20, 178, 221)$,
$(0, 38, 90, 166)$,\adfsplit
$(0, 28, 106, 135)$,
$(0, 18, 73, 194)$,
$(0, 34, 102, 146)$,
$(0, 69, 138, 207)$

%ADFvfyBlocksEnd
\adfLgap \noindent by the mapping:
$x \mapsto x + 3 j \adfmod{276}$ for $x < 276$,
$x \mapsto x$ for $x \ge 276$,
$0 \le j < 92$
 for the first 14 blocks;
$x \mapsto x + 2 j \adfmod{276}$ for $x < 276$,
$x \mapsto x$ for $x \ge 276$,
$0 \le j < 138$
 for the next 37 blocks,
$0 \le j < 69$
 for the last block.
\ADFvfyParStart{(290, ((14, 92, ((276, 3), (14, 14))), (37, 138, ((276, 2), (14, 14))), (1, 69, ((276, 2), (14, 14)))), ((23, 12), (14, 1)))} %ADFvfyParEnd
% End of 23^12 14^1
%%%%%%%%%%%%%%%%%%%%%%%%%%%%%%%%%%%%%%%%%%%%%%%%%%%%%%%%%%%%%%%%%%%%%%%%%%%%%%%%%%%%%%%%%%
%%%%%%%%%%%%%%%%%%%%%%%%%%%%%%%%%%%%%%%%%%%%%%%%%%%%%%%%%%%%%%%%%%%%%%%%%%%%%%%%%%%%%%%%%%

% Charlotte:GDD4-1-3-5-mod-6-TeX-gen-A:HITS-fun:4.10
\adfDgap
%ADFvfyBlocksStart {23,23,23,23,23,23,23,23,23,23,23,23,17}
\noindent{\boldmath $ 23^{12} 17^{1} $}~
With the point set $Z_{293}$ partitioned into
 residue classes modulo $12$ for $\{0, 1, \dots, 275\}$, and
 $\{276, 277, \dots, 292\}$,
 the design is generated from

\adfLgap %ADFvfyDesignStart
$(276, 0, 1, 2)$,
$(277, 0, 91, 185)$,
$(278, 0, 94, 92)$,
$(279, 0, 184, 275)$,\adfsplit
$(280, 0, 274, 182)$,
$(281, 0, 4, 11)$,
$(282, 0, 7, 272)$,
$(283, 0, 265, 269)$,\adfsplit
$(284, 0, 10, 23)$,
$(285, 0, 13, 266)$,
$(286, 0, 253, 263)$,
$(287, 0, 16, 35)$,\adfsplit
$(288, 0, 19, 260)$,
$(289, 0, 241, 257)$,
$(290, 0, 22, 5)$,
$(291, 0, 259, 254)$,\adfsplit
$(292, 0, 271, 17)$,
$(112, 211, 57, 193)$,
$(162, 233, 207, 130)$,
$(235, 226, 117, 263)$,\adfsplit
$(221, 255, 6, 215)$,
$(65, 186, 44, 225)$,
$(0, 3, 68, 83)$,
$(0, 8, 33, 157)$,\adfsplit
$(0, 14, 58, 78)$,
$(0, 29, 59, 170)$,
$(0, 46, 93, 163)$,
$(0, 38, 89, 175)$,\adfsplit
$(0, 50, 104, 201)$,
$(0, 52, 114, 202)$,
$(0, 49, 112, 197)$,
$(0, 41, 107, 194)$,\adfsplit
$(0, 31, 131, 174)$,
$(0, 53, 110, 200)$,
$(0, 42, 98, 203)$,
$(0, 69, 138, 207)$

%ADFvfyBlocksEnd
\adfLgap \noindent by the mapping:
$x \mapsto x + 3 j \adfmod{276}$ for $x < 276$,
$x \mapsto x$ for $x \ge 276$,
$0 \le j < 92$
 for the first 17 blocks;
$x \mapsto x +  j \adfmod{276}$ for $x < 276$,
$x \mapsto x$ for $x \ge 276$,
$0 \le j < 276$
 for the next 18 blocks,
$0 \le j < 69$
 for the last block.
\ADFvfyParStart{(293, ((17, 92, ((276, 3), (17, 17))), (18, 276, ((276, 1), (17, 17))), (1, 69, ((276, 1), (17, 17)))), ((23, 12), (17, 1)))} %ADFvfyParEnd
% End of 23^12 17^1
%%%%%%%%%%%%%%%%%%%%%%%%%%%%%%%%%%%%%%%%%%%%%%%%%%%%%%%%%%%%%%%%%%%%%%%%%%%%%%%%%%%%%%%%%%
%%%%%%%%%%%%%%%%%%%%%%%%%%%%%%%%%%%%%%%%%%%%%%%%%%%%%%%%%%%%%%%%%%%%%%%%%%%%%%%%%%%%%%%%%%

% Charlotte:GDD4-1-3-5-mod-6-TeX-gen-A:HITS-fun:4.10
\adfDgap
%ADFvfyBlocksStart {23,23,23,23,23,23,23,23,23,23,23,23,20}
\noindent{\boldmath $ 23^{12} 20^{1} $}~
With the point set $Z_{296}$ partitioned into
 residue classes modulo $12$ for $\{0, 1, \dots, 275\}$, and
 $\{276, 277, \dots, 295\}$,
 the design is generated from

\adfLgap %ADFvfyDesignStart
$(276, 0, 1, 2)$,
$(277, 0, 91, 185)$,
$(278, 0, 94, 92)$,
$(279, 0, 184, 275)$,\adfsplit
$(280, 0, 274, 182)$,
$(281, 0, 4, 11)$,
$(282, 0, 7, 272)$,
$(283, 0, 265, 269)$,\adfsplit
$(284, 0, 10, 23)$,
$(285, 0, 13, 266)$,
$(286, 0, 253, 263)$,
$(287, 0, 16, 35)$,\adfsplit
$(288, 0, 19, 260)$,
$(289, 0, 241, 257)$,
$(290, 0, 22, 5)$,
$(291, 0, 259, 254)$,\adfsplit
$(292, 0, 271, 17)$,
$(293, 0, 25, 53)$,
$(294, 0, 28, 251)$,
$(295, 0, 223, 248)$,\adfsplit
$(26, 84, 90, 177)$,
$(168, 199, 39, 125)$,
$(137, 204, 267, 100)$,
$(217, 43, 179, 34)$,\adfsplit
$(30, 106, 56, 135)$,
$(261, 96, 197, 158)$,
$(36, 158, 57, 172)$,
$(111, 50, 12, 189)$,\adfsplit
$(248, 37, 114, 40)$,
$(83, 206, 141, 12)$,
$(153, 120, 35, 266)$,
$(118, 66, 228, 29)$,\adfsplit
$(159, 246, 238, 45)$,
$(221, 112, 270, 167)$,
$(102, 217, 36, 272)$,
$(98, 67, 18, 219)$,\adfsplit
$(235, 255, 42, 165)$,
$(119, 60, 133, 185)$,
$(0, 15, 150, 206)$,
$(0, 20, 54, 98)$,\adfsplit
$(0, 3, 18, 247)$,
$(0, 32, 159, 235)$,
$(0, 30, 116, 249)$,
$(0, 47, 81, 152)$,\adfsplit
$(0, 27, 67, 186)$,
$(0, 41, 102, 237)$,
$(0, 43, 131, 176)$,
$(0, 57, 155, 164)$,\adfsplit
$(0, 111, 119, 137)$,
$(0, 42, 149, 230)$,
$(0, 113, 225, 255)$,
$(0, 51, 97, 179)$,\adfsplit
$(0, 89, 148, 243)$,
$(1, 7, 69, 111)$,
$(1, 45, 101, 151)$,
$(0, 69, 138, 207)$

%ADFvfyBlocksEnd
\adfLgap \noindent by the mapping:
$x \mapsto x + 3 j \adfmod{276}$ for $x < 276$,
$x \mapsto x$ for $x \ge 276$,
$0 \le j < 92$
 for the first 20 blocks;
$x \mapsto x + 2 j \adfmod{276}$ for $x < 276$,
$x \mapsto x$ for $x \ge 276$,
$0 \le j < 138$
 for the next 35 blocks,
$0 \le j < 69$
 for the last block.
\ADFvfyParStart{(296, ((20, 92, ((276, 3), (20, 20))), (35, 138, ((276, 2), (20, 20))), (1, 69, ((276, 2), (20, 20)))), ((23, 12), (20, 1)))} %ADFvfyParEnd
% End of 23^12 20^1
%%%%%%%%%%%%%%%%%%%%%%%%%%%%%%%%%%%%%%%%%%%%%%%%%%%%%%%%%%%%%%%%%%%%%%%%%%%%%%%%%%%%%%%%%%
%%%%%%%%%%%%%%%%%%%%%%%%%%%%%%%%%%%%%%%%%%%%%%%%%%%%%%%%%%%%%%%%%%%%%%%%%%%%%%%%%%%%%%%%%%

% Charlotte:GDD4-1-3-5-mod-6-TeX-gen-A:HITS-fun:4.10
\adfDgap
%ADFvfyBlocksStart {23,23,23,23,23,23,23,23,23,23,23,23,23,23,23,17}
\noindent{\boldmath $ 23^{15} 17^{1} $}~
With the point set $Z_{362}$ partitioned into
 residue classes modulo $15$ for $\{0, 1, \dots, 344\}$, and
 $\{345, 346, \dots, 361\}$,
 the design is generated from

\adfLgap %ADFvfyDesignStart
$(345, 0, 1, 2)$,
$(346, 0, 115, 113)$,
$(347, 0, 229, 344)$,
$(348, 0, 232, 116)$,\adfsplit
$(349, 0, 343, 230)$,
$(350, 0, 4, 11)$,
$(351, 0, 7, 341)$,
$(352, 0, 334, 338)$,\adfsplit
$(353, 0, 10, 23)$,
$(354, 0, 13, 335)$,
$(355, 0, 322, 332)$,
$(356, 0, 16, 35)$,\adfsplit
$(357, 0, 19, 329)$,
$(358, 0, 310, 326)$,
$(359, 0, 22, 5)$,
$(360, 0, 328, 323)$,\adfsplit
$(361, 0, 340, 17)$,
$(251, 113, 283, 139)$,
$(37, 255, 191, 138)$,
$(118, 147, 249, 127)$,\adfsplit
$(92, 254, 292, 55)$,
$(114, 183, 311, 238)$,
$(149, 13, 304, 102)$,
$(72, 45, 319, 224)$,\adfsplit
$(170, 338, 227, 209)$,
$(163, 151, 107, 293)$,
$(155, 308, 205, 238)$,
$(23, 97, 29, 54)$,\adfsplit
$(80, 167, 219, 34)$,
$(0, 3, 61, 97)$,
$(0, 8, 59, 253)$,
$(0, 49, 137, 213)$,\adfsplit
$(0, 14, 80, 241)$,
$(0, 65, 171, 238)$,
$(0, 34, 96, 254)$,
$(0, 24, 109, 235)$,\adfsplit
$(0, 42, 114, 198)$,
$(0, 21, 99, 178)$,
$(0, 28, 121, 169)$,
$(0, 40, 81, 163)$,\adfsplit
$(0, 63, 140, 259)$

%ADFvfyBlocksEnd
\adfLgap \noindent by the mapping:
$x \mapsto x + 3 j \adfmod{345}$ for $x < 345$,
$x \mapsto x$ for $x \ge 345$,
$0 \le j < 115$
 for the first 17 blocks;
$x \mapsto x +  j \adfmod{345}$ for $x < 345$,
$x \mapsto x$ for $x \ge 345$,
$0 \le j < 345$
 for the last 24 blocks.
\ADFvfyParStart{(362, ((17, 115, ((345, 3), (17, 17))), (24, 345, ((345, 1), (17, 17)))), ((23, 15), (17, 1)))} %ADFvfyParEnd
% End of 23^15 17^1
%%%%%%%%%%%%%%%%%%%%%%%%%%%%%%%%%%%%%%%%%%%%%%%%%%%%%%%%%%%%%%%%%%%%%%%%%%%%%%%%%%%%%%%%%%
%%%%%%%%%%%%%%%%%%%%%%%%%%%%%%%%%%%%%%%%%%%%%%%%%%%%%%%%%%%%%%%%%%%%%%%%%%%%%%%%%%%%%%%%%%

% Charlotte:GDD4-1-3-5-mod-6-TeX-gen-A:HITS-fun:4.10
\adfDgap
%ADFvfyBlocksStart {23,23,23,23,23,23,23,23,23,14}
\noindent{\boldmath $ 23^{9} 14^{1} $}~
With the point set $Z_{221}$ partitioned into
 residue classes modulo $9$ for $\{0, 1, \dots, 206\}$, and
 $\{207, 208, \dots, 220\}$,
 the design is generated from

\adfLgap %ADFvfyDesignStart
$(219, 0, 1, 2)$,
$(220, 0, 70, 140)$,
$(207, 100, 110, 141)$,
$(208, 62, 154, 15)$,\adfsplit
$(209, 191, 22, 45)$,
$(210, 53, 148, 183)$,
$(194, 85, 26, 73)$,
$(179, 136, 42, 121)$,\adfsplit
$(0, 3, 7, 69)$,
$(0, 5, 11, 134)$,
$(0, 8, 19, 178)$,
$(0, 13, 30, 46)$,\adfsplit
$(0, 28, 59, 147)$,
$(0, 17, 97, 118)$,
$(0, 22, 74, 125)$,
$(0, 24, 53, 111)$,\adfsplit
$(0, 40, 93, 142)$,
$(0, 25, 67, 150)$,
$(0, 23, 55, 130)$

%ADFvfyBlocksEnd
\adfLgap \noindent by the mapping:
$x \mapsto x \oplus (3 j)$ for $x < 207$,
$x \mapsto (x - 207 + 4 j \adfmod{12}) + 207$ for $207 \le x < 219$,
$x \mapsto x$ for $x \ge 219$,
$0 \le j < 69$
 for the first two blocks;
$x \mapsto x \oplus j$ for $x < 207$,
$x \mapsto (x - 207 + 4 j \adfmod{12}) + 207$ for $207 \le x < 219$,
$x \mapsto x$ for $x \ge 219$,
$0 \le j < 207$
 for the last 17 blocks.
\ADFvfyParStart{(221, ((2, 69, ((207, 3, (69, 3)), (12, 4), (2, 2))), (17, 207, ((207, 1, (69, 3)), (12, 4), (2, 2)))), ((23, 9), (14, 1)))} %ADFvfyParEnd
% End of 23^9 14^1
%%%%%%%%%%%%%%%%%%%%%%%%%%%%%%%%%%%%%%%%%%%%%%%%%%%%%%%%%%%%%%%%%%%%%%%%%%%%%%%%%%%%%%%%%%
%%%%%%%%%%%%%%%%%%%%%%%%%%%%%%%%%%%%%%%%%%%%%%%%%%%%%%%%%%%%%%%%%%%%%%%%%%%%%%%%%%%%%%%%%%

% Charlotte:GDD4-1-3-5-mod-6-TeX-gen-A:HITS-fun:4.10
\adfDgap
%ADFvfyBlocksStart {23,23,23,23,23,23,23,23,23,20}
\noindent{\boldmath $ 23^{9} 20^{1} $}~
With the point set $Z_{227}$ partitioned into
 residue classes modulo $9$ for $\{0, 1, \dots, 206\}$, and
 $\{207, 208, \dots, 226\}$,
 the design is generated from

\adfLgap %ADFvfyDesignStart
$(225, 0, 1, 2)$,
$(226, 0, 70, 140)$,
$(207, 52, 35, 198)$,
$(208, 123, 146, 130)$,\adfsplit
$(209, 111, 109, 23)$,
$(210, 58, 117, 182)$,
$(211, 111, 151, 38)$,
$(212, 69, 50, 16)$,\adfsplit
$(188, 51, 41, 70)$,
$(15, 54, 66, 58)$,
$(180, 186, 194, 165)$,
$(0, 3, 55, 69)$,\adfsplit
$(0, 26, 101, 129)$,
$(0, 25, 58, 125)$,
$(0, 24, 77, 111)$,
$(0, 10, 74, 122)$,\adfsplit
$(0, 28, 84, 130)$,
$(0, 30, 71, 113)$,
$(0, 38, 95, 142)$,
$(0, 35, 79, 128)$

%ADFvfyBlocksEnd
\adfLgap \noindent by the mapping:
$x \mapsto x \oplus (3 j)$ for $x < 207$,
$x \mapsto (x - 207 + 6 j \adfmod{18}) + 207$ for $207 \le x < 225$,
$x \mapsto x$ for $x \ge 225$,
$0 \le j < 69$
 for the first two blocks;
$x \mapsto x \oplus j$ for $x < 207$,
$x \mapsto (x - 207 + 6 j \adfmod{18}) + 207$ for $207 \le x < 225$,
$x \mapsto x$ for $x \ge 225$,
$0 \le j < 207$
 for the last 18 blocks.
\ADFvfyParStart{(227, ((2, 69, ((207, 3, (69, 3)), (18, 6), (2, 2))), (18, 207, ((207, 1, (69, 3)), (18, 6), (2, 2)))), ((23, 9), (20, 1)))} %ADFvfyParEnd
% End of 23^9 20^1
%%%%%%%%%%%%%%%%%%%%%%%%%%%%%%%%%%%%%%%%%%%%%%%%%%%%%%%%%%%%%%%%%%%%%%%%%%%%%%%%%%%%%%%%%%
%%%%%%%%%%%%%%%%%%%%%%%%%%%%%%%%%%%%%%%%%%%%%%%%%%%%%%%%%%%%%%%%%%%%%%%%%%%%%%%%%%%%%%%%%%

% Charlotte:GDD4-1-3-5-mod-6-TeX-gen-A:HITS-fun:4.10
\adfDgap
%ADFvfyBlocksStart {23,23,23,23,23,23,23,23,23,26}
\noindent{\boldmath $ 23^{9} 26^{1} $}~
With the point set $Z_{233}$ partitioned into
 residue classes modulo $9$ for $\{0, 1, \dots, 206\}$, and
 $\{207, 208, \dots, 232\}$,
 the design is generated from

\adfLgap %ADFvfyDesignStart
$(231, 0, 1, 2)$,
$(232, 0, 70, 140)$,
$(207, 206, 171, 196)$,
$(208, 80, 4, 12)$,\adfsplit
$(209, 179, 81, 190)$,
$(210, 164, 25, 75)$,
$(211, 184, 78, 200)$,
$(212, 174, 16, 131)$,\adfsplit
$(213, 93, 2, 154)$,
$(214, 188, 169, 51)$,
$(158, 136, 13, 143)$,
$(83, 8, 163, 14)$,\adfsplit
$(0, 3, 32, 110)$,
$(0, 4, 24, 37)$,
$(0, 17, 38, 77)$,
$(0, 40, 93, 136)$,\adfsplit
$(0, 8, 50, 91)$,
$(0, 12, 59, 125)$,
$(0, 20, 48, 176)$,
$(0, 5, 30, 150)$,\adfsplit
$(0, 26, 64, 131)$

%ADFvfyBlocksEnd
\adfLgap \noindent by the mapping:
$x \mapsto x \oplus (3 j)$ for $x < 207$,
$x \mapsto (x - 207 + 8 j \adfmod{24}) + 207$ for $207 \le x < 231$,
$x \mapsto x$ for $x \ge 231$,
$0 \le j < 69$
 for the first two blocks;
$x \mapsto x \oplus j$ for $x < 207$,
$x \mapsto (x - 207 + 8 j \adfmod{24}) + 207$ for $207 \le x < 231$,
$x \mapsto x$ for $x \ge 231$,
$0 \le j < 207$
 for the last 19 blocks.
\ADFvfyParStart{(233, ((2, 69, ((207, 3, (69, 3)), (24, 8), (2, 2))), (19, 207, ((207, 1, (69, 3)), (24, 8), (2, 2)))), ((23, 9), (26, 1)))} %ADFvfyParEnd
% End of 23^9 26^1
%%%%%%%%%%%%%%%%%%%%%%%%%%%%%%%%%%%%%%%%%%%%%%%%%%%%%%%%%%%%%%%%%%%%%%%%%%%%%%%%%%%%%%%%%%
%%%%%%%%%%%%%%%%%%%%%%%%%%%%%%%%%%%%%%%%%%%%%%%%%%%%%%%%%%%%%%%%%%%%%%%%%%%%%%%%%%%%%%%%%%

% Charlotte:GDD4-1-3-5-mod-6-TeX-gen-A:HITS-fun:4.10
\adfDgap
%ADFvfyBlocksStart {23,23,23,23,23,23,23,23,23,32}
\noindent{\boldmath $ 23^{9} 32^{1} $}~
With the point set $Z_{239}$ partitioned into
 residue classes modulo $9$ for $\{0, 1, \dots, 206\}$, and
 $\{207, 208, \dots, 238\}$,
 the design is generated from

\adfLgap %ADFvfyDesignStart
$(237, 0, 1, 2)$,
$(238, 0, 70, 140)$,
$(207, 134, 54, 37)$,
$(208, 26, 73, 174)$,\adfsplit
$(209, 130, 41, 9)$,
$(210, 127, 125, 99)$,
$(211, 186, 8, 13)$,
$(212, 94, 80, 162)$,\adfsplit
$(213, 135, 14, 184)$,
$(214, 54, 19, 185)$,
$(215, 181, 122, 69)$,
$(216, 147, 194, 202)$,\adfsplit
$(0, 3, 7, 69)$,
$(0, 6, 16, 156)$,
$(0, 12, 60, 95)$,
$(0, 21, 52, 105)$,\adfsplit
$(0, 15, 40, 122)$,
$(0, 30, 73, 164)$,
$(0, 22, 64, 129)$,
$(0, 13, 24, 111)$,\adfsplit
$(0, 23, 61, 94)$,
$(0, 19, 58, 133)$

%ADFvfyBlocksEnd
\adfLgap \noindent by the mapping:
$x \mapsto x \oplus (3 j)$ for $x < 207$,
$x \mapsto (x - 207 + 10 j \adfmod{30}) + 207$ for $207 \le x < 237$,
$x \mapsto x$ for $x \ge 237$,
$0 \le j < 69$
 for the first two blocks;
$x \mapsto x \oplus j$ for $x < 207$,
$x \mapsto (x - 207 + 10 j \adfmod{30}) + 207$ for $207 \le x < 237$,
$x \mapsto x$ for $x \ge 237$,
$0 \le j < 207$
 for the last 20 blocks.
\ADFvfyParStart{(239, ((2, 69, ((207, 3, (69, 3)), (30, 10), (2, 2))), (20, 207, ((207, 1, (69, 3)), (30, 10), (2, 2)))), ((23, 9), (32, 1)))} %ADFvfyParEnd
% End of 23^9 32^1
%%%%%%%%%%%%%%%%%%%%%%%%%%%%%%%%%%%%%%%%%%%%%%%%%%%%%%%%%%%%%%%%%%%%%%%%%%%%%%%%%%%%%%%%%%
%%%%%%%%%%%%%%%%%%%%%%%%%%%%%%%%%%%%%%%%%%%%%%%%%%%%%%%%%%%%%%%%%%%%%%%%%%%%%%%%%%%%%%%%%%

% Charlotte:GDD4-1-3-5-mod-6-TeX-gen-A:HITS-fun:4.10
\adfDgap
%ADFvfyBlocksStart {23,23,23,23,23,23,23,23,23,38}
\noindent{\boldmath $ 23^{9} 38^{1} $}~
With the point set $Z_{245}$ partitioned into
 residue classes modulo $9$ for $\{0, 1, \dots, 206\}$, and
 $\{207, 208, \dots, 244\}$,
 the design is generated from

\adfLgap %ADFvfyDesignStart
$(243, 0, 1, 2)$,
$(244, 0, 70, 140)$,
$(207, 180, 53, 118)$,
$(208, 201, 128, 199)$,\adfsplit
$(209, 180, 127, 134)$,
$(210, 192, 191, 148)$,
$(211, 6, 61, 164)$,
$(212, 126, 37, 149)$,\adfsplit
$(213, 134, 76, 12)$,
$(214, 30, 38, 199)$,
$(215, 166, 156, 119)$,
$(216, 127, 44, 33)$,\adfsplit
$(20, 134, 185, 1)$,
$(0, 13, 29, 218)$,
$(0, 14, 28, 217)$,
$(0, 12, 33, 95)$,\adfsplit
$(0, 24, 91, 120)$,
$(0, 34, 73, 139)$,
$(0, 30, 61, 109)$,
$(0, 20, 89, 123)$,\adfsplit
$(0, 6, 44, 59)$,
$(0, 3, 60, 85)$,
$(0, 22, 97, 129)$

%ADFvfyBlocksEnd
\adfLgap \noindent by the mapping:
$x \mapsto x \oplus (3 j)$ for $x < 207$,
$x \mapsto (x - 207 + 12 j \adfmod{36}) + 207$ for $207 \le x < 243$,
$x \mapsto x$ for $x \ge 243$,
$0 \le j < 69$
 for the first two blocks;
$x \mapsto x \oplus j$ for $x < 207$,
$x \mapsto (x - 207 + 12 j \adfmod{36}) + 207$ for $207 \le x < 243$,
$x \mapsto x$ for $x \ge 243$,
$0 \le j < 207$
 for the last 21 blocks.
\ADFvfyParStart{(245, ((2, 69, ((207, 3, (69, 3)), (36, 12), (2, 2))), (21, 207, ((207, 1, (69, 3)), (36, 12), (2, 2)))), ((23, 9), (38, 1)))} %ADFvfyParEnd
% End of 23^9 38^1
%%%%%%%%%%%%%%%%%%%%%%%%%%%%%%%%%%%%%%%%%%%%%%%%%%%%%%%%%%%%%%%%%%%%%%%%%%%%%%%%%%%%%%%%%%
%%%%%%%%%%%%%%%%%%%%%%%%%%%%%%%%%%%%%%%%%%%%%%%%%%%%%%%%%%%%%%%%%%%%%%%%%%%%%%%%%%%%%%%%%%

% Charlotte:GDD4-1-3-5-mod-6-TeX-gen-A:HITS-fun:4.10
\adfDgap
%ADFvfyBlocksStart {23,23,23,23,23,23,23,23,23,44}
\noindent{\boldmath $ 23^{9} 44^{1} $}~
With the point set $Z_{251}$ partitioned into
 residue classes modulo $9$ for $\{0, 1, \dots, 206\}$, and
 $\{207, 208, \dots, 250\}$,
 the design is generated from

\adfLgap %ADFvfyDesignStart
$(249, 0, 1, 2)$,
$(250, 0, 70, 140)$,
$(207, 200, 64, 177)$,
$(208, 80, 24, 91)$,\adfsplit
$(209, 180, 65, 22)$,
$(210, 109, 171, 62)$,
$(211, 138, 206, 193)$,
$(212, 74, 37, 45)$,\adfsplit
$(213, 118, 78, 50)$,
$(214, 72, 167, 34)$,
$(215, 99, 160, 56)$,
$(216, 187, 39, 14)$,\adfsplit
$(217, 39, 20, 139)$,
$(0, 4, 20, 219)$,
$(0, 5, 10, 220)$,
$(0, 7, 26, 218)$,\adfsplit
$(0, 6, 66, 119)$,
$(0, 3, 76, 132)$,
$(0, 15, 39, 121)$,
$(0, 32, 83, 125)$,\adfsplit
$(0, 17, 38, 86)$,
$(0, 30, 87, 131)$,
$(0, 25, 58, 102)$,
$(0, 12, 64, 96)$

%ADFvfyBlocksEnd
\adfLgap \noindent by the mapping:
$x \mapsto x \oplus (3 j)$ for $x < 207$,
$x \mapsto (x - 207 + 14 j \adfmod{42}) + 207$ for $207 \le x < 249$,
$x \mapsto x$ for $x \ge 249$,
$0 \le j < 69$
 for the first two blocks;
$x \mapsto x \oplus j$ for $x < 207$,
$x \mapsto (x - 207 + 14 j \adfmod{42}) + 207$ for $207 \le x < 249$,
$x \mapsto x$ for $x \ge 249$,
$0 \le j < 207$
 for the last 22 blocks.
\ADFvfyParStart{(251, ((2, 69, ((207, 3, (69, 3)), (42, 14), (2, 2))), (22, 207, ((207, 1, (69, 3)), (42, 14), (2, 2)))), ((23, 9), (44, 1)))} %ADFvfyParEnd
% End of 23^9 44^1
%%%%%%%%%%%%%%%%%%%%%%%%%%%%%%%%%%%%%%%%%%%%%%%%%%%%%%%%%%%%%%%%%%%%%%%%%%%%%%%%%%%%%%%%%%
%%%%%%%%%%%%%%%%%%%%%%%%%%%%%%%%%%%%%%%%%%%%%%%%%%%%%%%%%%%%%%%%%%%%%%%%%%%%%%%%%%%%%%%%%%

% Charlotte:GDD4-1-3-5-mod-6-TeX-gen-A:HITS-fun:4.10
\adfDgap
%ADFvfyBlocksStart {23,23,23,23,23,23,23,23,23,50}
\noindent{\boldmath $ 23^{9} 50^{1} $}~
With the point set $Z_{257}$ partitioned into
 residue classes modulo $9$ for $\{0, 1, \dots, 206\}$, and
 $\{207, 208, \dots, 256\}$,
 the design is generated from

\adfLgap %ADFvfyDesignStart
$(255, 24, 148, 197)$,
$(256, 151, 182, 162)$,
$(207, 160, 1, 15)$,
$(207, 175, 27, 57)$,\adfsplit
$(207, 23, 152, 119)$,
$(208, 20, 33, 138)$,
$(208, 140, 145, 40)$,
$(208, 27, 143, 178)$,\adfsplit
$(209, 198, 24, 104)$,
$(209, 70, 92, 181)$,
$(209, 53, 57, 112)$,
$(210, 24, 71, 85)$,\adfsplit
$(210, 66, 124, 189)$,
$(210, 64, 68, 92)$,
$(211, 134, 140, 196)$,
$(211, 55, 15, 18)$,\adfsplit
$(211, 164, 148, 120)$,
$(212, 0, 164, 94)$,
$(212, 16, 39, 14)$,
$(212, 195, 179, 64)$,\adfsplit
$(213, 94, 90, 119)$,
$(213, 143, 82, 111)$,
$(213, 5, 34, 87)$,
$(214, 127, 105, 47)$,\adfsplit
$(214, 135, 57, 130)$,
$(214, 140, 206, 43)$,
$(215, 39, 195, 92)$,
$(215, 67, 80, 122)$,\adfsplit
$(215, 36, 70, 136)$,
$(216, 157, 131, 73)$,
$(216, 114, 155, 174)$,
$(216, 160, 44, 180)$,\adfsplit
$(217, 34, 81, 49)$,
$(217, 68, 6, 143)$,
$(217, 20, 163, 84)$,
$(218, 79, 156, 68)$,\adfsplit
$(218, 109, 197, 112)$,
$(218, 2, 0, 69)$,
$(219, 31, 178, 57)$,
$(219, 145, 137, 152)$,\adfsplit
$(219, 15, 54, 32)$,
$(220, 94, 77, 25)$,
$(220, 189, 154, 69)$,
$(220, 134, 12, 20)$,\adfsplit
$(221, 25, 13, 182)$,
$(221, 138, 161, 123)$,
$(221, 68, 180, 199)$,
$(222, 44, 121, 123)$,\adfsplit
$(222, 14, 156, 83)$,
$(222, 81, 82, 7)$,
$(0, 5, 12, 25)$,
$(0, 6, 16, 83)$,\adfsplit
$(0, 11, 66, 159)$,
$(0, 7, 111, 161)$,
$(0, 64, 101, 158)$,
$(0, 14, 26, 74)$,\adfsplit
$(0, 68, 98, 166)$,
$(0, 42, 139, 169)$,
$(0, 49, 57, 163)$,
$(0, 75, 167, 190)$,\adfsplit
$(0, 31, 82, 206)$,
$(0, 21, 56, 107)$,
$(0, 52, 109, 155)$,
$(0, 43, 176, 179)$,\adfsplit
$(0, 112, 136, 146)$,
$(0, 24, 70, 91)$,
$(0, 110, 157, 197)$,
$(1, 20, 122, 143)$,\adfsplit
$(0, 88, 131, 170)$,
$(1, 2, 43, 130)$,
$(1, 7, 40, 113)$

%ADFvfyBlocksEnd
\adfLgap \noindent by the mapping:
$x \mapsto x + 3 j \adfmod{207}$ for $x < 207$,
$x \mapsto (x - 207 + 16 j \adfmod{48}) + 207$ for $207 \le x < 255$,
$x \mapsto x$ for $x \ge 255$,
$0 \le j < 69$.
\ADFvfyParStart{(257, ((71, 69, ((207, 3), (48, 16), (2, 2)))), ((23, 9), (50, 1)))} %ADFvfyParEnd
% End of 23^9 50^1
%%%%%%%%%%%%%%%%%%%%%%%%%%%%%%%%%%%%%%%%%%%%%%%%%%%%%%%%%%%%%%%%%%%%%%%%%%%%%%%%%%%%%%%%%%
%%%%%%%%%%%%%%%%%%%%%%%%%%%%%%%%%%%%%%%%%%%%%%%%%%%%%%%%%%%%%%%%%%%%%%%%%%%%%%%%%%%%%%%%%%

% Charlotte:GDD4-1-3-5-mod-6-TeX-gen-A:HITS-fun:4.10
\adfDgap
%ADFvfyBlocksStart {23,23,23,23,23,23,23,23,23,56}
\noindent{\boldmath $ 23^{9} 56^{1} $}~
With the point set $Z_{263}$ partitioned into
 residue classes modulo $9$ for $\{0, 1, \dots, 206\}$, and
 $\{207, 208, \dots, 262\}$,
 the design is generated from

\adfLgap %ADFvfyDesignStart
$(261, 90, 121, 8)$,
$(262, 91, 87, 134)$,
$(207, 42, 30, 134)$,
$(207, 79, 45, 22)$,\adfsplit
$(207, 20, 77, 199)$,
$(208, 150, 157, 113)$,
$(208, 57, 73, 196)$,
$(208, 179, 2, 99)$,\adfsplit
$(209, 178, 112, 146)$,
$(209, 113, 174, 99)$,
$(209, 80, 6, 28)$,
$(210, 181, 169, 60)$,\adfsplit
$(210, 119, 44, 167)$,
$(210, 108, 31, 66)$,
$(211, 63, 178, 194)$,
$(211, 195, 38, 73)$,\adfsplit
$(211, 120, 116, 13)$,
$(212, 16, 95, 63)$,
$(212, 164, 190, 12)$,
$(212, 141, 17, 202)$,\adfsplit
$(213, 51, 80, 40)$,
$(213, 12, 52, 9)$,
$(213, 119, 23, 73)$,
$(214, 8, 203, 82)$,\adfsplit
$(214, 7, 165, 159)$,
$(214, 29, 198, 139)$,
$(215, 109, 38, 17)$,
$(215, 148, 78, 138)$,\adfsplit
$(215, 32, 70, 153)$,
$(216, 143, 110, 139)$,
$(216, 54, 192, 87)$,
$(216, 194, 28, 133)$,\adfsplit
$(217, 187, 104, 184)$,
$(217, 46, 189, 191)$,
$(217, 141, 107, 84)$,
$(218, 29, 26, 94)$,\adfsplit
$(218, 187, 167, 19)$,
$(218, 144, 165, 24)$,
$(219, 103, 117, 141)$,
$(219, 149, 62, 73)$,\adfsplit
$(219, 102, 169, 155)$,
$(220, 181, 105, 4)$,
$(220, 117, 152, 47)$,
$(220, 185, 66, 25)$,\adfsplit
$(221, 18, 1, 43)$,
$(221, 175, 78, 68)$,
$(221, 183, 53, 92)$,
$(222, 152, 153, 57)$,\adfsplit
$(222, 199, 191, 121)$,
$(222, 70, 95, 24)$,
$(223, 50, 143, 189)$,
$(223, 64, 65, 138)$,\adfsplit
$(223, 186, 124, 184)$,
$(0, 5, 154, 202)$,
$(0, 11, 62, 151)$,
$(0, 1, 156, 175)$,\adfsplit
$(0, 37, 88, 157)$,
$(1, 20, 76, 97)$,
$(0, 101, 199, 224)$,
$(0, 94, 185, 187)$,\adfsplit
$(0, 103, 118, 167)$,
$(0, 56, 73, 79)$,
$(0, 163, 200, 242)$,
$(0, 15, 122, 127)$,\adfsplit
$(0, 30, 128, 159)$,
$(0, 91, 164, 179)$,
$(0, 136, 191, 260)$,
$(0, 58, 82, 89)$,\adfsplit
$(0, 8, 149, 155)$,
$(0, 13, 26, 39)$,
$(0, 20, 114, 158)$,
$(0, 28, 140, 182)$,\adfsplit
$(0, 17, 41, 142)$,
$(0, 59, 123, 188)$

%ADFvfyBlocksEnd
\adfLgap \noindent by the mapping:
$x \mapsto x + 3 j \adfmod{207}$ for $x < 207$,
$x \mapsto (x - 207 + 18 j \adfmod{54}) + 207$ for $207 \le x < 261$,
$x \mapsto x$ for $x \ge 261$,
$0 \le j < 69$.
\ADFvfyParStart{(263, ((74, 69, ((207, 3), (54, 18), (2, 2)))), ((23, 9), (56, 1)))} %ADFvfyParEnd
% End of 23^9 56^1
%%%%%%%%%%%%%%%%%%%%%%%%%%%%%%%%%%%%%%%%%%%%%%%%%%%%%%%%%%%%%%%%%%%%%%%%%%%%%%%%%%%%%%%%%%
%%%%%%%%%%%%%%%%%%%%%%%%%%%%%%%%%%%%%%%%%%%%%%%%%%%%%%%%%%%%%%%%%%%%%%%%%%%%%%%%%%%%%%%%%%

% Charlotte:GDD4-1-3-5-mod-6-TeX-gen-A:HITS-fun:4.10
\adfDgap
%ADFvfyBlocksStart {23,23,23,23,23,23,23,23,23,62}
\noindent{\boldmath $ 23^{9} 62^{1} $}~
With the point set $Z_{269}$ partitioned into
 residue classes modulo $9$ for $\{0, 1, \dots, 206\}$, and
 $\{207, 208, \dots, 268\}$,
 the design is generated from

\adfLgap %ADFvfyDesignStart
$(267, 125, 117, 142)$,
$(268, 145, 146, 171)$,
$(207, 174, 105, 8)$,
$(207, 41, 154, 70)$,\adfsplit
$(207, 81, 74, 166)$,
$(208, 115, 55, 171)$,
$(208, 57, 53, 2)$,
$(208, 130, 104, 141)$,\adfsplit
$(209, 188, 4, 135)$,
$(209, 91, 196, 123)$,
$(209, 110, 12, 50)$,
$(210, 186, 21, 200)$,\adfsplit
$(210, 89, 54, 14)$,
$(210, 169, 145, 184)$,
$(211, 166, 155, 194)$,
$(211, 16, 90, 190)$,\adfsplit
$(211, 186, 147, 125)$,
$(212, 63, 181, 74)$,
$(212, 152, 183, 70)$,
$(212, 130, 24, 50)$,\adfsplit
$(213, 110, 162, 199)$,
$(213, 24, 5, 26)$,
$(213, 201, 151, 202)$,
$(214, 97, 19, 23)$,\adfsplit
$(214, 90, 110, 76)$,
$(214, 24, 120, 107)$,
$(215, 164, 3, 134)$,
$(215, 94, 189, 52)$,\adfsplit
$(215, 181, 77, 33)$,
$(216, 2, 184, 98)$,
$(216, 150, 198, 140)$,
$(216, 124, 21, 136)$,\adfsplit
$(217, 117, 194, 31)$,
$(217, 46, 39, 56)$,
$(217, 96, 188, 124)$,
$(218, 134, 32, 135)$,\adfsplit
$(218, 119, 73, 67)$,
$(218, 105, 12, 70)$,
$(219, 73, 16, 114)$,
$(219, 200, 158, 99)$,\adfsplit
$(219, 188, 148, 93)$,
$(220, 33, 165, 197)$,
$(220, 74, 115, 67)$,
$(220, 99, 46, 104)$,\adfsplit
$(221, 70, 138, 5)$,
$(221, 163, 49, 27)$,
$(221, 96, 62, 29)$,
$(222, 130, 143, 159)$,\adfsplit
$(222, 34, 122, 55)$,
$(222, 9, 129, 56)$,
$(223, 170, 122, 16)$,
$(223, 192, 28, 186)$,\adfsplit
$(223, 162, 139, 92)$,
$(224, 156, 136, 180)$,
$(224, 11, 202, 125)$,
$(0, 10, 71, 244)$,\adfsplit
$(0, 4, 21, 123)$,
$(0, 12, 52, 125)$,
$(0, 3, 33, 202)$,
$(0, 23, 29, 150)$,\adfsplit
$(0, 15, 61, 225)$,
$(0, 31, 68, 156)$,
$(0, 60, 122, 226)$,
$(0, 19, 66, 116)$,\adfsplit
$(0, 79, 128, 265)$,
$(0, 13, 16, 78)$,
$(0, 34, 130, 266)$,
$(0, 64, 107, 205)$,\adfsplit
$(0, 56, 124, 143)$,
$(0, 67, 89, 158)$,
$(0, 65, 97, 127)$,
$(1, 32, 56, 70)$,\adfsplit
$(1, 71, 86, 121)$,
$(1, 76, 188, 200)$,
$(1, 80, 137, 203)$,
$(1, 77, 206, 265)$,\adfsplit
$(1, 149, 152, 266)$

%ADFvfyBlocksEnd
\adfLgap \noindent by the mapping:
$x \mapsto x + 3 j \adfmod{207}$ for $x < 207$,
$x \mapsto (x - 207 + 20 j \adfmod{60}) + 207$ for $207 \le x < 267$,
$x \mapsto x$ for $x \ge 267$,
$0 \le j < 69$.
\ADFvfyParStart{(269, ((77, 69, ((207, 3), (60, 20), (2, 2)))), ((23, 9), (62, 1)))} %ADFvfyParEnd
% End of 23^9 62^1
%%%%%%%%%%%%%%%%%%%%%%%%%%%%%%%%%%%%%%%%%%%%%%%%%%%%%%%%%%%%%%%%%%%%%%%%%%%%%%%%%%%%%%%%%%
%%%%%%%%%%%%%%%%%%%%%%%%%%%%%%%%%%%%%%%%%%%%%%%%%%%%%%%%%%%%%%%%%%%%%%%%%%%%%%%%%%%%%%%%%%

% Charlotte:GDD4-1-3-5-mod-6-TeX-gen-A:HITS-fun:4.10
\adfDgap
%ADFvfyBlocksStart {23,23,23,23,23,23,23,23,23,68}
\noindent{\boldmath $ 23^{9} 68^{1} $}~
With the point set $Z_{275}$ partitioned into
 residue classes modulo $9$ for $\{0, 1, \dots, 206\}$, and
 $\{207, 208, \dots, 274\}$,
 the design is generated from

\adfLgap %ADFvfyDesignStart
$(273, 120, 85, 188)$,
$(274, 21, 106, 44)$,
$(207, 119, 183, 199)$,
$(207, 105, 197, 104)$,\adfsplit
$(207, 169, 36, 13)$,
$(208, 178, 199, 155)$,
$(208, 18, 166, 195)$,
$(208, 147, 149, 44)$,\adfsplit
$(209, 103, 150, 109)$,
$(209, 126, 173, 34)$,
$(209, 161, 75, 86)$,
$(210, 60, 97, 182)$,\adfsplit
$(210, 158, 154, 17)$,
$(210, 18, 112, 30)$,
$(211, 141, 148, 3)$,
$(211, 97, 62, 0)$,\adfsplit
$(211, 92, 14, 109)$,
$(212, 177, 94, 93)$,
$(212, 137, 188, 36)$,
$(212, 43, 59, 91)$,\adfsplit
$(213, 34, 6, 54)$,
$(213, 94, 138, 68)$,
$(213, 92, 179, 91)$,
$(214, 110, 156, 45)$,\adfsplit
$(214, 167, 125, 31)$,
$(214, 7, 127, 204)$,
$(215, 169, 4, 105)$,
$(215, 37, 102, 8)$,\adfsplit
$(215, 131, 110, 0)$,
$(216, 191, 42, 160)$,
$(216, 9, 4, 82)$,
$(216, 30, 71, 23)$,\adfsplit
$(217, 57, 186, 95)$,
$(217, 82, 191, 22)$,
$(217, 169, 53, 18)$,
$(218, 40, 7, 80)$,\adfsplit
$(218, 154, 186, 99)$,
$(218, 174, 65, 149)$,
$(219, 7, 199, 74)$,
$(219, 161, 167, 87)$,\adfsplit
$(219, 90, 183, 103)$,
$(220, 161, 198, 148)$,
$(220, 111, 11, 169)$,
$(220, 167, 46, 42)$,\adfsplit
$(221, 43, 114, 81)$,
$(221, 44, 10, 156)$,
$(221, 149, 31, 164)$,
$(222, 21, 206, 130)$,\adfsplit
$(222, 135, 60, 47)$,
$(222, 122, 61, 127)$,
$(223, 37, 61, 47)$,
$(223, 188, 156, 198)$,\adfsplit
$(223, 33, 122, 76)$,
$(224, 23, 197, 20)$,
$(224, 87, 157, 66)$,
$(0, 19, 88, 224)$,\adfsplit
$(0, 3, 52, 158)$,
$(0, 5, 141, 192)$,
$(0, 26, 105, 164)$,
$(0, 17, 77, 205)$,\adfsplit
$(0, 8, 24, 100)$,
$(0, 14, 67, 79)$,
$(0, 6, 173, 225)$,
$(0, 31, 34, 56)$,\adfsplit
$(0, 46, 60, 226)$,
$(0, 39, 83, 249)$,
$(0, 50, 146, 228)$,
$(0, 29, 40, 57)$,\adfsplit
$(0, 134, 154, 270)$,
$(0, 181, 188, 271)$,
$(0, 53, 103, 203)$,
$(0, 22, 199, 269)$,\adfsplit
$(0, 112, 140, 272)$,
$(1, 125, 149, 225)$,
$(1, 44, 85, 249)$,
$(0, 25, 121, 176)$,\adfsplit
$(1, 53, 65, 248)$,
$(0, 139, 196, 250)$,
$(1, 59, 98, 106)$,
$(1, 38, 40, 115)$

%ADFvfyBlocksEnd
\adfLgap \noindent by the mapping:
$x \mapsto x + 3 j \adfmod{207}$ for $x < 207$,
$x \mapsto (x - 207 + 22 j \adfmod{66}) + 207$ for $207 \le x < 273$,
$x \mapsto x$ for $x \ge 273$,
$0 \le j < 69$.
\ADFvfyParStart{(275, ((80, 69, ((207, 3), (66, 22), (2, 2)))), ((23, 9), (68, 1)))} %ADFvfyParEnd
% End of 23^9 68^1
%%%%%%%%%%%%%%%%%%%%%%%%%%%%%%%%%%%%%%%%%%%%%%%%%%%%%%%%%%%%%%%%%%%%%%%%%%%%%%%%%%%%%%%%%%
%%%%%%%%%%%%%%%%%%%%%%%%%%%%%%%%%%%%%%%%%%%%%%%%%%%%%%%%%%%%%%%%%%%%%%%%%%%%%%%%%%%%%%%%%%

% Charlotte:GDD4-1-3-5-mod-6-TeX-gen-A:HITS-fun:4.10
\adfDgap
%ADFvfyBlocksStart {23,23,23,23,23,23,23,23,23,74}
\noindent{\boldmath $ 23^{9} 74^{1} $}~
With the point set $Z_{281}$ partitioned into
 residue classes modulo $9$ for $\{0, 1, \dots, 206\}$, and
 $\{207, 208, \dots, 280\}$,
 the design is generated from

\adfLgap %ADFvfyDesignStart
$(279, 173, 111, 52)$,
$(280, 17, 18, 52)$,
$(207, 41, 180, 98)$,
$(207, 109, 13, 66)$,\adfsplit
$(207, 29, 160, 141)$,
$(208, 115, 32, 92)$,
$(208, 161, 103, 78)$,
$(208, 1, 102, 45)$,\adfsplit
$(209, 98, 82, 21)$,
$(209, 180, 141, 115)$,
$(209, 101, 5, 202)$,
$(210, 203, 25, 45)$,\adfsplit
$(210, 183, 35, 109)$,
$(210, 58, 20, 186)$,
$(211, 4, 117, 190)$,
$(211, 89, 204, 14)$,\adfsplit
$(211, 182, 129, 187)$,
$(212, 98, 203, 177)$,
$(212, 63, 148, 181)$,
$(212, 11, 178, 21)$,\adfsplit
$(213, 13, 95, 24)$,
$(213, 44, 90, 56)$,
$(213, 106, 82, 201)$,
$(214, 82, 38, 75)$,\adfsplit
$(214, 63, 113, 202)$,
$(214, 15, 71, 70)$,
$(215, 201, 182, 195)$,
$(215, 59, 16, 157)$,\adfsplit
$(215, 82, 36, 143)$,
$(216, 33, 126, 71)$,
$(216, 148, 83, 163)$,
$(216, 113, 160, 39)$,\adfsplit
$(217, 37, 31, 0)$,
$(217, 195, 129, 86)$,
$(217, 205, 155, 26)$,
$(218, 52, 4, 189)$,\adfsplit
$(218, 146, 179, 174)$,
$(218, 46, 185, 51)$,
$(219, 111, 160, 5)$,
$(219, 190, 45, 74)$,\adfsplit
$(219, 139, 206, 195)$,
$(220, 1, 47, 15)$,
$(220, 179, 25, 130)$,
$(220, 90, 66, 203)$,\adfsplit
$(221, 160, 68, 109)$,
$(221, 94, 17, 173)$,
$(221, 192, 87, 27)$,
$(222, 167, 18, 173)$,\adfsplit
$(222, 154, 142, 78)$,
$(222, 8, 111, 193)$,
$(223, 192, 112, 131)$,
$(223, 114, 36, 29)$,\adfsplit
$(223, 10, 8, 142)$,
$(224, 55, 125, 184)$,
$(224, 173, 165, 198)$,
$(224, 142, 42, 131)$,\adfsplit
$(225, 12, 123, 136)$,
$(225, 133, 36, 103)$,
$(2, 5, 26, 225)$,
$(0, 52, 65, 136)$,\adfsplit
$(0, 1, 176, 191)$,
$(0, 3, 159, 169)$,
$(0, 15, 35, 130)$,
$(0, 2, 123, 167)$,\adfsplit
$(0, 28, 30, 250)$,
$(0, 23, 175, 227)$,
$(0, 4, 21, 252)$,
$(0, 12, 131, 253)$,\adfsplit
$(0, 143, 199, 230)$,
$(0, 87, 203, 251)$,
$(0, 69, 172, 254)$,
$(0, 16, 110, 132)$,\adfsplit
$(0, 109, 178, 274)$,
$(0, 40, 160, 277)$,
$(0, 47, 140, 276)$,
$(0, 14, 80, 184)$,\adfsplit
$(2, 32, 71, 226)$,
$(1, 146, 194, 277)$,
$(1, 101, 188, 278)$,
$(1, 5, 89, 115)$,\adfsplit
$(1, 26, 148, 276)$,
$(1, 8, 151, 227)$,
$(1, 4, 35, 43)$

%ADFvfyBlocksEnd
\adfLgap \noindent by the mapping:
$x \mapsto x + 3 j \adfmod{207}$ for $x < 207$,
$x \mapsto (x - 207 + 24 j \adfmod{72}) + 207$ for $207 \le x < 279$,
$x \mapsto x$ for $x \ge 279$,
$0 \le j < 69$.
\ADFvfyParStart{(281, ((83, 69, ((207, 3), (72, 24), (2, 2)))), ((23, 9), (74, 1)))} %ADFvfyParEnd
% End of 23^9 74^1
%%%%%%%%%%%%%%%%%%%%%%%%%%%%%%%%%%%%%%%%%%%%%%%%%%%%%%%%%%%%%%%%%%%%%%%%%%%%%%%%%%%%%%%%%%
%%%%%%%%%%%%%%%%%%%%%%%%%%%%%%%%%%%%%%%%%%%%%%%%%%%%%%%%%%%%%%%%%%%%%%%%%%%%%%%%%%%%%%%%%%

% Charlotte:GDD4-1-3-5-mod-6-TeX-gen-A:HITS-fun:4.10
\adfDgap
%ADFvfyBlocksStart {23,23,23,23,23,23,23,23,23,80}
\noindent{\boldmath $ 23^{9} 80^{1} $}~
With the point set $Z_{287}$ partitioned into
 residue classes modulo $9$ for $\{0, 1, \dots, 206\}$, and
 $\{207, 208, \dots, 286\}$,
 the design is generated from

\adfLgap %ADFvfyDesignStart
$(285, 30, 146, 31)$,
$(286, 171, 70, 74)$,
$(207, 96, 166, 95)$,
$(207, 143, 18, 181)$,\adfsplit
$(207, 20, 142, 12)$,
$(208, 76, 100, 25)$,
$(208, 203, 21, 89)$,
$(208, 200, 0, 123)$,\adfsplit
$(209, 97, 170, 201)$,
$(209, 148, 69, 86)$,
$(209, 10, 20, 171)$,
$(210, 79, 147, 23)$,\adfsplit
$(210, 27, 127, 22)$,
$(210, 33, 179, 146)$,
$(211, 19, 110, 75)$,
$(211, 131, 76, 8)$,\adfsplit
$(211, 178, 51, 171)$,
$(212, 54, 106, 107)$,
$(212, 75, 11, 168)$,
$(212, 94, 28, 176)$,\adfsplit
$(213, 106, 9, 184)$,
$(213, 177, 19, 206)$,
$(213, 38, 185, 201)$,
$(214, 179, 121, 90)$,\adfsplit
$(214, 169, 21, 190)$,
$(214, 191, 195, 50)$,
$(215, 89, 174, 199)$,
$(215, 85, 18, 20)$,\adfsplit
$(215, 86, 15, 205)$,
$(216, 2, 197, 69)$,
$(216, 124, 111, 172)$,
$(216, 104, 76, 117)$,\adfsplit
$(217, 189, 155, 199)$,
$(217, 48, 52, 195)$,
$(217, 98, 149, 130)$,
$(218, 142, 21, 145)$,\adfsplit
$(218, 29, 194, 179)$,
$(218, 132, 135, 85)$,
$(219, 6, 139, 100)$,
$(219, 137, 167, 151)$,\adfsplit
$(219, 27, 57, 206)$,
$(220, 113, 123, 8)$,
$(220, 21, 76, 133)$,
$(220, 173, 108, 181)$,\adfsplit
$(221, 46, 164, 53)$,
$(221, 194, 162, 205)$,
$(221, 49, 33, 12)$,
$(222, 35, 74, 183)$,\adfsplit
$(222, 159, 133, 185)$,
$(222, 148, 90, 55)$,
$(223, 185, 79, 144)$,
$(223, 191, 82, 53)$,\adfsplit
$(223, 102, 159, 130)$,
$(224, 171, 151, 157)$,
$(224, 69, 107, 32)$,
$(224, 191, 154, 39)$,\adfsplit
$(225, 104, 91, 99)$,
$(225, 58, 89, 43)$,
$(225, 93, 168, 47)$,
$(226, 70, 119, 197)$,\adfsplit
$(226, 162, 202, 57)$,
$(0, 22, 134, 252)$,
$(0, 6, 20, 39)$,
$(0, 15, 34, 156)$,\adfsplit
$(0, 42, 111, 227)$,
$(0, 23, 183, 254)$,
$(0, 48, 136, 255)$,
$(0, 158, 164, 230)$,\adfsplit
$(0, 78, 185, 231)$,
$(0, 12, 131, 282)$,
$(1, 13, 97, 230)$,
$(0, 95, 118, 280)$,\adfsplit
$(0, 80, 82, 104)$,
$(0, 91, 167, 283)$,
$(0, 74, 109, 232)$,
$(0, 101, 205, 258)$,\adfsplit
$(0, 11, 154, 281)$,
$(0, 137, 184, 284)$,
$(0, 76, 155, 196)$,
$(1, 71, 191, 229)$,\adfsplit
$(1, 26, 166, 253)$,
$(1, 35, 148, 254)$,
$(1, 34, 155, 158)$,
$(1, 31, 131, 257)$,\adfsplit
$(1, 41, 62, 279)$,
$(1, 44, 70, 203)$

%ADFvfyBlocksEnd
\adfLgap \noindent by the mapping:
$x \mapsto x + 3 j \adfmod{207}$ for $x < 207$,
$x \mapsto (x - 207 + 26 j \adfmod{78}) + 207$ for $207 \le x < 285$,
$x \mapsto x$ for $x \ge 285$,
$0 \le j < 69$.
\ADFvfyParStart{(287, ((86, 69, ((207, 3), (78, 26), (2, 2)))), ((23, 9), (80, 1)))} %ADFvfyParEnd
% End of 23^9 80^1
%%%%%%%%%%%%%%%%%%%%%%%%%%%%%%%%%%%%%%%%%%%%%%%%%%%%%%%%%%%%%%%%%%%%%%%%%%%%%%%%%%%%%%%%%%
%%%%%%%%%%%%%%%%%%%%%%%%%%%%%%%%%%%%%%%%%%%%%%%%%%%%%%%%%%%%%%%%%%%%%%%%%%%%%%%%%%%%%%%%%%

% Charlotte:GDD4-1-3-5-mod-6-TeX-gen-A:HITS-fun:4.10
\adfDgap
%ADFvfyBlocksStart {23,23,23,23,23,23,23,23,23,86}
\noindent{\boldmath $ 23^{9} 86^{1} $}~
With the point set $Z_{293}$ partitioned into
 residue classes modulo $9$ for $\{0, 1, \dots, 206\}$, and
 $\{207, 208, \dots, 292\}$,
 the design is generated from

\adfLgap %ADFvfyDesignStart
$(291, 45, 199, 200)$,
$(292, 28, 95, 45)$,
$(207, 70, 176, 172)$,
$(207, 170, 12, 42)$,\adfsplit
$(207, 20, 148, 54)$,
$(208, 28, 88, 114)$,
$(208, 171, 179, 158)$,
$(208, 121, 2, 183)$,\adfsplit
$(209, 65, 95, 22)$,
$(209, 1, 69, 134)$,
$(209, 117, 7, 57)$,
$(210, 180, 22, 192)$,\adfsplit
$(210, 168, 191, 109)$,
$(210, 17, 52, 77)$,
$(211, 120, 113, 163)$,
$(211, 180, 52, 103)$,\adfsplit
$(211, 60, 164, 44)$,
$(212, 139, 47, 16)$,
$(212, 24, 183, 189)$,
$(212, 116, 100, 23)$,\adfsplit
$(213, 50, 55, 180)$,
$(213, 62, 173, 33)$,
$(213, 102, 103, 79)$,
$(214, 83, 100, 142)$,\adfsplit
$(214, 129, 53, 14)$,
$(214, 58, 198, 114)$,
$(215, 136, 12, 144)$,
$(215, 85, 106, 96)$,\adfsplit
$(215, 38, 194, 44)$,
$(216, 204, 73, 95)$,
$(216, 90, 47, 62)$,
$(216, 30, 58, 133)$,\adfsplit
$(217, 199, 130, 176)$,
$(217, 188, 105, 119)$,
$(217, 178, 9, 183)$,
$(218, 7, 173, 21)$,\adfsplit
$(218, 114, 171, 32)$,
$(218, 31, 28, 152)$,
$(219, 68, 114, 40)$,
$(219, 53, 147, 20)$,\adfsplit
$(219, 55, 205, 198)$,
$(220, 84, 125, 123)$,
$(220, 0, 86, 110)$,
$(220, 154, 184, 88)$,\adfsplit
$(221, 59, 150, 70)$,
$(221, 183, 45, 98)$,
$(221, 182, 13, 1)$,
$(222, 42, 53, 40)$,\adfsplit
$(222, 21, 45, 77)$,
$(222, 29, 37, 151)$,
$(223, 179, 160, 108)$,
$(223, 69, 100, 120)$,\adfsplit
$(223, 103, 113, 101)$,
$(224, 80, 150, 16)$,
$(224, 128, 129, 31)$,
$(224, 81, 50, 10)$,\adfsplit
$(225, 36, 129, 7)$,
$(225, 121, 33, 176)$,
$(225, 64, 155, 71)$,
$(226, 150, 77, 190)$,\adfsplit
$(226, 129, 107, 144)$,
$(226, 7, 40, 182)$,
$(227, 34, 9, 104)$,
$(227, 110, 76, 152)$,\adfsplit
$(0, 3, 16, 283)$,
$(0, 78, 167, 228)$,
$(0, 4, 19, 284)$,
$(1, 38, 179, 228)$,\adfsplit
$(0, 17, 106, 112)$,
$(0, 34, 141, 188)$,
$(0, 22, 61, 229)$,
$(0, 70, 105, 230)$,\adfsplit
$(0, 21, 59, 259)$,
$(0, 44, 142, 285)$,
$(0, 119, 197, 257)$,
$(0, 58, 203, 232)$,\adfsplit
$(0, 5, 111, 288)$,
$(1, 50, 121, 260)$,
$(0, 146, 160, 233)$,
$(0, 74, 118, 289)$,\adfsplit
$(0, 20, 166, 261)$,
$(0, 46, 175, 287)$,
$(1, 59, 161, 231)$,
$(0, 62, 120, 262)$,\adfsplit
$(0, 35, 91, 286)$,
$(1, 140, 188, 230)$,
$(1, 101, 104, 262)$,
$(0, 115, 163, 290)$,\adfsplit
$(0, 55, 107, 182)$

%ADFvfyBlocksEnd
\adfLgap \noindent by the mapping:
$x \mapsto x + 3 j \adfmod{207}$ for $x < 207$,
$x \mapsto (x - 207 + 28 j \adfmod{84}) + 207$ for $207 \le x < 291$,
$x \mapsto x$ for $x \ge 291$,
$0 \le j < 69$.
\ADFvfyParStart{(293, ((89, 69, ((207, 3), (84, 28), (2, 2)))), ((23, 9), (86, 1)))} %ADFvfyParEnd
% End of 23^9 86^1
%%%%%%%%%%%%%%%%%%%%%%%%%%%%%%%%%%%%%%%%%%%%%%%%%%%%%%%%%%%%%%%%%%%%%%%%%%%%%%%%%%%%%%%%%%
%%%%%%%%%%%%%%%%%%%%%%%%%%%%%%%%%%%%%%%%%%%%%%%%%%%%%%%%%%%%%%%%%%%%%%%%%%%%%%%%%%%%%%%%%%

%%%%%%%%%%%%%%%%%%%%%%%%%%%%%%%%%%%%%%%%%%%%%%%%%%%%%%%%%%%%%%%%%%%%%%%%%%%%%%%%%%%%%%%%%%
%%%%%%%%%%%%%%%%%%%%%%%%%%%%%%%%%%%%%%%%%%%%%%%%%%%%%%%%%%%%%%%%%%%%%%%%%%%%%%%%%%%%%%%%%%
\section{4-GDDs for the proof of Lemma \ref{lem:4-GDD 25^u m^1}}
\label{app:4-GDD 25^u m^1}
\adfnull{
$ 25^{12} 7^1 $,
$ 25^{12} 13^1 $,
$ 25^{12} 16^1 $,
$ 25^{12} 19^1 $,
$ 25^{12} 22^1 $,
$ 25^{15} 13^1 $,
$ 25^{15} 19^1 $,
$ 25^{27} 19^1 $,
$ 25^9 16^1 $,
$ 25^9 22^1 $,
$ 25^9 28^1 $,
$ 25^9 34^1 $,
$ 25^9 46^1 $,
$ 25^9 52^1 $,
$ 25^9 58^1 $,
$ 25^9 64^1 $,
$ 25^9 76^1 $,
$ 25^9 82^1 $,
$ 25^9 88^1 $,
$ 25^9 94^1 $,
$ 25^{21} 16^1 $ and
$ 25^{21} 22^1 $.
}

% Charlotte:GDD4-1-3-5-mod-6-TeX-gen-A:HITS-fun:4.10
\adfDgap
%ADFvfyBlocksStart {25,25,25,25,25,25,25,25,25,25,25,25,7}
\noindent{\boldmath $ 25^{12} 7^{1} $}~
With the point set $Z_{307}$ partitioned into
 residue classes modulo $12$ for $\{0, 1, \dots, 299\}$, and
 $\{300, 301, \dots, 306\}$,
 the design is generated from

\adfLgap %ADFvfyDesignStart
$(300, 262, 271, 260)$,
$(300, 109, 60, 174)$,
$(300, 207, 170, 11)$,
$(300, 41, 201, 244)$,\adfsplit
$(292, 66, 165, 35)$,
$(186, 149, 238, 61)$,
$(110, 89, 264, 202)$,
$(284, 135, 187, 198)$,\adfsplit
$(140, 174, 235, 117)$,
$(265, 178, 188, 209)$,
$(182, 103, 27, 169)$,
$(121, 151, 15, 86)$,\adfsplit
$(140, 211, 193, 252)$,
$(26, 80, 211, 297)$,
$(164, 214, 159, 240)$,
$(153, 6, 121, 74)$,\adfsplit
$(287, 220, 178, 122)$,
$(172, 108, 123, 125)$,
$(174, 267, 275, 229)$,
$(223, 196, 14, 287)$,\adfsplit
$(12, 251, 146, 267)$,
$(150, 8, 97, 191)$,
$(277, 208, 296, 69)$,
$(107, 272, 254, 161)$,\adfsplit
$(272, 295, 273, 276)$,
$(21, 183, 14, 209)$,
$(149, 207, 251, 258)$,
$(148, 7, 153, 181)$,\adfsplit
$(211, 40, 296, 177)$,
$(126, 256, 166, 248)$,
$(0, 3, 45, 85)$,
$(0, 6, 113, 291)$,\adfsplit
$(0, 81, 205, 267)$,
$(0, 103, 193, 213)$,
$(0, 261, 265, 275)$,
$(0, 51, 129, 199)$,\adfsplit
$(0, 117, 167, 283)$,
$(0, 20, 83, 163)$,
$(0, 22, 133, 201)$,
$(0, 78, 227, 233)$,\adfsplit
$(0, 59, 157, 231)$,
$(0, 25, 128, 198)$,
$(0, 28, 57, 58)$,
$(0, 13, 136, 174)$,\adfsplit
$(0, 16, 110, 176)$,
$(0, 39, 116, 220)$,
$(0, 14, 46, 152)$,
$(0, 75, 150, 225)$,\adfsplit
$(306, 0, 100, 200)$,
$(306, 1, 101, 201)$

%ADFvfyBlocksEnd
\adfLgap \noindent by the mapping:
$x \mapsto x + 2 j \adfmod{300}$ for $x < 300$,
$x \mapsto (x +  j \adfmod{6}) + 300$ for $300 \le x < 306$,
$306 \mapsto 306$,
$0 \le j < 150$
 for the first 47 blocks,
$0 \le j < 75$
 for the next block,
$0 \le j < 50$
 for the last two blocks.
\ADFvfyParStart{(307, ((47, 150, ((300, 2), (6, 1), (1, 1))), (1, 75, ((300, 2), (6, 1), (1, 1))), (2, 50, ((300, 2), (6, 1), (1, 1)))), ((25, 12), (7, 1)))} %ADFvfyParEnd
% End of 25^12 7^1
%%%%%%%%%%%%%%%%%%%%%%%%%%%%%%%%%%%%%%%%%%%%%%%%%%%%%%%%%%%%%%%%%%%%%%%%%%%%%%%%%%%%%%%%%%
%%%%%%%%%%%%%%%%%%%%%%%%%%%%%%%%%%%%%%%%%%%%%%%%%%%%%%%%%%%%%%%%%%%%%%%%%%%%%%%%%%%%%%%%%%

% Charlotte:GDD4-1-3-5-mod-6-TeX-gen-A:HITS-fun:4.10
\adfDgap
%ADFvfyBlocksStart {25,25,25,25,25,25,25,25,25,25,25,25,13}
\noindent{\boldmath $ 25^{12} 13^{1} $}~
With the point set $Z_{313}$ partitioned into
 residue classes modulo $12$ for $\{0, 1, \dots, 299\}$, and
 $\{300, 301, \dots, 312\}$,
 the design is generated from

\adfLgap %ADFvfyDesignStart
$(300, 243, 65, 134)$,
$(300, 299, 117, 32)$,
$(300, 211, 298, 156)$,
$(300, 280, 186, 37)$,\adfsplit
$(301, 199, 149, 204)$,
$(301, 297, 278, 46)$,
$(301, 265, 131, 51)$,
$(301, 282, 16, 8)$,\adfsplit
$(52, 143, 213, 226)$,
$(23, 112, 150, 231)$,
$(172, 6, 12, 110)$,
$(161, 144, 14, 237)$,\adfsplit
$(1, 197, 234, 285)$,
$(280, 189, 210, 151)$,
$(158, 100, 169, 68)$,
$(43, 85, 261, 280)$,\adfsplit
$(7, 10, 15, 186)$,
$(114, 15, 289, 127)$,
$(73, 262, 170, 231)$,
$(249, 13, 84, 56)$,\adfsplit
$(42, 79, 158, 57)$,
$(133, 126, 20, 238)$,
$(40, 252, 275, 29)$,
$(148, 10, 97, 281)$,\adfsplit
$(4, 39, 286, 8)$,
$(81, 194, 171, 142)$,
$(41, 297, 14, 267)$,
$(24, 25, 110, 119)$,\adfsplit
$(62, 185, 248, 112)$,
$(218, 179, 138, 72)$,
$(239, 20, 205, 91)$,
$(124, 170, 60, 15)$,\adfsplit
$(0, 3, 122, 224)$,
$(0, 10, 30, 182)$,
$(0, 21, 89, 222)$,
$(0, 247, 265, 293)$,\adfsplit
$(0, 44, 227, 233)$,
$(0, 25, 259, 273)$,
$(0, 65, 143, 163)$,
$(0, 14, 131, 141)$,\adfsplit
$(0, 16, 59, 169)$,
$(0, 49, 74, 177)$,
$(0, 42, 125, 157)$,
$(0, 149, 205, 207)$,\adfsplit
$(0, 2, 47, 56)$,
$(0, 155, 159, 285)$,
$(0, 97, 137, 243)$,
$(0, 39, 40, 179)$,\adfsplit
$(0, 33, 135, 197)$,
$(0, 75, 150, 225)$,
$(312, 0, 100, 200)$,
$(312, 1, 101, 201)$

%ADFvfyBlocksEnd
\adfLgap \noindent by the mapping:
$x \mapsto x + 2 j \adfmod{300}$ for $x < 300$,
$x \mapsto (x + 2 j \adfmod{12}) + 300$ for $300 \le x < 312$,
$312 \mapsto 312$,
$0 \le j < 150$
 for the first 49 blocks,
$0 \le j < 75$
 for the next block,
$0 \le j < 50$
 for the last two blocks.
\ADFvfyParStart{(313, ((49, 150, ((300, 2), (12, 2), (1, 1))), (1, 75, ((300, 2), (12, 2), (1, 1))), (2, 50, ((300, 2), (12, 2), (1, 1)))), ((25, 12), (13, 1)))} %ADFvfyParEnd
% End of 25^12 13^1
%%%%%%%%%%%%%%%%%%%%%%%%%%%%%%%%%%%%%%%%%%%%%%%%%%%%%%%%%%%%%%%%%%%%%%%%%%%%%%%%%%%%%%%%%%
%%%%%%%%%%%%%%%%%%%%%%%%%%%%%%%%%%%%%%%%%%%%%%%%%%%%%%%%%%%%%%%%%%%%%%%%%%%%%%%%%%%%%%%%%%

% Charlotte:GDD4-1-3-5-mod-6-TeX-gen-A:HITS-fun:4.10
\adfDgap
%ADFvfyBlocksStart {25,25,25,25,25,25,25,25,25,25,25,25,16}
\noindent{\boldmath $ 25^{12} 16^{1} $}~
With the point set $Z_{316}$ partitioned into
 residue classes modulo $12$ for $\{0, 1, \dots, 299\}$, and
 $\{300, 301, \dots, 315\}$,
 the design is generated from

\adfLgap %ADFvfyDesignStart
$(300, 273, 259, 122)$,
$(301, 80, 142, 108)$,
$(302, 260, 84, 229)$,
$(303, 117, 215, 298)$,\adfsplit
$(304, 226, 239, 12)$,
$(38, 45, 204, 97)$,
$(130, 172, 246, 252)$,
$(19, 121, 238, 80)$,\adfsplit
$(1, 128, 107, 206)$,
$(223, 228, 76, 98)$,
$(227, 198, 242, 231)$,
$(294, 60, 148, 125)$,\adfsplit
$(247, 41, 138, 256)$,
$(0, 1, 3, 113)$,
$(0, 8, 27, 43)$,
$(0, 10, 30, 55)$,\adfsplit
$(0, 32, 71, 209)$,
$(0, 17, 54, 253)$,
$(0, 51, 121, 208)$,
$(0, 18, 58, 111)$,\adfsplit
$(0, 38, 105, 174)$,
$(0, 49, 128, 210)$,
$(0, 26, 76, 211)$,
$(0, 56, 133, 196)$,\adfsplit
$(0, 46, 103, 232)$,
$(0, 75, 150, 225)$,
$(315, 0, 100, 200)$

%ADFvfyBlocksEnd
\adfLgap \noindent by the mapping:
$x \mapsto x +  j \adfmod{300}$ for $x < 300$,
$x \mapsto (x + 5 j \adfmod{15}) + 300$ for $300 \le x < 315$,
$315 \mapsto 315$,
$0 \le j < 300$
 for the first 25 blocks,
$0 \le j < 75$
 for the next block,
$0 \le j < 100$
 for the last block.
\ADFvfyParStart{(316, ((25, 300, ((300, 1), (15, 5), (1, 1))), (1, 75, ((300, 1), (15, 5), (1, 1))), (1, 100, ((300, 1), (15, 5), (1, 1)))), ((25, 12), (16, 1)))} %ADFvfyParEnd
% End of 25^12 16^1
%%%%%%%%%%%%%%%%%%%%%%%%%%%%%%%%%%%%%%%%%%%%%%%%%%%%%%%%%%%%%%%%%%%%%%%%%%%%%%%%%%%%%%%%%%
%%%%%%%%%%%%%%%%%%%%%%%%%%%%%%%%%%%%%%%%%%%%%%%%%%%%%%%%%%%%%%%%%%%%%%%%%%%%%%%%%%%%%%%%%%

% Charlotte:GDD4-1-3-5-mod-6-TeX-gen-A:HITS-fun:4.10
\adfDgap
%ADFvfyBlocksStart {25,25,25,25,25,25,25,25,25,25,25,25,19}
\noindent{\boldmath $ 25^{12} 19^{1} $}~
With the point set $Z_{319}$ partitioned into
 residue classes modulo $12$ for $\{0, 1, \dots, 299\}$, and
 $\{300, 301, \dots, 318\}$,
 the design is generated from

\adfLgap %ADFvfyDesignStart
$(300, 176, 186, 22)$,
$(300, 35, 72, 115)$,
$(300, 14, 193, 29)$,
$(300, 244, 141, 159)$,\adfsplit
$(301, 160, 286, 213)$,
$(301, 171, 198, 272)$,
$(301, 161, 216, 277)$,
$(301, 242, 191, 43)$,\adfsplit
$(302, 185, 250, 28)$,
$(302, 51, 206, 176)$,
$(302, 108, 90, 237)$,
$(302, 131, 127, 25)$,\adfsplit
$(208, 96, 217, 90)$,
$(115, 121, 101, 18)$,
$(252, 118, 155, 89)$,
$(0, 106, 269, 67)$,\adfsplit
$(162, 228, 188, 232)$,
$(126, 45, 299, 173)$,
$(59, 113, 174, 289)$,
$(1, 44, 112, 288)$,\adfsplit
$(21, 161, 116, 18)$,
$(164, 153, 156, 198)$,
$(30, 61, 123, 11)$,
$(99, 34, 139, 246)$,\adfsplit
$(79, 170, 193, 272)$,
$(9, 113, 275, 24)$,
$(246, 295, 59, 218)$,
$(41, 254, 94, 256)$,\adfsplit
$(30, 124, 92, 259)$,
$(231, 220, 278, 240)$,
$(215, 64, 174, 57)$,
$(0, 5, 214, 231)$,\adfsplit
$(0, 1, 114, 213)$,
$(0, 7, 236, 286)$,
$(0, 19, 148, 170)$,
$(0, 27, 57, 116)$,\adfsplit
$(0, 29, 54, 277)$,
$(0, 51, 196, 259)$,
$(0, 59, 92, 220)$,
$(0, 39, 95, 158)$,\adfsplit
$(0, 46, 119, 224)$,
$(0, 16, 177, 295)$,
$(0, 81, 82, 191)$,
$(0, 35, 52, 169)$,\adfsplit
$(0, 25, 69, 287)$,
$(1, 17, 95, 171)$,
$(0, 123, 125, 133)$,
$(0, 131, 159, 217)$,\adfsplit
$(0, 165, 207, 233)$,
$(0, 33, 55, 90)$,
$(0, 79, 111, 201)$,
$(0, 75, 150, 225)$,\adfsplit
$(318, 0, 100, 200)$,
$(318, 1, 101, 201)$

%ADFvfyBlocksEnd
\adfLgap \noindent by the mapping:
$x \mapsto x + 2 j \adfmod{300}$ for $x < 300$,
$x \mapsto (x - 300 + 3 j \adfmod{18}) + 300$ for $300 \le x < 318$,
$318 \mapsto 318$,
$0 \le j < 150$
 for the first 51 blocks,
$0 \le j < 75$
 for the next block,
$0 \le j < 50$
 for the last two blocks.
\ADFvfyParStart{(319, ((51, 150, ((300, 2), (18, 3), (1, 1))), (1, 75, ((300, 2), (18, 3), (1, 1))), (2, 50, ((300, 2), (18, 3), (1, 1)))), ((25, 12), (19, 1)))} %ADFvfyParEnd
% End of 25^12 19^1
%%%%%%%%%%%%%%%%%%%%%%%%%%%%%%%%%%%%%%%%%%%%%%%%%%%%%%%%%%%%%%%%%%%%%%%%%%%%%%%%%%%%%%%%%%
%%%%%%%%%%%%%%%%%%%%%%%%%%%%%%%%%%%%%%%%%%%%%%%%%%%%%%%%%%%%%%%%%%%%%%%%%%%%%%%%%%%%%%%%%%

% Charlotte:GDD4-1-3-5-mod-6-TeX-gen-A:HITS-fun:4.10
\adfDgap
%ADFvfyBlocksStart {25,25,25,25,25,25,25,25,25,25,25,25,22}
\noindent{\boldmath $ 25^{12} 22^{1} $}~
With the point set $Z_{322}$ partitioned into
 residue classes modulo $12$ for $\{0, 1, \dots, 299\}$, and
 $\{300, 301, \dots, 321\}$,
 the design is generated from

\adfLgap %ADFvfyDesignStart
$(300, 60, 164, 142)$,
$(301, 232, 248, 111)$,
$(302, 253, 161, 126)$,
$(303, 235, 42, 122)$,\adfsplit
$(304, 228, 220, 131)$,
$(305, 245, 70, 111)$,
$(306, 249, 61, 293)$,
$(209, 11, 274, 256)$,\adfsplit
$(197, 223, 268, 0)$,
$(168, 140, 298, 285)$,
$(214, 120, 127, 158)$,
$(129, 78, 288, 290)$,\adfsplit
$(26, 180, 177, 186)$,
$(65, 148, 212, 134)$,
$(0, 1, 5, 290)$,
$(0, 17, 74, 128)$,\adfsplit
$(0, 19, 39, 258)$,
$(0, 29, 63, 186)$,
$(0, 27, 133, 191)$,
$(0, 25, 91, 215)$,\adfsplit
$(0, 23, 53, 205)$,
$(0, 40, 99, 214)$,
$(0, 43, 93, 224)$,
$(0, 33, 79, 195)$,\adfsplit
$(0, 62, 129, 227)$,
$(0, 21, 70, 122)$,
$(0, 75, 150, 225)$,
$(321, 0, 100, 200)$

%ADFvfyBlocksEnd
\adfLgap \noindent by the mapping:
$x \mapsto x +  j \adfmod{300}$ for $x < 300$,
$x \mapsto (x - 300 + 7 j \adfmod{21}) + 300$ for $300 \le x < 321$,
$321 \mapsto 321$,
$0 \le j < 300$
 for the first 26 blocks,
$0 \le j < 75$
 for the next block,
$0 \le j < 100$
 for the last block.
\ADFvfyParStart{(322, ((26, 300, ((300, 1), (21, 7), (1, 1))), (1, 75, ((300, 1), (21, 7), (1, 1))), (1, 100, ((300, 1), (21, 7), (1, 1)))), ((25, 12), (22, 1)))} %ADFvfyParEnd
% End of 25^12 22^1
%%%%%%%%%%%%%%%%%%%%%%%%%%%%%%%%%%%%%%%%%%%%%%%%%%%%%%%%%%%%%%%%%%%%%%%%%%%%%%%%%%%%%%%%%%
%%%%%%%%%%%%%%%%%%%%%%%%%%%%%%%%%%%%%%%%%%%%%%%%%%%%%%%%%%%%%%%%%%%%%%%%%%%%%%%%%%%%%%%%%%

% Charlotte:GDD4-1-3-5-mod-6-TeX-gen-A:HITS-fun:4.10
\adfDgap
%ADFvfyBlocksStart {25,25,25,25,25,25,25,25,25,25,25,25,25,25,25,13}
\noindent{\boldmath $ 25^{15} 13^{1} $}~
With the point set $Z_{388}$ partitioned into
 residue classes modulo $15$ for $\{0, 1, \dots, 374\}$, and
 $\{375, 376, \dots, 387\}$,
 the design is generated from

\adfLgap %ADFvfyDesignStart
$(375, 141, 229, 215)$,
$(376, 86, 363, 4)$,
$(377, 159, 215, 214)$,
$(378, 145, 231, 275)$,\adfsplit
$(151, 298, 52, 356)$,
$(320, 62, 143, 139)$,
$(363, 284, 297, 151)$,
$(21, 88, 324, 224)$,\adfsplit
$(83, 126, 304, 239)$,
$(272, 41, 81, 18)$,
$(71, 282, 318, 202)$,
$(255, 205, 185, 224)$,\adfsplit
$(147, 64, 112, 271)$,
$(99, 7, 186, 213)$,
$(8, 200, 166, 0)$,
$(10, 352, 364, 300)$,\adfsplit
$(194, 320, 346, 236)$,
$(123, 366, 361, 218)$,
$(78, 190, 180, 227)$,
$(0, 3, 9, 41)$,\adfsplit
$(0, 7, 96, 237)$,
$(0, 17, 118, 241)$,
$(0, 59, 153, 214)$,
$(0, 28, 97, 201)$,\adfsplit
$(0, 18, 109, 185)$,
$(0, 11, 68, 187)$,
$(0, 46, 108, 235)$,
$(0, 53, 107, 213)$,\adfsplit
$(0, 2, 24, 326)$,
$(0, 25, 103, 218)$,
$(0, 29, 122, 233)$,
$(387, 0, 125, 250)$

%ADFvfyBlocksEnd
\adfLgap \noindent by the mapping:
$x \mapsto x +  j \adfmod{375}$ for $x < 375$,
$x \mapsto (x - 375 + 4 j \adfmod{12}) + 375$ for $375 \le x < 387$,
$387 \mapsto 387$,
$0 \le j < 375$
 for the first 31 blocks,
$0 \le j < 125$
 for the last block.
\ADFvfyParStart{(388, ((31, 375, ((375, 1), (12, 4), (1, 1))), (1, 125, ((375, 1), (12, 4), (1, 1)))), ((25, 15), (13, 1)))} %ADFvfyParEnd
% End of 25^15 13^1
%%%%%%%%%%%%%%%%%%%%%%%%%%%%%%%%%%%%%%%%%%%%%%%%%%%%%%%%%%%%%%%%%%%%%%%%%%%%%%%%%%%%%%%%%%
%%%%%%%%%%%%%%%%%%%%%%%%%%%%%%%%%%%%%%%%%%%%%%%%%%%%%%%%%%%%%%%%%%%%%%%%%%%%%%%%%%%%%%%%%%

% Charlotte:GDD4-1-3-5-mod-6-TeX-gen-A:HITS-fun:4.10
\adfDgap
%ADFvfyBlocksStart {25,25,25,25,25,25,25,25,25,25,25,25,25,25,25,19}
\noindent{\boldmath $ 25^{15} 19^{1} $}~
With the point set $Z_{394}$ partitioned into
 residue classes modulo $15$ for $\{0, 1, \dots, 374\}$, and
 $\{375, 376, \dots, 393\}$,
 the design is generated from

\adfLgap %ADFvfyDesignStart
$(375, 32, 18, 34)$,
$(376, 189, 151, 362)$,
$(377, 52, 185, 357)$,
$(378, 361, 329, 282)$,\adfsplit
$(379, 152, 228, 94)$,
$(380, 170, 268, 183)$,
$(88, 284, 108, 100)$,
$(183, 217, 326, 218)$,\adfsplit
$(171, 70, 270, 372)$,
$(214, 60, 43, 188)$,
$(110, 41, 0, 322)$,
$(217, 163, 156, 106)$,\adfsplit
$(82, 319, 200, 59)$,
$(21, 165, 292, 87)$,
$(255, 74, 211, 368)$,
$(346, 324, 232, 281)$,\adfsplit
$(321, 159, 170, 199)$,
$(191, 228, 349, 246)$,
$(201, 112, 12, 356)$,
$(306, 247, 79, 118)$,\adfsplit
$(0, 3, 28, 152)$,
$(0, 5, 24, 214)$,
$(0, 6, 48, 293)$,
$(0, 4, 68, 95)$,\adfsplit
$(0, 33, 140, 192)$,
$(0, 51, 147, 259)$,
$(0, 36, 123, 229)$,
$(0, 21, 84, 177)$,\adfsplit
$(0, 10, 72, 308)$,
$(0, 56, 136, 278)$,
$(0, 46, 117, 243)$,
$(0, 9, 83, 169)$,\adfsplit
$(393, 0, 125, 250)$

%ADFvfyBlocksEnd
\adfLgap \noindent by the mapping:
$x \mapsto x +  j \adfmod{375}$ for $x < 375$,
$x \mapsto (x - 375 + 6 j \adfmod{18}) + 375$ for $375 \le x < 393$,
$393 \mapsto 393$,
$0 \le j < 375$
 for the first 32 blocks,
$0 \le j < 125$
 for the last block.
\ADFvfyParStart{(394, ((32, 375, ((375, 1), (18, 6), (1, 1))), (1, 125, ((375, 1), (18, 6), (1, 1)))), ((25, 15), (19, 1)))} %ADFvfyParEnd
% End of 25^15 19^1
%%%%%%%%%%%%%%%%%%%%%%%%%%%%%%%%%%%%%%%%%%%%%%%%%%%%%%%%%%%%%%%%%%%%%%%%%%%%%%%%%%%%%%%%%%
%%%%%%%%%%%%%%%%%%%%%%%%%%%%%%%%%%%%%%%%%%%%%%%%%%%%%%%%%%%%%%%%%%%%%%%%%%%%%%%%%%%%%%%%%%

% Charlotte:GDD4-1-3-5-mod-6-TeX-gen-A:HITS-fun:4.10
\adfDgap
%ADFvfyBlocksStart {25,25,25,25,25,25,25,25,25,25,25,25,25,25,25,25,25,25,25,25,25,25,25,25,25,25,25,19}
\noindent{\boldmath $ 25^{27} 19^{1} $}~
With the point set $Z_{694}$ partitioned into
 residue classes modulo $27$ for $\{0, 1, \dots, 674\}$, and
 $\{675, 676, \dots, 693\}$,
 the design is generated from

\adfLgap %ADFvfyDesignStart
$(675, 206, 523, 468)$,
$(676, 9, 431, 292)$,
$(677, 255, 290, 268)$,
$(678, 31, 234, 20)$,\adfsplit
$(679, 32, 177, 214)$,
$(680, 558, 88, 236)$,
$(415, 144, 606, 3)$,
$(335, 532, 656, 156)$,\adfsplit
$(589, 571, 500, 621)$,
$(661, 597, 636, 667)$,
$(103, 649, 502, 542)$,
$(521, 300, 15, 320)$,\adfsplit
$(469, 90, 474, 323)$,
$(570, 614, 119, 314)$,
$(295, 657, 598, 666)$,
$(595, 23, 517, 367)$,\adfsplit
$(18, 195, 146, 10)$,
$(656, 196, 308, 213)$,
$(634, 641, 204, 118)$,
$(344, 2, 583, 390)$,\adfsplit
$(527, 164, 59, 62)$,
$(348, 282, 59, 1)$,
$(289, 457, 47, 487)$,
$(149, 218, 92, 475)$,\adfsplit
$(10, 164, 283, 454)$,
$(268, 515, 138, 311)$,
$(176, 454, 406, 370)$,
$(150, 370, 183, 25)$,\adfsplit
$(118, 181, 25, 147)$,
$(671, 353, 297, 430)$,
$(320, 532, 200, 186)$,
$(328, 17, 420, 263)$,\adfsplit
$(46, 97, 385, 452)$,
$(467, 457, 367, 183)$,
$(189, 188, 597, 114)$,
$(41, 640, 139, 391)$,\adfsplit
$(244, 341, 199, 560)$,
$(188, 474, 216, 304)$,
$(288, 437, 449, 625)$,
$(73, 617, 526, 153)$,\adfsplit
$(126, 521, 73, 186)$,
$(497, 99, 637, 663)$,
$(0, 4, 190, 424)$,
$(0, 19, 118, 400)$,\adfsplit
$(0, 52, 217, 467)$,
$(0, 73, 160, 414)$,
$(0, 16, 218, 382)$,
$(0, 2, 117, 155)$,\adfsplit
$(0, 21, 62, 497)$,
$(0, 23, 106, 167)$,
$(0, 101, 307, 416)$,
$(0, 42, 268, 411)$,\adfsplit
$(0, 79, 183, 279)$,
$(0, 47, 132, 295)$,
$(0, 24, 138, 503)$,
$(0, 15, 244, 367)$,\adfsplit
$(0, 82, 209, 319)$,
$(693, 0, 225, 450)$

%ADFvfyBlocksEnd
\adfLgap \noindent by the mapping:
$x \mapsto x +  j \adfmod{675}$ for $x < 675$,
$x \mapsto (x - 675 + 6 j \adfmod{18}) + 675$ for $675 \le x < 693$,
$693 \mapsto 693$,
$0 \le j < 675$
 for the first 57 blocks,
$0 \le j < 225$
 for the last block.
\ADFvfyParStart{(694, ((57, 675, ((675, 1), (18, 6), (1, 1))), (1, 225, ((675, 1), (18, 6), (1, 1)))), ((25, 27), (19, 1)))} %ADFvfyParEnd
% End of 25^27 19^1
%%%%%%%%%%%%%%%%%%%%%%%%%%%%%%%%%%%%%%%%%%%%%%%%%%%%%%%%%%%%%%%%%%%%%%%%%%%%%%%%%%%%%%%%%%
%%%%%%%%%%%%%%%%%%%%%%%%%%%%%%%%%%%%%%%%%%%%%%%%%%%%%%%%%%%%%%%%%%%%%%%%%%%%%%%%%%%%%%%%%%

% Charlotte:GDD4-1-3-5-mod-6-TeX-gen-A:HITS-fun:4.10
\adfDgap
%ADFvfyBlocksStart {25,25,25,25,25,25,25,25,25,16}
\noindent{\boldmath $ 25^{9} 16^{1} $}~
With the point set $Z_{241}$ partitioned into
 residue classes modulo $9$ for $\{0, 1, \dots, 224\}$, and
 $\{225, 226, \dots, 240\}$,
 the design is generated from

\adfLgap %ADFvfyDesignStart
$(225, 166, 102, 35)$,
$(226, 219, 197, 52)$,
$(227, 96, 23, 202)$,
$(228, 192, 50, 163)$,\adfsplit
$(229, 148, 143, 9)$,
$(55, 153, 5, 173)$,
$(147, 180, 124, 163)$,
$(129, 77, 91, 36)$,\adfsplit
$(73, 120, 162, 169)$,
$(28, 49, 52, 123)$,
$(0, 1, 31, 35)$,
$(0, 2, 12, 104)$,\adfsplit
$(0, 8, 48, 163)$,
$(0, 26, 79, 164)$,
$(0, 11, 43, 111)$,
$(0, 6, 66, 103)$,\adfsplit
$(0, 13, 28, 137)$,
$(0, 19, 78, 160)$,
$(0, 25, 69, 174)$,
$(240, 0, 75, 150)$

%ADFvfyBlocksEnd
\adfLgap \noindent by the mapping:
$x \mapsto x +  j \adfmod{225}$ for $x < 225$,
$x \mapsto (x + 5 j \adfmod{15}) + 225$ for $225 \le x < 240$,
$240 \mapsto 240$,
$0 \le j < 225$
 for the first 19 blocks,
$0 \le j < 75$
 for the last block.
\ADFvfyParStart{(241, ((19, 225, ((225, 1), (15, 5), (1, 1))), (1, 75, ((225, 1), (15, 5), (1, 1)))), ((25, 9), (16, 1)))} %ADFvfyParEnd
% End of 25^9 16^1
%%%%%%%%%%%%%%%%%%%%%%%%%%%%%%%%%%%%%%%%%%%%%%%%%%%%%%%%%%%%%%%%%%%%%%%%%%%%%%%%%%%%%%%%%%
%%%%%%%%%%%%%%%%%%%%%%%%%%%%%%%%%%%%%%%%%%%%%%%%%%%%%%%%%%%%%%%%%%%%%%%%%%%%%%%%%%%%%%%%%%

% Charlotte:GDD4-1-3-5-mod-6-TeX-gen-A:HITS-fun:4.10
\adfDgap
%ADFvfyBlocksStart {25,25,25,25,25,25,25,25,25,22}
\noindent{\boldmath $ 25^{9} 22^{1} $}~
With the point set $Z_{247}$ partitioned into
 residue classes modulo $9$ for $\{0, 1, \dots, 224\}$, and
 $\{225, 226, \dots, 246\}$,
 the design is generated from

\adfLgap %ADFvfyDesignStart
$(225, 206, 172, 153)$,
$(226, 190, 212, 129)$,
$(227, 67, 11, 27)$,
$(228, 74, 214, 6)$,\adfsplit
$(229, 180, 217, 137)$,
$(230, 55, 203, 9)$,
$(231, 96, 122, 52)$,
$(149, 15, 62, 117)$,\adfsplit
$(112, 71, 12, 41)$,
$(176, 67, 69, 44)$,
$(0, 1, 6, 39)$,
$(0, 3, 10, 130)$,\adfsplit
$(0, 4, 52, 119)$,
$(0, 14, 35, 111)$,
$(0, 28, 78, 129)$,
$(0, 8, 92, 112)$,\adfsplit
$(0, 24, 66, 160)$,
$(0, 12, 69, 151)$,
$(0, 15, 79, 137)$,
$(0, 11, 60, 73)$,\adfsplit
$(246, 0, 75, 150)$

%ADFvfyBlocksEnd
\adfLgap \noindent by the mapping:
$x \mapsto x +  j \adfmod{225}$ for $x < 225$,
$x \mapsto (x - 225 + 7 j \adfmod{21}) + 225$ for $225 \le x < 246$,
$246 \mapsto 246$,
$0 \le j < 225$
 for the first 20 blocks,
$0 \le j < 75$
 for the last block.
\ADFvfyParStart{(247, ((20, 225, ((225, 1), (21, 7), (1, 1))), (1, 75, ((225, 1), (21, 7), (1, 1)))), ((25, 9), (22, 1)))} %ADFvfyParEnd
% End of 25^9 22^1
%%%%%%%%%%%%%%%%%%%%%%%%%%%%%%%%%%%%%%%%%%%%%%%%%%%%%%%%%%%%%%%%%%%%%%%%%%%%%%%%%%%%%%%%%%
%%%%%%%%%%%%%%%%%%%%%%%%%%%%%%%%%%%%%%%%%%%%%%%%%%%%%%%%%%%%%%%%%%%%%%%%%%%%%%%%%%%%%%%%%%

% Charlotte:GDD4-1-3-5-mod-6-TeX-gen-A:HITS-fun:4.10
\adfDgap
%ADFvfyBlocksStart {25,25,25,25,25,25,25,25,25,28}
\noindent{\boldmath $ 25^{9} 28^{1} $}~
With the point set $Z_{253}$ partitioned into
 residue classes modulo $9$ for $\{0, 1, \dots, 224\}$, and
 $\{225, 226, \dots, 252\}$,
 the design is generated from

\adfLgap %ADFvfyDesignStart
$(225, 82, 27, 89)$,
$(226, 222, 154, 224)$,
$(227, 34, 9, 38)$,
$(228, 47, 10, 186)$,\adfsplit
$(229, 66, 2, 184)$,
$(230, 130, 191, 42)$,
$(231, 112, 86, 117)$,
$(232, 175, 35, 117)$,\adfsplit
$(233, 194, 15, 193)$,
$(78, 11, 184, 203)$,
$(0, 3, 11, 141)$,
$(0, 6, 16, 129)$,\adfsplit
$(0, 12, 92, 105)$,
$(0, 17, 57, 77)$,
$(0, 14, 35, 191)$,
$(0, 23, 51, 134)$,\adfsplit
$(0, 22, 101, 131)$,
$(0, 15, 53, 181)$,
$(0, 32, 73, 147)$,
$(0, 24, 66, 122)$,\adfsplit
$(0, 39, 89, 160)$,
$(252, 0, 75, 150)$

%ADFvfyBlocksEnd
\adfLgap \noindent by the mapping:
$x \mapsto x +  j \adfmod{225}$ for $x < 225$,
$x \mapsto (x - 225 + 9 j \adfmod{27}) + 225$ for $225 \le x < 252$,
$252 \mapsto 252$,
$0 \le j < 225$
 for the first 21 blocks,
$0 \le j < 75$
 for the last block.
\ADFvfyParStart{(253, ((21, 225, ((225, 1), (27, 9), (1, 1))), (1, 75, ((225, 1), (27, 9), (1, 1)))), ((25, 9), (28, 1)))} %ADFvfyParEnd
% End of 25^9 28^1
%%%%%%%%%%%%%%%%%%%%%%%%%%%%%%%%%%%%%%%%%%%%%%%%%%%%%%%%%%%%%%%%%%%%%%%%%%%%%%%%%%%%%%%%%%
%%%%%%%%%%%%%%%%%%%%%%%%%%%%%%%%%%%%%%%%%%%%%%%%%%%%%%%%%%%%%%%%%%%%%%%%%%%%%%%%%%%%%%%%%%

% Charlotte:GDD4-1-3-5-mod-6-TeX-gen-A:HITS-fun:4.10
\adfDgap
%ADFvfyBlocksStart {25,25,25,25,25,25,25,25,25,34}
\noindent{\boldmath $ 25^{9} 34^{1} $}~
With the point set $Z_{259}$ partitioned into
 residue classes modulo $9$ for $\{0, 1, \dots, 224\}$, and
 $\{225, 226, \dots, 258\}$,
 the design is generated from

\adfLgap %ADFvfyDesignStart
$(225, 155, 139, 186)$,
$(226, 115, 189, 71)$,
$(227, 154, 177, 20)$,
$(228, 57, 146, 85)$,\adfsplit
$(229, 156, 28, 161)$,
$(230, 86, 145, 33)$,
$(231, 168, 125, 112)$,
$(232, 48, 83, 46)$,\adfsplit
$(233, 65, 171, 211)$,
$(234, 175, 60, 197)$,
$(0, 1, 8, 235)$,
$(0, 3, 32, 42)$,\adfsplit
$(0, 4, 19, 52)$,
$(0, 6, 17, 102)$,
$(0, 12, 77, 159)$,
$(0, 30, 64, 168)$,\adfsplit
$(0, 21, 71, 122)$,
$(0, 26, 93, 142)$,
$(0, 20, 58, 120)$,
$(0, 25, 98, 139)$,\adfsplit
$(0, 14, 69, 145)$,
$(0, 24, 70, 165)$,
$(258, 0, 75, 150)$

%ADFvfyBlocksEnd
\adfLgap \noindent by the mapping:
$x \mapsto x +  j \adfmod{225}$ for $x < 225$,
$x \mapsto (x - 225 + 11 j \adfmod{33}) + 225$ for $225 \le x < 258$,
$258 \mapsto 258$,
$0 \le j < 225$
 for the first 22 blocks,
$0 \le j < 75$
 for the last block.
\ADFvfyParStart{(259, ((22, 225, ((225, 1), (33, 11), (1, 1))), (1, 75, ((225, 1), (33, 11), (1, 1)))), ((25, 9), (34, 1)))} %ADFvfyParEnd
% End of 25^9 34^1
%%%%%%%%%%%%%%%%%%%%%%%%%%%%%%%%%%%%%%%%%%%%%%%%%%%%%%%%%%%%%%%%%%%%%%%%%%%%%%%%%%%%%%%%%%
%%%%%%%%%%%%%%%%%%%%%%%%%%%%%%%%%%%%%%%%%%%%%%%%%%%%%%%%%%%%%%%%%%%%%%%%%%%%%%%%%%%%%%%%%%

% Charlotte:GDD4-1-3-5-mod-6-TeX-gen-A:HITS-fun:4.10
\adfDgap
%ADFvfyBlocksStart {25,25,25,25,25,25,25,25,25,46}
\noindent{\boldmath $ 25^{9} 46^{1} $}~
With the point set $Z_{271}$ partitioned into
 residue classes modulo $9$ for $\{0, 1, \dots, 224\}$, and
 $\{225, 226, \dots, 270\}$,
 the design is generated from

\adfLgap %ADFvfyDesignStart
$(225, 40, 50, 161)$,
$(225, 7, 53, 219)$,
$(225, 190, 4, 179)$,
$(225, 48, 11, 194)$,\adfsplit
$(225, 45, 137, 125)$,
$(225, 205, 136, 209)$,
$(225, 60, 82, 83)$,
$(225, 117, 69, 112)$,\adfsplit
$(225, 151, 167, 20)$,
$(225, 23, 165, 193)$,
$(225, 132, 216, 79)$,
$(225, 154, 41, 102)$,\adfsplit
$(0, 2, 31, 239)$,
$(0, 3, 33, 252)$,
$(0, 15, 206, 242)$,
$(0, 6, 20, 106)$,\adfsplit
$(0, 8, 101, 118)$,
$(0, 21, 68, 159)$,
$(0, 26, 97, 148)$,
$(0, 32, 70, 127)$,\adfsplit
$(0, 40, 96, 163)$,
$(0, 7, 65, 89)$,
$(0, 25, 60, 176)$,
$(0, 41, 85, 161)$,\adfsplit
$(270, 0, 75, 150)$

%ADFvfyBlocksEnd
\adfLgap \noindent by the mapping:
$x \mapsto x +  j \adfmod{225}$ for $x < 225$,
$x \mapsto (x +  j \adfmod{45}) + 225$ for $225 \le x < 270$,
$270 \mapsto 270$,
$0 \le j < 225$
 for the first 24 blocks,
$0 \le j < 75$
 for the last block.
\ADFvfyParStart{(271, ((24, 225, ((225, 1), (45, 1), (1, 1))), (1, 75, ((225, 1), (45, 1), (1, 1)))), ((25, 9), (46, 1)))} %ADFvfyParEnd
% End of 25^9 46^1
%%%%%%%%%%%%%%%%%%%%%%%%%%%%%%%%%%%%%%%%%%%%%%%%%%%%%%%%%%%%%%%%%%%%%%%%%%%%%%%%%%%%%%%%%%
%%%%%%%%%%%%%%%%%%%%%%%%%%%%%%%%%%%%%%%%%%%%%%%%%%%%%%%%%%%%%%%%%%%%%%%%%%%%%%%%%%%%%%%%%%

% Charlotte:GDD4-1-3-5-mod-6-TeX-gen-A:HITS-fun:4.10
\adfDgap
%ADFvfyBlocksStart {25,25,25,25,25,25,25,25,25,52}
\noindent{\boldmath $ 25^{9} 52^{1} $}~
With the point set $Z_{277}$ partitioned into
 residue classes modulo $9$ for $\{0, 1, \dots, 224\}$, and
 $\{225, 226, \dots, 276\}$,
 the design is generated from

\adfLgap %ADFvfyDesignStart
$(270, 25, 170, 63)$,
$(271, 147, 112, 197)$,
$(225, 13, 26, 150)$,
$(225, 88, 33, 111)$,\adfsplit
$(225, 98, 74, 204)$,
$(225, 121, 219, 70)$,
$(225, 199, 89, 32)$,
$(225, 125, 140, 145)$,\adfsplit
$(225, 182, 99, 22)$,
$(225, 51, 48, 217)$,
$(225, 147, 4, 73)$,
$(225, 27, 79, 220)$,\adfsplit
$(225, 45, 203, 177)$,
$(0, 1, 204, 253)$,
$(0, 2, 42, 256)$,
$(0, 6, 215, 240)$,\adfsplit
$(0, 8, 194, 232)$,
$(0, 11, 28, 133)$,
$(0, 7, 71, 104)$,
$(0, 19, 87, 131)$,\adfsplit
$(0, 12, 46, 60)$,
$(0, 25, 66, 136)$,
$(0, 29, 91, 152)$,
$(0, 30, 79, 116)$,\adfsplit
$(0, 4, 47, 100)$,
$(276, 0, 75, 150)$

%ADFvfyBlocksEnd
\adfLgap \noindent by the mapping:
$x \mapsto x +  j \adfmod{225}$ for $x < 225$,
$x \mapsto (x +  j \adfmod{45}) + 225$ for $225 \le x < 270$,
$x \mapsto (x + 2 j \adfmod{6}) + 270$ for $270 \le x < 276$,
$276 \mapsto 276$,
$0 \le j < 225$
 for the first 25 blocks,
$0 \le j < 75$
 for the last block.
\ADFvfyParStart{(277, ((25, 225, ((225, 1), (45, 1), (6, 2), (1, 1))), (1, 75, ((225, 1), (45, 1), (6, 2), (1, 1)))), ((25, 9), (52, 1)))} %ADFvfyParEnd
% End of 25^9 52^1
%%%%%%%%%%%%%%%%%%%%%%%%%%%%%%%%%%%%%%%%%%%%%%%%%%%%%%%%%%%%%%%%%%%%%%%%%%%%%%%%%%%%%%%%%%
%%%%%%%%%%%%%%%%%%%%%%%%%%%%%%%%%%%%%%%%%%%%%%%%%%%%%%%%%%%%%%%%%%%%%%%%%%%%%%%%%%%%%%%%%%

% Charlotte:GDD4-1-3-5-mod-6-TeX-gen-A:HITS-fun:4.10
\adfDgap
%ADFvfyBlocksStart {25,25,25,25,25,25,25,25,25,58}
\noindent{\boldmath $ 25^{9} 58^{1} $}~
With the point set $Z_{283}$ partitioned into
 residue classes modulo $9$ for $\{0, 1, \dots, 224\}$, and
 $\{225, 226, \dots, 282\}$,
 the design is generated from

\adfLgap %ADFvfyDesignStart
$(270, 133, 134, 69)$,
$(271, 208, 21, 221)$,
$(272, 104, 61, 78)$,
$(273, 224, 87, 202)$,\adfsplit
$(225, 152, 113, 84)$,
$(225, 207, 49, 196)$,
$(225, 12, 178, 59)$,
$(225, 32, 144, 172)$,\adfsplit
$(225, 50, 204, 80)$,
$(225, 97, 186, 93)$,
$(225, 169, 1, 78)$,
$(225, 182, 166, 90)$,\adfsplit
$(225, 30, 224, 103)$,
$(225, 154, 110, 206)$,
$(0, 2, 5, 234)$,
$(0, 7, 19, 249)$,\adfsplit
$(0, 14, 193, 228)$,
$(0, 42, 125, 245)$,
$(0, 24, 79, 241)$,
$(0, 48, 98, 164)$,\adfsplit
$(0, 34, 74, 156)$,
$(0, 33, 70, 172)$,
$(0, 20, 80, 131)$,
$(0, 8, 95, 105)$,\adfsplit
$(0, 23, 58, 107)$,
$(0, 6, 21, 62)$,
$(282, 0, 75, 150)$

%ADFvfyBlocksEnd
\adfLgap \noindent by the mapping:
$x \mapsto x +  j \adfmod{225}$ for $x < 225$,
$x \mapsto (x +  j \adfmod{45}) + 225$ for $225 \le x < 270$,
$x \mapsto (x - 270 + 4 j \adfmod{12}) + 270$ for $270 \le x < 282$,
$282 \mapsto 282$,
$0 \le j < 225$
 for the first 26 blocks,
$0 \le j < 75$
 for the last block.
\ADFvfyParStart{(283, ((26, 225, ((225, 1), (45, 1), (12, 4), (1, 1))), (1, 75, ((225, 1), (45, 1), (12, 4), (1, 1)))), ((25, 9), (58, 1)))} %ADFvfyParEnd
% End of 25^9 58^1
%%%%%%%%%%%%%%%%%%%%%%%%%%%%%%%%%%%%%%%%%%%%%%%%%%%%%%%%%%%%%%%%%%%%%%%%%%%%%%%%%%%%%%%%%%
%%%%%%%%%%%%%%%%%%%%%%%%%%%%%%%%%%%%%%%%%%%%%%%%%%%%%%%%%%%%%%%%%%%%%%%%%%%%%%%%%%%%%%%%%%

% Charlotte:GDD4-1-3-5-mod-6-TeX-gen-A:HITS-fun:4.10
\adfDgap
%ADFvfyBlocksStart {25,25,25,25,25,25,25,25,25,64}
\noindent{\boldmath $ 25^{9} 64^{1} $}~
With the point set $Z_{289}$ partitioned into
 residue classes modulo $9$ for $\{0, 1, \dots, 224\}$, and
 $\{225, 226, \dots, 288\}$,
 the design is generated from

\adfLgap %ADFvfyDesignStart
$(270, 95, 120, 35)$,
$(270, 7, 100, 2)$,
$(270, 63, 13, 186)$,
$(271, 72, 140, 60)$,\adfsplit
$(271, 151, 37, 152)$,
$(271, 112, 29, 48)$,
$(225, 47, 80, 210)$,
$(225, 220, 222, 21)$,\adfsplit
$(225, 217, 101, 88)$,
$(225, 204, 193, 190)$,
$(225, 173, 188, 5)$,
$(225, 41, 57, 196)$,\adfsplit
$(225, 208, 126, 164)$,
$(225, 91, 211, 60)$,
$(225, 200, 67, 59)$,
$(0, 4, 195, 267)$,\adfsplit
$(0, 6, 103, 113)$,
$(0, 20, 67, 154)$,
$(0, 7, 29, 166)$,
$(0, 43, 94, 147)$,\adfsplit
$(0, 28, 176, 266)$,
$(0, 23, 55, 179)$,
$(0, 17, 73, 236)$,
$(0, 41, 89, 243)$,\adfsplit
$(0, 21, 61, 167)$,
$(0, 35, 100, 261)$,
$(0, 37, 186, 237)$,
$(288, 0, 75, 150)$

%ADFvfyBlocksEnd
\adfLgap \noindent by the mapping:
$x \mapsto x +  j \adfmod{225}$ for $x < 225$,
$x \mapsto (x +  j \adfmod{45}) + 225$ for $225 \le x < 270$,
$x \mapsto (x + 2 j \adfmod{18}) + 270$ for $270 \le x < 288$,
$288 \mapsto 288$,
$0 \le j < 225$
 for the first 27 blocks,
$0 \le j < 75$
 for the last block.
\ADFvfyParStart{(289, ((27, 225, ((225, 1), (45, 1), (18, 2), (1, 1))), (1, 75, ((225, 1), (45, 1), (18, 2), (1, 1)))), ((25, 9), (64, 1)))} %ADFvfyParEnd
% End of 25^9 64^1
%%%%%%%%%%%%%%%%%%%%%%%%%%%%%%%%%%%%%%%%%%%%%%%%%%%%%%%%%%%%%%%%%%%%%%%%%%%%%%%%%%%%%%%%%%
%%%%%%%%%%%%%%%%%%%%%%%%%%%%%%%%%%%%%%%%%%%%%%%%%%%%%%%%%%%%%%%%%%%%%%%%%%%%%%%%%%%%%%%%%%

% Charlotte:GDD4-1-3-5-mod-6-TeX-gen-A:HITS-fun:4.10
\adfDgap
%ADFvfyBlocksStart {25,25,25,25,25,25,25,25,25,76}
\noindent{\boldmath $ 25^{9} 76^{1} $}~
With the point set $Z_{301}$ partitioned into
 residue classes modulo $9$ for $\{0, 1, \dots, 224\}$, and
 $\{225, 226, \dots, 300\}$,
 the design is generated from

\adfLgap %ADFvfyDesignStart
$(225, 222, 135, 7)$,
$(225, 78, 103, 106)$,
$(225, 11, 121, 192)$,
$(225, 68, 127, 188)$,\adfsplit
$(225, 207, 40, 200)$,
$(225, 154, 122, 120)$,
$(225, 62, 184, 99)$,
$(225, 49, 75, 167)$,\adfsplit
$(225, 53, 211, 65)$,
$(225, 198, 134, 33)$,
$(225, 81, 13, 30)$,
$(225, 214, 51, 89)$,\adfsplit
$(225, 94, 171, 129)$,
$(225, 175, 145, 194)$,
$(225, 35, 58, 74)$,
$(225, 116, 102, 159)$,\adfsplit
$(225, 186, 104, 173)$,
$(0, 4, 205, 282)$,
$(0, 8, 147, 231)$,
$(0, 31, 83, 288)$,\adfsplit
$(0, 6, 112, 274)$,
$(0, 11, 102, 244)$,
$(0, 22, 96, 149)$,
$(0, 21, 114, 284)$,\adfsplit
$(0, 5, 55, 121)$,
$(0, 29, 70, 290)$,
$(0, 56, 136, 285)$,
$(0, 15, 48, 88)$,\adfsplit
$(0, 1, 47, 131)$,
$(300, 0, 75, 150)$

%ADFvfyBlocksEnd
\adfLgap \noindent by the mapping:
$x \mapsto x +  j \adfmod{225}$ for $x < 225$,
$x \mapsto (x +  j \adfmod{75}) + 225$ for $225 \le x < 300$,
$300 \mapsto 300$,
$0 \le j < 225$
 for the first 29 blocks,
$0 \le j < 75$
 for the last block.
\ADFvfyParStart{(301, ((29, 225, ((225, 1), (75, 1), (1, 1))), (1, 75, ((225, 1), (75, 1), (1, 1)))), ((25, 9), (76, 1)))} %ADFvfyParEnd
% End of 25^9 76^1
%%%%%%%%%%%%%%%%%%%%%%%%%%%%%%%%%%%%%%%%%%%%%%%%%%%%%%%%%%%%%%%%%%%%%%%%%%%%%%%%%%%%%%%%%%
%%%%%%%%%%%%%%%%%%%%%%%%%%%%%%%%%%%%%%%%%%%%%%%%%%%%%%%%%%%%%%%%%%%%%%%%%%%%%%%%%%%%%%%%%%

% Charlotte:GDD4-1-3-5-mod-6-TeX-gen-A:HITS-fun:4.10
\adfDgap
%ADFvfyBlocksStart {25,25,25,25,25,25,25,25,25,82}
\noindent{\boldmath $ 25^{9} 82^{1} $}~
With the point set $Z_{307}$ partitioned into
 residue classes modulo $9$ for $\{0, 1, \dots, 224\}$, and
 $\{225, 226, \dots, 306\}$,
 the design is generated from

\adfLgap %ADFvfyDesignStart
$(300, 146, 93, 1)$,
$(301, 165, 44, 76)$,
$(225, 164, 59, 114)$,
$(225, 119, 165, 9)$,\adfsplit
$(225, 47, 3, 10)$,
$(225, 57, 74, 70)$,
$(225, 143, 105, 214)$,
$(225, 25, 191, 91)$,\adfsplit
$(225, 161, 200, 4)$,
$(225, 99, 193, 75)$,
$(225, 53, 183, 51)$,
$(225, 32, 21, 62)$,\adfsplit
$(225, 2, 140, 204)$,
$(225, 176, 7, 29)$,
$(225, 55, 113, 98)$,
$(225, 31, 36, 196)$,\adfsplit
$(225, 155, 147, 28)$,
$(225, 37, 138, 185)$,
$(0, 1, 74, 84)$,
$(0, 3, 51, 227)$,\adfsplit
$(0, 20, 122, 278)$,
$(0, 6, 34, 91)$,
$(0, 12, 31, 204)$,
$(0, 62, 129, 239)$,\adfsplit
$(0, 40, 82, 280)$,
$(0, 14, 190, 248)$,
$(0, 70, 149, 287)$,
$(0, 26, 114, 281)$,\adfsplit
$(0, 25, 86, 244)$,
$(0, 16, 128, 282)$,
$(306, 0, 75, 150)$

%ADFvfyBlocksEnd
\adfLgap \noindent by the mapping:
$x \mapsto x +  j \adfmod{225}$ for $x < 225$,
$x \mapsto (x +  j \adfmod{75}) + 225$ for $225 \le x < 300$,
$x \mapsto (x + 2 j \adfmod{6}) + 300$ for $300 \le x < 306$,
$306 \mapsto 306$,
$0 \le j < 225$
 for the first 30 blocks,
$0 \le j < 75$
 for the last block.
\ADFvfyParStart{(307, ((30, 225, ((225, 1), (75, 1), (6, 2), (1, 1))), (1, 75, ((225, 1), (75, 1), (6, 2), (1, 1)))), ((25, 9), (82, 1)))} %ADFvfyParEnd
% End of 25^9 82^1
%%%%%%%%%%%%%%%%%%%%%%%%%%%%%%%%%%%%%%%%%%%%%%%%%%%%%%%%%%%%%%%%%%%%%%%%%%%%%%%%%%%%%%%%%%
%%%%%%%%%%%%%%%%%%%%%%%%%%%%%%%%%%%%%%%%%%%%%%%%%%%%%%%%%%%%%%%%%%%%%%%%%%%%%%%%%%%%%%%%%%

% Charlotte:GDD4-1-3-5-mod-6-TeX-gen-A:HITS-fun:4.10
\adfDgap
%ADFvfyBlocksStart {25,25,25,25,25,25,25,25,25,88}
\noindent{\boldmath $ 25^{9} 88^{1} $}~
With the point set $Z_{313}$ partitioned into
 residue classes modulo $9$ for $\{0, 1, \dots, 224\}$, and
 $\{225, 226, \dots, 312\}$,
 the design is generated from

\adfLgap %ADFvfyDesignStart
$(300, 185, 85, 120)$,
$(301, 96, 70, 74)$,
$(302, 98, 42, 199)$,
$(303, 28, 117, 53)$,\adfsplit
$(225, 31, 60, 23)$,
$(225, 118, 15, 152)$,
$(225, 220, 171, 161)$,
$(225, 56, 117, 221)$,\adfsplit
$(225, 87, 102, 88)$,
$(225, 82, 200, 148)$,
$(225, 53, 74, 81)$,
$(225, 95, 138, 18)$,\adfsplit
$(225, 207, 62, 110)$,
$(225, 188, 61, 190)$,
$(225, 65, 198, 195)$,
$(225, 78, 201, 67)$,\adfsplit
$(225, 4, 183, 34)$,
$(225, 101, 194, 143)$,
$(225, 29, 204, 172)$,
$(0, 5, 38, 151)$,\adfsplit
$(0, 12, 170, 236)$,
$(0, 19, 39, 264)$,
$(0, 44, 155, 259)$,
$(0, 13, 86, 283)$,\adfsplit
$(0, 16, 40, 87)$,
$(0, 62, 156, 248)$,
$(0, 31, 140, 231)$,
$(0, 53, 168, 281)$,\adfsplit
$(0, 17, 58, 251)$,
$(0, 6, 147, 292)$,
$(0, 23, 142, 276)$,
$(312, 0, 75, 150)$

%ADFvfyBlocksEnd
\adfLgap \noindent by the mapping:
$x \mapsto x +  j \adfmod{225}$ for $x < 225$,
$x \mapsto (x +  j \adfmod{75}) + 225$ for $225 \le x < 300$,
$x \mapsto (x + 4 j \adfmod{12}) + 300$ for $300 \le x < 312$,
$312 \mapsto 312$,
$0 \le j < 225$
 for the first 31 blocks,
$0 \le j < 75$
 for the last block.
\ADFvfyParStart{(313, ((31, 225, ((225, 1), (75, 1), (12, 4), (1, 1))), (1, 75, ((225, 1), (75, 1), (12, 4), (1, 1)))), ((25, 9), (88, 1)))} %ADFvfyParEnd
% End of 25^9 88^1
%%%%%%%%%%%%%%%%%%%%%%%%%%%%%%%%%%%%%%%%%%%%%%%%%%%%%%%%%%%%%%%%%%%%%%%%%%%%%%%%%%%%%%%%%%
%%%%%%%%%%%%%%%%%%%%%%%%%%%%%%%%%%%%%%%%%%%%%%%%%%%%%%%%%%%%%%%%%%%%%%%%%%%%%%%%%%%%%%%%%%

% Charlotte:GDD4-1-3-5-mod-6-TeX-gen-A:HITS-fun:4.10
\adfDgap
%ADFvfyBlocksStart {25,25,25,25,25,25,25,25,25,94}
\noindent{\boldmath $ 25^{9} 94^{1} $}~
With the point set $Z_{319}$ partitioned into
 residue classes modulo $9$ for $\{0, 1, \dots, 224\}$, and
 $\{225, 226, \dots, 318\}$,
 the design is generated from

\adfLgap %ADFvfyDesignStart
$(300, 97, 188, 149)$,
$(300, 164, 123, 103)$,
$(300, 219, 154, 126)$,
$(301, 81, 155, 49)$,\adfsplit
$(301, 87, 212, 21)$,
$(301, 62, 7, 55)$,
$(225, 46, 4, 200)$,
$(225, 14, 182, 81)$,\adfsplit
$(225, 92, 77, 19)$,
$(225, 64, 174, 13)$,
$(225, 12, 158, 25)$,
$(225, 150, 149, 133)$,\adfsplit
$(225, 16, 51, 220)$,
$(225, 26, 29, 131)$,
$(225, 186, 45, 142)$,
$(225, 57, 199, 103)$,\adfsplit
$(225, 5, 216, 73)$,
$(225, 165, 61, 159)$,
$(225, 30, 117, 115)$,
$(225, 222, 172, 203)$,\adfsplit
$(0, 4, 111, 116)$,
$(0, 23, 145, 261)$,
$(0, 60, 136, 231)$,
$(0, 43, 137, 267)$,\adfsplit
$(0, 11, 187, 293)$,
$(0, 26, 192, 240)$,
$(0, 12, 200, 265)$,
$(0, 22, 69, 263)$,\adfsplit
$(0, 8, 155, 237)$,
$(0, 10, 40, 262)$,
$(0, 24, 86, 259)$,
$(0, 53, 130, 257)$,\adfsplit
$(318, 0, 75, 150)$

%ADFvfyBlocksEnd
\adfLgap \noindent by the mapping:
$x \mapsto x +  j \adfmod{225}$ for $x < 225$,
$x \mapsto (x +  j \adfmod{75}) + 225$ for $225 \le x < 300$,
$x \mapsto (x - 300 + 2 j \adfmod{18}) + 300$ for $300 \le x < 318$,
$318 \mapsto 318$,
$0 \le j < 225$
 for the first 32 blocks,
$0 \le j < 75$
 for the last block.
\ADFvfyParStart{(319, ((32, 225, ((225, 1), (75, 1), (18, 2), (1, 1))), (1, 75, ((225, 1), (75, 1), (18, 2), (1, 1)))), ((25, 9), (94, 1)))} %ADFvfyParEnd
% End of 25^9 94^1
%%%%%%%%%%%%%%%%%%%%%%%%%%%%%%%%%%%%%%%%%%%%%%%%%%%%%%%%%%%%%%%%%%%%%%%%%%%%%%%%%%%%%%%%%%
%%%%%%%%%%%%%%%%%%%%%%%%%%%%%%%%%%%%%%%%%%%%%%%%%%%%%%%%%%%%%%%%%%%%%%%%%%%%%%%%%%%%%%%%%%

% Charlotte:GDD4-1-3-5-mod-6-TeX-gen-A:HITS-fun:4.10
\adfDgap
%ADFvfyBlocksStart {25,25,25,25,25,25,25,25,25,25,25,25,25,25,25,25,25,25,25,25,25,16}
\noindent{\boldmath $ 25^{21} 16^{1} $}~
With the point set $Z_{541}$ partitioned into
 residue classes modulo $21$ for $\{0, 1, \dots, 524\}$, and
 $\{525, 526, \dots, 540\}$,
 the design is generated from

\adfLgap %ADFvfyDesignStart
$(525, 379, 372, 50)$,
$(526, 245, 402, 142)$,
$(527, 485, 312, 88)$,
$(528, 464, 285, 406)$,\adfsplit
$(529, 251, 57, 247)$,
$(317, 419, 421, 508)$,
$(309, 169, 345, 5)$,
$(465, 449, 303, 4)$,\adfsplit
$(212, 199, 139, 227)$,
$(377, 238, 155, 138)$,
$(474, 297, 18, 167)$,
$(413, 480, 90, 53)$,\adfsplit
$(111, 61, 221, 291)$,
$(264, 152, 237, 245)$,
$(367, 209, 98, 43)$,
$(443, 465, 490, 74)$,\adfsplit
$(258, 506, 509, 463)$,
$(441, 149, 225, 282)$,
$(71, 122, 216, 30)$,
$(495, 302, 376, 170)$,\adfsplit
$(493, 259, 139, 414)$,
$(434, 138, 24, 500)$,
$(245, 7, 470, 1)$,
$(15, 256, 131, 443)$,\adfsplit
$(293, 53, 417, 170)$,
$(495, 97, 83, 169)$,
$(281, 14, 228, 19)$,
$(402, 250, 467, 368)$,\adfsplit
$(142, 412, 243, 413)$,
$(0, 9, 20, 38)$,
$(0, 10, 153, 192)$,
$(0, 24, 54, 371)$,\adfsplit
$(0, 31, 219, 314)$,
$(0, 35, 142, 232)$,
$(0, 52, 183, 264)$,
$(0, 12, 45, 302)$,\adfsplit
$(0, 68, 150, 288)$,
$(0, 61, 198, 276)$,
$(0, 44, 195, 272)$,
$(0, 23, 129, 204)$,\adfsplit
$(0, 32, 91, 298)$,
$(0, 26, 148, 384)$,
$(0, 71, 243, 351)$,
$(0, 40, 136, 184)$,\adfsplit
$(540, 0, 175, 350)$

%ADFvfyBlocksEnd
\adfLgap \noindent by the mapping:
$x \mapsto x +  j \adfmod{525}$ for $x < 525$,
$x \mapsto (x + 5 j \adfmod{15}) + 525$ for $525 \le x < 540$,
$540 \mapsto 540$,
$0 \le j < 525$
 for the first 44 blocks,
$0 \le j < 175$
 for the last block.
\ADFvfyParStart{(541, ((44, 525, ((525, 1), (15, 5), (1, 1))), (1, 175, ((525, 1), (15, 5), (1, 1)))), ((25, 21), (16, 1)))} %ADFvfyParEnd
% End of 25^21 16^1
%%%%%%%%%%%%%%%%%%%%%%%%%%%%%%%%%%%%%%%%%%%%%%%%%%%%%%%%%%%%%%%%%%%%%%%%%%%%%%%%%%%%%%%%%%
%%%%%%%%%%%%%%%%%%%%%%%%%%%%%%%%%%%%%%%%%%%%%%%%%%%%%%%%%%%%%%%%%%%%%%%%%%%%%%%%%%%%%%%%%%

% Charlotte:GDD4-1-3-5-mod-6-TeX-gen-A:HITS-fun:4.10
\adfDgap
%ADFvfyBlocksStart {25,25,25,25,25,25,25,25,25,25,25,25,25,25,25,25,25,25,25,25,25,22}
\noindent{\boldmath $ 25^{21} 22^{1} $}~
With the point set $Z_{547}$ partitioned into
 residue classes modulo $21$ for $\{0, 1, \dots, 524\}$, and
 $\{525, 526, \dots, 546\}$,
 the design is generated from

\adfLgap %ADFvfyDesignStart
$(525, 466, 258, 215)$,
$(526, 349, 39, 401)$,
$(527, 177, 305, 310)$,
$(528, 109, 92, 174)$,\adfsplit
$(529, 501, 32, 370)$,
$(530, 269, 474, 247)$,
$(531, 213, 77, 457)$,
$(410, 460, 422, 344)$,\adfsplit
$(263, 295, 50, 226)$,
$(374, 154, 256, 97)$,
$(358, 397, 243, 55)$,
$(135, 223, 41, 150)$,\adfsplit
$(48, 283, 252, 407)$,
$(50, 309, 429, 319)$,
$(57, 169, 397, 493)$,
$(90, 392, 154, 344)$,\adfsplit
$(514, 374, 434, 66)$,
$(95, 229, 69, 88)$,
$(131, 373, 115, 401)$,
$(113, 18, 90, 381)$,\adfsplit
$(509, 126, 317, 323)$,
$(431, 475, 455, 215)$,
$(104, 430, 497, 205)$,
$(443, 87, 505, 132)$,\adfsplit
$(442, 381, 23, 467)$,
$(357, 450, 286, 201)$,
$(328, 508, 96, 500)$,
$(352, 323, 45, 4)$,\adfsplit
$(221, 358, 151, 450)$,
$(169, 470, 340, 129)$,
$(0, 1, 3, 492)$,
$(0, 4, 104, 386)$,\adfsplit
$(0, 9, 27, 476)$,
$(0, 30, 117, 381)$,
$(0, 54, 129, 203)$,
$(0, 35, 158, 360)$,\adfsplit
$(0, 14, 125, 346)$,
$(0, 46, 99, 352)$,
$(0, 55, 153, 250)$,
$(0, 68, 151, 347)$,\adfsplit
$(0, 51, 170, 364)$,
$(0, 13, 127, 390)$,
$(0, 59, 181, 289)$,
$(0, 11, 90, 327)$,\adfsplit
$(0, 47, 138, 422)$,
$(546, 0, 175, 350)$

%ADFvfyBlocksEnd
\adfLgap \noindent by the mapping:
$x \mapsto x +  j \adfmod{525}$ for $x < 525$,
$x \mapsto (x + 7 j \adfmod{21}) + 525$ for $525 \le x < 546$,
$546 \mapsto 546$,
$0 \le j < 525$
 for the first 45 blocks,
$0 \le j < 175$
 for the last block.
\ADFvfyParStart{(547, ((45, 525, ((525, 1), (21, 7), (1, 1))), (1, 175, ((525, 1), (21, 7), (1, 1)))), ((25, 21), (22, 1)))} %ADFvfyParEnd
% End of 25^21 22^1
%%%%%%%%%%%%%%%%%%%%%%%%%%%%%%%%%%%%%%%%%%%%%%%%%%%%%%%%%%%%%%%%%%%%%%%%%%%%%%%%%%%%%%%%%%
%%%%%%%%%%%%%%%%%%%%%%%%%%%%%%%%%%%%%%%%%%%%%%%%%%%%%%%%%%%%%%%%%%%%%%%%%%%%%%%%%%%%%%%%%%

%%%%%%%%%%%%%%%%%%%%%%%%%%%%%%%%%%%%%%%%%%%%%%%%%%%%%%%%%%%%%%%%%%%%%%%%%%%%%%%%%%%%%%%%%%
%%%%%%%%%%%%%%%%%%%%%%%%%%%%%%%%%%%%%%%%%%%%%%%%%%%%%%%%%%%%%%%%%%%%%%%%%%%%%%%%%%%%%%%%%%
\section{4-GDDs for the proof of Lemma \ref{lem:4-GDD 29^u m^1}}
\label{app:4-GDD 29^u m^1}
\adfnull{
$ 29^{12} 14^1 $,
$ 29^{12} 17^1 $,
$ 29^{12} 20^1 $,
$ 29^{12} 23^1 $,
$ 29^{12} 26^1 $,
$ 29^{24} 26^1 $,
$ 29^{15} 17^1 $,
$ 29^{15} 23^1 $,
$ 29^9 14^1 $,
$ 29^9 20^1 $,
$ 29^9 26^1 $,
$ 29^9 32^1 $,
$ 29^9 38^1 $,
$ 29^9 44^1 $,
$ 29^9 50^1 $,
$ 29^9 56^1 $,
$ 29^9 62^1 $,
$ 29^9 68^1 $,
$ 29^9 74^1 $,
$ 29^9 80^1 $,
$ 29^9 86^1 $,
$ 29^9 92^1 $,
$ 29^9 98^1 $,
$ 29^9 104^1 $,
$ 29^9 110^1 $ and
$ 29^{21} 26^1 $.
}

% Charlotte:GDD4-1-3-5-mod-6-TeX-gen-A:HITS-fun:4.10
\adfDgap
%ADFvfyBlocksStart {29,29,29,29,29,29,29,29,29,29,29,29,14}
\noindent{\boldmath $ 29^{12} 14^{1} $}~
With the point set $Z_{362}$ partitioned into
 residue classes modulo $12$ for $\{0, 1, \dots, 347\}$, and
 $\{348, 349, \dots, 361\}$,
 the design is generated from

\adfLgap %ADFvfyDesignStart
$(348, 0, 1, 2)$,
$(349, 0, 115, 233)$,
$(350, 0, 118, 116)$,
$(351, 0, 232, 347)$,\adfsplit
$(352, 0, 346, 230)$,
$(353, 0, 4, 11)$,
$(354, 0, 7, 344)$,
$(355, 0, 337, 341)$,\adfsplit
$(356, 0, 10, 23)$,
$(357, 0, 13, 338)$,
$(358, 0, 325, 335)$,
$(359, 0, 16, 35)$,\adfsplit
$(360, 0, 19, 332)$,
$(361, 0, 313, 329)$,
$(205, 319, 155, 233)$,
$(98, 187, 277, 285)$,\adfsplit
$(115, 82, 27, 318)$,
$(329, 131, 237, 134)$,
$(295, 5, 277, 86)$,
$(145, 274, 330, 220)$,\adfsplit
$(180, 260, 53, 67)$,
$(323, 180, 76, 218)$,
$(146, 155, 180, 277)$,
$(295, 87, 263, 324)$,\adfsplit
$(321, 214, 254, 161)$,
$(156, 334, 308, 161)$,
$(222, 39, 98, 175)$,
$(0, 15, 42, 83)$,\adfsplit
$(0, 31, 70, 189)$,
$(0, 49, 128, 263)$,
$(0, 6, 52, 283)$,
$(0, 22, 91, 121)$,\adfsplit
$(0, 51, 133, 199)$,
$(0, 17, 126, 171)$,
$(0, 37, 100, 162)$,
$(0, 21, 95, 218)$,\adfsplit
$(0, 43, 137, 210)$,
$(0, 20, 64, 166)$,
$(0, 87, 174, 261)$

%ADFvfyBlocksEnd
\adfLgap \noindent by the mapping:
$x \mapsto x + 3 j \adfmod{348}$ for $x < 348$,
$x \mapsto x$ for $x \ge 348$,
$0 \le j < 116$
 for the first 14 blocks;
$x \mapsto x +  j \adfmod{348}$ for $x < 348$,
$x \mapsto x$ for $x \ge 348$,
$0 \le j < 348$
 for the next 24 blocks,
$0 \le j < 87$
 for the last block.
\ADFvfyParStart{(362, ((14, 116, ((348, 3), (14, 14))), (24, 348, ((348, 1), (14, 14))), (1, 87, ((348, 1), (14, 14)))), ((29, 12), (14, 1)))} %ADFvfyParEnd
% End of 29^12 14^1
%%%%%%%%%%%%%%%%%%%%%%%%%%%%%%%%%%%%%%%%%%%%%%%%%%%%%%%%%%%%%%%%%%%%%%%%%%%%%%%%%%%%%%%%%%
%%%%%%%%%%%%%%%%%%%%%%%%%%%%%%%%%%%%%%%%%%%%%%%%%%%%%%%%%%%%%%%%%%%%%%%%%%%%%%%%%%%%%%%%%%

% Charlotte:GDD4-1-3-5-mod-6-TeX-gen-A:HITS-fun:4.10
\adfDgap
%ADFvfyBlocksStart {29,29,29,29,29,29,29,29,29,29,29,29,17}
\noindent{\boldmath $ 29^{12} 17^{1} $}~
With the point set $Z_{365}$ partitioned into
 residue classes modulo $12$ for $\{0, 1, \dots, 347\}$, and
 $\{348, 349, \dots, 364\}$,
 the design is generated from

\adfLgap %ADFvfyDesignStart
$(348, 0, 1, 2)$,
$(349, 0, 115, 233)$,
$(350, 0, 118, 116)$,
$(351, 0, 232, 347)$,\adfsplit
$(352, 0, 346, 230)$,
$(353, 0, 4, 11)$,
$(354, 0, 7, 344)$,
$(355, 0, 337, 341)$,\adfsplit
$(356, 0, 10, 23)$,
$(357, 0, 13, 338)$,
$(358, 0, 325, 335)$,
$(359, 0, 16, 35)$,\adfsplit
$(360, 0, 19, 332)$,
$(361, 0, 313, 329)$,
$(362, 0, 22, 5)$,
$(363, 0, 331, 326)$,\adfsplit
$(364, 0, 343, 17)$,
$(15, 308, 302, 66)$,
$(80, 170, 301, 275)$,
$(52, 155, 110, 319)$,\adfsplit
$(199, 145, 48, 318)$,
$(219, 210, 170, 262)$,
$(261, 6, 194, 279)$,
$(158, 306, 103, 321)$,\adfsplit
$(286, 328, 188, 81)$,
$(32, 292, 167, 266)$,
$(32, 187, 245, 275)$,
$(19, 88, 254, 284)$,\adfsplit
$(134, 96, 283, 243)$,
$(25, 114, 46, 96)$,
$(176, 220, 11, 120)$,
$(90, 288, 133, 164)$,\adfsplit
$(232, 41, 296, 106)$,
$(146, 318, 248, 89)$,
$(9, 151, 298, 123)$,
$(125, 0, 254, 215)$,\adfsplit
$(85, 29, 196, 275)$,
$(108, 11, 199, 325)$,
$(297, 121, 203, 266)$,
$(308, 6, 340, 113)$,\adfsplit
$(1, 225, 184, 320)$,
$(187, 257, 336, 87)$,
$(209, 134, 276, 301)$,
$(31, 299, 313, 293)$,\adfsplit
$(0, 3, 28, 294)$,
$(0, 21, 62, 89)$,
$(0, 8, 47, 226)$,
$(0, 33, 65, 202)$,\adfsplit
$(0, 76, 210, 333)$,
$(0, 34, 162, 231)$,
$(0, 71, 104, 184)$,
$(0, 20, 245, 321)$,\adfsplit
$(0, 127, 171, 275)$,
$(0, 37, 86, 263)$,
$(0, 137, 265, 303)$,
$(0, 77, 178, 289)$,\adfsplit
$(1, 9, 51, 163)$,
$(0, 141, 287, 339)$,
$(0, 133, 207, 271)$,
$(0, 73, 194, 311)$,\adfsplit
$(0, 51, 129, 191)$,
$(0, 53, 159, 205)$,
$(0, 63, 66, 161)$,
$(0, 81, 110, 295)$,\adfsplit
$(0, 87, 174, 261)$

%ADFvfyBlocksEnd
\adfLgap \noindent by the mapping:
$x \mapsto x + 3 j \adfmod{348}$ for $x < 348$,
$x \mapsto x$ for $x \ge 348$,
$0 \le j < 116$
 for the first 17 blocks;
$x \mapsto x + 2 j \adfmod{348}$ for $x < 348$,
$x \mapsto x$ for $x \ge 348$,
$0 \le j < 174$
 for the next 47 blocks,
$0 \le j < 87$
 for the last block.
\ADFvfyParStart{(365, ((17, 116, ((348, 3), (17, 17))), (47, 174, ((348, 2), (17, 17))), (1, 87, ((348, 2), (17, 17)))), ((29, 12), (17, 1)))} %ADFvfyParEnd
% End of 29^12 17^1
%%%%%%%%%%%%%%%%%%%%%%%%%%%%%%%%%%%%%%%%%%%%%%%%%%%%%%%%%%%%%%%%%%%%%%%%%%%%%%%%%%%%%%%%%%
%%%%%%%%%%%%%%%%%%%%%%%%%%%%%%%%%%%%%%%%%%%%%%%%%%%%%%%%%%%%%%%%%%%%%%%%%%%%%%%%%%%%%%%%%%

% Charlotte:GDD4-1-3-5-mod-6-TeX-gen-A:HITS-fun:4.10
\adfDgap
%ADFvfyBlocksStart {29,29,29,29,29,29,29,29,29,29,29,29,20}
\noindent{\boldmath $ 29^{12} 20^{1} $}~
With the point set $Z_{368}$ partitioned into
 residue classes modulo $12$ for $\{0, 1, \dots, 347\}$, and
 $\{348, 349, \dots, 367\}$,
 the design is generated from

\adfLgap %ADFvfyDesignStart
$(348, 0, 1, 2)$,
$(349, 0, 115, 233)$,
$(350, 0, 118, 116)$,
$(351, 0, 232, 347)$,\adfsplit
$(352, 0, 346, 230)$,
$(353, 0, 4, 11)$,
$(354, 0, 7, 344)$,
$(355, 0, 337, 341)$,\adfsplit
$(356, 0, 10, 23)$,
$(357, 0, 13, 338)$,
$(358, 0, 325, 335)$,
$(359, 0, 16, 35)$,\adfsplit
$(360, 0, 19, 332)$,
$(361, 0, 313, 329)$,
$(362, 0, 22, 5)$,
$(363, 0, 331, 326)$,\adfsplit
$(364, 0, 343, 17)$,
$(365, 0, 25, 53)$,
$(366, 0, 28, 323)$,
$(367, 0, 295, 320)$,\adfsplit
$(98, 271, 205, 304)$,
$(121, 48, 233, 171)$,
$(296, 252, 290, 13)$,
$(50, 24, 104, 202)$,\adfsplit
$(168, 78, 121, 44)$,
$(57, 150, 304, 75)$,
$(275, 15, 302, 126)$,
$(213, 70, 18, 73)$,\adfsplit
$(214, 50, 1, 183)$,
$(211, 86, 116, 287)$,
$(330, 245, 38, 95)$,
$(99, 136, 7, 85)$,\adfsplit
$(0, 8, 145, 191)$,
$(0, 15, 104, 136)$,
$(0, 45, 103, 231)$,
$(0, 20, 94, 134)$,\adfsplit
$(0, 69, 155, 238)$,
$(0, 41, 105, 251)$,
$(0, 79, 160, 266)$,
$(0, 29, 68, 131)$,\adfsplit
$(0, 9, 148, 190)$,
$(0, 21, 91, 218)$,
$(0, 59, 126, 248)$,
$(0, 87, 174, 261)$

%ADFvfyBlocksEnd
\adfLgap \noindent by the mapping:
$x \mapsto x + 3 j \adfmod{348}$ for $x < 348$,
$x \mapsto x$ for $x \ge 348$,
$0 \le j < 116$
 for the first 20 blocks;
$x \mapsto x +  j \adfmod{348}$ for $x < 348$,
$x \mapsto x$ for $x \ge 348$,
$0 \le j < 348$
 for the next 23 blocks,
$0 \le j < 87$
 for the last block.
\ADFvfyParStart{(368, ((20, 116, ((348, 3), (20, 20))), (23, 348, ((348, 1), (20, 20))), (1, 87, ((348, 1), (20, 20)))), ((29, 12), (20, 1)))} %ADFvfyParEnd
% End of 29^12 20^1
%%%%%%%%%%%%%%%%%%%%%%%%%%%%%%%%%%%%%%%%%%%%%%%%%%%%%%%%%%%%%%%%%%%%%%%%%%%%%%%%%%%%%%%%%%
%%%%%%%%%%%%%%%%%%%%%%%%%%%%%%%%%%%%%%%%%%%%%%%%%%%%%%%%%%%%%%%%%%%%%%%%%%%%%%%%%%%%%%%%%%

% Charlotte:GDD4-1-3-5-mod-6-TeX-gen-A:HITS-fun:4.10
\adfDgap
%ADFvfyBlocksStart {29,29,29,29,29,29,29,29,29,29,29,29,23}
\noindent{\boldmath $ 29^{12} 23^{1} $}~
With the point set $Z_{371}$ partitioned into
 residue classes modulo $12$ for $\{0, 1, \dots, 347\}$, and
 $\{348, 349, \dots, 370\}$,
 the design is generated from

\adfLgap %ADFvfyDesignStart
$(348, 0, 1, 2)$,
$(349, 0, 115, 233)$,
$(350, 0, 118, 116)$,
$(351, 0, 232, 347)$,\adfsplit
$(352, 0, 346, 230)$,
$(353, 0, 4, 11)$,
$(354, 0, 7, 344)$,
$(355, 0, 337, 341)$,\adfsplit
$(356, 0, 10, 23)$,
$(357, 0, 13, 338)$,
$(358, 0, 325, 335)$,
$(359, 0, 16, 35)$,\adfsplit
$(360, 0, 19, 332)$,
$(361, 0, 313, 329)$,
$(362, 0, 22, 5)$,
$(363, 0, 331, 326)$,\adfsplit
$(364, 0, 343, 17)$,
$(365, 0, 25, 53)$,
$(366, 0, 28, 323)$,
$(367, 0, 295, 320)$,\adfsplit
$(368, 0, 31, 65)$,
$(369, 0, 34, 317)$,
$(370, 0, 283, 314)$,
$(90, 123, 24, 253)$,\adfsplit
$(228, 63, 66, 115)$,
$(35, 126, 185, 15)$,
$(330, 143, 192, 218)$,
$(63, 216, 287, 19)$,\adfsplit
$(147, 162, 259, 335)$,
$(318, 220, 137, 288)$,
$(106, 126, 231, 216)$,
$(234, 167, 117, 296)$,\adfsplit
$(275, 116, 84, 153)$,
$(307, 44, 194, 136)$,
$(12, 146, 322, 152)$,
$(298, 218, 51, 69)$,\adfsplit
$(235, 278, 157, 69)$,
$(297, 289, 80, 208)$,
$(122, 335, 130, 48)$,
$(328, 134, 227, 125)$,\adfsplit
$(311, 308, 256, 96)$,
$(211, 280, 155, 217)$,
$(210, 334, 168, 80)$,
$(231, 276, 59, 258)$,\adfsplit
$(257, 109, 206, 51)$,
$(68, 223, 330, 21)$,
$(63, 297, 31, 120)$,
$(104, 174, 307, 269)$,\adfsplit
$(274, 40, 324, 49)$,
$(332, 223, 61, 132)$,
$(0, 27, 196, 271)$,
$(0, 21, 76, 222)$,\adfsplit
$(0, 57, 78, 184)$,
$(0, 46, 157, 315)$,
$(0, 40, 246, 307)$,
$(0, 45, 135, 319)$,\adfsplit
$(0, 14, 43, 141)$,
$(0, 44, 129, 255)$,
$(0, 47, 187, 292)$,
$(0, 83, 122, 297)$,\adfsplit
$(0, 117, 245, 311)$,
$(0, 54, 121, 244)$,
$(0, 63, 109, 201)$,
$(0, 41, 100, 253)$,\adfsplit
$(0, 79, 189, 275)$,
$(0, 143, 185, 249)$,
$(0, 107, 137, 207)$,
$(1, 15, 41, 295)$,\adfsplit
$(0, 87, 174, 261)$

%ADFvfyBlocksEnd
\adfLgap \noindent by the mapping:
$x \mapsto x + 3 j \adfmod{348}$ for $x < 348$,
$x \mapsto x$ for $x \ge 348$,
$0 \le j < 116$
 for the first 23 blocks;
$x \mapsto x + 2 j \adfmod{348}$ for $x < 348$,
$x \mapsto x$ for $x \ge 348$,
$0 \le j < 174$
 for the next 45 blocks,
$0 \le j < 87$
 for the last block.
\ADFvfyParStart{(371, ((23, 116, ((348, 3), (23, 23))), (45, 174, ((348, 2), (23, 23))), (1, 87, ((348, 2), (23, 23)))), ((29, 12), (23, 1)))} %ADFvfyParEnd
% End of 29^12 23^1
%%%%%%%%%%%%%%%%%%%%%%%%%%%%%%%%%%%%%%%%%%%%%%%%%%%%%%%%%%%%%%%%%%%%%%%%%%%%%%%%%%%%%%%%%%
%%%%%%%%%%%%%%%%%%%%%%%%%%%%%%%%%%%%%%%%%%%%%%%%%%%%%%%%%%%%%%%%%%%%%%%%%%%%%%%%%%%%%%%%%%

% Charlotte:GDD4-1-3-5-mod-6-TeX-gen-A:HITS-fun:4.10
\adfDgap
%ADFvfyBlocksStart {29,29,29,29,29,29,29,29,29,29,29,29,26}
\noindent{\boldmath $ 29^{12} 26^{1} $}~
With the point set $Z_{374}$ partitioned into
 residue classes modulo $12$ for $\{0, 1, \dots, 347\}$, and
 $\{348, 349, \dots, 373\}$,
 the design is generated from

\adfLgap %ADFvfyDesignStart
$(348, 0, 1, 2)$,
$(349, 0, 115, 233)$,
$(350, 0, 118, 116)$,
$(351, 0, 232, 347)$,\adfsplit
$(352, 0, 346, 230)$,
$(353, 0, 4, 11)$,
$(354, 0, 7, 344)$,
$(355, 0, 337, 341)$,\adfsplit
$(356, 0, 10, 23)$,
$(357, 0, 13, 338)$,
$(358, 0, 325, 335)$,
$(359, 0, 16, 35)$,\adfsplit
$(360, 0, 19, 332)$,
$(361, 0, 313, 329)$,
$(362, 0, 22, 5)$,
$(363, 0, 331, 326)$,\adfsplit
$(364, 0, 343, 17)$,
$(365, 0, 25, 53)$,
$(366, 0, 28, 323)$,
$(367, 0, 295, 320)$,\adfsplit
$(368, 0, 31, 65)$,
$(369, 0, 34, 317)$,
$(370, 0, 283, 314)$,
$(371, 0, 37, 8)$,\adfsplit
$(372, 0, 319, 311)$,
$(373, 0, 340, 29)$,
$(259, 136, 314, 154)$,
$(62, 188, 226, 35)$,\adfsplit
$(154, 96, 337, 194)$,
$(273, 140, 72, 173)$,
$(249, 206, 197, 95)$,
$(205, 254, 343, 34)$,\adfsplit
$(125, 215, 284, 105)$,
$(129, 264, 337, 185)$,
$(57, 268, 218, 339)$,
$(327, 55, 161, 286)$,\adfsplit
$(115, 68, 182, 24)$,
$(0, 3, 80, 86)$,
$(0, 14, 95, 226)$,
$(0, 45, 119, 218)$,\adfsplit
$(0, 57, 127, 219)$,
$(0, 26, 85, 198)$,
$(0, 32, 94, 239)$,
$(0, 64, 142, 266)$,\adfsplit
$(0, 21, 75, 214)$,
$(0, 30, 93, 197)$,
$(0, 42, 88, 287)$,
$(0, 15, 112, 163)$,\adfsplit
$(0, 87, 174, 261)$

%ADFvfyBlocksEnd
\adfLgap \noindent by the mapping:
$x \mapsto x + 3 j \adfmod{348}$ for $x < 348$,
$x \mapsto x$ for $x \ge 348$,
$0 \le j < 116$
 for the first 26 blocks;
$x \mapsto x +  j \adfmod{348}$ for $x < 348$,
$x \mapsto x$ for $x \ge 348$,
$0 \le j < 348$
 for the next 22 blocks,
$0 \le j < 87$
 for the last block.
\ADFvfyParStart{(374, ((26, 116, ((348, 3), (26, 26))), (22, 348, ((348, 1), (26, 26))), (1, 87, ((348, 1), (26, 26)))), ((29, 12), (26, 1)))} %ADFvfyParEnd
% End of 29^12 26^1
%%%%%%%%%%%%%%%%%%%%%%%%%%%%%%%%%%%%%%%%%%%%%%%%%%%%%%%%%%%%%%%%%%%%%%%%%%%%%%%%%%%%%%%%%%
%%%%%%%%%%%%%%%%%%%%%%%%%%%%%%%%%%%%%%%%%%%%%%%%%%%%%%%%%%%%%%%%%%%%%%%%%%%%%%%%%%%%%%%%%%

% Charlotte:GDD4-1-3-5-mod-6-TeX-gen-A:HITS-fun:4.10
\adfDgap
%ADFvfyBlocksStart {29,29,29,29,29,29,29,29,29,29,29,29,29,29,29,29,29,29,29,29,29,29,29,29,26}
\noindent{\boldmath $ 29^{24} 26^{1} $}~
With the point set $Z_{722}$ partitioned into
 residue classes modulo $24$ for $\{0, 1, \dots, 695\}$, and
 $\{696, 697, \dots, 721\}$,
 the design is generated from

\adfLgap %ADFvfyDesignStart
$(696, 0, 1, 2)$,
$(697, 0, 232, 230)$,
$(698, 0, 463, 695)$,
$(699, 0, 466, 233)$,\adfsplit
$(700, 0, 694, 464)$,
$(701, 0, 4, 11)$,
$(702, 0, 7, 692)$,
$(703, 0, 685, 689)$,\adfsplit
$(704, 0, 10, 23)$,
$(705, 0, 13, 686)$,
$(706, 0, 673, 683)$,
$(707, 0, 16, 35)$,\adfsplit
$(708, 0, 19, 680)$,
$(709, 0, 661, 677)$,
$(710, 0, 22, 5)$,
$(711, 0, 679, 674)$,\adfsplit
$(712, 0, 691, 17)$,
$(713, 0, 25, 53)$,
$(714, 0, 28, 671)$,
$(715, 0, 643, 668)$,\adfsplit
$(716, 0, 31, 65)$,
$(717, 0, 34, 665)$,
$(718, 0, 631, 662)$,
$(719, 0, 37, 8)$,\adfsplit
$(720, 0, 667, 659)$,
$(721, 0, 688, 29)$,
$(87, 310, 601, 424)$,
$(581, 231, 165, 368)$,\adfsplit
$(383, 31, 478, 569)$,
$(62, 425, 7, 496)$,
$(627, 21, 362, 406)$,
$(127, 1, 461, 630)$,\adfsplit
$(327, 250, 524, 157)$,
$(84, 176, 3, 223)$,
$(41, 356, 232, 675)$,
$(464, 410, 347, 678)$,\adfsplit
$(159, 607, 491, 440)$,
$(226, 414, 609, 114)$,
$(298, 166, 32, 23)$,
$(371, 128, 332, 573)$,\adfsplit
$(573, 629, 85, 659)$,
$(346, 28, 691, 352)$,
$(455, 9, 220, 487)$,
$(194, 598, 149, 321)$,\adfsplit
$(522, 205, 208, 70)$,
$(71, 685, 89, 173)$,
$(172, 78, 552, 578)$,
$(668, 458, 567, 157)$,\adfsplit
$(63, 131, 262, 533)$,
$(154, 234, 71, 196)$,
$(111, 338, 311, 660)$,
$(689, 480, 323, 208)$,\adfsplit
$(392, 542, 88, 313)$,
$(224, 377, 522, 328)$,
$(455, 419, 28, 146)$,
$(91, 76, 16, 419)$,\adfsplit
$(117, 488, 9, 378)$,
$(0, 12, 97, 171)$,
$(0, 14, 87, 175)$,
$(0, 20, 61, 307)$,\adfsplit
$(0, 21, 64, 320)$,
$(0, 40, 146, 259)$,
$(0, 57, 156, 589)$,
$(0, 89, 254, 530)$,\adfsplit
$(0, 70, 260, 393)$,
$(0, 105, 245, 400)$,
$(0, 59, 257, 399)$,
$(0, 69, 180, 321)$,\adfsplit
$(0, 121, 302, 491)$,
$(0, 119, 270, 457)$,
$(0, 46, 149, 566)$,
$(0, 58, 136, 342)$,\adfsplit
$(0, 52, 231, 459)$,
$(0, 50, 148, 308)$,
$(0, 123, 306, 484)$,
$(0, 128, 282, 411)$,\adfsplit
$(0, 76, 238, 462)$,
$(0, 174, 348, 522)$

%ADFvfyBlocksEnd
\adfLgap \noindent by the mapping:
$x \mapsto x + 3 j \adfmod{696}$ for $x < 696$,
$x \mapsto x$ for $x \ge 696$,
$0 \le j < 232$
 for the first 26 blocks;
$x \mapsto x +  j \adfmod{696}$ for $x < 696$,
$x \mapsto x$ for $x \ge 696$,
$0 \le j < 696$
 for the next 51 blocks,
$0 \le j < 174$
 for the last block.
\ADFvfyParStart{(722, ((26, 232, ((696, 3), (26, 26))), (51, 696, ((696, 1), (26, 26))), (1, 174, ((696, 1), (26, 26)))), ((29, 24), (26, 1)))} %ADFvfyParEnd
% End of 29^24 26^1
%%%%%%%%%%%%%%%%%%%%%%%%%%%%%%%%%%%%%%%%%%%%%%%%%%%%%%%%%%%%%%%%%%%%%%%%%%%%%%%%%%%%%%%%%%
%%%%%%%%%%%%%%%%%%%%%%%%%%%%%%%%%%%%%%%%%%%%%%%%%%%%%%%%%%%%%%%%%%%%%%%%%%%%%%%%%%%%%%%%%%

% Charlotte:GDD4-1-3-5-mod-6-TeX-gen-A:HITS-fun:4.10
\adfDgap
%ADFvfyBlocksStart {29,29,29,29,29,29,29,29,29,29,29,29,29,29,29,17}
\noindent{\boldmath $ 29^{15} 17^{1} $}~
With the point set $Z_{452}$ partitioned into
 residue classes modulo $15$ for $\{0, 1, \dots, 434\}$, and
 $\{435, 436, \dots, 451\}$,
 the design is generated from

\adfLgap %ADFvfyDesignStart
$(435, 0, 1, 2)$,
$(436, 0, 145, 143)$,
$(437, 0, 289, 434)$,
$(438, 0, 292, 146)$,\adfsplit
$(439, 0, 433, 290)$,
$(440, 0, 4, 11)$,
$(441, 0, 7, 431)$,
$(442, 0, 424, 428)$,\adfsplit
$(443, 0, 10, 23)$,
$(444, 0, 13, 425)$,
$(445, 0, 412, 422)$,
$(446, 0, 16, 35)$,\adfsplit
$(447, 0, 19, 419)$,
$(448, 0, 400, 416)$,
$(449, 0, 22, 5)$,
$(450, 0, 418, 413)$,\adfsplit
$(451, 0, 430, 17)$,
$(395, 193, 224, 321)$,
$(149, 222, 272, 429)$,
$(57, 128, 51, 301)$,\adfsplit
$(22, 376, 265, 162)$,
$(12, 277, 346, 317)$,
$(247, 78, 182, 70)$,
$(171, 40, 266, 238)$,\adfsplit
$(109, 362, 338, 91)$,
$(15, 177, 163, 64)$,
$(391, 190, 202, 48)$,
$(407, 176, 29, 354)$,\adfsplit
$(71, 367, 159, 180)$,
$(201, 234, 68, 176)$,
$(41, 238, 202, 324)$,
$(176, 39, 411, 350)$,\adfsplit
$(205, 116, 368, 196)$,
$(153, 0, 106, 27)$,
$(13, 236, 321, 104)$,
$(85, 281, 203, 82)$,\adfsplit
$(0, 20, 59, 96)$,
$(0, 46, 98, 291)$,
$(0, 38, 93, 306)$,
$(0, 34, 117, 159)$,\adfsplit
$(0, 32, 94, 353)$,
$(0, 48, 102, 251)$,
$(0, 64, 136, 220)$,
$(0, 43, 158, 224)$,\adfsplit
$(0, 51, 107, 348)$,
$(0, 68, 168, 284)$,
$(0, 41, 160, 301)$,
$(0, 26, 70, 256)$

%ADFvfyBlocksEnd
\adfLgap \noindent by the mapping:
$x \mapsto x + 3 j \adfmod{435}$ for $x < 435$,
$x \mapsto x$ for $x \ge 435$,
$0 \le j < 145$
 for the first 17 blocks;
$x \mapsto x +  j \adfmod{435}$ for $x < 435$,
$x \mapsto x$ for $x \ge 435$,
$0 \le j < 435$
 for the last 31 blocks.
\ADFvfyParStart{(452, ((17, 145, ((435, 3), (17, 17))), (31, 435, ((435, 1), (17, 17)))), ((29, 15), (17, 1)))} %ADFvfyParEnd
% End of 29^15 17^1
%%%%%%%%%%%%%%%%%%%%%%%%%%%%%%%%%%%%%%%%%%%%%%%%%%%%%%%%%%%%%%%%%%%%%%%%%%%%%%%%%%%%%%%%%%
%%%%%%%%%%%%%%%%%%%%%%%%%%%%%%%%%%%%%%%%%%%%%%%%%%%%%%%%%%%%%%%%%%%%%%%%%%%%%%%%%%%%%%%%%%

% Charlotte:GDD4-1-3-5-mod-6-TeX-gen-A:HITS-fun:4.10
\adfDgap
%ADFvfyBlocksStart {29,29,29,29,29,29,29,29,29,29,29,29,29,29,29,23}
\noindent{\boldmath $ 29^{15} 23^{1} $}~
With the point set $Z_{458}$ partitioned into
 residue classes modulo $15$ for $\{0, 1, \dots, 434\}$, and
 $\{435, 436, \dots, 457\}$,
 the design is generated from

\adfLgap %ADFvfyDesignStart
$(435, 0, 1, 2)$,
$(436, 0, 145, 143)$,
$(437, 0, 289, 434)$,
$(438, 0, 292, 146)$,\adfsplit
$(439, 0, 433, 290)$,
$(440, 0, 4, 11)$,
$(441, 0, 7, 431)$,
$(442, 0, 424, 428)$,\adfsplit
$(443, 0, 10, 23)$,
$(444, 0, 13, 425)$,
$(445, 0, 412, 422)$,
$(446, 0, 16, 35)$,\adfsplit
$(447, 0, 19, 419)$,
$(448, 0, 400, 416)$,
$(449, 0, 22, 5)$,
$(450, 0, 418, 413)$,\adfsplit
$(451, 0, 430, 17)$,
$(452, 0, 25, 53)$,
$(453, 0, 28, 410)$,
$(454, 0, 382, 407)$,\adfsplit
$(455, 0, 31, 65)$,
$(456, 0, 34, 404)$,
$(457, 0, 370, 401)$,
$(78, 416, 69, 117)$,\adfsplit
$(244, 164, 31, 396)$,
$(403, 179, 337, 57)$,
$(101, 411, 245, 288)$,
$(175, 273, 88, 345)$,\adfsplit
$(21, 127, 347, 188)$,
$(328, 199, 380, 217)$,
$(272, 411, 43, 217)$,
$(11, 339, 23, 80)$,\adfsplit
$(298, 81, 269, 295)$,
$(3, 40, 362, 194)$,
$(166, 73, 335, 409)$,
$(156, 54, 338, 358)$,\adfsplit
$(356, 12, 383, 62)$,
$(84, 33, 263, 180)$,
$(175, 35, 167, 135)$,
$(184, 240, 76, 396)$,\adfsplit
$(34, 150, 292, 144)$,
$(0, 14, 68, 376)$,
$(0, 24, 82, 275)$,
$(0, 62, 200, 263)$,\adfsplit
$(0, 21, 124, 273)$,
$(0, 42, 134, 228)$,
$(0, 41, 85, 171)$,
$(0, 36, 117, 245)$,\adfsplit
$(0, 47, 126, 204)$,
$(0, 95, 196, 314)$,
$(0, 38, 84, 237)$,
$(0, 33, 137, 208)$,\adfsplit
$(0, 49, 161, 238)$

%ADFvfyBlocksEnd
\adfLgap \noindent by the mapping:
$x \mapsto x + 3 j \adfmod{435}$ for $x < 435$,
$x \mapsto x$ for $x \ge 435$,
$0 \le j < 145$
 for the first 23 blocks;
$x \mapsto x +  j \adfmod{435}$ for $x < 435$,
$x \mapsto x$ for $x \ge 435$,
$0 \le j < 435$
 for the last 30 blocks.
\ADFvfyParStart{(458, ((23, 145, ((435, 3), (23, 23))), (30, 435, ((435, 1), (23, 23)))), ((29, 15), (23, 1)))} %ADFvfyParEnd
% End of 29^15 23^1
%%%%%%%%%%%%%%%%%%%%%%%%%%%%%%%%%%%%%%%%%%%%%%%%%%%%%%%%%%%%%%%%%%%%%%%%%%%%%%%%%%%%%%%%%%
%%%%%%%%%%%%%%%%%%%%%%%%%%%%%%%%%%%%%%%%%%%%%%%%%%%%%%%%%%%%%%%%%%%%%%%%%%%%%%%%%%%%%%%%%%

% Charlotte:GDD4-1-3-5-mod-6-TeX-gen-A:HITS-fun:4.10
\adfDgap
%ADFvfyBlocksStart {29,29,29,29,29,29,29,29,29,14}
\noindent{\boldmath $ 29^{9} 14^{1} $}~
With the point set $Z_{275}$ partitioned into
 residue classes modulo $9$ for $\{0, 1, \dots, 260\}$, and
 $\{261, 262, \dots, 274\}$,
 the design is generated from

\adfLgap %ADFvfyDesignStart
$(273, 0, 1, 2)$,
$(274, 0, 88, 176)$,
$(261, 45, 253, 260)$,
$(262, 11, 255, 112)$,\adfsplit
$(263, 172, 56, 3)$,
$(264, 44, 51, 205)$,
$(222, 256, 55, 212)$,
$(260, 60, 97, 239)$,\adfsplit
$(60, 158, 197, 91)$,
$(72, 100, 247, 142)$,
$(63, 247, 106, 228)$,
$(92, 94, 224, 232)$,\adfsplit
$(0, 3, 68, 69)$,
$(0, 6, 22, 235)$,
$(0, 8, 25, 83)$,
$(0, 33, 84, 124)$,\adfsplit
$(0, 26, 78, 151)$,
$(0, 24, 111, 172)$,
$(0, 46, 93, 205)$,
$(0, 52, 109, 185)$,\adfsplit
$(0, 14, 44, 82)$,
$(0, 12, 86, 115)$,
$(0, 15, 38, 202)$

%ADFvfyBlocksEnd
\adfLgap \noindent by the mapping:
$x \mapsto x \oplus (3 j)$ for $x < 261$,
$x \mapsto (x - 261 + 4 j \adfmod{12}) + 261$ for $261 \le x < 273$,
$x \mapsto x$ for $x \ge 273$,
$0 \le j < 87$
 for the first two blocks;
$x \mapsto x \oplus j$ for $x < 261$,
$x \mapsto (x - 261 + 4 j \adfmod{12}) + 261$ for $261 \le x < 273$,
$x \mapsto x$ for $x \ge 273$,
$0 \le j < 261$
 for the last 21 blocks.
\ADFvfyParStart{(275, ((2, 87, ((261, 3, (87, 3)), (12, 4), (2, 2))), (21, 261, ((261, 1, (87, 3)), (12, 4), (2, 2)))), ((29, 9), (14, 1)))} %ADFvfyParEnd
% End of 29^9 14^1
%%%%%%%%%%%%%%%%%%%%%%%%%%%%%%%%%%%%%%%%%%%%%%%%%%%%%%%%%%%%%%%%%%%%%%%%%%%%%%%%%%%%%%%%%%
%%%%%%%%%%%%%%%%%%%%%%%%%%%%%%%%%%%%%%%%%%%%%%%%%%%%%%%%%%%%%%%%%%%%%%%%%%%%%%%%%%%%%%%%%%

% Charlotte:GDD4-1-3-5-mod-6-TeX-gen-A:HITS-fun:4.10
\adfDgap
%ADFvfyBlocksStart {29,29,29,29,29,29,29,29,29,20}
\noindent{\boldmath $ 29^{9} 20^{1} $}~
With the point set $Z_{281}$ partitioned into
 residue classes modulo $9$ for $\{0, 1, \dots, 260\}$, and
 $\{261, 262, \dots, 280\}$,
 the design is generated from

\adfLgap %ADFvfyDesignStart
$(279, 0, 1, 2)$,
$(280, 0, 88, 176)$,
$(261, 179, 123, 22)$,
$(262, 122, 132, 10)$,\adfsplit
$(263, 217, 29, 207)$,
$(264, 87, 217, 92)$,
$(265, 57, 184, 125)$,
$(266, 169, 162, 161)$,\adfsplit
$(257, 202, 196, 18)$,
$(52, 137, 163, 67)$,
$(56, 9, 85, 177)$,
$(204, 58, 237, 97)$,\adfsplit
$(0, 3, 19, 24)$,
$(0, 12, 42, 158)$,
$(0, 14, 28, 75)$,
$(0, 26, 92, 156)$,\adfsplit
$(0, 41, 91, 193)$,
$(0, 22, 59, 100)$,
$(0, 35, 119, 174)$,
$(0, 31, 65, 213)$,\adfsplit
$(0, 20, 77, 120)$,
$(0, 38, 98, 167)$,
$(0, 23, 74, 114)$,
$(0, 40, 89, 138)$

%ADFvfyBlocksEnd
\adfLgap \noindent by the mapping:
$x \mapsto x \oplus (3 j)$ for $x < 261$,
$x \mapsto (x - 261 + 6 j \adfmod{18}) + 261$ for $261 \le x < 279$,
$x \mapsto x$ for $x \ge 279$,
$0 \le j < 87$
 for the first two blocks;
$x \mapsto x \oplus j$ for $x < 261$,
$x \mapsto (x - 261 + 6 j \adfmod{18}) + 261$ for $261 \le x < 279$,
$x \mapsto x$ for $x \ge 279$,
$0 \le j < 261$
 for the last 22 blocks.
\ADFvfyParStart{(281, ((2, 87, ((261, 3, (87, 3)), (18, 6), (2, 2))), (22, 261, ((261, 1, (87, 3)), (18, 6), (2, 2)))), ((29, 9), (20, 1)))} %ADFvfyParEnd
% End of 29^9 20^1
%%%%%%%%%%%%%%%%%%%%%%%%%%%%%%%%%%%%%%%%%%%%%%%%%%%%%%%%%%%%%%%%%%%%%%%%%%%%%%%%%%%%%%%%%%
%%%%%%%%%%%%%%%%%%%%%%%%%%%%%%%%%%%%%%%%%%%%%%%%%%%%%%%%%%%%%%%%%%%%%%%%%%%%%%%%%%%%%%%%%%

% Charlotte:GDD4-1-3-5-mod-6-TeX-gen-A:HITS-fun:4.10
\adfDgap
%ADFvfyBlocksStart {29,29,29,29,29,29,29,29,29,26}
\noindent{\boldmath $ 29^{9} 26^{1} $}~
With the point set $Z_{287}$ partitioned into
 residue classes modulo $9$ for $\{0, 1, \dots, 260\}$, and
 $\{261, 262, \dots, 286\}$,
 the design is generated from

\adfLgap %ADFvfyDesignStart
$(285, 0, 1, 2)$,
$(286, 0, 88, 176)$,
$(261, 25, 225, 32)$,
$(262, 242, 67, 237)$,\adfsplit
$(263, 146, 82, 87)$,
$(264, 17, 138, 256)$,
$(265, 5, 78, 136)$,
$(266, 28, 200, 69)$,\adfsplit
$(267, 155, 121, 42)$,
$(268, 140, 165, 178)$,
$(209, 140, 223, 217)$,
$(72, 200, 129, 154)$,\adfsplit
$(112, 226, 186, 81)$,
$(0, 3, 23, 24)$,
$(0, 10, 26, 186)$,
$(0, 12, 78, 107)$,\adfsplit
$(0, 14, 30, 155)$,
$(0, 35, 74, 167)$,
$(0, 43, 110, 177)$,
$(0, 37, 87, 160)$,\adfsplit
$(0, 38, 98, 150)$,
$(0, 42, 94, 142)$,
$(0, 15, 62, 158)$,
$(0, 49, 102, 152)$,\adfsplit
$(0, 32, 65, 210)$

%ADFvfyBlocksEnd
\adfLgap \noindent by the mapping:
$x \mapsto x \oplus (3 j)$ for $x < 261$,
$x \mapsto (x - 261 + 8 j \adfmod{24}) + 261$ for $261 \le x < 285$,
$x \mapsto x$ for $x \ge 285$,
$0 \le j < 87$
 for the first two blocks;
$x \mapsto x \oplus j$ for $x < 261$,
$x \mapsto (x - 261 + 8 j \adfmod{24}) + 261$ for $261 \le x < 285$,
$x \mapsto x$ for $x \ge 285$,
$0 \le j < 261$
 for the last 23 blocks.
\ADFvfyParStart{(287, ((2, 87, ((261, 3, (87, 3)), (24, 8), (2, 2))), (23, 261, ((261, 1, (87, 3)), (24, 8), (2, 2)))), ((29, 9), (26, 1)))} %ADFvfyParEnd
% End of 29^9 26^1
%%%%%%%%%%%%%%%%%%%%%%%%%%%%%%%%%%%%%%%%%%%%%%%%%%%%%%%%%%%%%%%%%%%%%%%%%%%%%%%%%%%%%%%%%%
%%%%%%%%%%%%%%%%%%%%%%%%%%%%%%%%%%%%%%%%%%%%%%%%%%%%%%%%%%%%%%%%%%%%%%%%%%%%%%%%%%%%%%%%%%

% Charlotte:GDD4-1-3-5-mod-6-TeX-gen-A:HITS-fun:4.10
\adfDgap
%ADFvfyBlocksStart {29,29,29,29,29,29,29,29,29,32}
\noindent{\boldmath $ 29^{9} 32^{1} $}~
With the point set $Z_{293}$ partitioned into
 residue classes modulo $9$ for $\{0, 1, \dots, 260\}$, and
 $\{261, 262, \dots, 292\}$,
 the design is generated from

\adfLgap %ADFvfyDesignStart
$(291, 0, 1, 2)$,
$(292, 0, 88, 176)$,
$(261, 213, 1, 56)$,
$(262, 162, 128, 136)$,\adfsplit
$(263, 83, 156, 175)$,
$(264, 176, 31, 27)$,
$(265, 40, 195, 170)$,
$(266, 230, 22, 258)$,\adfsplit
$(267, 47, 96, 169)$,
$(268, 227, 211, 6)$,
$(269, 188, 138, 247)$,
$(270, 117, 239, 127)$,\adfsplit
$(153, 33, 253, 211)$,
$(0, 3, 8, 146)$,
$(0, 6, 20, 201)$,
$(0, 7, 75, 87)$,\adfsplit
$(0, 13, 92, 107)$,
$(0, 21, 98, 137)$,
$(0, 22, 57, 124)$,
$(0, 33, 74, 136)$,\adfsplit
$(0, 46, 96, 210)$,
$(0, 35, 93, 163)$,
$(0, 32, 69, 182)$,
$(0, 23, 47, 179)$,\adfsplit
$(0, 26, 91, 177)$,
$(0, 17, 48, 78)$

%ADFvfyBlocksEnd
\adfLgap \noindent by the mapping:
$x \mapsto x \oplus (3 j)$ for $x < 261$,
$x \mapsto (x - 261 + 10 j \adfmod{30}) + 261$ for $261 \le x < 291$,
$x \mapsto x$ for $x \ge 291$,
$0 \le j < 87$
 for the first two blocks;
$x \mapsto x \oplus j$ for $x < 261$,
$x \mapsto (x - 261 + 10 j \adfmod{30}) + 261$ for $261 \le x < 291$,
$x \mapsto x$ for $x \ge 291$,
$0 \le j < 261$
 for the last 24 blocks.
\ADFvfyParStart{(293, ((2, 87, ((261, 3, (87, 3)), (30, 10), (2, 2))), (24, 261, ((261, 1, (87, 3)), (30, 10), (2, 2)))), ((29, 9), (32, 1)))} %ADFvfyParEnd
% End of 29^9 32^1
%%%%%%%%%%%%%%%%%%%%%%%%%%%%%%%%%%%%%%%%%%%%%%%%%%%%%%%%%%%%%%%%%%%%%%%%%%%%%%%%%%%%%%%%%%
%%%%%%%%%%%%%%%%%%%%%%%%%%%%%%%%%%%%%%%%%%%%%%%%%%%%%%%%%%%%%%%%%%%%%%%%%%%%%%%%%%%%%%%%%%

% Charlotte:GDD4-1-3-5-mod-6-TeX-gen-A:HITS-fun:4.10
\adfDgap
%ADFvfyBlocksStart {29,29,29,29,29,29,29,29,29,38}
\noindent{\boldmath $ 29^{9} 38^{1} $}~
With the point set $Z_{299}$ partitioned into
 residue classes modulo $9$ for $\{0, 1, \dots, 260\}$, and
 $\{261, 262, \dots, 298\}$,
 the design is generated from

\adfLgap %ADFvfyDesignStart
$(297, 0, 1, 2)$,
$(298, 0, 88, 176)$,
$(261, 250, 234, 47)$,
$(262, 183, 28, 56)$,\adfsplit
$(263, 11, 249, 169)$,
$(264, 68, 31, 126)$,
$(265, 15, 55, 110)$,
$(266, 109, 116, 165)$,\adfsplit
$(267, 249, 82, 203)$,
$(268, 30, 209, 1)$,
$(269, 157, 57, 56)$,
$(270, 39, 22, 164)$,\adfsplit
$(271, 147, 2, 43)$,
$(0, 5, 10, 272)$,
$(0, 3, 65, 71)$,
$(0, 11, 26, 42)$,\adfsplit
$(0, 12, 79, 109)$,
$(0, 13, 57, 123)$,
$(0, 14, 53, 75)$,
$(0, 33, 110, 183)$,\adfsplit
$(0, 21, 43, 145)$,
$(0, 17, 87, 137)$,
$(0, 31, 101, 136)$,
$(0, 24, 93, 173)$,\adfsplit
$(0, 41, 84, 213)$,
$(0, 29, 112, 172)$,
$(0, 35, 114, 165)$

%ADFvfyBlocksEnd
\adfLgap \noindent by the mapping:
$x \mapsto x \oplus (3 j)$ for $x < 261$,
$x \mapsto (x - 261 + 12 j \adfmod{36}) + 261$ for $261 \le x < 297$,
$x \mapsto x$ for $x \ge 297$,
$0 \le j < 87$
 for the first two blocks;
$x \mapsto x \oplus j$ for $x < 261$,
$x \mapsto (x - 261 + 12 j \adfmod{36}) + 261$ for $261 \le x < 297$,
$x \mapsto x$ for $x \ge 297$,
$0 \le j < 261$
 for the last 25 blocks.
\ADFvfyParStart{(299, ((2, 87, ((261, 3, (87, 3)), (36, 12), (2, 2))), (25, 261, ((261, 1, (87, 3)), (36, 12), (2, 2)))), ((29, 9), (38, 1)))} %ADFvfyParEnd
% End of 29^9 38^1
%%%%%%%%%%%%%%%%%%%%%%%%%%%%%%%%%%%%%%%%%%%%%%%%%%%%%%%%%%%%%%%%%%%%%%%%%%%%%%%%%%%%%%%%%%
%%%%%%%%%%%%%%%%%%%%%%%%%%%%%%%%%%%%%%%%%%%%%%%%%%%%%%%%%%%%%%%%%%%%%%%%%%%%%%%%%%%%%%%%%%

% Charlotte:GDD4-1-3-5-mod-6-TeX-gen-A:HITS-fun:4.10
\adfDgap
%ADFvfyBlocksStart {29,29,29,29,29,29,29,29,29,44}
\noindent{\boldmath $ 29^{9} 44^{1} $}~
With the point set $Z_{305}$ partitioned into
 residue classes modulo $9$ for $\{0, 1, \dots, 260\}$, and
 $\{261, 262, \dots, 304\}$,
 the design is generated from

\adfLgap %ADFvfyDesignStart
$(303, 0, 1, 2)$,
$(304, 0, 88, 176)$,
$(261, 47, 81, 1)$,
$(262, 249, 17, 217)$,\adfsplit
$(263, 150, 206, 88)$,
$(264, 140, 40, 132)$,
$(265, 93, 235, 224)$,
$(266, 57, 229, 107)$,\adfsplit
$(267, 133, 161, 117)$,
$(268, 11, 241, 165)$,
$(269, 251, 93, 7)$,
$(270, 206, 207, 142)$,\adfsplit
$(271, 211, 128, 51)$,
$(272, 130, 182, 201)$,
$(0, 5, 13, 273)$,
$(0, 7, 17, 274)$,\adfsplit
$(0, 3, 51, 70)$,
$(0, 6, 91, 111)$,
$(0, 12, 110, 132)$,
$(0, 15, 75, 128)$,\adfsplit
$(0, 24, 66, 121)$,
$(0, 57, 115, 184)$,
$(0, 34, 78, 165)$,
$(0, 21, 62, 168)$,\adfsplit
$(0, 32, 102, 145)$,
$(0, 26, 59, 94)$,
$(0, 30, 79, 182)$,
$(0, 39, 123, 163)$

%ADFvfyBlocksEnd
\adfLgap \noindent by the mapping:
$x \mapsto x \oplus (3 j)$ for $x < 261$,
$x \mapsto (x - 261 + 14 j \adfmod{42}) + 261$ for $261 \le x < 303$,
$x \mapsto x$ for $x \ge 303$,
$0 \le j < 87$
 for the first two blocks;
$x \mapsto x \oplus j$ for $x < 261$,
$x \mapsto (x - 261 + 14 j \adfmod{42}) + 261$ for $261 \le x < 303$,
$x \mapsto x$ for $x \ge 303$,
$0 \le j < 261$
 for the last 26 blocks.
\ADFvfyParStart{(305, ((2, 87, ((261, 3, (87, 3)), (42, 14), (2, 2))), (26, 261, ((261, 1, (87, 3)), (42, 14), (2, 2)))), ((29, 9), (44, 1)))} %ADFvfyParEnd
% End of 29^9 44^1
%%%%%%%%%%%%%%%%%%%%%%%%%%%%%%%%%%%%%%%%%%%%%%%%%%%%%%%%%%%%%%%%%%%%%%%%%%%%%%%%%%%%%%%%%%
%%%%%%%%%%%%%%%%%%%%%%%%%%%%%%%%%%%%%%%%%%%%%%%%%%%%%%%%%%%%%%%%%%%%%%%%%%%%%%%%%%%%%%%%%%

% Charlotte:GDD4-1-3-5-mod-6-TeX-gen-A:HITS-fun:4.10
\adfDgap
%ADFvfyBlocksStart {29,29,29,29,29,29,29,29,29,50}
\noindent{\boldmath $ 29^{9} 50^{1} $}~
With the point set $Z_{311}$ partitioned into
 residue classes modulo $9$ for $\{0, 1, \dots, 260\}$, and
 $\{261, 262, \dots, 310\}$,
 the design is generated from

\adfLgap %ADFvfyDesignStart
$(309, 0, 1, 2)$,
$(310, 0, 88, 176)$,
$(261, 193, 8, 108)$,
$(262, 78, 196, 47)$,\adfsplit
$(263, 26, 81, 124)$,
$(264, 60, 164, 16)$,
$(265, 49, 27, 224)$,
$(266, 232, 90, 86)$,\adfsplit
$(267, 177, 40, 167)$,
$(268, 150, 188, 211)$,
$(269, 132, 142, 173)$,
$(270, 128, 21, 220)$,\adfsplit
$(271, 82, 173, 243)$,
$(272, 69, 86, 214)$,
$(273, 180, 89, 70)$,
$(0, 4, 20, 276)$,\adfsplit
$(0, 5, 25, 274)$,
$(0, 8, 37, 275)$,
$(0, 3, 53, 78)$,
$(0, 6, 21, 35)$,\adfsplit
$(0, 48, 97, 150)$,
$(0, 33, 93, 172)$,
$(0, 11, 52, 98)$,
$(0, 40, 82, 187)$,\adfsplit
$(0, 30, 96, 158)$,
$(0, 24, 83, 134)$,
$(0, 39, 109, 177)$,
$(0, 12, 69, 190)$,\adfsplit
$(0, 55, 120, 184)$

%ADFvfyBlocksEnd
\adfLgap \noindent by the mapping:
$x \mapsto x \oplus (3 j)$ for $x < 261$,
$x \mapsto (x - 261 + 16 j \adfmod{48}) + 261$ for $261 \le x < 309$,
$x \mapsto x$ for $x \ge 309$,
$0 \le j < 87$
 for the first two blocks;
$x \mapsto x \oplus j$ for $x < 261$,
$x \mapsto (x - 261 + 16 j \adfmod{48}) + 261$ for $261 \le x < 309$,
$x \mapsto x$ for $x \ge 309$,
$0 \le j < 261$
 for the last 27 blocks.
\ADFvfyParStart{(311, ((2, 87, ((261, 3, (87, 3)), (48, 16), (2, 2))), (27, 261, ((261, 1, (87, 3)), (48, 16), (2, 2)))), ((29, 9), (50, 1)))} %ADFvfyParEnd
% End of 29^9 50^1
%%%%%%%%%%%%%%%%%%%%%%%%%%%%%%%%%%%%%%%%%%%%%%%%%%%%%%%%%%%%%%%%%%%%%%%%%%%%%%%%%%%%%%%%%%
%%%%%%%%%%%%%%%%%%%%%%%%%%%%%%%%%%%%%%%%%%%%%%%%%%%%%%%%%%%%%%%%%%%%%%%%%%%%%%%%%%%%%%%%%%

% Charlotte:GDD4-1-3-5-mod-6-TeX-gen-A:HITS-fun:4.10
\adfDgap
%ADFvfyBlocksStart {29,29,29,29,29,29,29,29,29,56}
\noindent{\boldmath $ 29^{9} 56^{1} $}~
With the point set $Z_{317}$ partitioned into
 residue classes modulo $9$ for $\{0, 1, \dots, 260\}$, and
 $\{261, 262, \dots, 316\}$,
 the design is generated from

\adfLgap %ADFvfyDesignStart
$(315, 0, 1, 2)$,
$(316, 0, 88, 176)$,
$(261, 233, 184, 156)$,
$(262, 171, 169, 14)$,\adfsplit
$(263, 70, 29, 99)$,
$(264, 72, 91, 80)$,
$(265, 243, 76, 20)$,
$(266, 163, 200, 108)$,\adfsplit
$(267, 213, 136, 74)$,
$(268, 9, 230, 205)$,
$(269, 20, 142, 126)$,
$(270, 18, 109, 71)$,\adfsplit
$(271, 188, 81, 73)$,
$(272, 141, 124, 251)$,
$(0, 4, 17, 278)$,
$(0, 7, 29, 273)$,\adfsplit
$(0, 10, 50, 276)$,
$(0, 23, 46, 275)$,
$(0, 31, 83, 277)$,
$(0, 35, 103, 274)$,\adfsplit
$(0, 6, 85, 129)$,
$(0, 3, 67, 150)$,
$(0, 34, 112, 165)$,
$(0, 39, 101, 159)$,\adfsplit
$(0, 15, 84, 182)$,
$(0, 24, 66, 140)$,
$(0, 48, 105, 194)$,
$(0, 21, 51, 164)$,\adfsplit
$(0, 33, 93, 169)$,
$(0, 12, 87, 145)$

%ADFvfyBlocksEnd
\adfLgap \noindent by the mapping:
$x \mapsto x \oplus (3 j)$ for $x < 261$,
$x \mapsto (x - 261 + 18 j \adfmod{54}) + 261$ for $261 \le x < 315$,
$x \mapsto x$ for $x \ge 315$,
$0 \le j < 87$
 for the first two blocks;
$x \mapsto x \oplus j$ for $x < 261$,
$x \mapsto (x - 261 + 18 j \adfmod{54}) + 261$ for $261 \le x < 315$,
$x \mapsto x$ for $x \ge 315$,
$0 \le j < 261$
 for the last 28 blocks.
\ADFvfyParStart{(317, ((2, 87, ((261, 3, (87, 3)), (54, 18), (2, 2))), (28, 261, ((261, 1, (87, 3)), (54, 18), (2, 2)))), ((29, 9), (56, 1)))} %ADFvfyParEnd
% End of 29^9 56^1
%%%%%%%%%%%%%%%%%%%%%%%%%%%%%%%%%%%%%%%%%%%%%%%%%%%%%%%%%%%%%%%%%%%%%%%%%%%%%%%%%%%%%%%%%%
%%%%%%%%%%%%%%%%%%%%%%%%%%%%%%%%%%%%%%%%%%%%%%%%%%%%%%%%%%%%%%%%%%%%%%%%%%%%%%%%%%%%%%%%%%

% Charlotte:GDD4-1-3-5-mod-6-TeX-gen-A:HITS-fun:4.10
\adfDgap
%ADFvfyBlocksStart {29,29,29,29,29,29,29,29,29,62}
\noindent{\boldmath $ 29^{9} 62^{1} $}~
With the point set $Z_{323}$ partitioned into
 residue classes modulo $9$ for $\{0, 1, \dots, 260\}$, and
 $\{261, 262, \dots, 322\}$,
 the design is generated from

\adfLgap %ADFvfyDesignStart
$(321, 146, 165, 31)$,
$(322, 234, 209, 103)$,
$(261, 13, 56, 216)$,
$(261, 33, 19, 3)$,\adfsplit
$(261, 113, 8, 106)$,
$(262, 152, 79, 190)$,
$(262, 164, 201, 225)$,
$(262, 186, 4, 32)$,\adfsplit
$(263, 248, 214, 47)$,
$(263, 183, 27, 55)$,
$(263, 260, 114, 76)$,
$(264, 153, 246, 175)$,\adfsplit
$(264, 105, 143, 253)$,
$(264, 239, 115, 110)$,
$(265, 92, 89, 88)$,
$(265, 193, 159, 156)$,\adfsplit
$(265, 73, 158, 171)$,
$(266, 253, 25, 21)$,
$(266, 38, 123, 252)$,
$(266, 32, 98, 121)$,\adfsplit
$(267, 135, 161, 200)$,
$(267, 194, 259, 159)$,
$(267, 57, 235, 31)$,
$(268, 86, 12, 172)$,\adfsplit
$(268, 119, 197, 238)$,
$(268, 135, 87, 43)$,
$(269, 6, 116, 187)$,
$(269, 254, 94, 46)$,\adfsplit
$(269, 255, 108, 158)$,
$(270, 8, 12, 100)$,
$(270, 38, 121, 140)$,
$(270, 16, 180, 186)$,\adfsplit
$(271, 35, 57, 198)$,
$(271, 202, 146, 124)$,
$(271, 109, 248, 96)$,
$(272, 192, 157, 181)$,\adfsplit
$(272, 239, 216, 206)$,
$(272, 164, 196, 87)$,
$(273, 206, 58, 174)$,
$(273, 105, 18, 43)$,\adfsplit
$(273, 128, 158, 208)$,
$(274, 16, 33, 47)$,
$(274, 193, 126, 242)$,
$(274, 165, 194, 118)$,\adfsplit
$(275, 78, 174, 235)$,
$(275, 7, 9, 58)$,
$(275, 80, 68, 164)$,
$(276, 164, 149, 27)$,\adfsplit
$(276, 213, 22, 12)$,
$(276, 109, 80, 106)$,
$(277, 89, 28, 69)$,
$(277, 198, 131, 31)$,\adfsplit
$(277, 79, 147, 11)$,
$(278, 32, 67, 189)$,
$(278, 213, 88, 170)$,
$(278, 109, 155, 183)$,\adfsplit
$(279, 58, 160, 44)$,
$(279, 3, 59, 54)$,
$(279, 155, 145, 114)$,
$(280, 141, 83, 210)$,\adfsplit
$(280, 167, 79, 13)$,
$(280, 136, 89, 144)$,
$(76, 81, 89, 253)$,
$(43, 64, 230, 103)$,\adfsplit
$(0, 1, 106, 121)$,
$(0, 2, 64, 211)$,
$(0, 7, 76, 222)$,
$(0, 11, 202, 208)$,\adfsplit
$(0, 17, 142, 238)$,
$(0, 43, 55, 229)$,
$(1, 122, 130, 260)$,
$(0, 19, 84, 195)$,\adfsplit
$(0, 71, 92, 241)$,
$(0, 52, 131, 175)$,
$(0, 112, 205, 230)$,
$(0, 154, 184, 221)$,\adfsplit
$(0, 124, 166, 182)$,
$(0, 53, 73, 128)$,
$(0, 140, 151, 260)$,
$(0, 40, 143, 219)$,\adfsplit
$(0, 15, 59, 118)$,
$(1, 134, 158, 245)$,
$(0, 78, 158, 227)$,
$(0, 12, 95, 209)$,\adfsplit
$(0, 167, 173, 215)$,
$(0, 57, 119, 170)$,
$(0, 75, 161, 254)$,
$(0, 21, 89, 159)$,\adfsplit
$(0, 33, 188, 245)$

%ADFvfyBlocksEnd
\adfLgap \noindent by the mapping:
$x \mapsto x + 3 j \adfmod{261}$ for $x < 261$,
$x \mapsto (x - 261 + 20 j \adfmod{60}) + 261$ for $261 \le x < 321$,
$x \mapsto x$ for $x \ge 321$,
$0 \le j < 87$.
\ADFvfyParStart{(323, ((89, 87, ((261, 3), (60, 20), (2, 2)))), ((29, 9), (62, 1)))} %ADFvfyParEnd
% End of 29^9 62^1
%%%%%%%%%%%%%%%%%%%%%%%%%%%%%%%%%%%%%%%%%%%%%%%%%%%%%%%%%%%%%%%%%%%%%%%%%%%%%%%%%%%%%%%%%%
%%%%%%%%%%%%%%%%%%%%%%%%%%%%%%%%%%%%%%%%%%%%%%%%%%%%%%%%%%%%%%%%%%%%%%%%%%%%%%%%%%%%%%%%%%

% Charlotte:GDD4-1-3-5-mod-6-TeX-gen-A:HITS-fun:4.10
\adfDgap
%ADFvfyBlocksStart {29,29,29,29,29,29,29,29,29,68}
\noindent{\boldmath $ 29^{9} 68^{1} $}~
With the point set $Z_{329}$ partitioned into
 residue classes modulo $9$ for $\{0, 1, \dots, 260\}$, and
 $\{261, 262, \dots, 328\}$,
 the design is generated from

\adfLgap %ADFvfyDesignStart
$(327, 129, 8, 37)$,
$(328, 182, 241, 108)$,
$(261, 2, 249, 214)$,
$(261, 3, 188, 72)$,\adfsplit
$(261, 122, 64, 193)$,
$(262, 248, 139, 133)$,
$(262, 249, 170, 1)$,
$(262, 12, 47, 63)$,\adfsplit
$(263, 42, 181, 41)$,
$(263, 227, 49, 246)$,
$(263, 25, 135, 260)$,
$(264, 88, 28, 255)$,\adfsplit
$(264, 20, 220, 14)$,
$(264, 80, 60, 108)$,
$(265, 255, 148, 51)$,
$(265, 109, 5, 29)$,\adfsplit
$(265, 160, 81, 143)$,
$(266, 43, 96, 3)$,
$(266, 80, 217, 254)$,
$(266, 176, 130, 234)$,\adfsplit
$(267, 106, 212, 98)$,
$(267, 184, 96, 2)$,
$(267, 216, 226, 156)$,
$(268, 146, 61, 219)$,\adfsplit
$(268, 238, 51, 136)$,
$(268, 234, 143, 131)$,
$(269, 43, 99, 75)$,
$(269, 33, 161, 145)$,\adfsplit
$(269, 101, 131, 31)$,
$(270, 141, 171, 133)$,
$(270, 212, 13, 155)$,
$(270, 242, 19, 156)$,\adfsplit
$(271, 45, 182, 193)$,
$(271, 115, 177, 165)$,
$(271, 212, 199, 152)$,
$(272, 69, 14, 70)$,\adfsplit
$(272, 94, 192, 73)$,
$(272, 251, 218, 153)$,
$(273, 6, 74, 187)$,
$(273, 149, 190, 26)$,\adfsplit
$(273, 153, 192, 148)$,
$(274, 93, 204, 215)$,
$(274, 235, 79, 5)$,
$(274, 18, 193, 20)$,\adfsplit
$(275, 172, 61, 152)$,
$(275, 140, 228, 153)$,
$(275, 177, 200, 193)$,
$(276, 225, 16, 147)$,\adfsplit
$(276, 110, 60, 212)$,
$(276, 154, 121, 26)$,
$(277, 11, 141, 108)$,
$(277, 55, 5, 215)$,\adfsplit
$(277, 169, 246, 22)$,
$(278, 114, 2, 180)$,
$(278, 148, 25, 199)$,
$(278, 176, 246, 251)$,\adfsplit
$(279, 240, 81, 66)$,
$(279, 239, 26, 136)$,
$(279, 214, 218, 22)$,
$(280, 238, 111, 71)$,\adfsplit
$(280, 145, 11, 140)$,
$(280, 168, 189, 196)$,
$(281, 144, 119, 80)$,
$(281, 178, 156, 13)$,\adfsplit
$(281, 86, 19, 204)$,
$(282, 49, 18, 124)$,
$(0, 4, 43, 46)$,
$(0, 3, 235, 250)$,\adfsplit
$(0, 6, 115, 172)$,
$(0, 19, 67, 160)$,
$(0, 17, 61, 238)$,
$(0, 25, 55, 107)$,\adfsplit
$(0, 8, 77, 193)$,
$(0, 44, 58, 241)$,
$(0, 113, 136, 202)$,
$(1, 13, 35, 176)$,\adfsplit
$(0, 121, 145, 146)$,
$(0, 49, 120, 179)$,
$(0, 176, 244, 254)$,
$(0, 47, 100, 212)$,\adfsplit
$(0, 119, 161, 196)$,
$(0, 26, 41, 259)$,
$(1, 119, 230, 282)$,
$(0, 71, 73, 227)$,\adfsplit
$(0, 56, 123, 224)$,
$(0, 91, 110, 132)$,
$(0, 89, 92, 214)$,
$(0, 38, 215, 219)$,\adfsplit
$(0, 82, 155, 165)$,
$(0, 32, 53, 84)$,
$(0, 104, 147, 326)$,
$(0, 29, 95, 156)$

%ADFvfyBlocksEnd
\adfLgap \noindent by the mapping:
$x \mapsto x + 3 j \adfmod{261}$ for $x < 261$,
$x \mapsto (x - 261 + 22 j \adfmod{66}) + 261$ for $261 \le x < 327$,
$x \mapsto x$ for $x \ge 327$,
$0 \le j < 87$.
\ADFvfyParStart{(329, ((92, 87, ((261, 3), (66, 22), (2, 2)))), ((29, 9), (68, 1)))} %ADFvfyParEnd
% End of 29^9 68^1
%%%%%%%%%%%%%%%%%%%%%%%%%%%%%%%%%%%%%%%%%%%%%%%%%%%%%%%%%%%%%%%%%%%%%%%%%%%%%%%%%%%%%%%%%%
%%%%%%%%%%%%%%%%%%%%%%%%%%%%%%%%%%%%%%%%%%%%%%%%%%%%%%%%%%%%%%%%%%%%%%%%%%%%%%%%%%%%%%%%%%

% Charlotte:GDD4-1-3-5-mod-6-TeX-gen-A:HITS-fun:4.10
\adfDgap
%ADFvfyBlocksStart {29,29,29,29,29,29,29,29,29,74}
\noindent{\boldmath $ 29^{9} 74^{1} $}~
With the point set $Z_{335}$ partitioned into
 residue classes modulo $9$ for $\{0, 1, \dots, 260\}$, and
 $\{261, 262, \dots, 334\}$,
 the design is generated from

\adfLgap %ADFvfyDesignStart
$(333, 28, 87, 116)$,
$(334, 51, 112, 128)$,
$(261, 170, 184, 183)$,
$(261, 78, 7, 74)$,\adfsplit
$(261, 100, 104, 216)$,
$(262, 242, 164, 49)$,
$(262, 169, 19, 192)$,
$(262, 132, 99, 149)$,\adfsplit
$(263, 86, 145, 210)$,
$(263, 232, 188, 238)$,
$(263, 150, 252, 191)$,
$(264, 20, 171, 247)$,\adfsplit
$(264, 115, 120, 253)$,
$(264, 194, 15, 134)$,
$(265, 186, 5, 72)$,
$(265, 74, 238, 79)$,\adfsplit
$(265, 183, 172, 143)$,
$(266, 233, 74, 3)$,
$(266, 166, 124, 100)$,
$(266, 195, 149, 171)$,\adfsplit
$(267, 102, 124, 101)$,
$(267, 33, 181, 206)$,
$(267, 230, 117, 256)$,
$(268, 253, 32, 153)$,\adfsplit
$(268, 101, 58, 210)$,
$(268, 133, 152, 96)$,
$(269, 244, 75, 204)$,
$(269, 198, 97, 146)$,\adfsplit
$(269, 161, 256, 50)$,
$(270, 197, 79, 204)$,
$(270, 209, 171, 5)$,
$(270, 219, 157, 136)$,\adfsplit
$(271, 185, 34, 226)$,
$(271, 256, 45, 89)$,
$(271, 120, 69, 128)$,
$(272, 113, 209, 108)$,\adfsplit
$(272, 97, 69, 37)$,
$(272, 80, 4, 12)$,
$(273, 63, 150, 111)$,
$(273, 179, 203, 128)$,\adfsplit
$(273, 217, 250, 247)$,
$(274, 19, 60, 70)$,
$(274, 59, 182, 237)$,
$(274, 44, 72, 202)$,\adfsplit
$(275, 216, 82, 80)$,
$(275, 23, 237, 148)$,
$(275, 115, 200, 258)$,
$(276, 93, 208, 59)$,\adfsplit
$(276, 35, 223, 168)$,
$(276, 90, 155, 148)$,
$(277, 236, 67, 80)$,
$(277, 7, 140, 117)$,\adfsplit
$(277, 258, 109, 57)$,
$(278, 47, 17, 82)$,
$(278, 203, 193, 97)$,
$(278, 93, 90, 24)$,\adfsplit
$(279, 246, 236, 116)$,
$(279, 132, 13, 196)$,
$(279, 149, 27, 118)$,
$(280, 235, 188, 204)$,\adfsplit
$(280, 139, 93, 243)$,
$(280, 142, 203, 29)$,
$(281, 5, 17, 124)$,
$(281, 105, 84, 49)$,\adfsplit
$(281, 180, 163, 11)$,
$(282, 63, 231, 48)$,
$(282, 250, 68, 161)$,
$(282, 211, 118, 155)$,\adfsplit
$(283, 228, 244, 70)$,
$(283, 234, 209, 139)$,
$(0, 2, 35, 307)$,
$(0, 7, 121, 219)$,\adfsplit
$(0, 4, 19, 186)$,
$(0, 20, 73, 259)$,
$(0, 12, 82, 187)$,
$(0, 11, 14, 105)$,\adfsplit
$(0, 13, 57, 177)$,
$(0, 25, 89, 155)$,
$(0, 67, 106, 188)$,
$(0, 30, 214, 242)$,\adfsplit
$(0, 6, 104, 241)$,
$(0, 116, 158, 232)$,
$(0, 185, 193, 284)$,
$(0, 74, 143, 154)$,\adfsplit
$(0, 164, 181, 308)$,
$(0, 107, 146, 208)$,
$(0, 62, 176, 191)$,
$(0, 32, 53, 96)$,\adfsplit
$(0, 43, 134, 182)$,
$(0, 26, 124, 332)$,
$(0, 85, 86, 123)$,
$(0, 34, 161, 247)$,\adfsplit
$(1, 47, 53, 85)$,
$(1, 58, 158, 242)$,
$(1, 13, 133, 191)$

%ADFvfyBlocksEnd
\adfLgap \noindent by the mapping:
$x \mapsto x + 3 j \adfmod{261}$ for $x < 261$,
$x \mapsto (x - 261 + 24 j \adfmod{72}) + 261$ for $261 \le x < 333$,
$x \mapsto x$ for $x \ge 333$,
$0 \le j < 87$.
\ADFvfyParStart{(335, ((95, 87, ((261, 3), (72, 24), (2, 2)))), ((29, 9), (74, 1)))} %ADFvfyParEnd
% End of 29^9 74^1
%%%%%%%%%%%%%%%%%%%%%%%%%%%%%%%%%%%%%%%%%%%%%%%%%%%%%%%%%%%%%%%%%%%%%%%%%%%%%%%%%%%%%%%%%%
%%%%%%%%%%%%%%%%%%%%%%%%%%%%%%%%%%%%%%%%%%%%%%%%%%%%%%%%%%%%%%%%%%%%%%%%%%%%%%%%%%%%%%%%%%

% Charlotte:GDD4-1-3-5-mod-6-TeX-gen-A:HITS-fun:4.10
\adfDgap
%ADFvfyBlocksStart {29,29,29,29,29,29,29,29,29,80}
\noindent{\boldmath $ 29^{9} 80^{1} $}~
With the point set $Z_{341}$ partitioned into
 residue classes modulo $9$ for $\{0, 1, \dots, 260\}$, and
 $\{261, 262, \dots, 340\}$,
 the design is generated from

\adfLgap %ADFvfyDesignStart
$(339, 6, 199, 47)$,
$(340, 46, 248, 147)$,
$(261, 43, 229, 89)$,
$(261, 222, 201, 11)$,\adfsplit
$(261, 113, 235, 180)$,
$(262, 189, 245, 3)$,
$(262, 226, 140, 124)$,
$(262, 161, 229, 123)$,\adfsplit
$(263, 2, 161, 175)$,
$(263, 203, 6, 18)$,
$(263, 111, 142, 19)$,
$(264, 51, 233, 67)$,\adfsplit
$(264, 120, 70, 212)$,
$(264, 208, 126, 137)$,
$(265, 53, 78, 36)$,
$(265, 181, 133, 50)$,\adfsplit
$(265, 75, 29, 40)$,
$(266, 150, 239, 156)$,
$(266, 58, 199, 223)$,
$(266, 242, 36, 101)$,\adfsplit
$(267, 190, 106, 251)$,
$(267, 50, 13, 110)$,
$(267, 141, 27, 228)$,
$(268, 238, 242, 114)$,\adfsplit
$(268, 61, 20, 165)$,
$(268, 117, 55, 212)$,
$(269, 170, 100, 24)$,
$(269, 171, 194, 22)$,\adfsplit
$(269, 218, 7, 66)$,
$(270, 182, 77, 235)$,
$(270, 12, 4, 222)$,
$(270, 160, 17, 243)$,\adfsplit
$(271, 138, 205, 31)$,
$(271, 59, 155, 8)$,
$(271, 91, 135, 231)$,
$(272, 148, 188, 163)$,\adfsplit
$(272, 18, 174, 43)$,
$(272, 110, 249, 32)$,
$(273, 214, 85, 78)$,
$(273, 63, 183, 62)$,\adfsplit
$(273, 200, 23, 235)$,
$(274, 82, 137, 177)$,
$(274, 160, 72, 170)$,
$(274, 75, 112, 239)$,\adfsplit
$(275, 113, 255, 121)$,
$(275, 180, 259, 146)$,
$(275, 197, 163, 177)$,
$(276, 134, 65, 172)$,\adfsplit
$(276, 70, 68, 36)$,
$(276, 237, 166, 168)$,
$(277, 10, 206, 254)$,
$(277, 15, 85, 39)$,\adfsplit
$(277, 54, 241, 113)$,
$(278, 243, 76, 107)$,
$(278, 60, 241, 23)$,
$(278, 10, 146, 156)$,\adfsplit
$(279, 12, 146, 258)$,
$(279, 63, 5, 233)$,
$(279, 187, 40, 118)$,
$(280, 32, 126, 48)$,\adfsplit
$(280, 19, 260, 169)$,
$(280, 164, 13, 78)$,
$(281, 19, 183, 257)$,
$(281, 76, 182, 186)$,\adfsplit
$(281, 170, 117, 97)$,
$(282, 55, 87, 76)$,
$(282, 98, 209, 72)$,
$(282, 212, 21, 160)$,\adfsplit
$(283, 30, 248, 44)$,
$(283, 241, 211, 74)$,
$(283, 90, 249, 208)$,
$(284, 17, 115, 41)$,\adfsplit
$(0, 4, 62, 104)$,
$(0, 1, 8, 230)$,
$(0, 5, 33, 143)$,
$(0, 2, 176, 188)$,\adfsplit
$(0, 19, 179, 200)$,
$(0, 47, 91, 113)$,
$(0, 49, 248, 254)$,
$(1, 14, 29, 211)$,\adfsplit
$(0, 30, 161, 336)$,
$(1, 20, 113, 205)$,
$(1, 185, 215, 285)$,
$(0, 80, 212, 286)$,\adfsplit
$(0, 155, 158, 184)$,
$(0, 66, 238, 239)$,
$(0, 209, 214, 285)$,
$(0, 10, 48, 77)$,\adfsplit
$(1, 34, 116, 338)$,
$(0, 22, 39, 107)$,
$(0, 13, 73, 129)$,
$(0, 28, 40, 133)$,\adfsplit
$(0, 109, 175, 284)$,
$(0, 58, 64, 177)$,
$(0, 61, 100, 103)$,
$(0, 52, 57, 150)$,\adfsplit
$(0, 3, 235, 337)$,
$(0, 85, 138, 338)$

%ADFvfyBlocksEnd
\adfLgap \noindent by the mapping:
$x \mapsto x + 3 j \adfmod{261}$ for $x < 261$,
$x \mapsto (x - 261 + 26 j \adfmod{78}) + 261$ for $261 \le x < 339$,
$x \mapsto x$ for $x \ge 339$,
$0 \le j < 87$.
\ADFvfyParStart{(341, ((98, 87, ((261, 3), (78, 26), (2, 2)))), ((29, 9), (80, 1)))} %ADFvfyParEnd
% End of 29^9 80^1
%%%%%%%%%%%%%%%%%%%%%%%%%%%%%%%%%%%%%%%%%%%%%%%%%%%%%%%%%%%%%%%%%%%%%%%%%%%%%%%%%%%%%%%%%%
%%%%%%%%%%%%%%%%%%%%%%%%%%%%%%%%%%%%%%%%%%%%%%%%%%%%%%%%%%%%%%%%%%%%%%%%%%%%%%%%%%%%%%%%%%

% Charlotte:GDD4-1-3-5-mod-6-TeX-gen-A:HITS-fun:4.10
\adfDgap
%ADFvfyBlocksStart {29,29,29,29,29,29,29,29,29,86}
\noindent{\boldmath $ 29^{9} 86^{1} $}~
With the point set $Z_{347}$ partitioned into
 residue classes modulo $9$ for $\{0, 1, \dots, 260\}$, and
 $\{261, 262, \dots, 346\}$,
 the design is generated from

\adfLgap %ADFvfyDesignStart
$(345, 59, 64, 177)$,
$(346, 170, 118, 186)$,
$(261, 129, 244, 180)$,
$(261, 238, 38, 159)$,\adfsplit
$(261, 44, 232, 167)$,
$(262, 49, 107, 199)$,
$(262, 75, 36, 131)$,
$(262, 61, 101, 186)$,\adfsplit
$(263, 260, 63, 245)$,
$(263, 10, 21, 239)$,
$(263, 69, 4, 34)$,
$(264, 192, 238, 214)$,\adfsplit
$(264, 249, 28, 11)$,
$(264, 212, 224, 198)$,
$(265, 185, 151, 192)$,
$(265, 177, 148, 227)$,\adfsplit
$(265, 1, 134, 81)$,
$(266, 147, 41, 134)$,
$(266, 108, 37, 151)$,
$(266, 240, 137, 67)$,\adfsplit
$(267, 39, 135, 59)$,
$(267, 164, 260, 96)$,
$(267, 244, 106, 184)$,
$(268, 202, 153, 17)$,\adfsplit
$(268, 227, 167, 168)$,
$(268, 199, 246, 124)$,
$(269, 39, 8, 22)$,
$(269, 180, 50, 208)$,\adfsplit
$(269, 101, 169, 60)$,
$(270, 9, 16, 5)$,
$(270, 179, 55, 254)$,
$(270, 87, 255, 139)$,\adfsplit
$(271, 5, 215, 186)$,
$(271, 52, 255, 164)$,
$(271, 256, 208, 0)$,
$(272, 117, 224, 202)$,\adfsplit
$(272, 7, 50, 19)$,
$(272, 150, 101, 129)$,
$(273, 183, 122, 181)$,
$(273, 125, 195, 58)$,\adfsplit
$(273, 189, 250, 164)$,
$(274, 113, 229, 74)$,
$(274, 232, 107, 195)$,
$(274, 189, 127, 219)$,\adfsplit
$(275, 81, 231, 208)$,
$(275, 211, 68, 170)$,
$(275, 259, 165, 128)$,
$(276, 200, 163, 162)$,\adfsplit
$(276, 178, 255, 131)$,
$(276, 132, 202, 251)$,
$(277, 119, 199, 241)$,
$(277, 176, 202, 107)$,\adfsplit
$(277, 12, 198, 240)$,
$(278, 23, 100, 81)$,
$(278, 182, 184, 97)$,
$(278, 107, 75, 141)$,\adfsplit
$(279, 125, 63, 78)$,
$(279, 56, 50, 31)$,
$(279, 61, 156, 46)$,
$(280, 197, 9, 226)$,\adfsplit
$(280, 176, 148, 250)$,
$(280, 186, 3, 164)$,
$(281, 9, 251, 138)$,
$(281, 95, 85, 79)$,\adfsplit
$(281, 186, 191, 37)$,
$(282, 121, 128, 54)$,
$(282, 244, 33, 104)$,
$(282, 170, 3, 7)$,\adfsplit
$(283, 260, 29, 156)$,
$(283, 234, 55, 94)$,
$(283, 69, 142, 104)$,
$(284, 143, 258, 99)$,\adfsplit
$(284, 155, 57, 190)$,
$(284, 221, 43, 94)$,
$(285, 203, 24, 138)$,
$(0, 6, 92, 116)$,\adfsplit
$(0, 11, 34, 89)$,
$(0, 77, 123, 251)$,
$(0, 17, 101, 149)$,
$(1, 2, 242, 285)$,\adfsplit
$(0, 13, 16, 285)$,
$(0, 2, 122, 286)$,
$(1, 70, 152, 218)$,
$(1, 101, 158, 287)$,\adfsplit
$(1, 14, 47, 166)$,
$(0, 152, 194, 288)$,
$(1, 5, 161, 169)$,
$(1, 92, 95, 206)$,\adfsplit
$(1, 89, 142, 286)$,
$(0, 83, 172, 229)$,
$(0, 3, 209, 315)$,
$(0, 8, 48, 205)$,\adfsplit
$(1, 22, 212, 316)$,
$(0, 25, 156, 253)$,
$(0, 10, 76, 343)$,
$(0, 12, 103, 187)$,\adfsplit
$(0, 55, 160, 174)$,
$(0, 31, 57, 163)$,
$(0, 24, 84, 202)$,
$(0, 69, 223, 316)$,\adfsplit
$(0, 100, 120, 342)$

%ADFvfyBlocksEnd
\adfLgap \noindent by the mapping:
$x \mapsto x + 3 j \adfmod{261}$ for $x < 261$,
$x \mapsto (x - 261 + 28 j \adfmod{84}) + 261$ for $261 \le x < 345$,
$x \mapsto x$ for $x \ge 345$,
$0 \le j < 87$.
\ADFvfyParStart{(347, ((101, 87, ((261, 3), (84, 28), (2, 2)))), ((29, 9), (86, 1)))} %ADFvfyParEnd
% End of 29^9 86^1
%%%%%%%%%%%%%%%%%%%%%%%%%%%%%%%%%%%%%%%%%%%%%%%%%%%%%%%%%%%%%%%%%%%%%%%%%%%%%%%%%%%%%%%%%%
%%%%%%%%%%%%%%%%%%%%%%%%%%%%%%%%%%%%%%%%%%%%%%%%%%%%%%%%%%%%%%%%%%%%%%%%%%%%%%%%%%%%%%%%%%

% Charlotte:GDD4-1-3-5-mod-6-TeX-gen-A:HITS-fun:4.10
\adfDgap
%ADFvfyBlocksStart {29,29,29,29,29,29,29,29,29,92}
\noindent{\boldmath $ 29^{9} 92^{1} $}~
With the point set $Z_{353}$ partitioned into
 residue classes modulo $9$ for $\{0, 1, \dots, 260\}$, and
 $\{261, 262, \dots, 352\}$,
 the design is generated from

\adfLgap %ADFvfyDesignStart
$(351, 54, 83, 34)$,
$(352, 198, 233, 196)$,
$(261, 140, 160, 90)$,
$(261, 195, 49, 156)$,\adfsplit
$(261, 146, 125, 82)$,
$(262, 116, 88, 37)$,
$(262, 131, 6, 192)$,
$(262, 112, 189, 101)$,\adfsplit
$(263, 131, 39, 145)$,
$(263, 98, 168, 252)$,
$(263, 88, 164, 193)$,
$(264, 236, 183, 35)$,\adfsplit
$(264, 118, 221, 103)$,
$(264, 117, 178, 159)$,
$(265, 203, 233, 0)$,
$(265, 2, 109, 42)$,\adfsplit
$(265, 223, 174, 40)$,
$(266, 189, 14, 156)$,
$(266, 152, 196, 46)$,
$(266, 83, 240, 58)$,\adfsplit
$(267, 79, 198, 24)$,
$(267, 119, 185, 228)$,
$(267, 55, 62, 157)$,
$(268, 160, 56, 50)$,\adfsplit
$(268, 15, 202, 1)$,
$(268, 93, 216, 215)$,
$(269, 9, 255, 112)$,
$(269, 179, 187, 78)$,\adfsplit
$(269, 239, 1, 2)$,
$(270, 235, 230, 8)$,
$(270, 213, 202, 128)$,
$(270, 21, 106, 189)$,\adfsplit
$(271, 72, 222, 97)$,
$(271, 247, 253, 197)$,
$(271, 84, 5, 173)$,
$(272, 78, 212, 62)$,\adfsplit
$(272, 172, 175, 133)$,
$(272, 192, 45, 11)$,
$(273, 45, 42, 206)$,
$(273, 109, 210, 32)$,\adfsplit
$(273, 92, 25, 94)$,
$(274, 62, 186, 30)$,
$(274, 70, 103, 191)$,
$(274, 217, 126, 14)$,\adfsplit
$(275, 165, 213, 181)$,
$(275, 50, 121, 137)$,
$(275, 70, 44, 63)$,
$(276, 128, 87, 163)$,\adfsplit
$(276, 52, 179, 111)$,
$(276, 32, 13, 63)$,
$(277, 89, 28, 138)$,
$(277, 221, 2, 15)$,\adfsplit
$(277, 135, 106, 130)$,
$(278, 227, 60, 180)$,
$(278, 160, 206, 201)$,
$(278, 157, 127, 95)$,\adfsplit
$(279, 120, 82, 18)$,
$(279, 143, 141, 158)$,
$(279, 97, 56, 4)$,
$(280, 182, 105, 156)$,\adfsplit
$(280, 98, 0, 1)$,
$(280, 241, 157, 140)$,
$(281, 172, 120, 85)$,
$(281, 225, 134, 59)$,\adfsplit
$(281, 186, 245, 214)$,
$(282, 121, 236, 243)$,
$(282, 86, 82, 17)$,
$(282, 57, 123, 34)$,\adfsplit
$(283, 181, 49, 189)$,
$(283, 187, 57, 119)$,
$(283, 176, 240, 260)$,
$(284, 198, 209, 219)$,\adfsplit
$(284, 249, 62, 232)$,
$(284, 185, 46, 103)$,
$(285, 98, 227, 33)$,
$(285, 109, 247, 189)$,\adfsplit
$(285, 3, 196, 158)$,
$(286, 77, 236, 166)$,
$(0, 8, 131, 236)$,
$(0, 10, 110, 158)$,\adfsplit
$(0, 22, 146, 224)$,
$(0, 14, 71, 163)$,
$(0, 38, 185, 287)$,
$(0, 140, 143, 288)$,\adfsplit
$(0, 34, 44, 56)$,
$(0, 23, 148, 188)$,
$(1, 95, 146, 289)$,
$(1, 74, 215, 290)$,\adfsplit
$(1, 22, 131, 347)$,
$(1, 14, 67, 142)$,
$(0, 73, 128, 316)$,
$(0, 37, 215, 320)$,\adfsplit
$(0, 97, 239, 318)$,
$(0, 31, 116, 196)$,
$(0, 46, 209, 319)$,
$(0, 4, 179, 183)$,\adfsplit
$(0, 12, 145, 286)$,
$(0, 40, 88, 132)$,
$(0, 112, 124, 348)$,
$(0, 43, 96, 190)$,\adfsplit
$(0, 6, 30, 205)$,
$(0, 13, 60, 317)$,
$(0, 57, 157, 289)$,
$(0, 69, 235, 350)$

%ADFvfyBlocksEnd
\adfLgap \noindent by the mapping:
$x \mapsto x + 3 j \adfmod{261}$ for $x < 261$,
$x \mapsto (x - 261 + 30 j \adfmod{90}) + 261$ for $261 \le x < 351$,
$x \mapsto x$ for $x \ge 351$,
$0 \le j < 87$.
\ADFvfyParStart{(353, ((104, 87, ((261, 3), (90, 30), (2, 2)))), ((29, 9), (92, 1)))} %ADFvfyParEnd
% End of 29^9 92^1
%%%%%%%%%%%%%%%%%%%%%%%%%%%%%%%%%%%%%%%%%%%%%%%%%%%%%%%%%%%%%%%%%%%%%%%%%%%%%%%%%%%%%%%%%%
%%%%%%%%%%%%%%%%%%%%%%%%%%%%%%%%%%%%%%%%%%%%%%%%%%%%%%%%%%%%%%%%%%%%%%%%%%%%%%%%%%%%%%%%%%

% Charlotte:GDD4-1-3-5-mod-6-TeX-gen-A:HITS-fun:4.10
\adfDgap
%ADFvfyBlocksStart {29,29,29,29,29,29,29,29,29,98}
\noindent{\boldmath $ 29^{9} 98^{1} $}~
With the point set $Z_{359}$ partitioned into
 residue classes modulo $9$ for $\{0, 1, \dots, 260\}$, and
 $\{261, 262, \dots, 358\}$,
 the design is generated from

\adfLgap %ADFvfyDesignStart
$(357, 7, 179, 63)$,
$(358, 172, 50, 69)$,
$(261, 73, 65, 183)$,
$(261, 43, 125, 211)$,\adfsplit
$(261, 234, 86, 222)$,
$(262, 0, 141, 93)$,
$(262, 95, 202, 127)$,
$(262, 227, 152, 205)$,\adfsplit
$(263, 156, 121, 99)$,
$(263, 253, 254, 197)$,
$(263, 33, 196, 185)$,
$(264, 148, 83, 174)$,\adfsplit
$(264, 26, 77, 141)$,
$(264, 28, 7, 207)$,
$(265, 171, 250, 190)$,
$(265, 147, 158, 98)$,\adfsplit
$(265, 218, 132, 256)$,
$(266, 218, 31, 131)$,
$(266, 48, 234, 78)$,
$(266, 233, 187, 181)$,\adfsplit
$(267, 131, 100, 72)$,
$(267, 229, 53, 39)$,
$(267, 47, 141, 34)$,
$(268, 254, 170, 198)$,\adfsplit
$(268, 84, 176, 46)$,
$(268, 69, 52, 49)$,
$(269, 45, 65, 35)$,
$(269, 150, 64, 196)$,\adfsplit
$(269, 203, 112, 39)$,
$(270, 55, 102, 236)$,
$(270, 168, 108, 26)$,
$(270, 70, 184, 86)$,\adfsplit
$(271, 254, 78, 9)$,
$(271, 210, 212, 251)$,
$(271, 70, 154, 193)$,
$(272, 40, 119, 90)$,\adfsplit
$(272, 145, 96, 35)$,
$(272, 219, 59, 196)$,
$(273, 34, 243, 122)$,
$(273, 100, 6, 53)$,\adfsplit
$(273, 191, 174, 22)$,
$(274, 208, 117, 106)$,
$(274, 157, 95, 2)$,
$(274, 105, 102, 35)$,\adfsplit
$(275, 259, 131, 144)$,
$(275, 244, 53, 186)$,
$(275, 147, 245, 211)$,
$(276, 109, 146, 72)$,\adfsplit
$(276, 231, 58, 169)$,
$(276, 237, 113, 107)$,
$(277, 31, 197, 165)$,
$(277, 64, 59, 47)$,\adfsplit
$(277, 16, 96, 81)$,
$(278, 250, 257, 108)$,
$(278, 128, 118, 159)$,
$(278, 40, 89, 39)$,\adfsplit
$(279, 120, 125, 204)$,
$(279, 136, 79, 140)$,
$(279, 92, 135, 211)$,
$(280, 6, 116, 45)$,\adfsplit
$(280, 254, 48, 113)$,
$(280, 250, 220, 1)$,
$(281, 85, 201, 241)$,
$(281, 32, 36, 186)$,\adfsplit
$(281, 209, 71, 235)$,
$(282, 236, 178, 109)$,
$(282, 211, 95, 26)$,
$(282, 201, 87, 225)$,\adfsplit
$(283, 221, 138, 206)$,
$(283, 69, 173, 175)$,
$(283, 70, 153, 253)$,
$(284, 82, 155, 115)$,\adfsplit
$(284, 9, 96, 256)$,
$(284, 66, 89, 41)$,
$(285, 187, 172, 104)$,
$(285, 66, 224, 29)$,\adfsplit
$(285, 42, 58, 207)$,
$(286, 27, 193, 26)$,
$(286, 77, 155, 129)$,
$(286, 253, 222, 79)$,\adfsplit
$(1, 26, 137, 158)$,
$(0, 35, 193, 287)$,
$(0, 4, 240, 288)$,
$(0, 10, 102, 321)$,\adfsplit
$(0, 6, 13, 161)$,
$(0, 44, 77, 121)$,
$(0, 42, 97, 215)$,
$(0, 8, 148, 290)$,\adfsplit
$(0, 34, 78, 210)$,
$(0, 66, 136, 323)$,
$(0, 33, 172, 319)$,
$(1, 56, 203, 287)$,\adfsplit
$(0, 122, 157, 322)$,
$(0, 227, 256, 354)$,
$(0, 202, 221, 352)$,
$(1, 65, 221, 288)$,\adfsplit
$(0, 80, 130, 239)$,
$(0, 133, 229, 253)$,
$(0, 43, 203, 292)$,
$(0, 53, 67, 95)$,\adfsplit
$(1, 44, 67, 289)$,
$(0, 38, 62, 353)$,
$(0, 184, 232, 324)$,
$(0, 187, 254, 355)$,\adfsplit
$(1, 113, 116, 355)$,
$(0, 89, 185, 356)$,
$(0, 188, 208, 259)$

%ADFvfyBlocksEnd
\adfLgap \noindent by the mapping:
$x \mapsto x + 3 j \adfmod{261}$ for $x < 261$,
$x \mapsto (x - 261 + 32 j \adfmod{96}) + 261$ for $261 \le x < 357$,
$x \mapsto x$ for $x \ge 357$,
$0 \le j < 87$.
\ADFvfyParStart{(359, ((107, 87, ((261, 3), (96, 32), (2, 2)))), ((29, 9), (98, 1)))} %ADFvfyParEnd
% End of 29^9 98^1
%%%%%%%%%%%%%%%%%%%%%%%%%%%%%%%%%%%%%%%%%%%%%%%%%%%%%%%%%%%%%%%%%%%%%%%%%%%%%%%%%%%%%%%%%%
%%%%%%%%%%%%%%%%%%%%%%%%%%%%%%%%%%%%%%%%%%%%%%%%%%%%%%%%%%%%%%%%%%%%%%%%%%%%%%%%%%%%%%%%%%

% Charlotte:GDD4-1-3-5-mod-6-TeX-gen-A:HITS-fun:4.10
\adfDgap
%ADFvfyBlocksStart {29,29,29,29,29,29,29,29,29,104}
\noindent{\boldmath $ 29^{9} 104^{1} $}~
With the point set $Z_{365}$ partitioned into
 residue classes modulo $8$ for $\{0, 1, \dots, 231\}$,
 $\{232, 233, \dots, 260\}$, and
 $\{261, 262, \dots, 364\}$,
 the design is generated from

\adfLgap %ADFvfyDesignStart
$(261, 146, 29, 193)$,
$(262, 254, 110, 93)$,
$(263, 96, 155, 220)$,
$(264, 204, 103, 22)$,\adfsplit
$(265, 178, 15, 107)$,
$(266, 73, 151, 203)$,
$(267, 234, 159, 26)$,
$(268, 143, 204, 99)$,\adfsplit
$(269, 113, 180, 27)$,
$(270, 260, 31, 150)$,
$(271, 32, 121, 154)$,
$(272, 178, 213, 174)$,\adfsplit
$(273, 45, 237, 96)$,
$(261, 182, 87, 228)$,
$(262, 162, 180, 51)$,
$(263, 13, 183, 247)$,\adfsplit
$(264, 179, 1, 210)$,
$(265, 254, 148, 225)$,
$(266, 212, 86, 239)$,
$(267, 164, 198, 185)$,\adfsplit
$(268, 57, 64, 45)$,
$(269, 162, 250, 102)$,
$(270, 8, 220, 122)$,
$(0, 28, 241, 294)$,\adfsplit
$(0, 2, 38, 335)$,
$(0, 1, 10, 275)$,
$(0, 3, 76, 289)$,
$(0, 5, 11, 305)$,\adfsplit
$(0, 14, 43, 306)$,
$(0, 27, 243, 261)$,
$(0, 37, 244, 264)$,
$(0, 15, 157, 269)$,\adfsplit
$(0, 26, 175, 291)$,
$(0, 49, 158, 297)$,
$(0, 22, 246, 310)$,
$(0, 53, 138, 247)$,\adfsplit
$(0, 70, 239, 363)$,
$(0, 45, 145, 337)$,
$(0, 30, 169, 299)$,
$(0, 25, 66, 150)$,\adfsplit
$(0, 42, 97, 364)$,
$(0, 58, 116, 174)$

%ADFvfyBlocksEnd
\adfLgap \noindent by the mapping:
$x \mapsto x +  j \adfmod{232}$ for $x < 232$,
$x \mapsto (x +  j \adfmod{29}) + 232$ for $232 \le x < 261$,
$x \mapsto (x - 261 + 13 j \adfmod{104}) + 261$ for $x \ge 261$,
$0 \le j < 232$
 for the first 41 blocks,
$0 \le j < 58$
 for the last block.
\ADFvfyParStart{(365, ((41, 232, ((232, 1), (29, 1), (104, 13))), (1, 58, ((232, 1), (29, 1), (104, 13)))), ((29, 8), (29, 1), (104, 1)))} %ADFvfyParEnd
% End of 29^9 104^1
%%%%%%%%%%%%%%%%%%%%%%%%%%%%%%%%%%%%%%%%%%%%%%%%%%%%%%%%%%%%%%%%%%%%%%%%%%%%%%%%%%%%%%%%%%
%%%%%%%%%%%%%%%%%%%%%%%%%%%%%%%%%%%%%%%%%%%%%%%%%%%%%%%%%%%%%%%%%%%%%%%%%%%%%%%%%%%%%%%%%%

% Charlotte:GDD4-1-3-5-mod-6-TeX-gen-A:HITS-fun:4.10
\adfDgap
%ADFvfyBlocksStart {29,29,29,29,29,29,29,29,29,110}
\noindent{\boldmath $ 29^{9} 110^{1} $}~
With the point set $Z_{371}$ partitioned into
 residue classes modulo $9$ for $\{0, 1, \dots, 260\}$, and
 $\{261, 262, \dots, 370\}$,
 the design is generated from

\adfLgap %ADFvfyDesignStart
$(261, 62, 103, 69)$,
$(273, 63, 104, 70)$,
$(285, 64, 105, 71)$,
$(261, 39, 223, 54)$,\adfsplit
$(273, 40, 224, 55)$,
$(285, 41, 225, 56)$,
$(261, 59, 109, 146)$,
$(273, 60, 110, 147)$,\adfsplit
$(285, 61, 111, 148)$,
$(262, 15, 124, 121)$,
$(274, 16, 125, 122)$,
$(286, 17, 126, 123)$,\adfsplit
$(262, 118, 8, 90)$,
$(274, 119, 9, 91)$,
$(286, 120, 10, 92)$,
$(262, 221, 119, 183)$,\adfsplit
$(274, 222, 120, 184)$,
$(286, 223, 121, 185)$,
$(263, 168, 39, 53)$,
$(275, 169, 40, 54)$,\adfsplit
$(287, 170, 41, 55)$,
$(263, 254, 244, 86)$,
$(275, 255, 245, 87)$,
$(287, 256, 246, 88)$,\adfsplit
$(263, 157, 133, 117)$,
$(275, 158, 134, 118)$,
$(287, 159, 135, 119)$,
$(264, 75, 258, 133)$,\adfsplit
$(276, 76, 259, 134)$,
$(288, 77, 260, 135)$,
$(264, 99, 28, 26)$,
$(276, 100, 29, 27)$,\adfsplit
$(288, 101, 30, 28)$,
$(264, 230, 83, 4)$,
$(276, 231, 84, 5)$,
$(288, 232, 85, 6)$,\adfsplit
$(265, 109, 240, 193)$,
$(277, 110, 241, 194)$,
$(289, 111, 242, 195)$,
$(265, 137, 48, 54)$,\adfsplit
$(277, 138, 49, 55)$,
$(289, 139, 50, 56)$,
$(265, 233, 133, 122)$,
$(277, 234, 134, 123)$,\adfsplit
$(289, 235, 135, 124)$,
$(266, 247, 87, 113)$,
$(278, 248, 88, 114)$,
$(290, 249, 89, 115)$,\adfsplit
$(266, 226, 107, 255)$,
$(278, 227, 108, 256)$,
$(290, 228, 109, 257)$,
$(266, 79, 83, 207)$,\adfsplit
$(278, 80, 84, 208)$,
$(290, 81, 85, 209)$,
$(267, 241, 74, 120)$,
$(279, 242, 75, 121)$,\adfsplit
$(291, 243, 76, 122)$,
$(267, 212, 251, 94)$,
$(279, 213, 252, 95)$,
$(291, 214, 253, 96)$,\adfsplit
$(267, 141, 45, 46)$,
$(279, 142, 46, 47)$,
$(291, 143, 47, 48)$,
$(268, 69, 256, 176)$,\adfsplit
$(280, 70, 257, 177)$,
$(292, 71, 258, 178)$,
$(268, 142, 89, 19)$,
$(280, 143, 90, 20)$,\adfsplit
$(292, 144, 91, 21)$,
$(268, 83, 246, 180)$,
$(280, 84, 247, 181)$,
$(292, 85, 248, 182)$,\adfsplit
$(269, 19, 88, 131)$,
$(281, 20, 89, 132)$,
$(293, 21, 90, 133)$,
$(269, 96, 29, 84)$,\adfsplit
$(281, 97, 30, 85)$,
$(293, 98, 31, 86)$,
$(269, 99, 175, 116)$,
$(281, 100, 176, 117)$,\adfsplit
$(293, 101, 177, 118)$,
$(270, 166, 196, 217)$,
$(282, 167, 197, 218)$,
$(294, 168, 198, 219)$,\adfsplit
$(0, 22, 44, 271)$,
$(0, 20, 42, 270)$,
$(0, 5, 65, 342)$,
$(0, 25, 86, 242)$,\adfsplit
$(0, 61, 217, 236)$,
$(0, 19, 105, 272)$,
$(0, 48, 248, 253)$,
$(0, 49, 62, 283)$,\adfsplit
$(0, 23, 209, 284)$,
$(0, 13, 186, 295)$,
$(0, 33, 229, 296)$,
$(0, 32, 120, 356)$,\adfsplit
$(1, 34, 76, 284)$,
$(0, 8, 57, 343)$,
$(1, 53, 121, 307)$,
$(0, 52, 256, 355)$,\adfsplit
$(0, 68, 116, 319)$,
$(0, 176, 241, 367)$,
$(0, 31, 56, 368)$,
$(1, 200, 242, 296)$,\adfsplit
$(1, 89, 230, 331)$,
$(1, 206, 239, 308)$,
$(0, 145, 193, 282)$,
$(0, 60, 199, 318)$,\adfsplit
$(0, 91, 230, 344)$,
$(0, 85, 170, 369)$,
$(0, 122, 238, 370)$,
$(1, 32, 202, 330)$,\adfsplit
$(1, 50, 254, 366)$

%ADFvfyBlocksEnd
\adfLgap \noindent by the mapping:
$x \mapsto x + 3 j \adfmod{261}$ for $x < 261$,
$x \mapsto (x - 261 + 36 j \adfmod{108}) + 261$ for $261 \le x < 369$,
$x \mapsto x$ for $x \ge 369$,
$0 \le j < 87$.
\ADFvfyParStart{(371, ((113, 87, ((261, 3), (108, 36), (2, 2)))), ((29, 9), (110, 1)))} %ADFvfyParEnd
% End of 29^9 110^1
%%%%%%%%%%%%%%%%%%%%%%%%%%%%%%%%%%%%%%%%%%%%%%%%%%%%%%%%%%%%%%%%%%%%%%%%%%%%%%%%%%%%%%%%%%
%%%%%%%%%%%%%%%%%%%%%%%%%%%%%%%%%%%%%%%%%%%%%%%%%%%%%%%%%%%%%%%%%%%%%%%%%%%%%%%%%%%%%%%%%%

% Charlotte:GDD4-1-3-5-mod-6-TeX-gen-A:HITS-fun:4.10
\adfDgap
%ADFvfyBlocksStart {29,29,29,29,29,29,29,29,29,29,29,29,29,29,29,29,29,29,29,29,29,26}
\noindent{\boldmath $ 29^{21} 26^{1} $}~
With the point set $Z_{635}$ partitioned into
 residue classes modulo $21$ for $\{0, 1, \dots, 608\}$, and
 $\{609, 610, \dots, 634\}$,
 the design is generated from

\adfLgap %ADFvfyDesignStart
$(609, 0, 1, 2)$,
$(610, 0, 202, 407)$,
$(611, 0, 205, 203)$,
$(612, 0, 406, 608)$,\adfsplit
$(613, 0, 607, 404)$,
$(614, 0, 4, 11)$,
$(615, 0, 7, 605)$,
$(616, 0, 598, 602)$,\adfsplit
$(617, 0, 10, 23)$,
$(618, 0, 13, 599)$,
$(619, 0, 586, 596)$,
$(620, 0, 16, 35)$,\adfsplit
$(621, 0, 19, 593)$,
$(622, 0, 574, 590)$,
$(623, 0, 22, 5)$,
$(624, 0, 592, 587)$,\adfsplit
$(625, 0, 604, 17)$,
$(626, 0, 25, 53)$,
$(627, 0, 28, 584)$,
$(628, 0, 556, 581)$,\adfsplit
$(629, 0, 31, 65)$,
$(630, 0, 34, 578)$,
$(631, 0, 544, 575)$,
$(632, 0, 37, 8)$,\adfsplit
$(633, 0, 580, 572)$,
$(634, 0, 601, 29)$,
$(53, 260, 165, 321)$,
$(80, 567, 471, 558)$,\adfsplit
$(606, 66, 48, 115)$,
$(480, 540, 121, 176)$,
$(23, 285, 70, 8)$,
$(528, 453, 356, 72)$,\adfsplit
$(148, 511, 204, 84)$,
$(209, 445, 182, 254)$,
$(577, 190, 41, 342)$,
$(361, 459, 505, 9)$,\adfsplit
$(479, 252, 346, 390)$,
$(227, 151, 593, 58)$,
$(46, 375, 234, 89)$,
$(311, 184, 511, 392)$,\adfsplit
$(270, 170, 570, 310)$,
$(419, 433, 109, 595)$,
$(84, 25, 353, 557)$,
$(411, 276, 385, 27)$,\adfsplit
$(175, 338, 567, 509)$,
$(396, 0, 52, 524)$,
$(365, 62, 513, 101)$,
$(433, 376, 115, 308)$,\adfsplit
$(19, 341, 52, 233)$,
$(606, 292, 112, 182)$,
$(324, 545, 497, 207)$,
$(43, 411, 287, 173)$,\adfsplit
$(161, 73, 303, 575)$,
$(0, 3, 157, 296)$,
$(0, 6, 30, 253)$,
$(0, 12, 83, 342)$,\adfsplit
$(0, 80, 183, 317)$,
$(0, 32, 233, 383)$,
$(0, 91, 255, 430)$,
$(0, 20, 212, 355)$,\adfsplit
$(0, 90, 196, 493)$,
$(0, 41, 211, 431)$,
$(0, 36, 86, 165)$,
$(0, 78, 177, 411)$,\adfsplit
$(0, 38, 199, 415)$,
$(0, 132, 283, 443)$,
$(0, 54, 146, 385)$,
$(0, 66, 187, 435)$,\adfsplit
$(0, 101, 256, 367)$,
$(0, 82, 184, 505)$

%ADFvfyBlocksEnd
\adfLgap \noindent by the mapping:
$x \mapsto x + 3 j \adfmod{609}$ for $x < 609$,
$x \mapsto x$ for $x \ge 609$,
$0 \le j < 203$
 for the first 26 blocks;
$x \mapsto x +  j \adfmod{609}$ for $x < 609$,
$x \mapsto x$ for $x \ge 609$,
$0 \le j < 609$
 for the last 44 blocks.
\ADFvfyParStart{(635, ((26, 203, ((609, 3), (26, 26))), (44, 609, ((609, 1), (26, 26)))), ((29, 21), (26, 1)))} %ADFvfyParEnd
% End of 29^21 26^1
%%%%%%%%%%%%%%%%%%%%%%%%%%%%%%%%%%%%%%%%%%%%%%%%%%%%%%%%%%%%%%%%%%%%%%%%%%%%%%%%%%%%%%%%%%
%%%%%%%%%%%%%%%%%%%%%%%%%%%%%%%%%%%%%%%%%%%%%%%%%%%%%%%%%%%%%%%%%%%%%%%%%%%%%%%%%%%%%%%%%%

%%%%%%%%%%%%%%%%%%%%%%%%%%%%%%%%%%%%%%%%%%%%%%%%%%%%%%%%%%%%%%%%%%%%%%%%%%%%%%%%%%%%%%%%%%
%%%%%%%%%%%%%%%%%%%%%%%%%%%%%%%%%%%%%%%%%%%%%%%%%%%%%%%%%%%%%%%%%%%%%%%%%%%%%%%%%%%%%%%%%%
\section{4-GDDs for the proof of Lemma \ref{lem:4-GDD 31^u m^1}}
\label{app:4-GDD 31^u m^1}
\adfnull{
$ 31^{12} 25^1 $,
$ 31^{12} 28^1 $,
$ 31^9 22^1 $,
$ 31^9 28^1 $,
$ 31^9 34^1 $,
$ 31^9 40^1 $,
$ 31^9 46^1 $,
$ 31^9 52^1 $,
$ 31^9 58^1 $,
$ 31^9 64^1 $,
$ 31^9 70^1 $,
$ 31^9 76^1 $,
$ 31^9 82^1 $,
$ 31^9 88^1 $,
$ 31^9 94^1 $,
$ 31^9 100^1 $,
$ 31^9 106^1 $,
$ 31^9 112^1 $ and
$ 31^9 118^1 $.
}

% Charlotte:GDD4-1-3-5-mod-6-TeX-gen-A:HITS-fun:4.10
\adfDgap
%ADFvfyBlocksStart {31,31,31,31,31,31,31,31,31,31,31,31,25}
\noindent{\boldmath $ 31^{12} 25^{1} $}~
With the point set $Z_{397}$ partitioned into
 residue classes modulo $12$ for $\{0, 1, \dots, 371\}$, and
 $\{372, 373, \dots, 396\}$,
 the design is generated from

\adfLgap %ADFvfyDesignStart
$(372, 267, 187, 257)$,
$(373, 137, 196, 153)$,
$(374, 179, 103, 276)$,
$(375, 99, 265, 113)$,\adfsplit
$(376, 367, 137, 63)$,
$(377, 254, 13, 120)$,
$(378, 351, 212, 355)$,
$(379, 5, 321, 205)$,\adfsplit
$(282, 277, 19, 195)$,
$(363, 140, 364, 217)$,
$(155, 49, 216, 52)$,
$(186, 285, 233, 111)$,\adfsplit
$(107, 125, 261, 344)$,
$(332, 173, 58, 45)$,
$(114, 88, 121, 153)$,
$(62, 172, 152, 203)$,\adfsplit
$(349, 249, 191, 360)$,
$(25, 343, 210, 165)$,
$(114, 305, 51, 28)$,
$(308, 158, 143, 196)$,\adfsplit
$(24, 68, 213, 194)$,
$(0, 6, 40, 177)$,
$(0, 9, 78, 268)$,
$(0, 42, 130, 253)$,\adfsplit
$(0, 8, 29, 299)$,
$(0, 17, 66, 155)$,
$(0, 30, 71, 121)$,
$(0, 22, 79, 267)$,\adfsplit
$(0, 46, 101, 163)$,
$(0, 2, 27, 94)$,
$(0, 35, 160, 197)$,
$(0, 28, 157, 221)$,\adfsplit
$(0, 93, 186, 279)$,
$(396, 0, 124, 248)$

%ADFvfyBlocksEnd
\adfLgap \noindent by the mapping:
$x \mapsto x +  j \adfmod{372}$ for $x < 372$,
$x \mapsto (x - 372 + 8 j \adfmod{24}) + 372$ for $372 \le x < 396$,
$396 \mapsto 396$,
$0 \le j < 372$
 for the first 32 blocks,
$0 \le j < 93$
 for the next block,
$0 \le j < 124$
 for the last block.
\ADFvfyParStart{(397, ((32, 372, ((372, 1), (24, 8), (1, 1))), (1, 93, ((372, 1), (24, 8), (1, 1))), (1, 124, ((372, 1), (24, 8), (1, 1)))), ((31, 12), (25, 1)))} %ADFvfyParEnd
% End of 31^12 25^1
%%%%%%%%%%%%%%%%%%%%%%%%%%%%%%%%%%%%%%%%%%%%%%%%%%%%%%%%%%%%%%%%%%%%%%%%%%%%%%%%%%%%%%%%%%
%%%%%%%%%%%%%%%%%%%%%%%%%%%%%%%%%%%%%%%%%%%%%%%%%%%%%%%%%%%%%%%%%%%%%%%%%%%%%%%%%%%%%%%%%%

% Charlotte:GDD4-1-3-5-mod-6-TeX-gen-A:HITS-fun:4.10
\adfDgap
%ADFvfyBlocksStart {31,31,31,31,31,31,31,31,31,31,31,31,28}
\noindent{\boldmath $ 31^{12} 28^{1} $}~
With the point set $Z_{400}$ partitioned into
 residue classes modulo $12$ for $\{0, 1, \dots, 371\}$, and
 $\{372, 373, \dots, 399\}$,
 the design is generated from

\adfLgap %ADFvfyDesignStart
$(372, 9, 106, 342)$,
$(372, 319, 305, 188)$,
$(373, 193, 248, 84)$,
$(373, 33, 317, 148)$,\adfsplit
$(374, 66, 75, 1)$,
$(374, 124, 218, 149)$,
$(375, 69, 173, 192)$,
$(375, 133, 320, 118)$,\adfsplit
$(376, 83, 99, 268)$,
$(376, 162, 290, 181)$,
$(377, 0, 242, 256)$,
$(377, 61, 93, 53)$,\adfsplit
$(378, 131, 248, 331)$,
$(378, 190, 324, 117)$,
$(379, 236, 180, 160)$,
$(379, 299, 205, 351)$,\adfsplit
$(380, 77, 112, 1)$,
$(380, 312, 213, 170)$,
$(54, 340, 301, 274)$,
$(334, 294, 266, 123)$,\adfsplit
$(369, 60, 119, 157)$,
$(325, 2, 94, 100)$,
$(75, 29, 157, 0)$,
$(366, 271, 315, 325)$,\adfsplit
$(142, 187, 207, 288)$,
$(117, 140, 50, 265)$,
$(351, 208, 158, 156)$,
$(307, 311, 157, 320)$,\adfsplit
$(127, 22, 311, 145)$,
$(331, 157, 46, 50)$,
$(75, 42, 296, 77)$,
$(289, 66, 24, 113)$,\adfsplit
$(180, 149, 70, 211)$,
$(125, 320, 95, 172)$,
$(335, 208, 42, 201)$,
$(269, 235, 0, 272)$,\adfsplit
$(71, 82, 29, 348)$,
$(23, 1, 281, 182)$,
$(307, 65, 110, 165)$,
$(298, 145, 234, 283)$,\adfsplit
$(273, 124, 329, 323)$,
$(156, 134, 271, 155)$,
$(22, 104, 161, 280)$,
$(193, 326, 252, 214)$,\adfsplit
$(281, 73, 204, 255)$,
$(0, 1, 69, 179)$,
$(0, 3, 81, 347)$,
$(0, 8, 181, 267)$,\adfsplit
$(0, 61, 297, 367)$,
$(0, 10, 101, 227)$,
$(0, 41, 251, 315)$,
$(0, 63, 121, 233)$,\adfsplit
$(0, 85, 175, 246)$,
$(0, 7, 221, 250)$,
$(0, 16, 46, 119)$,
$(0, 44, 207, 287)$,\adfsplit
$(0, 17, 135, 237)$,
$(0, 5, 194, 354)$,
$(0, 113, 174, 345)$,
$(0, 13, 140, 218)$,\adfsplit
$(0, 58, 145, 162)$,
$(0, 62, 133, 200)$,
$(0, 11, 54, 346)$,
$(0, 34, 189, 222)$,\adfsplit
$(0, 32, 102, 190)$,
$(0, 93, 186, 279)$,
$(399, 0, 124, 248)$,
$(399, 1, 125, 249)$

%ADFvfyBlocksEnd
\adfLgap \noindent by the mapping:
$x \mapsto x + 2 j \adfmod{372}$ for $x < 372$,
$x \mapsto (x - 372 + 9 j \adfmod{27}) + 372$ for $372 \le x < 399$,
$399 \mapsto 399$,
$0 \le j < 186$
 for the first 65 blocks,
$0 \le j < 93$
 for the next block,
$0 \le j < 62$
 for the last two blocks.
\ADFvfyParStart{(400, ((65, 186, ((372, 2), (27, 9), (1, 1))), (1, 93, ((372, 2), (27, 9), (1, 1))), (2, 62, ((372, 2), (27, 9), (1, 1)))), ((31, 12), (28, 1)))} %ADFvfyParEnd
% End of 31^12 28^1
%%%%%%%%%%%%%%%%%%%%%%%%%%%%%%%%%%%%%%%%%%%%%%%%%%%%%%%%%%%%%%%%%%%%%%%%%%%%%%%%%%%%%%%%%%
%%%%%%%%%%%%%%%%%%%%%%%%%%%%%%%%%%%%%%%%%%%%%%%%%%%%%%%%%%%%%%%%%%%%%%%%%%%%%%%%%%%%%%%%%%

% Charlotte:GDD4-1-3-5-mod-6-TeX-gen-A:HITS-fun:4.10
\adfDgap
%ADFvfyBlocksStart {31,31,31,31,31,31,31,31,31,22}
\noindent{\boldmath $ 31^{9} 22^{1} $}~
With the point set $Z_{301}$ partitioned into
 residue classes modulo $9$ for $\{0, 1, \dots, 278\}$, and
 $\{279, 280, \dots, 300\}$,
 the design is generated from

\adfLgap %ADFvfyDesignStart
$(297, 219, 116, 118)$,
$(279, 169, 72, 253)$,
$(279, 139, 138, 114)$,
$(279, 203, 35, 191)$,\adfsplit
$(280, 191, 0, 8)$,
$(280, 31, 149, 138)$,
$(280, 181, 115, 42)$,
$(218, 90, 3, 124)$,\adfsplit
$(267, 254, 140, 264)$,
$(95, 166, 126, 44)$,
$(185, 138, 223, 166)$,
$(241, 137, 132, 206)$,\adfsplit
$(174, 220, 79, 118)$,
$(56, 133, 147, 117)$,
$(0, 4, 52, 59)$,
$(0, 15, 37, 115)$,\adfsplit
$(0, 32, 65, 199)$,
$(0, 26, 68, 173)$,
$(0, 21, 83, 169)$,
$(0, 20, 70, 149)$,\adfsplit
$(0, 29, 89, 142)$,
$(0, 6, 49, 125)$,
$(0, 23, 67, 159)$,
$(0, 17, 58, 133)$,\adfsplit
$(300, 0, 93, 186)$

%ADFvfyBlocksEnd
\adfLgap \noindent by the mapping:
$x \mapsto x +  j \adfmod{279}$ for $x < 279$,
$x \mapsto (x - 279 + 2 j \adfmod{18}) + 279$ for $279 \le x < 297$,
$x \mapsto (x +  j \adfmod{3}) + 297$ for $297 \le x < 300$,
$300 \mapsto 300$,
$0 \le j < 279$
 for the first 24 blocks,
$0 \le j < 93$
 for the last block.
\ADFvfyParStart{(301, ((24, 279, ((279, 1), (18, 2), (3, 1), (1, 1))), (1, 93, ((279, 1), (18, 2), (3, 1), (1, 1)))), ((31, 9), (22, 1)))} %ADFvfyParEnd
% End of 31^9 22^1
%%%%%%%%%%%%%%%%%%%%%%%%%%%%%%%%%%%%%%%%%%%%%%%%%%%%%%%%%%%%%%%%%%%%%%%%%%%%%%%%%%%%%%%%%%
%%%%%%%%%%%%%%%%%%%%%%%%%%%%%%%%%%%%%%%%%%%%%%%%%%%%%%%%%%%%%%%%%%%%%%%%%%%%%%%%%%%%%%%%%%

% Charlotte:GDD4-1-3-5-mod-6-TeX-gen-A:HITS-fun:4.10
\adfDgap
%ADFvfyBlocksStart {31,31,31,31,31,31,31,31,31,28}
\noindent{\boldmath $ 31^{9} 28^{1} $}~
With the point set $Z_{307}$ partitioned into
 residue classes modulo $9$ for $\{0, 1, \dots, 278\}$, and
 $\{279, 280, \dots, 306\}$,
 the design is generated from

\adfLgap %ADFvfyDesignStart
$(279, 154, 146, 159)$,
$(279, 57, 225, 256)$,
$(279, 34, 230, 17)$,
$(280, 53, 51, 50)$,\adfsplit
$(280, 238, 57, 97)$,
$(280, 108, 244, 29)$,
$(281, 85, 172, 129)$,
$(281, 242, 86, 204)$,\adfsplit
$(281, 47, 115, 180)$,
$(54, 107, 50, 166)$,
$(75, 190, 6, 166)$,
$(198, 127, 92, 69)$,\adfsplit
$(167, 42, 179, 94)$,
$(131, 259, 152, 120)$,
$(131, 112, 105, 146)$,
$(0, 6, 16, 67)$,\adfsplit
$(0, 14, 88, 193)$,
$(0, 42, 89, 166)$,
$(0, 25, 121, 170)$,
$(0, 30, 78, 182)$,\adfsplit
$(0, 22, 50, 219)$,
$(0, 29, 75, 159)$,
$(0, 20, 76, 177)$,
$(0, 37, 92, 131)$,\adfsplit
$(0, 33, 103, 165)$,
$(306, 0, 93, 186)$

%ADFvfyBlocksEnd
\adfLgap \noindent by the mapping:
$x \mapsto x +  j \adfmod{279}$ for $x < 279$,
$x \mapsto (x - 279 + 3 j \adfmod{27}) + 279$ for $279 \le x < 306$,
$306 \mapsto 306$,
$0 \le j < 279$
 for the first 25 blocks,
$0 \le j < 93$
 for the last block.
\ADFvfyParStart{(307, ((25, 279, ((279, 1), (27, 3), (1, 1))), (1, 93, ((279, 1), (27, 3), (1, 1)))), ((31, 9), (28, 1)))} %ADFvfyParEnd
% End of 31^9 28^1
%%%%%%%%%%%%%%%%%%%%%%%%%%%%%%%%%%%%%%%%%%%%%%%%%%%%%%%%%%%%%%%%%%%%%%%%%%%%%%%%%%%%%%%%%%
%%%%%%%%%%%%%%%%%%%%%%%%%%%%%%%%%%%%%%%%%%%%%%%%%%%%%%%%%%%%%%%%%%%%%%%%%%%%%%%%%%%%%%%%%%

% Charlotte:GDD4-1-3-5-mod-6-TeX-gen-A:HITS-fun:4.10
\adfDgap
%ADFvfyBlocksStart {31,31,31,31,31,31,31,31,31,34}
\noindent{\boldmath $ 31^{9} 34^{1} $}~
With the point set $Z_{313}$ partitioned into
 residue classes modulo $9$ for $\{0, 1, \dots, 278\}$, and
 $\{279, 280, \dots, 312\}$,
 the design is generated from

\adfLgap %ADFvfyDesignStart
$(306, 73, 110, 78)$,
$(307, 129, 224, 148)$,
$(279, 135, 1, 197)$,
$(279, 239, 101, 178)$,\adfsplit
$(279, 273, 193, 141)$,
$(280, 141, 38, 129)$,
$(280, 43, 157, 32)$,
$(280, 127, 161, 18)$,\adfsplit
$(281, 141, 253, 101)$,
$(281, 121, 174, 178)$,
$(281, 143, 189, 248)$,
$(139, 260, 24, 97)$,\adfsplit
$(242, 213, 30, 85)$,
$(250, 153, 104, 65)$,
$(189, 169, 227, 275)$,
$(127, 112, 212, 111)$,\adfsplit
$(0, 2, 8, 131)$,
$(0, 10, 35, 210)$,
$(0, 22, 92, 120)$,
$(0, 13, 64, 124)$,\adfsplit
$(0, 7, 31, 208)$,
$(0, 21, 65, 140)$,
$(0, 50, 118, 192)$,
$(0, 23, 89, 130)$,\adfsplit
$(0, 3, 17, 249)$,
$(0, 26, 82, 195)$,
$(312, 0, 93, 186)$

%ADFvfyBlocksEnd
\adfLgap \noindent by the mapping:
$x \mapsto x +  j \adfmod{279}$ for $x < 279$,
$x \mapsto (x - 279 + 3 j \adfmod{27}) + 279$ for $279 \le x < 306$,
$x \mapsto (x + 2 j \adfmod{6}) + 306$ for $306 \le x < 312$,
$312 \mapsto 312$,
$0 \le j < 279$
 for the first 26 blocks,
$0 \le j < 93$
 for the last block.
\ADFvfyParStart{(313, ((26, 279, ((279, 1), (27, 3), (6, 2), (1, 1))), (1, 93, ((279, 1), (27, 3), (6, 2), (1, 1)))), ((31, 9), (34, 1)))} %ADFvfyParEnd
% End of 31^9 34^1
%%%%%%%%%%%%%%%%%%%%%%%%%%%%%%%%%%%%%%%%%%%%%%%%%%%%%%%%%%%%%%%%%%%%%%%%%%%%%%%%%%%%%%%%%%
%%%%%%%%%%%%%%%%%%%%%%%%%%%%%%%%%%%%%%%%%%%%%%%%%%%%%%%%%%%%%%%%%%%%%%%%%%%%%%%%%%%%%%%%%%

% Charlotte:GDD4-1-3-5-mod-6-TeX-gen-A:HITS-fun:4.10
\adfDgap
%ADFvfyBlocksStart {31,31,31,31,31,31,31,31,31,40}
\noindent{\boldmath $ 31^{9} 40^{1} $}~
With the point set $Z_{319}$ partitioned into
 residue classes modulo $9$ for $\{0, 1, \dots, 278\}$, and
 $\{279, 280, \dots, 318\}$,
 the design is generated from

\adfLgap %ADFvfyDesignStart
$(315, 101, 124, 258)$,
$(279, 167, 26, 78)$,
$(279, 81, 25, 172)$,
$(279, 112, 263, 102)$,\adfsplit
$(280, 50, 85, 210)$,
$(280, 206, 136, 243)$,
$(280, 133, 263, 195)$,
$(281, 84, 125, 50)$,\adfsplit
$(281, 37, 151, 96)$,
$(281, 27, 121, 137)$,
$(282, 252, 273, 227)$,
$(282, 266, 1, 34)$,\adfsplit
$(282, 168, 13, 89)$,
$(209, 3, 105, 7)$,
$(266, 227, 126, 124)$,
$(8, 176, 181, 123)$,\adfsplit
$(116, 131, 11, 0)$,
$(68, 142, 110, 111)$,
$(0, 13, 30, 196)$,
$(0, 8, 65, 201)$,\adfsplit
$(0, 6, 50, 88)$,
$(0, 7, 71, 194)$,
$(0, 12, 109, 133)$,
$(0, 28, 95, 195)$,\adfsplit
$(0, 20, 69, 219)$,
$(0, 26, 66, 218)$,
$(0, 3, 22, 51)$,
$(318, 0, 93, 186)$

%ADFvfyBlocksEnd
\adfLgap \noindent by the mapping:
$x \mapsto x +  j \adfmod{279}$ for $x < 279$,
$x \mapsto (x - 279 + 4 j \adfmod{36}) + 279$ for $279 \le x < 315$,
$x \mapsto (x +  j \adfmod{3}) + 315$ for $315 \le x < 318$,
$318 \mapsto 318$,
$0 \le j < 279$
 for the first 27 blocks,
$0 \le j < 93$
 for the last block.
\ADFvfyParStart{(319, ((27, 279, ((279, 1), (36, 4), (3, 1), (1, 1))), (1, 93, ((279, 1), (36, 4), (3, 1), (1, 1)))), ((31, 9), (40, 1)))} %ADFvfyParEnd
% End of 31^9 40^1
%%%%%%%%%%%%%%%%%%%%%%%%%%%%%%%%%%%%%%%%%%%%%%%%%%%%%%%%%%%%%%%%%%%%%%%%%%%%%%%%%%%%%%%%%%
%%%%%%%%%%%%%%%%%%%%%%%%%%%%%%%%%%%%%%%%%%%%%%%%%%%%%%%%%%%%%%%%%%%%%%%%%%%%%%%%%%%%%%%%%%

% Charlotte:GDD4-1-3-5-mod-6-TeX-gen-A:HITS-fun:4.10
\adfDgap
%ADFvfyBlocksStart {31,31,31,31,31,31,31,31,31,46}
\noindent{\boldmath $ 31^{9} 46^{1} $}~
With the point set $Z_{325}$ partitioned into
 residue classes modulo $9$ for $\{0, 1, \dots, 278\}$, and
 $\{279, 280, \dots, 324\}$,
 the design is generated from

\adfLgap %ADFvfyDesignStart
$(279, 182, 150, 75)$,
$(279, 208, 266, 238)$,
$(279, 206, 187, 252)$,
$(280, 267, 134, 75)$,\adfsplit
$(280, 235, 162, 115)$,
$(280, 95, 110, 40)$,
$(281, 103, 60, 214)$,
$(281, 131, 118, 125)$,\adfsplit
$(281, 191, 255, 252)$,
$(282, 51, 225, 104)$,
$(282, 109, 146, 43)$,
$(282, 197, 157, 246)$,\adfsplit
$(283, 255, 70, 95)$,
$(283, 91, 173, 274)$,
$(283, 242, 42, 252)$,
$(179, 277, 167, 40)$,\adfsplit
$(234, 97, 213, 229)$,
$(162, 11, 206, 49)$,
$(77, 63, 3, 139)$,
$(0, 20, 77, 106)$,\adfsplit
$(0, 17, 39, 184)$,
$(0, 4, 35, 191)$,
$(0, 8, 41, 196)$,
$(0, 24, 80, 194)$,\adfsplit
$(0, 1, 51, 212)$,
$(0, 2, 102, 150)$,
$(0, 11, 34, 175)$,
$(0, 26, 78, 208)$,\adfsplit
$(324, 0, 93, 186)$

%ADFvfyBlocksEnd
\adfLgap \noindent by the mapping:
$x \mapsto x +  j \adfmod{279}$ for $x < 279$,
$x \mapsto (x - 279 + 5 j \adfmod{45}) + 279$ for $279 \le x < 324$,
$324 \mapsto 324$,
$0 \le j < 279$
 for the first 28 blocks,
$0 \le j < 93$
 for the last block.
\ADFvfyParStart{(325, ((28, 279, ((279, 1), (45, 5), (1, 1))), (1, 93, ((279, 1), (45, 5), (1, 1)))), ((31, 9), (46, 1)))} %ADFvfyParEnd
% End of 31^9 46^1
%%%%%%%%%%%%%%%%%%%%%%%%%%%%%%%%%%%%%%%%%%%%%%%%%%%%%%%%%%%%%%%%%%%%%%%%%%%%%%%%%%%%%%%%%%
%%%%%%%%%%%%%%%%%%%%%%%%%%%%%%%%%%%%%%%%%%%%%%%%%%%%%%%%%%%%%%%%%%%%%%%%%%%%%%%%%%%%%%%%%%

% Charlotte:GDD4-1-3-5-mod-6-TeX-gen-A:HITS-fun:4.10
\adfDgap
%ADFvfyBlocksStart {31,31,31,31,31,31,31,31,31,52}
\noindent{\boldmath $ 31^{9} 52^{1} $}~
With the point set $Z_{331}$ partitioned into
 residue classes modulo $9$ for $\{0, 1, \dots, 278\}$, and
 $\{279, 280, \dots, 330\}$,
 the design is generated from

\adfLgap %ADFvfyDesignStart
$(324, 128, 76, 21)$,
$(325, 106, 204, 251)$,
$(279, 203, 61, 102)$,
$(279, 127, 54, 11)$,\adfsplit
$(279, 69, 130, 161)$,
$(280, 249, 121, 273)$,
$(280, 137, 171, 142)$,
$(280, 242, 86, 28)$,\adfsplit
$(281, 156, 179, 239)$,
$(281, 136, 88, 175)$,
$(281, 110, 225, 96)$,
$(282, 174, 68, 162)$,\adfsplit
$(282, 209, 213, 19)$,
$(282, 223, 76, 269)$,
$(283, 36, 141, 167)$,
$(283, 100, 191, 219)$,\adfsplit
$(283, 49, 115, 98)$,
$(141, 133, 215, 131)$,
$(217, 76, 60, 180)$,
$(0, 1, 21, 51)$,\adfsplit
$(0, 3, 35, 79)$,
$(0, 6, 103, 167)$,
$(0, 25, 125, 165)$,
$(0, 13, 70, 183)$,\adfsplit
$(0, 19, 78, 155)$,
$(0, 33, 71, 217)$,
$(0, 7, 22, 75)$,
$(0, 11, 67, 169)$,\adfsplit
$(0, 42, 111, 191)$,
$(330, 0, 93, 186)$

%ADFvfyBlocksEnd
\adfLgap \noindent by the mapping:
$x \mapsto x +  j \adfmod{279}$ for $x < 279$,
$x \mapsto (x - 279 + 5 j \adfmod{45}) + 279$ for $279 \le x < 324$,
$x \mapsto (x + 2 j \adfmod{6}) + 324$ for $324 \le x < 330$,
$330 \mapsto 330$,
$0 \le j < 279$
 for the first 29 blocks,
$0 \le j < 93$
 for the last block.
\ADFvfyParStart{(331, ((29, 279, ((279, 1), (45, 5), (6, 2), (1, 1))), (1, 93, ((279, 1), (45, 5), (6, 2), (1, 1)))), ((31, 9), (52, 1)))} %ADFvfyParEnd
% End of 31^9 52^1
%%%%%%%%%%%%%%%%%%%%%%%%%%%%%%%%%%%%%%%%%%%%%%%%%%%%%%%%%%%%%%%%%%%%%%%%%%%%%%%%%%%%%%%%%%
%%%%%%%%%%%%%%%%%%%%%%%%%%%%%%%%%%%%%%%%%%%%%%%%%%%%%%%%%%%%%%%%%%%%%%%%%%%%%%%%%%%%%%%%%%

% Charlotte:GDD4-1-3-5-mod-6-TeX-gen-A:HITS-fun:4.10
\adfDgap
%ADFvfyBlocksStart {31,31,31,31,31,31,31,31,31,58}
\noindent{\boldmath $ 31^{9} 58^{1} $}~
With the point set $Z_{337}$ partitioned into
 residue classes modulo $9$ for $\{0, 1, \dots, 278\}$, and
 $\{279, 280, \dots, 336\}$,
 the design is generated from

\adfLgap %ADFvfyDesignStart
$(333, 24, 158, 28)$,
$(279, 240, 130, 74)$,
$(279, 233, 253, 176)$,
$(279, 61, 147, 243)$,\adfsplit
$(280, 269, 217, 90)$,
$(280, 157, 231, 142)$,
$(280, 32, 191, 66)$,
$(281, 249, 38, 113)$,\adfsplit
$(281, 112, 138, 107)$,
$(281, 73, 241, 90)$,
$(282, 254, 276, 0)$,
$(282, 278, 154, 104)$,\adfsplit
$(282, 112, 39, 88)$,
$(283, 70, 116, 28)$,
$(283, 246, 51, 18)$,
$(283, 50, 38, 85)$,\adfsplit
$(284, 35, 196, 166)$,
$(284, 180, 140, 181)$,
$(284, 39, 2, 177)$,
$(217, 138, 94, 236)$,\adfsplit
$(0, 2, 8, 67)$,
$(0, 7, 21, 116)$,
$(0, 16, 103, 173)$,
$(0, 23, 62, 226)$,\adfsplit
$(0, 29, 112, 194)$,
$(0, 11, 66, 199)$,
$(0, 28, 60, 129)$,
$(0, 43, 107, 201)$,\adfsplit
$(0, 10, 48, 187)$,
$(0, 13, 71, 132)$,
$(336, 0, 93, 186)$

%ADFvfyBlocksEnd
\adfLgap \noindent by the mapping:
$x \mapsto x +  j \adfmod{279}$ for $x < 279$,
$x \mapsto (x - 279 + 6 j \adfmod{54}) + 279$ for $279 \le x < 333$,
$x \mapsto (x +  j \adfmod{3}) + 333$ for $333 \le x < 336$,
$336 \mapsto 336$,
$0 \le j < 279$
 for the first 30 blocks,
$0 \le j < 93$
 for the last block.
\ADFvfyParStart{(337, ((30, 279, ((279, 1), (54, 6), (3, 1), (1, 1))), (1, 93, ((279, 1), (54, 6), (3, 1), (1, 1)))), ((31, 9), (58, 1)))} %ADFvfyParEnd
% End of 31^9 58^1
%%%%%%%%%%%%%%%%%%%%%%%%%%%%%%%%%%%%%%%%%%%%%%%%%%%%%%%%%%%%%%%%%%%%%%%%%%%%%%%%%%%%%%%%%%
%%%%%%%%%%%%%%%%%%%%%%%%%%%%%%%%%%%%%%%%%%%%%%%%%%%%%%%%%%%%%%%%%%%%%%%%%%%%%%%%%%%%%%%%%%

% Charlotte:GDD4-1-3-5-mod-6-TeX-gen-A:HITS-fun:4.10
\adfDgap
%ADFvfyBlocksStart {31,31,31,31,31,31,31,31,31,64}
\noindent{\boldmath $ 31^{9} 64^{1} $}~
With the point set $Z_{343}$ partitioned into
 residue classes modulo $9$ for $\{0, 1, \dots, 278\}$, and
 $\{279, 280, \dots, 342\}$,
 the design is generated from

\adfLgap %ADFvfyDesignStart
$(279, 205, 67, 269)$,
$(279, 207, 275, 65)$,
$(279, 75, 262, 204)$,
$(280, 211, 16, 176)$,\adfsplit
$(280, 135, 267, 19)$,
$(280, 245, 174, 260)$,
$(281, 92, 105, 166)$,
$(281, 262, 205, 39)$,\adfsplit
$(281, 243, 95, 215)$,
$(282, 127, 106, 227)$,
$(282, 230, 39, 0)$,
$(282, 42, 148, 251)$,\adfsplit
$(283, 96, 257, 58)$,
$(283, 243, 262, 8)$,
$(283, 192, 155, 88)$,
$(284, 78, 25, 68)$,\adfsplit
$(284, 253, 164, 211)$,
$(284, 156, 44, 126)$,
$(285, 208, 14, 31)$,
$(285, 132, 9, 183)$,\adfsplit
$(285, 128, 88, 224)$,
$(0, 1, 4, 115)$,
$(0, 5, 34, 46)$,
$(0, 11, 33, 157)$,\adfsplit
$(0, 23, 98, 130)$,
$(0, 16, 110, 170)$,
$(0, 6, 20, 220)$,
$(0, 8, 91, 192)$,\adfsplit
$(0, 7, 55, 152)$,
$(0, 62, 128, 201)$,
$(0, 2, 26, 229)$,
$(342, 0, 93, 186)$

%ADFvfyBlocksEnd
\adfLgap \noindent by the mapping:
$x \mapsto x +  j \adfmod{279}$ for $x < 279$,
$x \mapsto (x - 279 + 7 j \adfmod{63}) + 279$ for $279 \le x < 342$,
$342 \mapsto 342$,
$0 \le j < 279$
 for the first 31 blocks,
$0 \le j < 93$
 for the last block.
\ADFvfyParStart{(343, ((31, 279, ((279, 1), (63, 7), (1, 1))), (1, 93, ((279, 1), (63, 7), (1, 1)))), ((31, 9), (64, 1)))} %ADFvfyParEnd
% End of 31^9 64^1
%%%%%%%%%%%%%%%%%%%%%%%%%%%%%%%%%%%%%%%%%%%%%%%%%%%%%%%%%%%%%%%%%%%%%%%%%%%%%%%%%%%%%%%%%%
%%%%%%%%%%%%%%%%%%%%%%%%%%%%%%%%%%%%%%%%%%%%%%%%%%%%%%%%%%%%%%%%%%%%%%%%%%%%%%%%%%%%%%%%%%

% Charlotte:GDD4-1-3-5-mod-6-TeX-gen-A:HITS-fun:4.10
\adfDgap
%ADFvfyBlocksStart {31,31,31,31,31,31,31,31,31,70}
\noindent{\boldmath $ 31^{9} 70^{1} $}~
With the point set $Z_{349}$ partitioned into
 residue classes modulo $9$ for $\{0, 1, \dots, 278\}$, and
 $\{279, 280, \dots, 348\}$,
 the design is generated from

\adfLgap %ADFvfyDesignStart
$(342, 62, 148, 12)$,
$(343, 121, 6, 53)$,
$(279, 200, 157, 223)$,
$(279, 123, 206, 201)$,\adfsplit
$(279, 126, 266, 172)$,
$(280, 127, 202, 272)$,
$(280, 233, 138, 270)$,
$(280, 230, 177, 97)$,\adfsplit
$(281, 131, 228, 241)$,
$(281, 186, 20, 112)$,
$(281, 107, 28, 90)$,
$(282, 139, 39, 98)$,\adfsplit
$(282, 158, 190, 33)$,
$(282, 250, 234, 173)$,
$(283, 31, 64, 42)$,
$(283, 34, 92, 8)$,\adfsplit
$(283, 257, 264, 243)$,
$(284, 14, 66, 2)$,
$(284, 6, 107, 178)$,
$(284, 55, 94, 243)$,\adfsplit
$(285, 179, 151, 63)$,
$(285, 200, 86, 192)$,
$(0, 1, 3, 231)$,
$(0, 6, 104, 161)$,\adfsplit
$(0, 34, 76, 306)$,
$(0, 4, 24, 109)$,
$(0, 30, 119, 150)$,
$(0, 15, 82, 142)$,\adfsplit
$(0, 35, 73, 138)$,
$(0, 40, 96, 151)$,
$(0, 10, 29, 131)$,
$(0, 25, 69, 192)$,\adfsplit
$(348, 0, 93, 186)$

%ADFvfyBlocksEnd
\adfLgap \noindent by the mapping:
$x \mapsto x +  j \adfmod{279}$ for $x < 279$,
$x \mapsto (x - 279 + 7 j \adfmod{63}) + 279$ for $279 \le x < 342$,
$x \mapsto (x + 2 j \adfmod{6}) + 342$ for $342 \le x < 348$,
$348 \mapsto 348$,
$0 \le j < 279$
 for the first 32 blocks,
$0 \le j < 93$
 for the last block.
\ADFvfyParStart{(349, ((32, 279, ((279, 1), (63, 7), (6, 2), (1, 1))), (1, 93, ((279, 1), (63, 7), (6, 2), (1, 1)))), ((31, 9), (70, 1)))} %ADFvfyParEnd
% End of 31^9 70^1
%%%%%%%%%%%%%%%%%%%%%%%%%%%%%%%%%%%%%%%%%%%%%%%%%%%%%%%%%%%%%%%%%%%%%%%%%%%%%%%%%%%%%%%%%%
%%%%%%%%%%%%%%%%%%%%%%%%%%%%%%%%%%%%%%%%%%%%%%%%%%%%%%%%%%%%%%%%%%%%%%%%%%%%%%%%%%%%%%%%%%

% Charlotte:GDD4-1-3-5-mod-6-TeX-gen-A:HITS-fun:4.10
\adfDgap
%ADFvfyBlocksStart {31,31,31,31,31,31,31,31,31,76}
\noindent{\boldmath $ 31^{9} 76^{1} $}~
With the point set $Z_{355}$ partitioned into
 residue classes modulo $9$ for $\{0, 1, \dots, 278\}$, and
 $\{279, 280, \dots, 354\}$,
 the design is generated from

\adfLgap %ADFvfyDesignStart
$(351, 272, 243, 127)$,
$(279, 197, 277, 108)$,
$(279, 209, 138, 94)$,
$(279, 208, 159, 131)$,\adfsplit
$(280, 147, 61, 242)$,
$(280, 132, 234, 110)$,
$(280, 181, 40, 140)$,
$(281, 201, 239, 232)$,\adfsplit
$(281, 182, 168, 117)$,
$(281, 40, 244, 44)$,
$(282, 50, 276, 144)$,
$(282, 247, 187, 272)$,\adfsplit
$(282, 10, 219, 143)$,
$(283, 251, 228, 167)$,
$(283, 227, 172, 69)$,
$(283, 277, 0, 166)$,\adfsplit
$(284, 255, 125, 261)$,
$(284, 52, 157, 149)$,
$(284, 213, 154, 2)$,
$(285, 201, 236, 180)$,\adfsplit
$(285, 134, 230, 42)$,
$(285, 49, 271, 61)$,
$(286, 201, 168, 50)$,
$(0, 1, 20, 157)$,\adfsplit
$(0, 3, 50, 67)$,
$(0, 15, 58, 165)$,
$(0, 11, 131, 342)$,
$(0, 39, 87, 206)$,\adfsplit
$(0, 13, 37, 191)$,
$(0, 5, 83, 109)$,
$(0, 10, 40, 237)$,
$(0, 32, 66, 205)$,\adfsplit
$(0, 16, 62, 350)$,
$(354, 0, 93, 186)$

%ADFvfyBlocksEnd
\adfLgap \noindent by the mapping:
$x \mapsto x +  j \adfmod{279}$ for $x < 279$,
$x \mapsto (x - 279 + 8 j \adfmod{72}) + 279$ for $279 \le x < 351$,
$x \mapsto (x +  j \adfmod{3}) + 351$ for $351 \le x < 354$,
$354 \mapsto 354$,
$0 \le j < 279$
 for the first 33 blocks,
$0 \le j < 93$
 for the last block.
\ADFvfyParStart{(355, ((33, 279, ((279, 1), (72, 8), (3, 1), (1, 1))), (1, 93, ((279, 1), (72, 8), (3, 1), (1, 1)))), ((31, 9), (76, 1)))} %ADFvfyParEnd
% End of 31^9 76^1
%%%%%%%%%%%%%%%%%%%%%%%%%%%%%%%%%%%%%%%%%%%%%%%%%%%%%%%%%%%%%%%%%%%%%%%%%%%%%%%%%%%%%%%%%%
%%%%%%%%%%%%%%%%%%%%%%%%%%%%%%%%%%%%%%%%%%%%%%%%%%%%%%%%%%%%%%%%%%%%%%%%%%%%%%%%%%%%%%%%%%

% Charlotte:GDD4-1-3-5-mod-6-TeX-gen-A:HITS-fun:4.10
\adfDgap
%ADFvfyBlocksStart {31,31,31,31,31,31,31,31,31,82}
\noindent{\boldmath $ 31^{9} 82^{1} $}~
With the point set $Z_{361}$ partitioned into
 residue classes modulo $9$ for $\{0, 1, \dots, 278\}$, and
 $\{279, 280, \dots, 360\}$,
 the design is generated from

\adfLgap %ADFvfyDesignStart
$(279, 232, 141, 209)$,
$(279, 10, 58, 251)$,
$(279, 153, 147, 59)$,
$(280, 64, 210, 166)$,\adfsplit
$(280, 87, 223, 167)$,
$(280, 29, 36, 188)$,
$(281, 8, 241, 65)$,
$(281, 54, 33, 22)$,\adfsplit
$(281, 140, 75, 199)$,
$(282, 76, 36, 245)$,
$(282, 136, 114, 61)$,
$(282, 251, 30, 32)$,\adfsplit
$(283, 62, 75, 265)$,
$(283, 266, 70, 235)$,
$(283, 11, 6, 216)$,
$(284, 26, 127, 99)$,\adfsplit
$(284, 202, 61, 41)$,
$(284, 29, 78, 129)$,
$(285, 32, 7, 74)$,
$(285, 265, 27, 134)$,\adfsplit
$(285, 186, 57, 73)$,
$(286, 19, 116, 0)$,
$(286, 33, 178, 50)$,
$(286, 211, 92, 156)$,\adfsplit
$(0, 1, 35, 85)$,
$(0, 8, 104, 181)$,
$(0, 12, 26, 213)$,
$(0, 62, 149, 287)$,\adfsplit
$(0, 15, 52, 157)$,
$(0, 24, 95, 156)$,
$(0, 4, 47, 158)$,
$(0, 10, 39, 341)$,\adfsplit
$(0, 30, 139, 305)$,
$(0, 3, 82, 115)$,
$(360, 0, 93, 186)$

%ADFvfyBlocksEnd
\adfLgap \noindent by the mapping:
$x \mapsto x +  j \adfmod{279}$ for $x < 279$,
$x \mapsto (x - 279 + 9 j \adfmod{81}) + 279$ for $279 \le x < 360$,
$360 \mapsto 360$,
$0 \le j < 279$
 for the first 34 blocks,
$0 \le j < 93$
 for the last block.
\ADFvfyParStart{(361, ((34, 279, ((279, 1), (81, 9), (1, 1))), (1, 93, ((279, 1), (81, 9), (1, 1)))), ((31, 9), (82, 1)))} %ADFvfyParEnd
% End of 31^9 82^1
%%%%%%%%%%%%%%%%%%%%%%%%%%%%%%%%%%%%%%%%%%%%%%%%%%%%%%%%%%%%%%%%%%%%%%%%%%%%%%%%%%%%%%%%%%
%%%%%%%%%%%%%%%%%%%%%%%%%%%%%%%%%%%%%%%%%%%%%%%%%%%%%%%%%%%%%%%%%%%%%%%%%%%%%%%%%%%%%%%%%%

% Charlotte:GDD4-1-3-5-mod-6-TeX-gen-A:HITS-fun:4.10
\adfDgap
%ADFvfyBlocksStart {31,31,31,31,31,31,31,31,31,88}
\noindent{\boldmath $ 31^{9} 88^{1} $}~
With the point set $Z_{367}$ partitioned into
 residue classes modulo $9$ for $\{0, 1, \dots, 278\}$, and
 $\{279, 280, \dots, 366\}$,
 the design is generated from

\adfLgap %ADFvfyDesignStart
$(360, 240, 115, 203)$,
$(361, 156, 191, 109)$,
$(279, 130, 154, 15)$,
$(279, 104, 52, 119)$,\adfsplit
$(279, 53, 180, 102)$,
$(280, 173, 212, 151)$,
$(280, 6, 180, 228)$,
$(280, 148, 260, 73)$,\adfsplit
$(281, 255, 106, 117)$,
$(281, 217, 67, 107)$,
$(281, 222, 2, 275)$,
$(282, 164, 28, 66)$,\adfsplit
$(282, 116, 97, 54)$,
$(282, 130, 150, 230)$,
$(283, 126, 258, 29)$,
$(283, 253, 93, 248)$,\adfsplit
$(283, 241, 125, 49)$,
$(284, 229, 59, 84)$,
$(284, 141, 252, 226)$,
$(284, 88, 74, 71)$,\adfsplit
$(285, 62, 208, 237)$,
$(285, 196, 31, 240)$,
$(285, 236, 99, 257)$,
$(286, 98, 200, 166)$,\adfsplit
$(286, 252, 178, 186)$,
$(0, 2, 58, 358)$,
$(0, 4, 107, 287)$,
$(0, 16, 71, 122)$,\adfsplit
$(0, 7, 91, 323)$,
$(0, 28, 69, 148)$,
$(0, 12, 77, 190)$,
$(0, 13, 246, 305)$,\adfsplit
$(0, 23, 83, 206)$,
$(0, 10, 42, 128)$,
$(0, 1, 31, 95)$,
$(366, 0, 93, 186)$

%ADFvfyBlocksEnd
\adfLgap \noindent by the mapping:
$x \mapsto x +  j \adfmod{279}$ for $x < 279$,
$x \mapsto (x - 279 + 9 j \adfmod{81}) + 279$ for $279 \le x < 360$,
$x \mapsto (x + 2 j \adfmod{6}) + 360$ for $360 \le x < 366$,
$366 \mapsto 366$,
$0 \le j < 279$
 for the first 35 blocks,
$0 \le j < 93$
 for the last block.
\ADFvfyParStart{(367, ((35, 279, ((279, 1), (81, 9), (6, 2), (1, 1))), (1, 93, ((279, 1), (81, 9), (6, 2), (1, 1)))), ((31, 9), (88, 1)))} %ADFvfyParEnd
% End of 31^9 88^1
%%%%%%%%%%%%%%%%%%%%%%%%%%%%%%%%%%%%%%%%%%%%%%%%%%%%%%%%%%%%%%%%%%%%%%%%%%%%%%%%%%%%%%%%%%
%%%%%%%%%%%%%%%%%%%%%%%%%%%%%%%%%%%%%%%%%%%%%%%%%%%%%%%%%%%%%%%%%%%%%%%%%%%%%%%%%%%%%%%%%%

% Charlotte:GDD4-1-3-5-mod-6-TeX-gen-A:HITS-fun:4.10
\adfDgap
%ADFvfyBlocksStart {31,31,31,31,31,31,31,31,31,94}
\noindent{\boldmath $ 31^{9} 94^{1} $}~
With the point set $Z_{373}$ partitioned into
 residue classes modulo $9$ for $\{0, 1, \dots, 278\}$, and
 $\{279, 280, \dots, 372\}$,
 the design is generated from

\adfLgap %ADFvfyDesignStart
$(369, 146, 136, 186)$,
$(279, 24, 16, 29)$,
$(279, 91, 215, 189)$,
$(279, 165, 149, 211)$,\adfsplit
$(280, 141, 29, 175)$,
$(280, 176, 237, 55)$,
$(280, 270, 233, 61)$,
$(281, 59, 164, 195)$,\adfsplit
$(281, 16, 35, 144)$,
$(281, 121, 100, 48)$,
$(282, 213, 126, 34)$,
$(282, 66, 218, 211)$,\adfsplit
$(282, 109, 278, 176)$,
$(283, 248, 197, 236)$,
$(283, 210, 252, 79)$,
$(283, 109, 40, 15)$,\adfsplit
$(284, 1, 69, 84)$,
$(284, 128, 85, 108)$,
$(284, 167, 169, 134)$,
$(285, 182, 257, 64)$,\adfsplit
$(285, 180, 102, 184)$,
$(285, 143, 160, 105)$,
$(286, 159, 250, 29)$,
$(286, 40, 153, 156)$,\adfsplit
$(286, 53, 127, 5)$,
$(287, 61, 1, 72)$,
$(287, 274, 194, 71)$,
$(0, 64, 202, 357)$,\adfsplit
$(0, 1, 29, 176)$,
$(0, 6, 120, 288)$,
$(0, 47, 96, 184)$,
$(0, 24, 65, 298)$,\adfsplit
$(0, 56, 140, 358)$,
$(0, 53, 119, 178)$,
$(0, 22, 79, 111)$,
$(0, 14, 44, 129)$,\adfsplit
$(372, 0, 93, 186)$

%ADFvfyBlocksEnd
\adfLgap \noindent by the mapping:
$x \mapsto x +  j \adfmod{279}$ for $x < 279$,
$x \mapsto (x - 279 + 10 j \adfmod{90}) + 279$ for $279 \le x < 369$,
$x \mapsto (x +  j \adfmod{3}) + 369$ for $369 \le x < 372$,
$372 \mapsto 372$,
$0 \le j < 279$
 for the first 36 blocks,
$0 \le j < 93$
 for the last block.
\ADFvfyParStart{(373, ((36, 279, ((279, 1), (90, 10), (3, 1), (1, 1))), (1, 93, ((279, 1), (90, 10), (3, 1), (1, 1)))), ((31, 9), (94, 1)))} %ADFvfyParEnd
% End of 31^9 94^1
%%%%%%%%%%%%%%%%%%%%%%%%%%%%%%%%%%%%%%%%%%%%%%%%%%%%%%%%%%%%%%%%%%%%%%%%%%%%%%%%%%%%%%%%%%
%%%%%%%%%%%%%%%%%%%%%%%%%%%%%%%%%%%%%%%%%%%%%%%%%%%%%%%%%%%%%%%%%%%%%%%%%%%%%%%%%%%%%%%%%%

% Charlotte:GDD4-1-3-5-mod-6-TeX-gen-A:HITS-fun:4.10
\adfDgap
%ADFvfyBlocksStart {31,31,31,31,31,31,31,31,31,100}
\noindent{\boldmath $ 31^{9} 100^{1} $}~
With the point set $Z_{379}$ partitioned into
 residue classes modulo $9$ for $\{0, 1, \dots, 278\}$, and
 $\{279, 280, \dots, 378\}$,
 the design is generated from

\adfLgap %ADFvfyDesignStart
$(279, 27, 83, 197)$,
$(279, 223, 68, 10)$,
$(279, 183, 213, 103)$,
$(280, 82, 213, 72)$,\adfsplit
$(280, 161, 70, 65)$,
$(280, 228, 230, 265)$,
$(281, 0, 164, 224)$,
$(281, 40, 111, 250)$,\adfsplit
$(281, 226, 59, 276)$,
$(282, 171, 98, 82)$,
$(282, 236, 102, 43)$,
$(282, 23, 184, 177)$,\adfsplit
$(283, 161, 48, 169)$,
$(283, 202, 2, 234)$,
$(283, 41, 42, 82)$,
$(284, 166, 46, 146)$,\adfsplit
$(284, 264, 52, 153)$,
$(284, 60, 8, 86)$,
$(285, 209, 156, 125)$,
$(285, 216, 168, 73)$,\adfsplit
$(285, 193, 95, 187)$,
$(286, 36, 163, 33)$,
$(286, 130, 143, 66)$,
$(286, 176, 263, 241)$,\adfsplit
$(287, 234, 222, 115)$,
$(287, 121, 17, 191)$,
$(287, 257, 66, 172)$,
$(288, 110, 266, 8)$,\adfsplit
$(0, 17, 51, 133)$,
$(0, 23, 97, 289)$,
$(0, 28, 150, 288)$,
$(0, 24, 68, 310)$,\adfsplit
$(0, 11, 49, 196)$,
$(0, 14, 39, 142)$,
$(0, 29, 75, 377)$,
$(0, 33, 76, 366)$,\adfsplit
$(0, 4, 19, 61)$,
$(378, 0, 93, 186)$

%ADFvfyBlocksEnd
\adfLgap \noindent by the mapping:
$x \mapsto x +  j \adfmod{279}$ for $x < 279$,
$x \mapsto (x - 279 + 11 j \adfmod{99}) + 279$ for $279 \le x < 378$,
$378 \mapsto 378$,
$0 \le j < 279$
 for the first 37 blocks,
$0 \le j < 93$
 for the last block.
\ADFvfyParStart{(379, ((37, 279, ((279, 1), (99, 11), (1, 1))), (1, 93, ((279, 1), (99, 11), (1, 1)))), ((31, 9), (100, 1)))} %ADFvfyParEnd
% End of 31^9 100^1
%%%%%%%%%%%%%%%%%%%%%%%%%%%%%%%%%%%%%%%%%%%%%%%%%%%%%%%%%%%%%%%%%%%%%%%%%%%%%%%%%%%%%%%%%%
%%%%%%%%%%%%%%%%%%%%%%%%%%%%%%%%%%%%%%%%%%%%%%%%%%%%%%%%%%%%%%%%%%%%%%%%%%%%%%%%%%%%%%%%%%

% Charlotte:GDD4-1-3-5-mod-6-TeX-gen-A:HITS-fun:4.10
\adfDgap
%ADFvfyBlocksStart {31,31,31,31,31,31,31,31,31,106}
\noindent{\boldmath $ 31^{9} 106^{1} $}~
With the point set $Z_{385}$ partitioned into
 residue classes modulo $9$ for $\{0, 1, \dots, 278\}$, and
 $\{279, 280, \dots, 384\}$,
 the design is generated from

\adfLgap %ADFvfyDesignStart
$(378, 219, 85, 260)$,
$(379, 6, 262, 224)$,
$(279, 38, 237, 86)$,
$(279, 168, 234, 43)$,\adfsplit
$(279, 46, 242, 130)$,
$(280, 32, 85, 11)$,
$(280, 1, 276, 169)$,
$(280, 198, 17, 147)$,\adfsplit
$(281, 7, 267, 93)$,
$(281, 50, 247, 100)$,
$(281, 171, 146, 44)$,
$(282, 224, 218, 34)$,\adfsplit
$(282, 96, 129, 50)$,
$(282, 262, 198, 22)$,
$(283, 227, 112, 44)$,
$(283, 52, 253, 84)$,\adfsplit
$(283, 50, 150, 108)$,
$(284, 143, 150, 256)$,
$(284, 236, 127, 248)$,
$(284, 273, 88, 18)$,\adfsplit
$(285, 64, 2, 221)$,
$(285, 261, 205, 274)$,
$(285, 219, 80, 60)$,
$(286, 268, 5, 64)$,\adfsplit
$(286, 119, 238, 120)$,
$(286, 197, 252, 249)$,
$(287, 110, 239, 147)$,
$(287, 100, 223, 60)$,\adfsplit
$(0, 5, 31, 342)$,
$(0, 11, 87, 310)$,
$(0, 8, 65, 202)$,
$(0, 14, 49, 344)$,\adfsplit
$(0, 47, 138, 299)$,
$(0, 10, 124, 146)$,
$(0, 2, 264, 354)$,
$(0, 29, 73, 333)$,\adfsplit
$(0, 28, 71, 366)$,
$(0, 30, 97, 131)$,
$(384, 0, 93, 186)$

%ADFvfyBlocksEnd
\adfLgap \noindent by the mapping:
$x \mapsto x +  j \adfmod{279}$ for $x < 279$,
$x \mapsto (x - 279 + 11 j \adfmod{99}) + 279$ for $279 \le x < 378$,
$x \mapsto (x + 2 j \adfmod{6}) + 378$ for $378 \le x < 384$,
$384 \mapsto 384$,
$0 \le j < 279$
 for the first 38 blocks,
$0 \le j < 93$
 for the last block.
\ADFvfyParStart{(385, ((38, 279, ((279, 1), (99, 11), (6, 2), (1, 1))), (1, 93, ((279, 1), (99, 11), (6, 2), (1, 1)))), ((31, 9), (106, 1)))} %ADFvfyParEnd
% End of 31^9 106^1
%%%%%%%%%%%%%%%%%%%%%%%%%%%%%%%%%%%%%%%%%%%%%%%%%%%%%%%%%%%%%%%%%%%%%%%%%%%%%%%%%%%%%%%%%%
%%%%%%%%%%%%%%%%%%%%%%%%%%%%%%%%%%%%%%%%%%%%%%%%%%%%%%%%%%%%%%%%%%%%%%%%%%%%%%%%%%%%%%%%%%

% Charlotte:GDD4-1-3-5-mod-6-TeX-gen-A:HITS-fun:4.10
\adfDgap
%ADFvfyBlocksStart {31,31,31,31,31,31,31,31,31,112}
\noindent{\boldmath $ 31^{9} 112^{1} $}~
With the point set $Z_{391}$ partitioned into
 residue classes modulo $9$ for $\{0, 1, \dots, 278\}$, and
 $\{279, 280, \dots, 390\}$,
 the design is generated from

\adfLgap %ADFvfyDesignStart
$(387, 102, 103, 155)$,
$(279, 223, 59, 89)$,
$(279, 42, 119, 40)$,
$(279, 271, 120, 198)$,\adfsplit
$(280, 212, 81, 43)$,
$(280, 21, 87, 130)$,
$(280, 161, 154, 182)$,
$(281, 219, 118, 2)$,\adfsplit
$(281, 171, 95, 142)$,
$(281, 44, 229, 249)$,
$(282, 135, 67, 77)$,
$(282, 114, 147, 268)$,\adfsplit
$(282, 47, 127, 152)$,
$(283, 24, 243, 110)$,
$(283, 25, 37, 112)$,
$(283, 8, 21, 32)$,\adfsplit
$(284, 157, 205, 68)$,
$(284, 172, 267, 138)$,
$(284, 209, 224, 117)$,
$(285, 132, 193, 201)$,\adfsplit
$(285, 88, 65, 104)$,
$(285, 116, 180, 271)$,
$(286, 196, 5, 19)$,
$(286, 54, 76, 206)$,\adfsplit
$(286, 110, 210, 51)$,
$(287, 178, 138, 55)$,
$(287, 38, 170, 266)$,
$(287, 229, 90, 6)$,\adfsplit
$(288, 182, 224, 127)$,
$(0, 17, 49, 214)$,
$(0, 3, 122, 288)$,
$(0, 70, 181, 324)$,\adfsplit
$(0, 19, 104, 301)$,
$(0, 6, 235, 302)$,
$(0, 5, 118, 337)$,
$(0, 46, 103, 326)$,\adfsplit
$(0, 26, 67, 138)$,
$(0, 37, 143, 362)$,
$(0, 4, 35, 361)$,
$(390, 0, 93, 186)$

%ADFvfyBlocksEnd
\adfLgap \noindent by the mapping:
$x \mapsto x +  j \adfmod{279}$ for $x < 279$,
$x \mapsto (x - 279 + 12 j \adfmod{108}) + 279$ for $279 \le x < 387$,
$x \mapsto (x +  j \adfmod{3}) + 387$ for $387 \le x < 390$,
$390 \mapsto 390$,
$0 \le j < 279$
 for the first 39 blocks,
$0 \le j < 93$
 for the last block.
\ADFvfyParStart{(391, ((39, 279, ((279, 1), (108, 12), (3, 1), (1, 1))), (1, 93, ((279, 1), (108, 12), (3, 1), (1, 1)))), ((31, 9), (112, 1)))} %ADFvfyParEnd
% End of 31^9 112^1
%%%%%%%%%%%%%%%%%%%%%%%%%%%%%%%%%%%%%%%%%%%%%%%%%%%%%%%%%%%%%%%%%%%%%%%%%%%%%%%%%%%%%%%%%%
%%%%%%%%%%%%%%%%%%%%%%%%%%%%%%%%%%%%%%%%%%%%%%%%%%%%%%%%%%%%%%%%%%%%%%%%%%%%%%%%%%%%%%%%%%

% Charlotte:GDD4-1-3-5-mod-6-TeX-gen-A:HITS-fun:4.10
\adfDgap
%ADFvfyBlocksStart {31,31,31,31,31,31,31,31,31,118}
\noindent{\boldmath $ 31^{9} 118^{1} $}~
With the point set $Z_{397}$ partitioned into
 residue classes modulo $9$ for $\{0, 1, \dots, 278\}$, and
 $\{279, 280, \dots, 396\}$,
 the design is generated from

\adfLgap %ADFvfyDesignStart
$(279, 102, 34, 132)$,
$(279, 137, 152, 58)$,
$(279, 180, 235, 239)$,
$(280, 186, 32, 115)$,\adfsplit
$(280, 119, 235, 103)$,
$(280, 36, 17, 102)$,
$(281, 112, 54, 192)$,
$(281, 105, 134, 182)$,\adfsplit
$(281, 277, 109, 86)$,
$(282, 139, 259, 26)$,
$(282, 140, 83, 45)$,
$(282, 253, 141, 66)$,\adfsplit
$(283, 204, 190, 184)$,
$(283, 228, 155, 179)$,
$(283, 169, 41, 63)$,
$(284, 35, 230, 42)$,\adfsplit
$(284, 217, 75, 189)$,
$(284, 52, 157, 182)$,
$(285, 238, 56, 77)$,
$(285, 205, 162, 188)$,\adfsplit
$(285, 37, 78, 156)$,
$(286, 137, 224, 14)$,
$(286, 19, 198, 237)$,
$(286, 124, 42, 274)$,\adfsplit
$(287, 50, 240, 196)$,
$(287, 63, 64, 31)$,
$(287, 110, 107, 120)$,
$(288, 215, 159, 164)$,\adfsplit
$(288, 133, 257, 135)$,
$(288, 130, 264, 190)$,
$(0, 34, 86, 289)$,
$(0, 40, 110, 354)$,\adfsplit
$(0, 35, 136, 380)$,
$(0, 64, 131, 303)$,
$(0, 12, 121, 343)$,
$(0, 37, 140, 381)$,\adfsplit
$(0, 31, 107, 329)$,
$(0, 50, 115, 177)$,
$(0, 8, 104, 382)$,
$(0, 11, 53, 356)$,\adfsplit
$(396, 0, 93, 186)$

%ADFvfyBlocksEnd
\adfLgap \noindent by the mapping:
$x \mapsto x +  j \adfmod{279}$ for $x < 279$,
$x \mapsto (x - 279 + 13 j \adfmod{117}) + 279$ for $279 \le x < 396$,
$396 \mapsto 396$,
$0 \le j < 279$
 for the first 40 blocks,
$0 \le j < 93$
 for the last block.
\ADFvfyParStart{(397, ((40, 279, ((279, 1), (117, 13), (1, 1))), (1, 93, ((279, 1), (117, 13), (1, 1)))), ((31, 9), (118, 1)))} %ADFvfyParEnd
% End of 31^9 118^1
%%%%%%%%%%%%%%%%%%%%%%%%%%%%%%%%%%%%%%%%%%%%%%%%%%%%%%%%%%%%%%%%%%%%%%%%%%%%%%%%%%%%%%%%%%
%%%%%%%%%%%%%%%%%%%%%%%%%%%%%%%%%%%%%%%%%%%%%%%%%%%%%%%%%%%%%%%%%%%%%%%%%%%%%%%%%%%%%%%%%%

%%%%%%%%%%%%%%%%%%%%%%%%%%%%%%%%%%%%%%%%%%%%%%%%%%%%%%%%%%%%%%%%%%%%%%%%%%%%%%%%%%%%%%%%%%
%%%%%%%%%%%%%%%%%%%%%%%%%%%%%%%%%%%%%%%%%%%%%%%%%%%%%%%%%%%%%%%%%%%%%%%%%%%%%%%%%%%%%%%%%%
\section{4-GDDs for the proof of Lemma \ref{lem:4-GDD 35^u m^1}}
\label{app:4-GDD 35^u m^1}
\adfnull{
$ 35^{12} 29^1 $,
$ 35^{12} 32^1 $,
$ 35^9 26^1 $,
$ 35^9 32^1 $,
$ 35^9 38^1 $,
$ 35^9 44^1 $,
$ 35^9 62^1 $,
$ 35^9 68^1 $,
$ 35^9 74^1 $,
$ 35^9 86^1 $,
$ 35^9 92^1 $,
$ 35^9 104^1 $,
$ 35^9 116^1 $,
$ 35^9 122^1 $,
$ 35^9 128^1 $ and
$ 35^9 134^1 $.
}

% Charlotte:GDD4-1-3-5-mod-6-TeX-gen-A:HITS-fun:4.10
\adfDgap
%ADFvfyBlocksStart {35,35,35,35,35,35,35,35,35,35,35,35,29}
\noindent{\boldmath $ 35^{12} 29^{1} $}~
With the point set $Z_{449}$ partitioned into
 residue classes modulo $12$ for $\{0, 1, \dots, 419\}$, and
 $\{420, 421, \dots, 448\}$,
 the design is generated from

\adfLgap %ADFvfyDesignStart
$(420, 0, 1, 2)$,
$(421, 0, 139, 281)$,
$(422, 0, 142, 140)$,
$(423, 0, 280, 419)$,\adfsplit
$(424, 0, 418, 278)$,
$(425, 0, 4, 11)$,
$(426, 0, 7, 416)$,
$(427, 0, 409, 413)$,\adfsplit
$(428, 0, 10, 23)$,
$(429, 0, 13, 410)$,
$(430, 0, 397, 407)$,
$(431, 0, 16, 35)$,\adfsplit
$(432, 0, 19, 404)$,
$(433, 0, 385, 401)$,
$(434, 0, 22, 5)$,
$(435, 0, 403, 398)$,\adfsplit
$(436, 0, 415, 17)$,
$(437, 0, 25, 53)$,
$(438, 0, 28, 395)$,
$(439, 0, 367, 392)$,\adfsplit
$(440, 0, 31, 65)$,
$(441, 0, 34, 389)$,
$(442, 0, 355, 386)$,
$(443, 0, 37, 8)$,\adfsplit
$(444, 0, 391, 383)$,
$(445, 0, 412, 29)$,
$(446, 0, 40, 14)$,
$(447, 0, 394, 380)$,\adfsplit
$(448, 0, 406, 26)$,
$(136, 182, 75, 13)$,
$(267, 54, 383, 34)$,
$(106, 385, 176, 290)$,\adfsplit
$(414, 159, 277, 367)$,
$(219, 286, 408, 390)$,
$(54, 244, 3, 253)$,
$(313, 416, 41, 214)$,\adfsplit
$(221, 143, 331, 298)$,
$(233, 120, 405, 362)$,
$(220, 363, 307, 275)$,
$(172, 151, 214, 289)$,\adfsplit
$(150, 303, 385, 223)$,
$(35, 272, 322, 76)$,
$(64, 353, 391, 194)$,
$(0, 3, 57, 295)$,\adfsplit
$(0, 6, 152, 167)$,
$(0, 30, 106, 335)$,
$(0, 52, 164, 266)$,
$(0, 64, 150, 239)$,\adfsplit
$(0, 68, 149, 294)$,
$(0, 39, 186, 244)$,
$(0, 49, 160, 269)$,
$(0, 59, 193, 262)$,\adfsplit
$(0, 66, 163, 320)$,
$(0, 44, 127, 328)$,
$(0, 79, 177, 301)$,
$(0, 27, 101, 326)$,\adfsplit
$(0, 105, 210, 315)$

%ADFvfyBlocksEnd
\adfLgap \noindent by the mapping:
$x \mapsto x + 3 j \adfmod{420}$ for $x < 420$,
$x \mapsto x$ for $x \ge 420$,
$0 \le j < 140$
 for the first 29 blocks;
$x \mapsto x +  j \adfmod{420}$ for $x < 420$,
$x \mapsto x$ for $x \ge 420$,
$0 \le j < 420$
 for the next 27 blocks;
$x \mapsto x + 2 j \adfmod{420}$ for $x < 420$,
$x \mapsto x$ for $x \ge 420$,
$0 \le j < 105$
 for the last block.
\ADFvfyParStart{(449, ((29, 140, ((420, 3), (29, 29))), (27, 420, ((420, 1), (29, 29))), (1, 105, ((420, 2), (29, 29)))), ((35, 12), (29, 1)))} %ADFvfyParEnd
% End of 35^12 29^1
%%%%%%%%%%%%%%%%%%%%%%%%%%%%%%%%%%%%%%%%%%%%%%%%%%%%%%%%%%%%%%%%%%%%%%%%%%%%%%%%%%%%%%%%%%
%%%%%%%%%%%%%%%%%%%%%%%%%%%%%%%%%%%%%%%%%%%%%%%%%%%%%%%%%%%%%%%%%%%%%%%%%%%%%%%%%%%%%%%%%%

% Charlotte:GDD4-1-3-5-mod-6-TeX-gen-A:HITS-fun:4.10
\adfDgap
%ADFvfyBlocksStart {35,35,35,35,35,35,35,35,35,35,35,35,32}
\noindent{\boldmath $ 35^{12} 32^{1} $}~
With the point set $Z_{452}$ partitioned into
 residue classes modulo $12$ for $\{0, 1, \dots, 419\}$, and
 $\{420, 421, \dots, 451\}$,
 the design is generated from

\adfLgap %ADFvfyDesignStart
$(420, 0, 1, 2)$,
$(421, 0, 139, 281)$,
$(422, 0, 142, 140)$,
$(423, 0, 280, 419)$,\adfsplit
$(424, 0, 418, 278)$,
$(425, 0, 4, 11)$,
$(426, 0, 7, 416)$,
$(427, 0, 409, 413)$,\adfsplit
$(428, 0, 10, 23)$,
$(429, 0, 13, 410)$,
$(430, 0, 397, 407)$,
$(431, 0, 16, 35)$,\adfsplit
$(432, 0, 19, 404)$,
$(433, 0, 385, 401)$,
$(434, 0, 22, 5)$,
$(435, 0, 403, 398)$,\adfsplit
$(436, 0, 415, 17)$,
$(437, 0, 25, 53)$,
$(438, 0, 28, 395)$,
$(439, 0, 367, 392)$,\adfsplit
$(440, 0, 31, 65)$,
$(441, 0, 34, 389)$,
$(442, 0, 355, 386)$,
$(443, 0, 37, 8)$,\adfsplit
$(444, 0, 391, 383)$,
$(445, 0, 412, 29)$,
$(446, 0, 40, 14)$,
$(447, 0, 394, 380)$,\adfsplit
$(448, 0, 406, 26)$,
$(449, 0, 43, 89)$,
$(450, 0, 46, 377)$,
$(451, 0, 331, 374)$,\adfsplit
$(214, 138, 292, 361)$,
$(202, 107, 69, 312)$,
$(324, 187, 69, 262)$,
$(417, 322, 148, 307)$,\adfsplit
$(375, 33, 149, 191)$,
$(107, 283, 410, 51)$,
$(155, 202, 18, 380)$,
$(242, 210, 93, 360)$,\adfsplit
$(357, 256, 37, 214)$,
$(392, 216, 153, 85)$,
$(394, 291, 412, 119)$,
$(326, 306, 257, 154)$,\adfsplit
$(398, 299, 68, 365)$,
$(383, 254, 109, 161)$,
$(236, 182, 288, 403)$,
$(302, 87, 246, 376)$,\adfsplit
$(385, 298, 216, 63)$,
$(245, 90, 364, 134)$,
$(210, 119, 267, 160)$,
$(240, 106, 279, 49)$,\adfsplit
$(42, 175, 93, 327)$,
$(182, 0, 284, 220)$,
$(359, 33, 15, 344)$,
$(150, 37, 220, 403)$,\adfsplit
$(311, 188, 375, 129)$,
$(408, 328, 221, 325)$,
$(415, 155, 222, 97)$,
$(60, 211, 33, 325)$,\adfsplit
$(402, 178, 276, 237)$,
$(222, 106, 397, 347)$,
$(327, 114, 189, 336)$,
$(26, 120, 205, 359)$,\adfsplit
$(105, 137, 167, 368)$,
$(0, 9, 124, 349)$,
$(0, 3, 104, 203)$,
$(0, 71, 122, 292)$,\adfsplit
$(0, 63, 160, 272)$,
$(0, 86, 207, 277)$,
$(0, 77, 186, 347)$,
$(0, 67, 73, 352)$,\adfsplit
$(1, 21, 113, 315)$,
$(0, 27, 163, 251)$,
$(0, 21, 66, 232)$,
$(0, 45, 257, 343)$,\adfsplit
$(0, 30, 149, 279)$,
$(0, 49, 93, 263)$,
$(0, 33, 114, 206)$,
$(0, 6, 47, 208)$,\adfsplit
$(0, 97, 158, 320)$,
$(0, 235, 309, 399)$,
$(0, 79, 164, 245)$,
$(0, 83, 138, 332)$,\adfsplit
$(0, 55, 217, 341)$,
$(0, 105, 210, 315)$

%ADFvfyBlocksEnd
\adfLgap \noindent by the mapping:
$x \mapsto x + 3 j \adfmod{420}$ for $x < 420$,
$x \mapsto x$ for $x \ge 420$,
$0 \le j < 140$
 for the first 32 blocks;
$x \mapsto x + 2 j \adfmod{420}$ for $x < 420$,
$x \mapsto x$ for $x \ge 420$,
$0 \le j < 210$
 for the next 53 blocks,
$0 \le j < 105$
 for the last block.
\ADFvfyParStart{(452, ((32, 140, ((420, 3), (32, 32))), (53, 210, ((420, 2), (32, 32))), (1, 105, ((420, 2), (32, 32)))), ((35, 12), (32, 1)))} %ADFvfyParEnd
% End of 35^12 32^1
%%%%%%%%%%%%%%%%%%%%%%%%%%%%%%%%%%%%%%%%%%%%%%%%%%%%%%%%%%%%%%%%%%%%%%%%%%%%%%%%%%%%%%%%%%
%%%%%%%%%%%%%%%%%%%%%%%%%%%%%%%%%%%%%%%%%%%%%%%%%%%%%%%%%%%%%%%%%%%%%%%%%%%%%%%%%%%%%%%%%%

% Charlotte:GDD4-1-3-5-mod-6-TeX-gen-A:HITS-fun:4.10
\adfDgap
%ADFvfyBlocksStart {35,35,35,35,35,35,35,35,35,26}
\noindent{\boldmath $ 35^{9} 26^{1} $}~
With the point set $Z_{341}$ partitioned into
 residue classes modulo $9$ for $\{0, 1, \dots, 314\}$, and
 $\{315, 316, \dots, 340\}$,
 the design is generated from

\adfLgap %ADFvfyDesignStart
$(339, 0, 1, 2)$,
$(340, 0, 106, 212)$,
$(315, 237, 232, 155)$,
$(316, 276, 109, 74)$,\adfsplit
$(317, 294, 175, 290)$,
$(318, 226, 87, 194)$,
$(319, 147, 106, 158)$,
$(320, 176, 213, 160)$,\adfsplit
$(321, 242, 229, 6)$,
$(322, 307, 264, 65)$,
$(116, 86, 259, 213)$,
$(218, 152, 50, 84)$,\adfsplit
$(305, 219, 284, 139)$,
$(164, 141, 7, 238)$,
$(59, 8, 164, 187)$,
$(248, 0, 110, 49)$,\adfsplit
$(268, 131, 10, 237)$,
$(0, 3, 22, 28)$,
$(0, 4, 98, 124)$,
$(0, 10, 34, 78)$,\adfsplit
$(0, 12, 67, 123)$,
$(0, 32, 128, 183)$,
$(0, 20, 89, 150)$,
$(0, 5, 76, 227)$,\adfsplit
$(0, 41, 101, 215)$,
$(0, 17, 104, 163)$,
$(0, 14, 53, 129)$,
$(0, 33, 75, 125)$,\adfsplit
$(0, 15, 127, 175)$

%ADFvfyBlocksEnd
\adfLgap \noindent by the mapping:
$x \mapsto x \oplus (3 j)$ for $x < 315$,
$x \mapsto (x - 315 + 8 j \adfmod{24}) + 315$ for $315 \le x < 339$,
$x \mapsto x$ for $x \ge 339$,
$0 \le j < 105$
 for the first two blocks;
$x \mapsto x \oplus j$ for $x < 315$,
$x \mapsto (x - 315 + 8 j \adfmod{24}) + 315$ for $315 \le x < 339$,
$x \mapsto x$ for $x \ge 339$,
$0 \le j < 315$
 for the last 27 blocks.
\ADFvfyParStart{(341, ((2, 105, ((315, 3, (105, 3)), (24, 8), (2, 2))), (27, 315, ((315, 1, (105, 3)), (24, 8), (2, 2)))), ((35, 9), (26, 1)))} %ADFvfyParEnd
% End of 35^9 26^1
%%%%%%%%%%%%%%%%%%%%%%%%%%%%%%%%%%%%%%%%%%%%%%%%%%%%%%%%%%%%%%%%%%%%%%%%%%%%%%%%%%%%%%%%%%
%%%%%%%%%%%%%%%%%%%%%%%%%%%%%%%%%%%%%%%%%%%%%%%%%%%%%%%%%%%%%%%%%%%%%%%%%%%%%%%%%%%%%%%%%%

% Charlotte:GDD4-1-3-5-mod-6-TeX-gen-A:HITS-fun:4.10
\adfDgap
%ADFvfyBlocksStart {35,35,35,35,35,35,35,35,35,32}
\noindent{\boldmath $ 35^{9} 32^{1} $}~
With the point set $Z_{347}$ partitioned into
 residue classes modulo $9$ for $\{0, 1, \dots, 314\}$, and
 $\{315, 316, \dots, 346\}$,
 the design is generated from

\adfLgap %ADFvfyDesignStart
$(345, 0, 1, 2)$,
$(346, 0, 106, 212)$,
$(315, 184, 179, 108)$,
$(316, 108, 274, 56)$,\adfsplit
$(317, 283, 147, 164)$,
$(318, 117, 136, 56)$,
$(319, 129, 95, 283)$,
$(320, 120, 100, 251)$,\adfsplit
$(321, 20, 103, 219)$,
$(322, 107, 261, 238)$,
$(323, 203, 111, 10)$,
$(324, 159, 269, 229)$,\adfsplit
$(246, 150, 71, 207)$,
$(230, 188, 164, 136)$,
$(266, 124, 146, 189)$,
$(54, 65, 165, 212)$,\adfsplit
$(302, 303, 271, 297)$,
$(248, 164, 139, 195)$,
$(0, 7, 93, 105)$,
$(0, 3, 16, 141)$,\adfsplit
$(0, 10, 61, 230)$,
$(0, 30, 78, 163)$,
$(0, 50, 112, 165)$,
$(0, 20, 41, 256)$,\adfsplit
$(0, 32, 73, 240)$,
$(0, 38, 102, 197)$,
$(0, 15, 58, 145)$,
$(0, 14, 137, 183)$,\adfsplit
$(0, 33, 113, 227)$,
$(0, 60, 129, 203)$

%ADFvfyBlocksEnd
\adfLgap \noindent by the mapping:
$x \mapsto x \oplus (3 j)$ for $x < 315$,
$x \mapsto (x - 315 + 10 j \adfmod{30}) + 315$ for $315 \le x < 345$,
$x \mapsto x$ for $x \ge 345$,
$0 \le j < 105$
 for the first two blocks;
$x \mapsto x \oplus j$ for $x < 315$,
$x \mapsto (x - 315 + 10 j \adfmod{30}) + 315$ for $315 \le x < 345$,
$x \mapsto x$ for $x \ge 345$,
$0 \le j < 315$
 for the last 28 blocks.
\ADFvfyParStart{(347, ((2, 105, ((315, 3, (105, 3)), (30, 10), (2, 2))), (28, 315, ((315, 1, (105, 3)), (30, 10), (2, 2)))), ((35, 9), (32, 1)))} %ADFvfyParEnd
% End of 35^9 32^1
%%%%%%%%%%%%%%%%%%%%%%%%%%%%%%%%%%%%%%%%%%%%%%%%%%%%%%%%%%%%%%%%%%%%%%%%%%%%%%%%%%%%%%%%%%
%%%%%%%%%%%%%%%%%%%%%%%%%%%%%%%%%%%%%%%%%%%%%%%%%%%%%%%%%%%%%%%%%%%%%%%%%%%%%%%%%%%%%%%%%%

% Charlotte:GDD4-1-3-5-mod-6-TeX-gen-A:HITS-fun:4.10
\adfDgap
%ADFvfyBlocksStart {35,35,35,35,35,35,35,35,35,38}
\noindent{\boldmath $ 35^{9} 38^{1} $}~
With the point set $Z_{353}$ partitioned into
 residue classes modulo $9$ for $\{0, 1, \dots, 314\}$, and
 $\{315, 316, \dots, 352\}$,
 the design is generated from

\adfLgap %ADFvfyDesignStart
$(351, 0, 1, 2)$,
$(352, 0, 106, 212)$,
$(315, 186, 200, 85)$,
$(316, 243, 4, 92)$,\adfsplit
$(317, 75, 245, 214)$,
$(318, 219, 59, 112)$,
$(319, 162, 304, 140)$,
$(320, 103, 201, 200)$,\adfsplit
$(321, 128, 67, 147)$,
$(322, 291, 97, 41)$,
$(323, 14, 45, 229)$,
$(324, 257, 1, 9)$,\adfsplit
$(325, 11, 79, 228)$,
$(326, 3, 23, 253)$,
$(230, 157, 181, 134)$,
$(267, 27, 129, 109)$,\adfsplit
$(228, 50, 288, 238)$,
$(300, 312, 27, 2)$,
$(221, 292, 20, 188)$,
$(0, 3, 32, 69)$,\adfsplit
$(0, 6, 58, 178)$,
$(0, 8, 118, 133)$,
$(0, 21, 67, 186)$,
$(0, 35, 123, 209)$,\adfsplit
$(0, 41, 89, 173)$,
$(0, 19, 112, 204)$,
$(0, 7, 44, 228)$,
$(0, 38, 143, 194)$,\adfsplit
$(0, 47, 125, 202)$,
$(0, 16, 55, 119)$,
$(0, 13, 92, 149)$

%ADFvfyBlocksEnd
\adfLgap \noindent by the mapping:
$x \mapsto x \oplus (3 j)$ for $x < 315$,
$x \mapsto (x - 315 + 12 j \adfmod{36}) + 315$ for $315 \le x < 351$,
$x \mapsto x$ for $x \ge 351$,
$0 \le j < 105$
 for the first two blocks;
$x \mapsto x \oplus j$ for $x < 315$,
$x \mapsto (x - 315 + 12 j \adfmod{36}) + 315$ for $315 \le x < 351$,
$x \mapsto x$ for $x \ge 351$,
$0 \le j < 315$
 for the last 29 blocks.
\ADFvfyParStart{(353, ((2, 105, ((315, 3, (105, 3)), (36, 12), (2, 2))), (29, 315, ((315, 1, (105, 3)), (36, 12), (2, 2)))), ((35, 9), (38, 1)))} %ADFvfyParEnd
% End of 35^9 38^1
%%%%%%%%%%%%%%%%%%%%%%%%%%%%%%%%%%%%%%%%%%%%%%%%%%%%%%%%%%%%%%%%%%%%%%%%%%%%%%%%%%%%%%%%%%
%%%%%%%%%%%%%%%%%%%%%%%%%%%%%%%%%%%%%%%%%%%%%%%%%%%%%%%%%%%%%%%%%%%%%%%%%%%%%%%%%%%%%%%%%%

% Charlotte:GDD4-1-3-5-mod-6-TeX-gen-A:HITS-fun:4.10
\adfDgap
%ADFvfyBlocksStart {35,35,35,35,35,35,35,35,35,44}
\noindent{\boldmath $ 35^{9} 44^{1} $}~
With the point set $Z_{359}$ partitioned into
 residue classes modulo $9$ for $\{0, 1, \dots, 314\}$, and
 $\{315, 316, \dots, 358\}$,
 the design is generated from

\adfLgap %ADFvfyDesignStart
$(357, 0, 1, 2)$,
$(358, 0, 106, 212)$,
$(315, 246, 49, 311)$,
$(316, 187, 72, 119)$,\adfsplit
$(317, 249, 100, 308)$,
$(318, 145, 155, 18)$,
$(319, 169, 203, 108)$,
$(320, 206, 310, 207)$,\adfsplit
$(321, 221, 273, 142)$,
$(322, 85, 242, 222)$,
$(323, 137, 216, 292)$,
$(324, 119, 69, 28)$,\adfsplit
$(325, 81, 151, 92)$,
$(326, 57, 244, 17)$,
$(327, 59, 157, 231)$,
$(328, 0, 205, 92)$,\adfsplit
$(239, 236, 265, 46)$,
$(122, 1, 25, 157)$,
$(168, 67, 309, 7)$,
$(93, 139, 232, 263)$,\adfsplit
$(0, 5, 26, 155)$,
$(0, 6, 74, 111)$,
$(0, 7, 22, 145)$,
$(0, 33, 84, 142)$,\adfsplit
$(0, 42, 94, 237)$,
$(0, 12, 98, 114)$,
$(0, 8, 28, 85)$,
$(0, 14, 48, 89)$,\adfsplit
$(0, 25, 64, 228)$,
$(0, 30, 130, 199)$,
$(0, 17, 49, 196)$,
$(0, 16, 66, 251)$

%ADFvfyBlocksEnd
\adfLgap \noindent by the mapping:
$x \mapsto x \oplus (3 j)$ for $x < 315$,
$x \mapsto (x - 315 + 14 j \adfmod{42}) + 315$ for $315 \le x < 357$,
$x \mapsto x$ for $x \ge 357$,
$0 \le j < 105$
 for the first two blocks;
$x \mapsto x \oplus j$ for $x < 315$,
$x \mapsto (x - 315 + 14 j \adfmod{42}) + 315$ for $315 \le x < 357$,
$x \mapsto x$ for $x \ge 357$,
$0 \le j < 315$
 for the last 30 blocks.
\ADFvfyParStart{(359, ((2, 105, ((315, 3, (105, 3)), (42, 14), (2, 2))), (30, 315, ((315, 1, (105, 3)), (42, 14), (2, 2)))), ((35, 9), (44, 1)))} %ADFvfyParEnd
% End of 35^9 44^1
%%%%%%%%%%%%%%%%%%%%%%%%%%%%%%%%%%%%%%%%%%%%%%%%%%%%%%%%%%%%%%%%%%%%%%%%%%%%%%%%%%%%%%%%%%
%%%%%%%%%%%%%%%%%%%%%%%%%%%%%%%%%%%%%%%%%%%%%%%%%%%%%%%%%%%%%%%%%%%%%%%%%%%%%%%%%%%%%%%%%%

% Charlotte:GDD4-1-3-5-mod-6-TeX-gen-A:HITS-fun:4.10
\adfDgap
%ADFvfyBlocksStart {35,35,35,35,35,35,35,35,35,62}
\noindent{\boldmath $ 35^{9} 62^{1} $}~
With the point set $Z_{377}$ partitioned into
 residue classes modulo $9$ for $\{0, 1, \dots, 314\}$, and
 $\{315, 316, \dots, 376\}$,
 the design is generated from

\adfLgap %ADFvfyDesignStart
$(375, 0, 1, 2)$,
$(376, 0, 106, 212)$,
$(315, 107, 247, 54)$,
$(316, 15, 20, 265)$,\adfsplit
$(317, 138, 152, 130)$,
$(318, 60, 46, 164)$,
$(319, 124, 156, 275)$,
$(320, 234, 260, 28)$,\adfsplit
$(321, 87, 44, 115)$,
$(322, 48, 142, 47)$,
$(323, 314, 187, 192)$,
$(324, 215, 175, 312)$,\adfsplit
$(325, 133, 293, 48)$,
$(326, 190, 252, 134)$,
$(327, 87, 308, 118)$,
$(328, 10, 138, 215)$,\adfsplit
$(329, 236, 312, 100)$,
$(330, 295, 213, 47)$,
$(331, 127, 168, 110)$,
$(332, 66, 11, 91)$,\adfsplit
$(333, 281, 201, 109)$,
$(334, 133, 78, 272)$,
$(13, 242, 308, 284)$,
$(0, 3, 19, 87)$,\adfsplit
$(0, 6, 13, 264)$,
$(0, 30, 132, 184)$,
$(0, 12, 88, 177)$,
$(0, 37, 148, 186)$,\adfsplit
$(0, 43, 96, 201)$,
$(0, 10, 101, 130)$,
$(0, 33, 145, 192)$,
$(0, 15, 49, 75)$,\adfsplit
$(0, 21, 137, 176)$,
$(0, 48, 110, 222)$,
$(0, 38, 107, 185)$

%ADFvfyBlocksEnd
\adfLgap \noindent by the mapping:
$x \mapsto x \oplus (3 j)$ for $x < 315$,
$x \mapsto (x - 315 + 20 j \adfmod{60}) + 315$ for $315 \le x < 375$,
$x \mapsto x$ for $x \ge 375$,
$0 \le j < 105$
 for the first two blocks;
$x \mapsto x \oplus j$ for $x < 315$,
$x \mapsto (x - 315 + 20 j \adfmod{60}) + 315$ for $315 \le x < 375$,
$x \mapsto x$ for $x \ge 375$,
$0 \le j < 315$
 for the last 33 blocks.
\ADFvfyParStart{(377, ((2, 105, ((315, 3, (105, 3)), (60, 20), (2, 2))), (33, 315, ((315, 1, (105, 3)), (60, 20), (2, 2)))), ((35, 9), (62, 1)))} %ADFvfyParEnd
% End of 35^9 62^1
%%%%%%%%%%%%%%%%%%%%%%%%%%%%%%%%%%%%%%%%%%%%%%%%%%%%%%%%%%%%%%%%%%%%%%%%%%%%%%%%%%%%%%%%%%
%%%%%%%%%%%%%%%%%%%%%%%%%%%%%%%%%%%%%%%%%%%%%%%%%%%%%%%%%%%%%%%%%%%%%%%%%%%%%%%%%%%%%%%%%%

% Charlotte:GDD4-1-3-5-mod-6-TeX-gen-A:HITS-fun:4.10
\adfDgap
%ADFvfyBlocksStart {35,35,35,35,35,35,35,35,35,68}
\noindent{\boldmath $ 35^{9} 68^{1} $}~
With the point set $Z_{383}$ partitioned into
 residue classes modulo $9$ for $\{0, 1, \dots, 314\}$, and
 $\{315, 316, \dots, 382\}$,
 the design is generated from

\adfLgap %ADFvfyDesignStart
$(381, 0, 1, 2)$,
$(382, 0, 106, 212)$,
$(315, 47, 6, 109)$,
$(316, 53, 234, 250)$,\adfsplit
$(317, 14, 151, 279)$,
$(318, 34, 92, 135)$,
$(319, 179, 154, 168)$,
$(320, 264, 23, 40)$,\adfsplit
$(321, 93, 77, 106)$,
$(322, 38, 103, 72)$,
$(323, 18, 256, 137)$,
$(324, 81, 301, 95)$,\adfsplit
$(325, 19, 266, 123)$,
$(326, 56, 309, 16)$,
$(327, 254, 112, 219)$,
$(328, 313, 243, 191)$,\adfsplit
$(329, 55, 249, 128)$,
$(330, 267, 301, 137)$,
$(331, 309, 163, 158)$,
$(332, 140, 226, 219)$,\adfsplit
$(333, 24, 256, 203)$,
$(334, 248, 43, 132)$,
$(335, 251, 70, 150)$,
$(336, 55, 29, 75)$,\adfsplit
$(0, 3, 47, 141)$,
$(0, 4, 60, 192)$,
$(0, 5, 30, 69)$,
$(0, 10, 115, 166)$,\adfsplit
$(0, 38, 95, 188)$,
$(0, 33, 76, 201)$,
$(0, 6, 48, 100)$,
$(0, 21, 85, 240)$,\adfsplit
$(0, 24, 77, 206)$,
$(0, 12, 78, 157)$,
$(0, 15, 102, 163)$,
$(0, 26, 111, 195)$

%ADFvfyBlocksEnd
\adfLgap \noindent by the mapping:
$x \mapsto x \oplus (3 j)$ for $x < 315$,
$x \mapsto (x - 315 + 22 j \adfmod{66}) + 315$ for $315 \le x < 381$,
$x \mapsto x$ for $x \ge 381$,
$0 \le j < 105$
 for the first two blocks;
$x \mapsto x \oplus j$ for $x < 315$,
$x \mapsto (x - 315 + 22 j \adfmod{66}) + 315$ for $315 \le x < 381$,
$x \mapsto x$ for $x \ge 381$,
$0 \le j < 315$
 for the last 34 blocks.
\ADFvfyParStart{(383, ((2, 105, ((315, 3, (105, 3)), (66, 22), (2, 2))), (34, 315, ((315, 1, (105, 3)), (66, 22), (2, 2)))), ((35, 9), (68, 1)))} %ADFvfyParEnd
% End of 35^9 68^1
%%%%%%%%%%%%%%%%%%%%%%%%%%%%%%%%%%%%%%%%%%%%%%%%%%%%%%%%%%%%%%%%%%%%%%%%%%%%%%%%%%%%%%%%%%
%%%%%%%%%%%%%%%%%%%%%%%%%%%%%%%%%%%%%%%%%%%%%%%%%%%%%%%%%%%%%%%%%%%%%%%%%%%%%%%%%%%%%%%%%%

% Charlotte:GDD4-1-3-5-mod-6-TeX-gen-A:HITS-fun:4.10
\adfDgap
%ADFvfyBlocksStart {35,35,35,35,35,35,35,35,35,74}
\noindent{\boldmath $ 35^{9} 74^{1} $}~
With the point set $Z_{389}$ partitioned into
 residue classes modulo $9$ for $\{0, 1, \dots, 314\}$, and
 $\{315, 316, \dots, 388\}$,
 the design is generated from

\adfLgap %ADFvfyDesignStart
$(387, 202, 179, 6)$,
$(388, 302, 220, 303)$,
$(315, 43, 176, 57)$,
$(315, 65, 130, 278)$,\adfsplit
$(315, 225, 82, 42)$,
$(316, 262, 69, 44)$,
$(316, 180, 157, 25)$,
$(316, 200, 230, 3)$,\adfsplit
$(317, 261, 4, 308)$,
$(317, 55, 178, 57)$,
$(317, 152, 293, 312)$,
$(318, 70, 99, 82)$,\adfsplit
$(318, 87, 138, 254)$,
$(318, 8, 41, 310)$,
$(319, 65, 60, 262)$,
$(319, 241, 239, 282)$,\adfsplit
$(319, 99, 35, 148)$,
$(320, 216, 4, 239)$,
$(320, 200, 115, 186)$,
$(320, 93, 8, 271)$,\adfsplit
$(321, 181, 150, 206)$,
$(321, 290, 207, 202)$,
$(321, 30, 178, 248)$,
$(322, 89, 121, 217)$,\adfsplit
$(322, 83, 131, 286)$,
$(322, 222, 165, 9)$,
$(323, 307, 98, 78)$,
$(323, 27, 283, 212)$,\adfsplit
$(323, 300, 209, 223)$,
$(324, 127, 9, 84)$,
$(324, 247, 123, 71)$,
$(324, 34, 110, 239)$,\adfsplit
$(325, 176, 288, 226)$,
$(325, 308, 178, 292)$,
$(325, 98, 66, 159)$,
$(326, 70, 154, 66)$,\adfsplit
$(326, 107, 112, 252)$,
$(326, 168, 290, 311)$,
$(327, 62, 154, 300)$,
$(327, 52, 110, 203)$,\adfsplit
$(327, 175, 99, 105)$,
$(328, 7, 275, 220)$,
$(328, 24, 210, 269)$,
$(328, 128, 99, 91)$,\adfsplit
$(329, 249, 212, 268)$,
$(329, 307, 155, 279)$,
$(329, 228, 314, 94)$,
$(330, 289, 228, 200)$,\adfsplit
$(330, 215, 207, 259)$,
$(330, 68, 114, 229)$,
$(331, 211, 195, 290)$,
$(331, 255, 122, 44)$,\adfsplit
$(331, 36, 298, 169)$,
$(332, 213, 131, 35)$,
$(332, 241, 235, 92)$,
$(332, 120, 270, 40)$,\adfsplit
$(333, 87, 109, 174)$,
$(333, 266, 225, 175)$,
$(333, 286, 38, 314)$,
$(334, 25, 183, 123)$,\adfsplit
$(334, 104, 101, 139)$,
$(334, 26, 153, 46)$,
$(335, 206, 228, 202)$,
$(335, 90, 69, 136)$,\adfsplit
$(335, 313, 5, 47)$,
$(336, 6, 155, 45)$,
$(336, 163, 193, 95)$,
$(336, 152, 219, 88)$,\adfsplit
$(337, 280, 74, 221)$,
$(337, 204, 0, 192)$,
$(337, 232, 247, 206)$,
$(338, 271, 263, 192)$,\adfsplit
$(338, 212, 177, 139)$,
$(338, 297, 80, 205)$,
$(270, 55, 236, 21)$,
$(179, 191, 208, 122)$,\adfsplit
$(176, 142, 49, 36)$,
$(0, 1, 25, 112)$,
$(0, 3, 94, 145)$,
$(0, 7, 73, 168)$,\adfsplit
$(0, 2, 55, 273)$,
$(0, 11, 130, 295)$,
$(0, 37, 141, 205)$,
$(0, 65, 139, 187)$,\adfsplit
$(1, 11, 142, 211)$,
$(1, 2, 256, 277)$,
$(0, 109, 128, 268)$,
$(0, 127, 247, 304)$,\adfsplit
$(1, 4, 104, 241)$,
$(0, 80, 226, 259)$,
$(0, 10, 84, 283)$,
$(0, 89, 176, 280)$,\adfsplit
$(0, 82, 125, 239)$,
$(0, 107, 113, 214)$,
$(0, 194, 209, 271)$,
$(0, 15, 68, 151)$,\adfsplit
$(0, 163, 284, 308)$,
$(0, 50, 166, 260)$,
$(0, 17, 190, 212)$,
$(0, 69, 161, 311)$,\adfsplit
$(0, 26, 177, 210)$,
$(0, 62, 200, 266)$,
$(0, 24, 120, 299)$,
$(0, 38, 201, 302)$,\adfsplit
$(0, 74, 134, 257)$,
$(0, 30, 78, 236)$,
$(0, 146, 221, 305)$

%ADFvfyBlocksEnd
\adfLgap \noindent by the mapping:
$x \mapsto x + 3 j \adfmod{315}$ for $x < 315$,
$x \mapsto (x - 315 + 24 j \adfmod{72}) + 315$ for $315 \le x < 387$,
$x \mapsto x$ for $x \ge 387$,
$0 \le j < 105$.
\ADFvfyParStart{(389, ((107, 105, ((315, 3), (72, 24), (2, 2)))), ((35, 9), (74, 1)))} %ADFvfyParEnd
% End of 35^9 74^1
%%%%%%%%%%%%%%%%%%%%%%%%%%%%%%%%%%%%%%%%%%%%%%%%%%%%%%%%%%%%%%%%%%%%%%%%%%%%%%%%%%%%%%%%%%
%%%%%%%%%%%%%%%%%%%%%%%%%%%%%%%%%%%%%%%%%%%%%%%%%%%%%%%%%%%%%%%%%%%%%%%%%%%%%%%%%%%%%%%%%%

% Charlotte:GDD4-1-3-5-mod-6-TeX-gen-A:HITS-fun:4.10
\adfDgap
%ADFvfyBlocksStart {35,35,35,35,35,35,35,35,35,86}
\noindent{\boldmath $ 35^{9} 86^{1} $}~
With the point set $Z_{401}$ partitioned into
 residue classes modulo $9$ for $\{0, 1, \dots, 314\}$, and
 $\{315, 316, \dots, 400\}$,
 the design is generated from

\adfLgap %ADFvfyDesignStart
$(396, 87, 131, 58)$,
$(397, 88, 132, 59)$,
$(398, 89, 133, 60)$,
$(315, 295, 131, 99)$,\adfsplit
$(324, 296, 132, 100)$,
$(333, 297, 133, 101)$,
$(315, 218, 46, 228)$,
$(324, 219, 47, 229)$,\adfsplit
$(333, 220, 48, 230)$,
$(315, 269, 294, 274)$,
$(324, 270, 295, 275)$,
$(333, 271, 296, 276)$,\adfsplit
$(316, 288, 227, 163)$,
$(325, 289, 228, 164)$,
$(334, 290, 229, 165)$,
$(316, 106, 301, 273)$,\adfsplit
$(325, 107, 302, 274)$,
$(334, 108, 303, 275)$,
$(316, 224, 213, 77)$,
$(325, 225, 214, 78)$,\adfsplit
$(334, 226, 215, 79)$,
$(317, 109, 8, 135)$,
$(326, 110, 9, 136)$,
$(335, 111, 10, 137)$,\adfsplit
$(317, 157, 237, 275)$,
$(326, 158, 238, 276)$,
$(335, 159, 239, 277)$,
$(317, 285, 268, 47)$,\adfsplit
$(326, 286, 269, 48)$,
$(335, 287, 270, 49)$,
$(318, 101, 58, 295)$,
$(327, 102, 59, 296)$,\adfsplit
$(336, 103, 60, 297)$,
$(318, 246, 197, 86)$,
$(327, 247, 198, 87)$,
$(336, 248, 199, 88)$,\adfsplit
$(318, 96, 0, 1)$,
$(327, 97, 1, 2)$,
$(336, 98, 2, 3)$,
$(319, 71, 64, 111)$,\adfsplit
$(328, 72, 65, 112)$,
$(337, 73, 66, 113)$,
$(319, 200, 166, 196)$,
$(328, 201, 167, 197)$,\adfsplit
$(337, 202, 168, 198)$,
$(319, 135, 284, 123)$,
$(328, 136, 285, 124)$,
$(337, 137, 286, 125)$,\adfsplit
$(320, 183, 38, 296)$,
$(329, 184, 39, 297)$,
$(338, 185, 40, 298)$,
$(320, 257, 159, 54)$,\adfsplit
$(329, 258, 160, 55)$,
$(338, 259, 161, 56)$,
$(320, 166, 289, 124)$,
$(329, 167, 290, 125)$,\adfsplit
$(338, 168, 291, 126)$,
$(321, 291, 306, 250)$,
$(330, 292, 307, 251)$,
$(339, 293, 308, 252)$,\adfsplit
$(321, 248, 114, 145)$,
$(330, 249, 115, 146)$,
$(339, 250, 116, 147)$,
$(321, 200, 310, 206)$,\adfsplit
$(330, 201, 311, 207)$,
$(339, 202, 312, 208)$,
$(322, 158, 228, 28)$,
$(331, 159, 229, 29)$,\adfsplit
$(340, 160, 230, 30)$,
$(322, 74, 16, 8)$,
$(331, 75, 17, 9)$,
$(340, 76, 18, 10)$,\adfsplit
$(322, 261, 177, 85)$,
$(331, 262, 178, 86)$,
$(340, 263, 179, 87)$,
$(323, 308, 285, 156)$,\adfsplit
$(332, 309, 286, 157)$,
$(341, 310, 287, 158)$,
$(323, 99, 215, 283)$,
$(1, 4, 17, 350)$,\adfsplit
$(0, 3, 62, 263)$,
$(0, 13, 75, 276)$,
$(0, 14, 60, 93)$,
$(0, 16, 76, 294)$,\adfsplit
$(0, 22, 55, 291)$,
$(0, 35, 132, 241)$,
$(0, 48, 131, 206)$,
$(0, 68, 82, 128)$,\adfsplit
$(1, 22, 70, 254)$,
$(0, 67, 199, 296)$,
$(0, 65, 100, 301)$,
$(0, 86, 173, 232)$,\adfsplit
$(0, 91, 156, 313)$,
$(1, 83, 176, 332)$,
$(0, 137, 239, 278)$,
$(1, 68, 175, 251)$,\adfsplit
$(0, 169, 193, 244)$,
$(0, 106, 208, 293)$,
$(0, 140, 256, 368)$,
$(0, 187, 224, 226)$,\adfsplit
$(0, 2, 50, 265)$,
$(1, 242, 245, 266)$,
$(0, 19, 87, 175)$,
$(0, 51, 229, 395)$,\adfsplit
$(1, 56, 107, 229)$,
$(0, 142, 233, 399)$,
$(0, 74, 107, 174)$,
$(0, 89, 213, 332)$,\adfsplit
$(1, 20, 89, 178)$,
$(0, 53, 230, 246)$,
$(0, 124, 146, 302)$,
$(0, 227, 280, 400)$,\adfsplit
$(0, 71, 177, 262)$

%ADFvfyBlocksEnd
\adfLgap \noindent by the mapping:
$x \mapsto x + 3 j \adfmod{315}$ for $x < 315$,
$x \mapsto (x - 315 + 27 j \adfmod{81}) + 315$ for $315 \le x < 396$,
$x \mapsto x$ for $x \ge 396$,
$0 \le j < 105$.
\ADFvfyParStart{(401, ((113, 105, ((315, 3), (81, 27), (5, 5)))), ((35, 9), (86, 1)))} %ADFvfyParEnd
% End of 35^9 86^1
%%%%%%%%%%%%%%%%%%%%%%%%%%%%%%%%%%%%%%%%%%%%%%%%%%%%%%%%%%%%%%%%%%%%%%%%%%%%%%%%%%%%%%%%%%
%%%%%%%%%%%%%%%%%%%%%%%%%%%%%%%%%%%%%%%%%%%%%%%%%%%%%%%%%%%%%%%%%%%%%%%%%%%%%%%%%%%%%%%%%%

% Charlotte:GDD4-1-3-5-mod-6-TeX-gen-A:HITS-fun:4.10
\adfDgap
%ADFvfyBlocksStart {35,35,35,35,35,35,35,35,35,92}
\noindent{\boldmath $ 35^{9} 92^{1} $}~
With the point set $Z_{407}$ partitioned into
 residue classes modulo $9$ for $\{0, 1, \dots, 314\}$, and
 $\{315, 316, \dots, 406\}$,
 the design is generated from

\adfLgap %ADFvfyDesignStart
$(315, 105, 183, 200)$,
$(325, 106, 184, 201)$,
$(335, 107, 185, 202)$,
$(315, 296, 46, 4)$,\adfsplit
$(325, 297, 47, 5)$,
$(335, 298, 48, 6)$,
$(315, 54, 250, 68)$,
$(325, 55, 251, 69)$,\adfsplit
$(335, 56, 252, 70)$,
$(316, 301, 215, 226)$,
$(326, 302, 216, 227)$,
$(336, 303, 217, 228)$,\adfsplit
$(316, 218, 167, 147)$,
$(326, 219, 168, 148)$,
$(336, 220, 169, 149)$,
$(316, 52, 162, 177)$,\adfsplit
$(326, 53, 163, 178)$,
$(336, 54, 164, 179)$,
$(317, 51, 1, 272)$,
$(327, 52, 2, 273)$,\adfsplit
$(337, 53, 3, 274)$,
$(317, 264, 296, 304)$,
$(327, 265, 297, 305)$,
$(337, 266, 298, 306)$,\adfsplit
$(317, 243, 256, 68)$,
$(327, 244, 257, 69)$,
$(337, 245, 258, 70)$,
$(318, 288, 239, 291)$,\adfsplit
$(328, 289, 240, 292)$,
$(338, 290, 241, 293)$,
$(318, 305, 67, 223)$,
$(328, 306, 68, 224)$,\adfsplit
$(338, 307, 69, 225)$,
$(318, 42, 281, 253)$,
$(328, 43, 282, 254)$,
$(338, 44, 283, 255)$,\adfsplit
$(319, 264, 176, 115)$,
$(329, 265, 177, 116)$,
$(339, 266, 178, 117)$,
$(319, 247, 136, 89)$,\adfsplit
$(329, 248, 137, 90)$,
$(339, 249, 138, 91)$,
$(319, 87, 281, 126)$,
$(329, 88, 282, 127)$,\adfsplit
$(339, 89, 283, 128)$,
$(320, 181, 108, 250)$,
$(330, 182, 109, 251)$,
$(340, 183, 110, 252)$,\adfsplit
$(320, 42, 251, 166)$,
$(330, 43, 252, 167)$,
$(340, 44, 253, 168)$,
$(320, 111, 2, 194)$,\adfsplit
$(330, 112, 3, 195)$,
$(340, 113, 4, 196)$,
$(321, 71, 87, 185)$,
$(331, 72, 88, 186)$,\adfsplit
$(341, 73, 89, 187)$,
$(321, 192, 38, 108)$,
$(331, 193, 39, 109)$,
$(341, 194, 40, 110)$,\adfsplit
$(321, 310, 91, 178)$,
$(331, 311, 92, 179)$,
$(341, 312, 93, 180)$,
$(322, 42, 214, 1)$,\adfsplit
$(332, 43, 215, 2)$,
$(342, 44, 216, 3)$,
$(322, 68, 148, 47)$,
$(332, 69, 149, 48)$,\adfsplit
$(342, 70, 150, 49)$,
$(322, 224, 39, 72)$,
$(332, 225, 40, 73)$,
$(342, 226, 41, 74)$,\adfsplit
$(323, 105, 197, 196)$,
$(333, 106, 198, 197)$,
$(343, 107, 199, 198)$,
$(323, 131, 127, 65)$,\adfsplit
$(333, 132, 128, 66)$,
$(343, 133, 129, 67)$,
$(323, 93, 4, 261)$,
$(333, 94, 5, 262)$,\adfsplit
$(343, 95, 6, 263)$,
$(324, 251, 106, 275)$,
$(334, 252, 107, 276)$,
$(344, 253, 108, 277)$,\adfsplit
$(324, 28, 141, 135)$,
$(0, 10, 35, 384)$,
$(0, 2, 141, 178)$,
$(0, 19, 25, 137)$,\adfsplit
$(0, 5, 112, 290)$,
$(0, 26, 38, 285)$,
$(0, 46, 151, 197)$,
$(0, 29, 129, 247)$,\adfsplit
$(0, 34, 64, 93)$,
$(0, 113, 148, 150)$,
$(0, 22, 57, 269)$,
$(0, 53, 120, 187)$,\adfsplit
$(0, 60, 139, 259)$,
$(0, 12, 128, 134)$,
$(1, 38, 94, 151)$,
$(0, 103, 218, 277)$,\adfsplit
$(0, 115, 179, 394)$,
$(1, 11, 188, 304)$,
$(1, 68, 263, 404)$,
$(0, 48, 251, 404)$,\adfsplit
$(0, 59, 200, 210)$,
$(0, 7, 138, 193)$,
$(0, 97, 281, 405)$,
$(0, 236, 262, 296)$,\adfsplit
$(0, 100, 107, 293)$,
$(1, 56, 149, 314)$,
$(1, 101, 139, 394)$,
$(0, 260, 289, 308)$,\adfsplit
$(0, 43, 74, 310)$,
$(0, 31, 167, 272)$,
$(0, 131, 136, 406)$,
$(0, 181, 241, 284)$

%ADFvfyBlocksEnd
\adfLgap \noindent by the mapping:
$x \mapsto x + 3 j \adfmod{315}$ for $x < 315$,
$x \mapsto (x - 315 + 30 j \adfmod{90}) + 315$ for $315 \le x < 405$,
$x \mapsto x$ for $x \ge 405$,
$0 \le j < 105$.
\ADFvfyParStart{(407, ((116, 105, ((315, 3), (90, 30), (2, 2)))), ((35, 9), (92, 1)))} %ADFvfyParEnd
% End of 35^9 92^1
%%%%%%%%%%%%%%%%%%%%%%%%%%%%%%%%%%%%%%%%%%%%%%%%%%%%%%%%%%%%%%%%%%%%%%%%%%%%%%%%%%%%%%%%%%
%%%%%%%%%%%%%%%%%%%%%%%%%%%%%%%%%%%%%%%%%%%%%%%%%%%%%%%%%%%%%%%%%%%%%%%%%%%%%%%%%%%%%%%%%%

% Charlotte:GDD4-1-3-5-mod-6-TeX-gen-A:HITS-fun:4.10
\adfDgap
%ADFvfyBlocksStart {35,35,35,35,35,35,35,35,35,104}
\noindent{\boldmath $ 35^{9} 104^{1} $}~
With the point set $Z_{419}$ partitioned into
 residue classes modulo $9$ for $\{0, 1, \dots, 314\}$, and
 $\{315, 316, \dots, 418\}$,
 the design is generated from

\adfLgap %ADFvfyDesignStart
$(414, 150, 23, 181)$,
$(415, 151, 24, 182)$,
$(416, 152, 25, 183)$,
$(315, 239, 301, 168)$,\adfsplit
$(326, 240, 302, 169)$,
$(337, 241, 303, 170)$,
$(315, 45, 133, 84)$,
$(326, 46, 134, 85)$,\adfsplit
$(337, 47, 135, 86)$,
$(315, 28, 164, 197)$,
$(326, 29, 165, 198)$,
$(337, 30, 166, 199)$,\adfsplit
$(316, 265, 284, 182)$,
$(327, 266, 285, 183)$,
$(338, 267, 286, 184)$,
$(316, 295, 24, 84)$,\adfsplit
$(327, 296, 25, 85)$,
$(338, 297, 26, 86)$,
$(316, 206, 19, 126)$,
$(327, 207, 20, 127)$,\adfsplit
$(338, 208, 21, 128)$,
$(317, 89, 75, 47)$,
$(328, 90, 76, 48)$,
$(339, 91, 77, 49)$,\adfsplit
$(317, 199, 267, 158)$,
$(328, 200, 268, 159)$,
$(339, 201, 269, 160)$,
$(317, 112, 268, 189)$,\adfsplit
$(328, 113, 269, 190)$,
$(339, 114, 270, 191)$,
$(318, 95, 88, 209)$,
$(329, 96, 89, 210)$,\adfsplit
$(340, 97, 90, 211)$,
$(318, 62, 156, 253)$,
$(329, 63, 157, 254)$,
$(340, 64, 158, 255)$,\adfsplit
$(318, 202, 213, 162)$,
$(329, 203, 214, 163)$,
$(340, 204, 215, 164)$,
$(319, 258, 2, 0)$,\adfsplit
$(330, 259, 3, 1)$,
$(341, 260, 4, 2)$,
$(319, 179, 67, 232)$,
$(330, 180, 68, 233)$,\adfsplit
$(341, 181, 69, 234)$,
$(319, 73, 23, 282)$,
$(330, 74, 24, 283)$,
$(341, 75, 25, 284)$,\adfsplit
$(320, 286, 168, 29)$,
$(331, 287, 169, 30)$,
$(342, 288, 170, 31)$,
$(320, 63, 309, 233)$,\adfsplit
$(331, 64, 310, 234)$,
$(342, 65, 311, 235)$,
$(320, 37, 184, 41)$,
$(331, 38, 185, 42)$,\adfsplit
$(342, 39, 186, 43)$,
$(321, 238, 5, 312)$,
$(332, 239, 6, 313)$,
$(343, 240, 7, 314)$,\adfsplit
$(321, 180, 115, 255)$,
$(332, 181, 116, 256)$,
$(343, 182, 117, 257)$,
$(321, 188, 307, 47)$,\adfsplit
$(332, 189, 308, 48)$,
$(343, 190, 309, 49)$,
$(322, 39, 267, 220)$,
$(333, 40, 268, 221)$,\adfsplit
$(344, 41, 269, 222)$,
$(322, 70, 287, 262)$,
$(333, 71, 288, 263)$,
$(344, 72, 289, 264)$,\adfsplit
$(322, 128, 36, 122)$,
$(333, 129, 37, 123)$,
$(344, 130, 38, 124)$,
$(323, 89, 137, 189)$,\adfsplit
$(334, 90, 138, 190)$,
$(345, 91, 139, 191)$,
$(323, 311, 159, 22)$,
$(334, 312, 160, 23)$,\adfsplit
$(345, 313, 161, 24)$,
$(323, 280, 210, 295)$,
$(334, 281, 211, 296)$,
$(345, 282, 212, 297)$,\adfsplit
$(324, 161, 183, 199)$,
$(335, 162, 184, 200)$,
$(346, 163, 185, 201)$,
$(324, 144, 157, 187)$,\adfsplit
$(335, 145, 158, 188)$,
$(346, 146, 159, 189)$,
$(324, 294, 68, 2)$,
$(335, 295, 69, 3)$,\adfsplit
$(346, 296, 70, 4)$,
$(325, 165, 101, 268)$,
$(0, 1, 96, 295)$,
$(0, 10, 177, 214)$,\adfsplit
$(0, 3, 24, 125)$,
$(0, 5, 204, 224)$,
$(0, 17, 84, 130)$,
$(0, 12, 154, 369)$,\adfsplit
$(0, 32, 35, 186)$,
$(0, 29, 113, 298)$,
$(0, 61, 212, 347)$,
$(0, 67, 91, 222)$,\adfsplit
$(0, 73, 166, 195)$,
$(0, 78, 280, 283)$,
$(0, 95, 210, 413)$,
$(0, 110, 183, 358)$,\adfsplit
$(0, 34, 314, 417)$,
$(0, 64, 155, 167)$,
$(0, 131, 269, 402)$,
$(0, 190, 254, 418)$,\adfsplit
$(0, 149, 173, 278)$,
$(0, 151, 185, 305)$,
$(0, 115, 281, 310)$,
$(1, 62, 187, 413)$,\adfsplit
$(1, 68, 161, 178)$,
$(1, 13, 97, 296)$,
$(1, 11, 194, 215)$,
$(1, 38, 116, 211)$,\adfsplit
$(1, 74, 79, 391)$,
$(1, 2, 133, 402)$

%ADFvfyBlocksEnd
\adfLgap \noindent by the mapping:
$x \mapsto x + 3 j \adfmod{315}$ for $x < 315$,
$x \mapsto (x - 315 + 33 j \adfmod{99}) + 315$ for $315 \le x < 414$,
$x \mapsto x$ for $x \ge 414$,
$0 \le j < 105$.
\ADFvfyParStart{(419, ((122, 105, ((315, 3), (99, 33), (5, 5)))), ((35, 9), (104, 1)))} %ADFvfyParEnd
% End of 35^9 104^1
%%%%%%%%%%%%%%%%%%%%%%%%%%%%%%%%%%%%%%%%%%%%%%%%%%%%%%%%%%%%%%%%%%%%%%%%%%%%%%%%%%%%%%%%%%
%%%%%%%%%%%%%%%%%%%%%%%%%%%%%%%%%%%%%%%%%%%%%%%%%%%%%%%%%%%%%%%%%%%%%%%%%%%%%%%%%%%%%%%%%%

% Charlotte:GDD4-1-3-5-mod-6-TeX-gen-A:HITS-fun:4.10
\adfDgap
%ADFvfyBlocksStart {35,35,35,35,35,35,35,35,35,116}
\noindent{\boldmath $ 35^{9} 116^{1} $}~
With the point set $Z_{431}$ partitioned into
 residue classes modulo $9$ for $\{0, 1, \dots, 314\}$, and
 $\{315, 316, \dots, 430\}$,
 the design is generated from

\adfLgap %ADFvfyDesignStart
$(423, 70, 246, 86)$,
$(425, 71, 247, 87)$,
$(427, 72, 248, 88)$,
$(424, 266, 283, 198)$,\adfsplit
$(426, 267, 284, 199)$,
$(428, 268, 285, 200)$,
$(315, 135, 51, 197)$,
$(327, 136, 52, 198)$,\adfsplit
$(339, 137, 53, 199)$,
$(315, 272, 136, 201)$,
$(327, 273, 137, 202)$,
$(339, 274, 138, 203)$,\adfsplit
$(315, 52, 248, 202)$,
$(327, 53, 249, 203)$,
$(339, 54, 250, 204)$,
$(316, 7, 65, 109)$,\adfsplit
$(328, 8, 66, 110)$,
$(340, 9, 67, 111)$,
$(316, 159, 239, 72)$,
$(328, 160, 240, 73)$,\adfsplit
$(340, 161, 241, 74)$,
$(316, 282, 98, 22)$,
$(328, 283, 99, 23)$,
$(340, 284, 100, 24)$,\adfsplit
$(317, 132, 191, 298)$,
$(329, 133, 192, 299)$,
$(341, 134, 193, 300)$,
$(317, 192, 301, 98)$,\adfsplit
$(329, 193, 302, 99)$,
$(341, 194, 303, 100)$,
$(317, 212, 79, 18)$,
$(329, 213, 80, 19)$,\adfsplit
$(341, 214, 81, 20)$,
$(318, 170, 201, 56)$,
$(330, 171, 202, 57)$,
$(342, 172, 203, 58)$,\adfsplit
$(318, 145, 212, 115)$,
$(330, 146, 213, 116)$,
$(342, 147, 214, 117)$,
$(318, 159, 76, 27)$,\adfsplit
$(330, 160, 77, 28)$,
$(342, 161, 78, 29)$,
$(319, 115, 310, 192)$,
$(331, 116, 311, 193)$,\adfsplit
$(343, 117, 312, 194)$,
$(319, 27, 204, 122)$,
$(331, 28, 205, 123)$,
$(343, 29, 206, 124)$,\adfsplit
$(319, 116, 28, 155)$,
$(331, 117, 29, 156)$,
$(343, 118, 30, 157)$,
$(320, 171, 196, 15)$,\adfsplit
$(332, 172, 197, 16)$,
$(344, 173, 198, 17)$,
$(320, 188, 221, 236)$,
$(332, 189, 222, 237)$,\adfsplit
$(344, 190, 223, 238)$,
$(320, 192, 199, 139)$,
$(332, 193, 200, 140)$,
$(344, 194, 201, 141)$,\adfsplit
$(321, 55, 95, 124)$,
$(333, 56, 96, 125)$,
$(345, 57, 97, 126)$,
$(321, 189, 98, 47)$,\adfsplit
$(333, 190, 99, 48)$,
$(345, 191, 100, 49)$,
$(321, 112, 168, 147)$,
$(333, 113, 169, 148)$,\adfsplit
$(345, 114, 170, 149)$,
$(322, 154, 43, 117)$,
$(334, 155, 44, 118)$,
$(346, 156, 45, 119)$,\adfsplit
$(322, 267, 67, 192)$,
$(334, 268, 68, 193)$,
$(346, 269, 69, 194)$,
$(322, 77, 89, 182)$,\adfsplit
$(334, 78, 90, 183)$,
$(346, 79, 91, 184)$,
$(323, 69, 230, 192)$,
$(335, 70, 231, 193)$,\adfsplit
$(347, 71, 232, 194)$,
$(323, 136, 155, 179)$,
$(335, 137, 156, 180)$,
$(347, 138, 157, 181)$,\adfsplit
$(323, 16, 252, 274)$,
$(335, 17, 253, 275)$,
$(347, 18, 254, 276)$,
$(324, 113, 300, 253)$,\adfsplit
$(336, 114, 301, 254)$,
$(348, 115, 302, 255)$,
$(324, 58, 35, 108)$,
$(336, 59, 36, 109)$,\adfsplit
$(348, 60, 37, 110)$,
$(324, 61, 69, 56)$,
$(336, 62, 70, 57)$,
$(348, 63, 71, 58)$,\adfsplit
$(325, 306, 77, 57)$,
$(337, 307, 78, 58)$,
$(349, 308, 79, 59)$,
$(0, 2, 6, 143)$,\adfsplit
$(0, 4, 10, 147)$,
$(0, 14, 28, 157)$,
$(0, 41, 42, 193)$,
$(0, 1, 11, 274)$,\adfsplit
$(0, 52, 152, 158)$,
$(0, 26, 129, 163)$,
$(0, 100, 164, 326)$,
$(0, 78, 263, 304)$,\adfsplit
$(0, 64, 174, 374)$,
$(0, 130, 281, 361)$,
$(1, 35, 148, 397)$,
$(1, 79, 175, 305)$,\adfsplit
$(1, 158, 200, 410)$,
$(0, 215, 313, 410)$,
$(0, 98, 245, 350)$,
$(0, 70, 96, 421)$,\adfsplit
$(0, 106, 209, 385)$,
$(0, 301, 305, 386)$,
$(0, 199, 202, 422)$,
$(0, 191, 217, 287)$,\adfsplit
$(0, 3, 116, 337)$,
$(0, 122, 223, 251)$,
$(1, 206, 284, 373)$,
$(0, 32, 124, 429)$,\adfsplit
$(0, 89, 211, 430)$,
$(0, 110, 214, 362)$,
$(0, 101, 103, 104)$,
$(0, 92, 283, 398)$

%ADFvfyBlocksEnd
\adfLgap \noindent by the mapping:
$x \mapsto x + 3 j \adfmod{315}$ for $x < 315$,
$x \mapsto (x - 315 + 36 j \adfmod{108}) + 315$ for $315 \le x < 423$,
$x \mapsto x$ for $x \ge 423$,
$0 \le j < 105$.
\ADFvfyParStart{(431, ((128, 105, ((315, 3), (108, 36), (8, 8)))), ((35, 9), (116, 1)))} %ADFvfyParEnd
% End of 35^9 116^1
%%%%%%%%%%%%%%%%%%%%%%%%%%%%%%%%%%%%%%%%%%%%%%%%%%%%%%%%%%%%%%%%%%%%%%%%%%%%%%%%%%%%%%%%%%
%%%%%%%%%%%%%%%%%%%%%%%%%%%%%%%%%%%%%%%%%%%%%%%%%%%%%%%%%%%%%%%%%%%%%%%%%%%%%%%%%%%%%%%%%%

% Charlotte:GDD4-1-3-5-mod-6-TeX-gen-A:HITS-fun:4.10
\adfDgap
%ADFvfyBlocksStart {35,35,35,35,35,35,35,35,35,122}
\noindent{\boldmath $ 35^{9} 122^{1} $}~
With the point set $Z_{437}$ partitioned into
 residue classes modulo $9$ for $\{0, 1, \dots, 314\}$, and
 $\{315, 316, \dots, 436\}$,
 the design is generated from

\adfLgap %ADFvfyDesignStart
$(432, 112, 272, 66)$,
$(433, 113, 273, 67)$,
$(434, 114, 274, 68)$,
$(315, 182, 102, 68)$,\adfsplit
$(328, 183, 103, 69)$,
$(341, 184, 104, 70)$,
$(315, 241, 163, 17)$,
$(328, 242, 164, 18)$,\adfsplit
$(341, 243, 165, 19)$,
$(315, 148, 222, 54)$,
$(328, 149, 223, 55)$,
$(341, 150, 224, 56)$,\adfsplit
$(316, 202, 286, 261)$,
$(329, 203, 287, 262)$,
$(342, 204, 288, 263)$,
$(316, 249, 305, 19)$,\adfsplit
$(329, 250, 306, 20)$,
$(342, 251, 307, 21)$,
$(316, 41, 164, 219)$,
$(329, 42, 165, 220)$,\adfsplit
$(342, 43, 166, 221)$,
$(317, 17, 141, 34)$,
$(330, 18, 142, 35)$,
$(343, 19, 143, 36)$,\adfsplit
$(317, 136, 219, 270)$,
$(330, 137, 220, 271)$,
$(343, 138, 221, 272)$,
$(317, 182, 184, 86)$,\adfsplit
$(330, 183, 185, 87)$,
$(343, 184, 186, 88)$,
$(318, 33, 59, 247)$,
$(331, 34, 60, 248)$,\adfsplit
$(344, 35, 61, 249)$,
$(318, 98, 128, 133)$,
$(331, 99, 129, 134)$,
$(344, 100, 130, 135)$,\adfsplit
$(318, 156, 63, 145)$,
$(331, 157, 64, 146)$,
$(344, 158, 65, 147)$,
$(319, 31, 308, 50)$,\adfsplit
$(332, 32, 309, 51)$,
$(345, 33, 310, 52)$,
$(319, 0, 268, 199)$,
$(332, 1, 269, 200)$,\adfsplit
$(345, 2, 270, 201)$,
$(319, 152, 264, 267)$,
$(332, 153, 265, 268)$,
$(345, 154, 266, 269)$,\adfsplit
$(320, 99, 79, 283)$,
$(333, 100, 80, 284)$,
$(346, 101, 81, 285)$,
$(320, 303, 200, 185)$,\adfsplit
$(333, 304, 201, 186)$,
$(346, 305, 202, 187)$,
$(320, 134, 55, 174)$,
$(333, 135, 56, 175)$,\adfsplit
$(346, 136, 57, 176)$,
$(321, 173, 165, 141)$,
$(334, 174, 166, 142)$,
$(347, 175, 167, 143)$,\adfsplit
$(321, 253, 86, 63)$,
$(334, 254, 87, 64)$,
$(347, 255, 88, 65)$,
$(321, 247, 313, 134)$,\adfsplit
$(334, 248, 314, 135)$,
$(347, 249, 0, 136)$,
$(322, 14, 308, 267)$,
$(335, 15, 309, 268)$,\adfsplit
$(348, 16, 310, 269)$,
$(322, 94, 161, 171)$,
$(335, 95, 162, 172)$,
$(348, 96, 163, 173)$,\adfsplit
$(322, 57, 100, 61)$,
$(335, 58, 101, 62)$,
$(348, 59, 102, 63)$,
$(323, 101, 259, 36)$,\adfsplit
$(336, 102, 260, 37)$,
$(349, 103, 261, 38)$,
$(323, 183, 122, 51)$,
$(336, 184, 123, 52)$,\adfsplit
$(349, 185, 124, 53)$,
$(323, 310, 136, 188)$,
$(336, 311, 137, 189)$,
$(349, 312, 138, 190)$,\adfsplit
$(324, 173, 19, 204)$,
$(337, 174, 20, 205)$,
$(350, 175, 21, 206)$,
$(324, 106, 234, 57)$,\adfsplit
$(337, 107, 235, 58)$,
$(350, 108, 236, 59)$,
$(324, 179, 121, 266)$,
$(337, 180, 122, 267)$,\adfsplit
$(350, 181, 123, 268)$,
$(325, 141, 129, 81)$,
$(338, 142, 130, 82)$,
$(351, 143, 131, 83)$,\adfsplit
$(325, 82, 158, 211)$,
$(338, 83, 159, 212)$,
$(351, 84, 160, 213)$,
$(325, 110, 305, 304)$,\adfsplit
$(338, 111, 306, 305)$,
$(351, 112, 307, 306)$,
$(326, 162, 11, 312)$,
$(0, 75, 293, 378)$,\adfsplit
$(0, 88, 227, 340)$,
$(1, 65, 151, 404)$,
$(0, 73, 173, 392)$,
$(1, 38, 71, 211)$,\adfsplit
$(0, 13, 86, 404)$,
$(0, 151, 226, 327)$,
$(0, 50, 64, 215)$,
$(1, 23, 98, 339)$,\adfsplit
$(0, 163, 265, 417)$,
$(0, 16, 44, 435)$,
$(0, 6, 28, 352)$,
$(0, 100, 106, 418)$,\adfsplit
$(0, 97, 139, 302)$,
$(0, 37, 70, 210)$,
$(0, 42, 152, 247)$,
$(0, 301, 308, 436)$,\adfsplit
$(0, 7, 102, 278)$,
$(1, 107, 157, 430)$,
$(0, 68, 172, 353)$,
$(0, 143, 211, 299)$,\adfsplit
$(0, 33, 166, 182)$,
$(0, 89, 271, 431)$,
$(0, 140, 229, 242)$,
$(0, 104, 209, 430)$,\adfsplit
$(1, 167, 173, 405)$,
$(0, 245, 287, 379)$,
$(0, 95, 159, 366)$

%ADFvfyBlocksEnd
\adfLgap \noindent by the mapping:
$x \mapsto x + 3 j \adfmod{315}$ for $x < 315$,
$x \mapsto (x - 315 + 39 j \adfmod{117}) + 315$ for $315 \le x < 432$,
$x \mapsto x$ for $x \ge 432$,
$0 \le j < 105$.
\ADFvfyParStart{(437, ((131, 105, ((315, 3), (117, 39), (5, 5)))), ((35, 9), (122, 1)))} %ADFvfyParEnd
% End of 35^9 122^1
%%%%%%%%%%%%%%%%%%%%%%%%%%%%%%%%%%%%%%%%%%%%%%%%%%%%%%%%%%%%%%%%%%%%%%%%%%%%%%%%%%%%%%%%%%
%%%%%%%%%%%%%%%%%%%%%%%%%%%%%%%%%%%%%%%%%%%%%%%%%%%%%%%%%%%%%%%%%%%%%%%%%%%%%%%%%%%%%%%%%%

% Charlotte:GDD4-1-3-5-mod-6-TeX-gen-A:HITS-fun:4.10
\adfDgap
%ADFvfyBlocksStart {35,35,35,35,35,35,35,35,35,128}
\noindent{\boldmath $ 35^{9} 128^{1} $}~
With the point set $Z_{443}$ partitioned into
 residue classes modulo $9$ for $\{0, 1, \dots, 314\}$, and
 $\{315, 316, \dots, 442\}$,
 the design is generated from

\adfLgap %ADFvfyDesignStart
$(315, 145, 69, 256)$,
$(329, 146, 70, 257)$,
$(343, 147, 71, 258)$,
$(315, 70, 243, 156)$,\adfsplit
$(329, 71, 244, 157)$,
$(343, 72, 245, 158)$,
$(315, 32, 56, 224)$,
$(329, 33, 57, 225)$,\adfsplit
$(343, 34, 58, 226)$,
$(316, 134, 105, 137)$,
$(330, 135, 106, 138)$,
$(344, 136, 107, 139)$,\adfsplit
$(316, 219, 221, 37)$,
$(330, 220, 222, 38)$,
$(344, 221, 223, 39)$,
$(316, 135, 286, 292)$,\adfsplit
$(330, 136, 287, 293)$,
$(344, 137, 288, 294)$,
$(317, 164, 163, 186)$,
$(331, 165, 164, 187)$,\adfsplit
$(345, 166, 165, 188)$,
$(317, 205, 9, 89)$,
$(331, 206, 10, 90)$,
$(345, 207, 11, 91)$,\adfsplit
$(317, 301, 113, 237)$,
$(331, 302, 114, 238)$,
$(345, 303, 115, 239)$,
$(318, 281, 184, 163)$,\adfsplit
$(332, 282, 185, 164)$,
$(346, 283, 186, 165)$,
$(318, 210, 98, 135)$,
$(332, 211, 99, 136)$,\adfsplit
$(346, 212, 100, 137)$,
$(318, 285, 311, 295)$,
$(332, 286, 312, 296)$,
$(346, 287, 313, 297)$,\adfsplit
$(319, 143, 87, 160)$,
$(333, 144, 88, 161)$,
$(347, 145, 89, 162)$,
$(319, 111, 180, 122)$,\adfsplit
$(333, 112, 181, 123)$,
$(347, 113, 182, 124)$,
$(319, 218, 211, 55)$,
$(333, 219, 212, 56)$,\adfsplit
$(347, 220, 213, 57)$,
$(320, 255, 242, 230)$,
$(334, 256, 243, 231)$,
$(348, 257, 244, 232)$,\adfsplit
$(320, 85, 280, 189)$,
$(334, 86, 281, 190)$,
$(348, 87, 282, 191)$,
$(320, 299, 123, 151)$,\adfsplit
$(334, 300, 124, 152)$,
$(348, 301, 125, 153)$,
$(321, 100, 196, 210)$,
$(335, 101, 197, 211)$,\adfsplit
$(349, 102, 198, 212)$,
$(321, 301, 263, 179)$,
$(335, 302, 264, 180)$,
$(349, 303, 265, 181)$,\adfsplit
$(321, 285, 153, 68)$,
$(335, 286, 154, 69)$,
$(349, 287, 155, 70)$,
$(322, 248, 282, 208)$,\adfsplit
$(336, 249, 283, 209)$,
$(350, 250, 284, 210)$,
$(322, 177, 144, 259)$,
$(336, 178, 145, 260)$,\adfsplit
$(350, 179, 146, 261)$,
$(322, 211, 26, 65)$,
$(336, 212, 27, 66)$,
$(350, 213, 28, 67)$,\adfsplit
$(323, 104, 112, 33)$,
$(337, 105, 113, 34)$,
$(351, 106, 114, 35)$,
$(323, 99, 192, 37)$,\adfsplit
$(337, 100, 193, 38)$,
$(351, 101, 194, 39)$,
$(323, 125, 232, 263)$,
$(337, 126, 233, 264)$,\adfsplit
$(351, 127, 234, 265)$,
$(324, 48, 5, 177)$,
$(338, 49, 6, 178)$,
$(352, 50, 7, 179)$,\adfsplit
$(324, 26, 92, 271)$,
$(338, 27, 93, 272)$,
$(352, 28, 94, 273)$,
$(324, 238, 72, 160)$,\adfsplit
$(338, 239, 73, 161)$,
$(352, 240, 74, 162)$,
$(325, 63, 113, 8)$,
$(339, 64, 114, 9)$,\adfsplit
$(353, 65, 115, 10)$,
$(325, 240, 236, 75)$,
$(339, 241, 237, 76)$,
$(353, 242, 238, 77)$,\adfsplit
$(325, 265, 307, 205)$,
$(339, 266, 308, 206)$,
$(353, 267, 309, 207)$,
$(326, 243, 311, 202)$,\adfsplit
$(340, 244, 312, 203)$,
$(354, 245, 313, 204)$,
$(326, 271, 236, 206)$,
$(0, 15, 190, 368)$,\adfsplit
$(0, 46, 61, 113)$,
$(0, 52, 254, 327)$,
$(0, 202, 248, 263)$,
$(0, 30, 65, 280)$,\adfsplit
$(0, 95, 114, 220)$,
$(0, 59, 103, 328)$,
$(0, 209, 268, 340)$,
$(0, 5, 100, 424)$,\adfsplit
$(0, 121, 140, 341)$,
$(0, 44, 256, 342)$,
$(0, 51, 134, 355)$,
$(0, 101, 215, 356)$,\adfsplit
$(0, 170, 221, 384)$,
$(0, 49, 223, 426)$,
$(0, 20, 141, 438)$,
$(1, 50, 259, 438)$,\adfsplit
$(0, 53, 145, 398)$,
$(0, 181, 295, 440)$,
$(0, 57, 269, 412)$,
$(0, 19, 67, 439)$,\adfsplit
$(1, 52, 233, 328)$,
$(1, 122, 263, 355)$,
$(1, 239, 296, 327)$,
$(0, 48, 226, 411)$,\adfsplit
$(0, 92, 232, 262)$,
$(0, 47, 94, 383)$,
$(0, 89, 310, 425)$,
$(0, 77, 125, 214)$,\adfsplit
$(0, 266, 271, 441)$,
$(0, 137, 238, 442)$

%ADFvfyBlocksEnd
\adfLgap \noindent by the mapping:
$x \mapsto x + 3 j \adfmod{315}$ for $x < 315$,
$x \mapsto (x - 315 + 42 j \adfmod{126}) + 315$ for $315 \le x < 441$,
$x \mapsto x$ for $x \ge 441$,
$0 \le j < 105$.
\ADFvfyParStart{(443, ((134, 105, ((315, 3), (126, 42), (2, 2)))), ((35, 9), (128, 1)))} %ADFvfyParEnd
% End of 35^9 128^1
%%%%%%%%%%%%%%%%%%%%%%%%%%%%%%%%%%%%%%%%%%%%%%%%%%%%%%%%%%%%%%%%%%%%%%%%%%%%%%%%%%%%%%%%%%
%%%%%%%%%%%%%%%%%%%%%%%%%%%%%%%%%%%%%%%%%%%%%%%%%%%%%%%%%%%%%%%%%%%%%%%%%%%%%%%%%%%%%%%%%%

% Charlotte:GDD4-1-3-5-mod-6-TeX-gen-A:HITS-fun:4.10
\adfDgap
%ADFvfyBlocksStart {35,35,35,35,35,35,35,35,35,134}
\noindent{\boldmath $ 35^{9} 134^{1} $}~
With the point set $Z_{449}$ partitioned into
 residue classes modulo $9$ for $\{0, 1, \dots, 314\}$, and
 $\{315, 316, \dots, 448\}$,
 the design is generated from

\adfLgap %ADFvfyDesignStart
$(441, 215, 207, 130)$,
$(443, 216, 208, 131)$,
$(445, 217, 209, 132)$,
$(442, 187, 42, 158)$,\adfsplit
$(444, 188, 43, 159)$,
$(446, 189, 44, 160)$,
$(315, 172, 245, 88)$,
$(329, 173, 246, 89)$,\adfsplit
$(343, 174, 247, 90)$,
$(315, 189, 76, 123)$,
$(329, 190, 77, 124)$,
$(343, 191, 78, 125)$,\adfsplit
$(315, 140, 269, 201)$,
$(329, 141, 270, 202)$,
$(343, 142, 271, 203)$,
$(316, 113, 237, 262)$,\adfsplit
$(330, 114, 238, 263)$,
$(344, 115, 239, 264)$,
$(316, 125, 69, 211)$,
$(330, 126, 70, 212)$,\adfsplit
$(344, 127, 71, 213)$,
$(316, 216, 47, 34)$,
$(330, 217, 48, 35)$,
$(344, 218, 49, 36)$,\adfsplit
$(317, 144, 280, 215)$,
$(331, 145, 281, 216)$,
$(345, 146, 282, 217)$,
$(317, 194, 210, 268)$,\adfsplit
$(331, 195, 211, 269)$,
$(345, 196, 212, 270)$,
$(317, 173, 274, 240)$,
$(331, 174, 275, 241)$,\adfsplit
$(345, 175, 276, 242)$,
$(318, 52, 307, 30)$,
$(332, 53, 308, 31)$,
$(346, 54, 309, 32)$,\adfsplit
$(318, 231, 310, 92)$,
$(332, 232, 311, 93)$,
$(346, 233, 312, 94)$,
$(318, 225, 269, 284)$,\adfsplit
$(332, 226, 270, 285)$,
$(346, 227, 271, 286)$,
$(319, 208, 153, 313)$,
$(333, 209, 154, 314)$,\adfsplit
$(347, 210, 155, 0)$,
$(319, 89, 203, 47)$,
$(333, 90, 204, 48)$,
$(347, 91, 205, 49)$,\adfsplit
$(319, 24, 220, 120)$,
$(333, 25, 221, 121)$,
$(347, 26, 222, 122)$,
$(320, 208, 133, 83)$,\adfsplit
$(334, 209, 134, 84)$,
$(348, 210, 135, 85)$,
$(320, 78, 219, 76)$,
$(334, 79, 220, 77)$,\adfsplit
$(348, 80, 221, 78)$,
$(320, 80, 50, 198)$,
$(334, 81, 51, 199)$,
$(348, 82, 52, 200)$,\adfsplit
$(321, 39, 76, 0)$,
$(335, 40, 77, 1)$,
$(349, 41, 78, 2)$,
$(321, 150, 47, 212)$,\adfsplit
$(335, 151, 48, 213)$,
$(349, 152, 49, 214)$,
$(321, 127, 205, 314)$,
$(335, 128, 206, 0)$,\adfsplit
$(349, 129, 207, 1)$,
$(322, 306, 14, 2)$,
$(336, 307, 15, 3)$,
$(350, 308, 16, 4)$,\adfsplit
$(322, 133, 102, 69)$,
$(336, 134, 103, 70)$,
$(350, 135, 104, 71)$,
$(322, 280, 58, 53)$,\adfsplit
$(336, 281, 59, 54)$,
$(350, 282, 60, 55)$,
$(323, 122, 2, 260)$,
$(337, 123, 3, 261)$,\adfsplit
$(351, 124, 4, 262)$,
$(323, 150, 115, 219)$,
$(337, 151, 116, 220)$,
$(351, 152, 117, 221)$,\adfsplit
$(323, 306, 112, 10)$,
$(337, 307, 113, 11)$,
$(351, 308, 114, 12)$,
$(324, 280, 126, 102)$,\adfsplit
$(338, 281, 127, 103)$,
$(352, 282, 128, 104)$,
$(324, 130, 35, 182)$,
$(338, 131, 36, 183)$,\adfsplit
$(352, 132, 37, 184)$,
$(324, 97, 123, 140)$,
$(338, 98, 124, 141)$,
$(352, 99, 125, 142)$,\adfsplit
$(325, 81, 309, 203)$,
$(339, 82, 310, 204)$,
$(353, 83, 311, 205)$,
$(325, 202, 188, 209)$,\adfsplit
$(339, 203, 189, 210)$,
$(353, 204, 190, 211)$,
$(325, 294, 205, 154)$,
$(339, 295, 206, 155)$,\adfsplit
$(353, 296, 207, 156)$,
$(326, 35, 129, 310)$,
$(340, 36, 130, 311)$,
$(354, 37, 131, 312)$,\adfsplit
$(0, 4, 309, 368)$,
$(1, 5, 116, 368)$,
$(0, 217, 223, 327)$,
$(0, 46, 92, 328)$,\adfsplit
$(0, 3, 269, 383)$,
$(0, 305, 311, 342)$,
$(0, 70, 80, 354)$,
$(0, 1, 152, 438)$,\adfsplit
$(0, 20, 127, 130)$,
$(0, 185, 188, 355)$,
$(0, 110, 208, 356)$,
$(0, 107, 295, 370)$,\adfsplit
$(0, 112, 203, 412)$,
$(1, 112, 224, 341)$,
$(0, 49, 233, 425)$,
$(0, 111, 262, 426)$,\adfsplit
$(1, 83, 124, 342)$,
$(0, 28, 98, 439)$,
$(0, 91, 274, 397)$,
$(0, 164, 287, 411)$,\adfsplit
$(0, 41, 132, 382)$,
$(1, 29, 49, 382)$,
$(0, 200, 232, 369)$,
$(1, 2, 50, 440)$,\adfsplit
$(0, 48, 283, 440)$,
$(0, 131, 184, 314)$,
$(0, 53, 123, 205)$,
$(0, 32, 115, 447)$,\adfsplit
$(0, 83, 163, 448)$

%ADFvfyBlocksEnd
\adfLgap \noindent by the mapping:
$x \mapsto x + 3 j \adfmod{315}$ for $x < 315$,
$x \mapsto (x - 315 + 42 j \adfmod{126}) + 315$ for $315 \le x < 441$,
$x \mapsto x$ for $x \ge 441$,
$0 \le j < 105$.
\ADFvfyParStart{(449, ((137, 105, ((315, 3), (126, 42), (8, 8)))), ((35, 9), (134, 1)))} %ADFvfyParEnd
% End of 35^9 134^1
%%%%%%%%%%%%%%%%%%%%%%%%%%%%%%%%%%%%%%%%%%%%%%%%%%%%%%%%%%%%%%%%%%%%%%%%%%%%%%%%%%%%%%%%%%
%%%%%%%%%%%%%%%%%%%%%%%%%%%%%%%%%%%%%%%%%%%%%%%%%%%%%%%%%%%%%%%%%%%%%%%%%%%%%%%%%%%%%%%%%%

%ADFvfySectionEnd 

\end{document}